\documentclass[12pt,reqno]{amsart}

\usepackage{amsmath}
\usepackage{amsfonts}
\usepackage{amsbsy}
\usepackage{amssymb}
\usepackage[normalem]{ulem}
\usepackage{times}
\usepackage{mathrsfs}
\usepackage{exscale}
\usepackage{graphicx}
\usepackage[colorlinks,citecolor=blue,linkcolor=rouge,
            bookmarksopen,
            bookmarksnumbered
           ]{hyperref}

\allowdisplaybreaks

\usepackage{times}
\usepackage[nomessages]{fp}

\ifx\figforTeXisloaded\relax \else\global\let\figforTeXisloaded=\relax\fi
\catcode`\@=11
\ifx\ctr@ln@m\undefined\else%
    \immediate\write16{*** Fig4TeX WARNING : \string\ctr@ln@m\space already defined.}\fi
\def\ctr@ln@m#1{\ifx#1\undefined\else%
    \immediate\write16{*** Fig4TeX WARNING : \string#1 already defined.}\fi}
\ctr@ln@m\ctr@ld@f
\def\ctr@ld@f#1#2{\ctr@ln@m#2#1#2}
\ctr@ld@f\def\ctr@ln@w#1#2{\ctr@ln@m#2\csname#1\endcsname#2}
{\catcode`\/=0 \catcode`/\=12 /ctr@ld@f/gdef/BS@{\}}
\ctr@ld@f\def\ctr@lcsn@m#1{\expandafter\ifx\csname#1\endcsname\relax\else%
    \immediate\write16{*** Fig4TeX WARNING : \BS@\expandafter\string#1\space already defined.}\fi}
\ctr@ld@f\edef\colonc@tcode{\the\catcode`\:}
\ctr@ld@f\edef\semicolonc@tcode{\the\catcode`\;}
\ctr@ld@f\def\t@stc@tcodech@nge{{\let\c@tcodech@nged=\z@%
    \ifnum\colonc@tcode=\the\catcode`\:\else\let\c@tcodech@nged=\@ne\fi%
    \ifnum\semicolonc@tcode=\the\catcode`\;\else\let\c@tcodech@nged=\@ne\fi%
    \ifx\c@tcodech@nged\@ne%
    \immediate\write16{}
    \immediate\write16{!!!=============================================================!!!}
    \immediate\write16{ Fig4TeX WARNING:}
    \immediate\write16{ The category code of some characters has been changed, which will}
    \immediate\write16{ result in an error (message "Runaway argument?").}
    \immediate\write16{ This probably comes from another package that changed the category}
    \immediate\write16{ code after Fig4TeX was loaded. If that proves to be exact, the}
    \immediate\write16{ solution is to exchange the loading commands on top of your file}
    \immediate\write16{ so that Fig4TeX is loaded last. For example, in LaTeX, we should}
    \immediate\write16{ say :}
    \immediate\write16{\BS@ usepackage[french]{babel}}
    \immediate\write16{\BS@ usepackage{fig4tex}}
    \immediate\write16{!!!=============================================================!!!}
    \immediate\write16{}
    \fi}}
\ctr@ld@f\def\FigforTeX{F\kern-.05em i\kern-.05em g\kern-.1em\raise-.14em\hbox{4}\kern-.19em\TeX}
\ctr@ld@f\def\W@rnmesoldA#1{\W@rnmesold}
\ctr@ld@f\def\W@rnmesoldAB#1(#2){\W@rnmesold}
\ctr@ld@f\def\W@rnmesold{%
    \immediate\write16{}
    \immediate\write16{!!!=============================================================!!!}
    \immediate\write16{ Fig4TeX WARNING:}
    \immediate\write16{ The file to be compiled is not compatible with the current version}
    \immediate\write16{ of Fig4TeX. To fix that, upgrade the source file (mainly change \BS@ ps*}
    \immediate\write16{ macros by \BS@ fig* macros), or use fig4tex184.tex instead (\BS@ input fig4tex184}
    \immediate\write16{ or \BS@ usepackage{fig4tex184}).}
    \immediate\write16{!!!=============================================================!!!}
    \immediate\write16{}}
\ctr@ln@m\psbeginfig\let\psbeginfig\W@rnmesoldA
\ctr@ln@m\psset\let\psset\W@rnmesoldAB
\ctr@ln@m\pssetdefault\let\pssetdefault\W@rnmesoldAB
\ctr@ln@m\pssetupdate\let\pssetupdate\W@rnmesoldA
\ctr@ln@w{newdimen}\epsil@n\epsil@n=0.00005pt
\ctr@ln@w{newdimen}\Cepsil@n\Cepsil@n=0.005pt
\ctr@ln@w{newdimen}\dcq@\dcq@=254pt
\ctr@ln@w{newdimen}\PI@\PI@=3.141592pt
\ctr@ln@w{newdimen}\DemiPI@deg\DemiPI@deg=90pt
\ctr@ln@w{newdimen}\PI@deg\PI@deg=180pt
\ctr@ln@w{newdimen}\DePI@deg\DePI@deg=360pt
\ctr@ld@f\chardef\t@n=10
\ctr@ld@f\chardef\c@nt=100
\ctr@ld@f\chardef\@lxxiv=74
\ctr@ld@f\chardef\@xci=91
\ctr@ld@f\mathchardef\@nMnCQn=9949
\ctr@ld@f\chardef\@vi=6
\ctr@ld@f\chardef\@xxx=30
\ctr@ld@f\chardef\@lvi=56
\ctr@ld@f\chardef\@@lxxi=71
\ctr@ld@f\chardef\@lxxxv=85
\ctr@ld@f\mathchardef\@@mmmmlxviii=4068
\ctr@ld@f\mathchardef\@ccclx=360
\ctr@ld@f\mathchardef\@dccxx=720
\ctr@ln@w{newcount}\p@rtent \ctr@ln@w{newcount}\f@ctech \ctr@ln@w{newcount}\result@tent
\ctr@ln@w{newdimen}\v@lmin \ctr@ln@w{newdimen}\v@lmax \ctr@ln@w{newdimen}\v@leur
\ctr@ln@w{newdimen}\result@t\ctr@ln@w{newdimen}\result@@t
\ctr@ln@w{newdimen}\mili@u \ctr@ln@w{newdimen}\c@rre \ctr@ln@w{newdimen}\delt@
\ctr@ld@f\def\degT@rd{0.017453 }  
\ctr@ld@f\def\rdT@deg{57.295779 } 
\ctr@ln@m\v@leurseule
{\catcode`p=12 \catcode`t=12 \gdef\v@leurseule#1pt{#1}}
\ctr@ld@f\def\repdecn@mb#1{\expandafter\v@leurseule\the#1\space}
\ctr@ld@f\def\arct@n#1(#2,#3){{\v@lmin=#2\v@lmax=#3%
    \maxim@m{\mili@u}{-\v@lmin}{\v@lmin}\maxim@m{\c@rre}{-\v@lmax}{\v@lmax}%
    \delt@=\mili@u\m@ech\mili@u%
    \ifdim\c@rre>\@nMnCQn\mili@u\divide\v@lmax\tw@\c@lATAN\v@leur(\z@,\v@lmax)
    \else%
    \maxim@m{\mili@u}{-\v@lmin}{\v@lmin}\maxim@m{\c@rre}{-\v@lmax}{\v@lmax}%
    \m@ech\c@rre%
    \ifdim\mili@u>\@nMnCQn\c@rre\divide\v@lmin\tw@
    \maxim@m{\mili@u}{-\v@lmin}{\v@lmin}\c@lATAN\v@leur(\mili@u,\z@)%
    \else\c@lATAN\v@leur(\delt@,\v@lmax)\fi\fi%
    \ifdim\v@lmin<\z@\v@leur=-\v@leur\ifdim\v@lmax<\z@\advance\v@leur-\PI@%
    \else\advance\v@leur\PI@\fi\fi%
    \global\result@t=\v@leur}#1=\result@t}
\ctr@ld@f\def\m@ech#1{\ifdim#1>1.646pt\divide\mili@u\t@n\divide\c@rre\t@n\m@ech#1\fi}
\ctr@ld@f\def\c@lATAN#1(#2,#3){{\v@lmin=#2\v@lmax=#3\v@leur=\z@\delt@=\tw@ pt%
    \un@iter{0.785398}{\v@lmax<}%
    \un@iter{0.463648}{\v@lmax<}%
    \un@iter{0.244979}{\v@lmax<}%
    \un@iter{0.124355}{\v@lmax<}%
    \un@iter{0.062419}{\v@lmax<}%
    \un@iter{0.031240}{\v@lmax<}%
    \un@iter{0.015624}{\v@lmax<}%
    \un@iter{0.007812}{\v@lmax<}%
    \un@iter{0.003906}{\v@lmax<}%
    \un@iter{0.001953}{\v@lmax<}%
    \un@iter{0.000976}{\v@lmax<}%
    \un@iter{0.000488}{\v@lmax<}%
    \un@iter{0.000244}{\v@lmax<}%
    \un@iter{0.000122}{\v@lmax<}%
    \un@iter{0.000061}{\v@lmax<}%
    \un@iter{0.000030}{\v@lmax<}%
    \un@iter{0.000015}{\v@lmax<}%
    \global\result@t=\v@leur}#1=\result@t}
\ctr@ld@f\def\un@iter#1#2{%
    \divide\delt@\tw@\edef\dpmn@{\repdecn@mb{\delt@}}%
    \mili@u=\v@lmin%
    \ifdim#2\z@%
      \advance\v@lmin-\dpmn@\v@lmax\advance\v@lmax\dpmn@\mili@u%
      \advance\v@leur-#1pt%
    \else%
      \advance\v@lmin\dpmn@\v@lmax\advance\v@lmax-\dpmn@\mili@u%
      \advance\v@leur#1pt%
    \fi}
\ctr@ld@f\def\c@ssin#1#2#3{\expandafter\ifx\csname COS@\number#3\endcsname\relax\c@lCS{#3pt}%
    \expandafter\xdef\csname COS@\number#3\endcsname{\repdecn@mb\result@t}%
    \expandafter\xdef\csname SIN@\number#3\endcsname{\repdecn@mb\result@@t}\fi%
    \edef#1{\csname COS@\number#3\endcsname}\edef#2{\csname SIN@\number#3\endcsname}}
\ctr@ld@f\def\c@lCS#1{{\mili@u=#1\p@rtent=\@ne%
    \relax\ifdim\mili@u<\z@\red@ng<-\else\red@ng>+\fi\f@ctech=\p@rtent%
    \relax\ifdim\mili@u<\z@\mili@u=-\mili@u\f@ctech=-\f@ctech\fi\c@@lCS}}
\ctr@ld@f\def\c@@lCS{\v@lmin=\mili@u\c@rre=-\mili@u\advance\c@rre\DemiPI@deg\v@lmax=\c@rre%
    \mili@u\@@lxxi\mili@u\divide\mili@u\@@mmmmlxviii%
    \edef\v@larg{\repdecn@mb{\mili@u}}\mili@u=-\v@larg\mili@u%
    \edef\v@lmxde{\repdecn@mb{\mili@u}}%
    \c@rre\@@lxxi\c@rre\divide\c@rre\@@mmmmlxviii%
    \edef\v@largC{\repdecn@mb{\c@rre}}\c@rre=-\v@largC\c@rre%
    \edef\v@lmxdeC{\repdecn@mb{\c@rre}}%
    \fctc@s\mili@u\v@lmin\global\result@t\p@rtent\v@leur%
    \let\t@mp=\v@larg\let\v@larg=\v@largC\let\v@largC=\t@mp%
    \let\t@mp=\v@lmxde\let\v@lmxde=\v@lmxdeC\let\v@lmxdeC=\t@mp%
    \fctc@s\c@rre\v@lmax\global\result@@t\f@ctech\v@leur}
\ctr@ld@f\def\fctc@s#1#2{\v@leur=#1\relax\ifdim#2<\@lxxxv\p@\cosser@h\else\sinser@t\fi}
\ctr@ld@f\def\cosser@h{\advance\v@leur\@lvi\p@\divide\v@leur\@lvi%
    \v@leur=\v@lmxde\v@leur\advance\v@leur\@xxx\p@%
    \v@leur=\v@lmxde\v@leur\advance\v@leur\@ccclx\p@%
    \v@leur=\v@lmxde\v@leur\advance\v@leur\@dccxx\p@\divide\v@leur\@dccxx}
\ctr@ld@f\def\sinser@t{\v@leur=\v@lmxdeC\p@\advance\v@leur\@vi\p@%
    \v@leur=\v@largC\v@leur\divide\v@leur\@vi}
\ctr@ld@f\def\red@ng#1#2{\relax\ifdim\mili@u#1#2\DemiPI@deg\advance\mili@u#2-\PI@deg%
    \p@rtent=-\p@rtent\red@ng#1#2\fi}
\ctr@ld@f\def\pr@c@lCS#1#2#3{\ctr@lcsn@m{COS@\number#3 }%
    \expandafter\xdef\csname COS@\number#3\endcsname{#1}%
    \expandafter\xdef\csname SIN@\number#3\endcsname{#2}}
\pr@c@lCS{1}{0}{0}
\pr@c@lCS{0.7071}{0.7071}{45}\pr@c@lCS{0.7071}{-0.7071}{-45}
\pr@c@lCS{0}{1}{90}          \pr@c@lCS{0}{-1}{-90}
\pr@c@lCS{-1}{0}{180}        \pr@c@lCS{-1}{0}{-180}
\pr@c@lCS{0}{-1}{270}        \pr@c@lCS{0}{1}{-270}
\ctr@ld@f\def\invers@#1#2{{\v@leur=#2\maxim@m{\v@lmax}{-\v@leur}{\v@leur}%
    \f@ctech=\@ne\m@inv@rs%
    \multiply\v@leur\f@ctech\edef\v@lv@leur{\repdecn@mb{\v@leur}}%
    \p@rtentiere{\p@rtent}{\v@leur}\v@lmin=\p@\divide\v@lmin\p@rtent%
    \inv@rs@\multiply\v@lmax\f@ctech\global\result@t=\v@lmax}#1=\result@t}
\ctr@ld@f\def\m@inv@rs{\ifdim\v@lmax<\p@\multiply\v@lmax\t@n\multiply\f@ctech\t@n\m@inv@rs\fi}
\ctr@ld@f\def\inv@rs@{\v@lmax=-\v@lmin\v@lmax=\v@lv@leur\v@lmax%
    \advance\v@lmax\tw@ pt\v@lmax=\repdecn@mb{\v@lmin}\v@lmax%
    \delt@=\v@lmax\advance\delt@-\v@lmin\ifdim\delt@<\z@\delt@=-\delt@\fi%
    \ifdim\delt@>\epsil@n\v@lmin=\v@lmax\inv@rs@\fi}
\ctr@ld@f\def\minim@m#1#2#3{\relax\ifdim#2<#3#1=#2\else#1=#3\fi}
\ctr@ld@f\def\maxim@m#1#2#3{\relax\ifdim#2>#3#1=#2\else#1=#3\fi}
\ctr@ld@f\def\p@rtentiere#1#2{#1=#2\divide#1by65536 }
\ctr@ld@f\def\r@undint#1#2{{\v@leur=#2\divide\v@leur\t@n\p@rtentiere{\p@rtent}{\v@leur}%
    \v@leur=\p@rtent pt\global\result@t=\t@n\v@leur}#1=\result@t}
\ctr@ld@f\def\sqrt@#1#2{{\v@leur=#2%
    \minim@m{\v@lmin}{\p@}{\v@leur}\maxim@m{\v@lmax}{\p@}{\v@leur}%
    \f@ctech=\@ne\m@sqrt@\sqrt@@%
    \mili@u=\v@lmin\advance\mili@u\v@lmax\divide\mili@u\tw@\multiply\mili@u\f@ctech%
    \global\result@t=\mili@u}#1=\result@t}
\ctr@ld@f\def\m@sqrt@{\ifdim\v@leur>\dcq@\divide\v@leur\c@nt\v@lmax=\v@leur%
    \multiply\f@ctech\t@n\m@sqrt@\fi}
\ctr@ld@f\def\sqrt@@{\mili@u=\v@lmin\advance\mili@u\v@lmax\divide\mili@u\tw@%
    \c@rre=\repdecn@mb{\mili@u}\mili@u%
    \ifdim\c@rre<\v@leur\v@lmin=\mili@u\else\v@lmax=\mili@u\fi%
    \delt@=\v@lmax\advance\delt@-\v@lmin\ifdim\delt@>\epsil@n\sqrt@@\fi}
\ctr@ld@f\def\extrairelepremi@r#1\de#2{\expandafter\lepremi@r#2@#1#2}
\ctr@ld@f\def\lepremi@r#1,#2@#3#4{\def#3{#1}\def#4{#2}\ignorespaces}
\ctr@ld@f\def\@cfor#1:=#2\do#3{%
  \edef\@fortemp{#2}%
  \ifx\@fortemp\empty\else\@cforloop#2,\@nil,\@nil\@@#1{#3}\fi}
\ctr@ln@m\@nextwhile
\ctr@ld@f\def\@cforloop#1,#2\@@#3#4{%
  \def#3{#1}%
  \ifx#3\Fig@nnil\let\@nextwhile=\Fig@fornoop\else#4\relax\let\@nextwhile=\@cforloop\fi%
  \@nextwhile#2\@@#3{#4}}

\ctr@ld@f\def\@ecfor#1:=#2\do#3{%
  \def\@@cfor{\@cfor#1:=}%
  \edef\@@@cfor{#2}%
  \expandafter\@@cfor\@@@cfor\do{#3}}
\ctr@ld@f\def\Fig@nnil{\@nil}
\ctr@ld@f\def\Fig@fornoop#1\@@#2#3{}
\ctr@ln@m\list@@rg
\ctr@ld@f\def\trtlis@rg#1#2{\def\list@@rg{#1}%
    \@ecfor\p@rv@l:=\list@@rg\do{\expandafter#2\p@rv@l|}}
\ctr@ld@f\def\trtlis@rgtok#1{\let@xte={}\let\n@xt\addt@t@xt\addt@t@xt #1}
\ctr@ln@m\M@cro
\ctr@ln@m\n@xt
\ctr@ld@f\def\addt@t@xt#1{\if#1|\let\n@xt\relax\else%
    \if#1,\expandafter\M@cro\the\let@xte|\let@xte={}%
    \else\let@xte=\expandafter{\the\let@xte #1}\fi\fi\n@xt}
\ctr@ln@w{newbox}\b@xvisu
\ctr@ln@w{newtoks}\let@xte
\ctr@ln@w{newif}\ifitis@K
\ctr@ln@w{newcount}\s@mme
\ctr@ln@w{newcount}\l@mbd@un \ctr@ln@w{newcount}\l@mbd@de
\ctr@ln@w{newcount}\superc@ntr@l\superc@ntr@l=\@ne        
\ctr@ln@w{newcount}\typec@ntr@l\typec@ntr@l=\superc@ntr@l 
\ctr@ln@w{newdimen}\v@lX  \ctr@ln@w{newdimen}\v@lY  \ctr@ln@w{newdimen}\v@lZ
\ctr@ln@w{newdimen}\v@lXa \ctr@ln@w{newdimen}\v@lYa \ctr@ln@w{newdimen}\v@lZa
\ctr@ln@w{newdimen}\unit@\unit@=\p@ 
\ctr@ld@f\def\unit@util{pt}
\ctr@ld@f\def\ptT@ptps{0.996264 }
\ctr@ld@f\def\ptpsT@pt{1.00375 }
\ctr@ld@f\def\ptT@unit@{1} 
\ctr@ld@f\def\setunit@#1{\def\unit@util{#1}\setunit@@#1:\invers@{\result@t}{\unit@}%
    \edef\ptT@unit@{\repdecn@mb\result@t}}
\ctr@ld@f\def\setunit@@#1#2:{\ifcat#1a\unit@=\@ne#1#2\else\unit@=#1#2\fi}
\ctr@ld@f\def\d@fm@cdim#1#2{{\v@leur=#2\v@leur=\ptT@unit@\v@leur\xdef#1{\repdecn@mb\v@leur}}}
\ctr@ln@w{newif}\ifBdingB@x\BdingB@xtrue
\ctr@ln@w{newdimen}\c@@rdXmin \ctr@ln@w{newdimen}\c@@rdYmin  
\ctr@ln@w{newdimen}\c@@rdXmax \ctr@ln@w{newdimen}\c@@rdYmax
\ctr@ld@f\def\b@undb@x#1#2{\ifBdingB@x%
    \relax\ifdim#1<\c@@rdXmin\global\c@@rdXmin=#1\fi%
    \relax\ifdim#2<\c@@rdYmin\global\c@@rdYmin=#2\fi%
    \relax\ifdim#1>\c@@rdXmax\global\c@@rdXmax=#1\fi%
    \relax\ifdim#2>\c@@rdYmax\global\c@@rdYmax=#2\fi\fi}
\ctr@ld@f\def\b@undb@xP#1{{\Figg@tXY{#1}\b@undb@x{\v@lX}{\v@lY}}}
\ctr@ld@f\def\ellBB@x#1;#2,#3(#4,#5,#6){{\s@uvc@ntr@l\et@tellBB@x%
    \setc@ntr@l{2}\figptell-2::#1;#2,#3(#4,#6)\b@undb@xP{-2}%
    \figptell-2::#1;#2,#3(#5,#6)\b@undb@xP{-2}%
    \c@ssin{\C@}{\S@}{#6}\v@lmin=\C@ pt\v@lmax=\S@ pt%
    \mili@u=#3\v@lmin\delt@=#2\v@lmax\arct@n\v@leur(\delt@,\mili@u)%
    \mili@u=-#3\v@lmax\delt@=#2\v@lmin\arct@n\c@rre(\delt@,\mili@u)%
    \v@leur=\rdT@deg\v@leur\advance\v@leur-\DePI@deg%
    \c@rre=\rdT@deg\c@rre\advance\c@rre-\DePI@deg%
    \v@lmin=#4pt\v@lmax=#5pt%
    \loop\ifdim\v@leur<\v@lmax\ifdim\v@leur>\v@lmin%
    \edef\@ngle{\repdecn@mb\v@leur}\figptell-2::#1;#2,#3(\@ngle,#6)%
    \b@undb@xP{-2}\fi\advance\v@leur\PI@deg\repeat%
    \loop\ifdim\c@rre<\v@lmax\ifdim\c@rre>\v@lmin%
    \edef\@ngle{\repdecn@mb\c@rre}\figptell-2::#1;#2,#3(\@ngle,#6)%
    \b@undb@xP{-2}\fi\advance\c@rre\PI@deg\repeat%
    \resetc@ntr@l\et@tellBB@x}\ignorespaces}
\ctr@ld@f\def\initb@undb@x{\c@@rdXmin=\maxdimen\c@@rdYmin=\maxdimen%
    \c@@rdXmax=-\maxdimen\c@@rdYmax=-\maxdimen}
\ctr@ld@f\def\c@ntr@lnum#1{%
    \relax\ifnum\typec@ntr@l=\@ne%
    \ifnum#1<\z@%
    \immediate\write16{*** Forbidden point number (#1). Abort.}\end\fi\fi%
    \set@bjc@de{#1}}
\ctr@ln@m\objc@de
\ctr@ld@f\def\set@bjc@de#1{\edef\objc@de{@BJ\ifnum#1<\z@ M\romannumeral-#1\else\romannumeral#1\fi}}
\s@mme=\m@ne\loop\ifnum\s@mme>-19
  \set@bjc@de{\s@mme}\ctr@lcsn@m\objc@de\ctr@lcsn@m{\objc@de T}
\advance\s@mme\m@ne\repeat
\s@mme=\@ne\loop\ifnum\s@mme<6
  \set@bjc@de{\s@mme}\ctr@lcsn@m\objc@de\ctr@lcsn@m{\objc@de T}
\advance\s@mme\@ne\repeat
\ctr@ld@f\def\setc@ntr@l#1{\ifnum\superc@ntr@l>#1\typec@ntr@l=\superc@ntr@l%
    \else\typec@ntr@l=#1\fi}
\ctr@ld@f\def\resetc@ntr@l#1{\global\superc@ntr@l=#1\setc@ntr@l{#1}}
\ctr@ld@f\def\s@uvc@ntr@l#1{\edef#1{\the\superc@ntr@l}}
\ctr@ln@m\c@lproscal
\ctr@ld@f\def\c@lproscalDD#1[#2,#3]{{\Figg@tXY{#2}%
    \edef\Xu@{\repdecn@mb{\v@lX}}\edef\Yu@{\repdecn@mb{\v@lY}}\Figg@tXY{#3}%
    \global\result@t=\Xu@\v@lX\global\advance\result@t\Yu@\v@lY}#1=\result@t}
\ctr@ld@f\def\c@lproscalTD#1[#2,#3]{{\Figg@tXY{#2}\edef\Xu@{\repdecn@mb{\v@lX}}%
    \edef\Yu@{\repdecn@mb{\v@lY}}\edef\Zu@{\repdecn@mb{\v@lZ}}%
    \Figg@tXY{#3}\global\result@t=\Xu@\v@lX\global\advance\result@t\Yu@\v@lY%
    \global\advance\result@t\Zu@\v@lZ}#1=\result@t}
\ctr@ld@f\def\c@lprovec#1{%
    \det@rmC\v@lZa(\v@lX,\v@lY,\v@lmin,\v@lmax)%
    \det@rmC\v@lXa(\v@lY,\v@lZ,\v@lmax,\v@leur)%
    \det@rmC\v@lYa(\v@lZ,\v@lX,\v@leur,\v@lmin)%
    \Figv@ctCreg#1(\v@lXa,\v@lYa,\v@lZa)}
\ctr@ld@f\def\det@rm#1[#2,#3]{{\Figg@tXY{#2}\Figg@tXYa{#3}%
    \delt@=\repdecn@mb{\v@lX}\v@lYa\advance\delt@-\repdecn@mb{\v@lY}\v@lXa%
    \global\result@t=\delt@}#1=\result@t}
\ctr@ld@f\def\det@rmC#1(#2,#3,#4,#5){{\global\result@t=\repdecn@mb{#2}#5%
    \global\advance\result@t-\repdecn@mb{#3}#4}#1=\result@t}
\ctr@ld@f\def\getredf@ctDD#1(#2,#3){{\maxim@m{\v@lXa}{-#2}{#2}\maxim@m{\v@lYa}{-#3}{#3}%
    \maxim@m{\v@lXa}{\v@lXa}{\v@lYa}
    \ifdim\v@lXa>\@xci pt\divide\v@lXa\@xci%
    \p@rtentiere{\p@rtent}{\v@lXa}\advance\p@rtent\@ne\else\p@rtent=\@ne\fi%
    \global\result@tent=\p@rtent}#1=\result@tent\ignorespaces}
\ctr@ld@f\def\getredf@ctTD#1(#2,#3,#4){{\maxim@m{\v@lXa}{-#2}{#2}\maxim@m{\v@lYa}{-#3}{#3}%
    \maxim@m{\v@lZa}{-#4}{#4}\maxim@m{\v@lXa}{\v@lXa}{\v@lYa}%
    \maxim@m{\v@lXa}{\v@lXa}{\v@lZa}
    \ifdim\v@lXa>\@lxxiv pt\divide\v@lXa\@lxxiv%
    \p@rtentiere{\p@rtent}{\v@lXa}\advance\p@rtent\@ne\else\p@rtent=\@ne\fi%
    \global\result@tent=\p@rtent}#1=\result@tent\ignorespaces}
\ctr@ln@m\getredf@ctB
\ctr@ld@f\def\getredf@ctBDD#1{\getredf@ctDD#1(\v@lX,\v@lY)}
\ctr@ld@f\def\getredf@ctBTD#1{\getredf@ctTD#1(\v@lX,\v@lY,\v@lZ)}
\ctr@ld@f\def\FigptintercircB@zDD#1:#2:#3,#4[#5,#6,#7,#8]{{\s@uvc@ntr@l\et@tfigptintercircB@zDD%
    \setc@ntr@l{2}\figvectPDD-1[#5,#8]\Figg@tXY{-1}\getredf@ctDD\f@ctech(\v@lX,\v@lY)%
    \mili@u=#4\unit@\divide\mili@u\f@ctech\c@rre=\repdecn@mb{\mili@u}\mili@u%
    \figptBezierDD-5::#3[#5,#6,#7,#8]%
    \v@lmin=#3\p@\v@lmax=\v@lmin\advance\v@lmax0.1\p@%
    \loop\edef\T@{\repdecn@mb{\v@lmax}}\figptBezierDD-2::\T@[#5,#6,#7,#8]%
    \figvectPDD-1[-5,-2]\n@rmeucCDD{\delt@}{-1}\ifdim\delt@<\c@rre\v@lmin=\v@lmax%
    \advance\v@lmax0.1\p@\repeat%
    \loop\mili@u=\v@lmin\advance\mili@u\v@lmax%
    \divide\mili@u\tw@\edef\T@{\repdecn@mb{\mili@u}}\figptBezierDD-2::\T@[#5,#6,#7,#8]%
    \figvectPDD-1[-5,-2]\n@rmeucCDD{\delt@}{-1}\ifdim\delt@>\c@rre\v@lmax=\mili@u%
    \else\v@lmin=\mili@u\fi\v@leur=\v@lmax\advance\v@leur-\v@lmin%
    \ifdim\v@leur>\epsil@n\repeat\figptcopyDD#1:#2/-2/%
    \resetc@ntr@l\et@tfigptintercircB@zDD}\ignorespaces}
\ctr@ln@m\figptinterlines
\ctr@ld@f\def\inters@cDD#1:#2[#3,#4;#5,#6]{{\s@uvc@ntr@l\et@tinters@cDD%
    \setc@ntr@l{2}\vecunit@{-1}{#4}\vecunit@{-2}{#6}%
    \Figg@tXY{-1}\setc@ntr@l{1}\Figg@tXYa{#3}%
    \edef\A@{\repdecn@mb{\v@lX}}\edef\B@{\repdecn@mb{\v@lY}}%
    \v@lmin=\B@\v@lXa\advance\v@lmin-\A@\v@lYa%
    \Figg@tXYa{#5}\setc@ntr@l{2}\Figg@tXY{-2}%
    \edef\C@{\repdecn@mb{\v@lX}}\edef\D@{\repdecn@mb{\v@lY}}%
    \v@lmax=\D@\v@lXa\advance\v@lmax-\C@\v@lYa%
    \delt@=\A@\v@lY\advance\delt@-\B@\v@lX%
    \invers@{\v@leur}{\delt@}\edef\v@ldelta{\repdecn@mb{\v@leur}}%
    \v@lXa=\A@\v@lmax\advance\v@lXa-\C@\v@lmin%
    \v@lYa=\B@\v@lmax\advance\v@lYa-\D@\v@lmin%
    \v@lXa=\v@ldelta\v@lXa\v@lYa=\v@ldelta\v@lYa%
    \setc@ntr@l{1}\Figp@intregDD#1:{#2}(\v@lXa,\v@lYa)%
    \resetc@ntr@l\et@tinters@cDD}\ignorespaces}
\ctr@ld@f\def\inters@cTD#1:#2[#3,#4;#5,#6]{{\s@uvc@ntr@l\et@tinters@cTD%
    \setc@ntr@l{2}\figvectNVTD-1[#4,#6]\figvectNVTD-2[#6,-1]\figvectPTD-1[#3,#5]%
    \r@pPSTD\v@leur[-2,-1,#4]\edef\v@lcoef{\repdecn@mb{\v@leur}}%
    \figpttraTD#1:{#2}=#3/\v@lcoef,#4/\resetc@ntr@l\et@tinters@cTD}\ignorespaces}
\ctr@ld@f\def\r@pPSTD#1[#2,#3,#4]{{\Figg@tXY{#2}\edef\Xu@{\repdecn@mb{\v@lX}}%
    \edef\Yu@{\repdecn@mb{\v@lY}}\edef\Zu@{\repdecn@mb{\v@lZ}}%
    \Figg@tXY{#3}\v@lmin=\Xu@\v@lX\advance\v@lmin\Yu@\v@lY\advance\v@lmin\Zu@\v@lZ%
    \Figg@tXY{#4}\v@lmax=\Xu@\v@lX\advance\v@lmax\Yu@\v@lY\advance\v@lmax\Zu@\v@lZ%
    \invers@{\v@leur}{\v@lmax}\global\result@t=\repdecn@mb{\v@leur}\v@lmin}%
    #1=\result@t}
\ctr@ln@m\n@rminf
\ctr@ld@f\def\n@rminfDD#1#2{{\Figg@tXY{#2}\maxim@m{\v@lX}{\v@lX}{-\v@lX}%
    \maxim@m{\v@lY}{\v@lY}{-\v@lY}\maxim@m{\global\result@t}{\v@lX}{\v@lY}}%
    #1=\result@t}
\ctr@ld@f\def\n@rminfTD#1#2{{\Figg@tXY{#2}\maxim@m{\v@lX}{\v@lX}{-\v@lX}%
    \maxim@m{\v@lY}{\v@lY}{-\v@lY}\maxim@m{\v@lZ}{\v@lZ}{-\v@lZ}%
    \maxim@m{\v@lX}{\v@lX}{\v@lY}\maxim@m{\global\result@t}{\v@lX}{\v@lZ}}%
    #1=\result@t}
\ctr@ln@m\n@rmeucC
\ctr@ld@f\def\n@rmeucCDD#1#2{\Figg@tXY{#2}\divide\v@lX\f@ctech\divide\v@lY\f@ctech%
    #1=\repdecn@mb{\v@lX}\v@lX\v@lX=\repdecn@mb{\v@lY}\v@lY\advance#1\v@lX}
\ctr@ld@f\def\n@rmeucCTD#1#2{\Figg@tXY{#2}%
    \divide\v@lX\f@ctech\divide\v@lY\f@ctech\divide\v@lZ\f@ctech%
    #1=\repdecn@mb{\v@lX}\v@lX\v@lX=\repdecn@mb{\v@lY}\v@lY\advance#1\v@lX%
    \v@lX=\repdecn@mb{\v@lZ}\v@lZ\advance#1\v@lX}
\ctr@ln@m\n@rmeucSV
\ctr@ld@f\def\n@rmeucSVDD#1#2{{\Figg@tXY{#2}%
    \v@lXa=\repdecn@mb{\v@lX}\v@lX\v@lYa=\repdecn@mb{\v@lY}\v@lY%
    \advance\v@lXa\v@lYa\sqrt@{\global\result@t}{\v@lXa}}#1=\result@t}
\ctr@ld@f\def\n@rmeucSVTD#1#2{{\Figg@tXY{#2}\v@lXa=\repdecn@mb{\v@lX}\v@lX%
    \v@lYa=\repdecn@mb{\v@lY}\v@lY\v@lZa=\repdecn@mb{\v@lZ}\v@lZ%
    \advance\v@lXa\v@lYa\advance\v@lXa\v@lZa\sqrt@{\global\result@t}{\v@lXa}}#1=\result@t}
\ctr@ln@m\n@rmeuc
\ctr@ld@f\def\n@rmeucDD#1#2{{\Figg@tXY{#2}\getredf@ctDD\f@ctech(\v@lX,\v@lY)%
    \divide\v@lX\f@ctech\divide\v@lY\f@ctech%
    \v@lXa=\repdecn@mb{\v@lX}\v@lX\v@lYa=\repdecn@mb{\v@lY}\v@lY%
    \advance\v@lXa\v@lYa\sqrt@{\global\result@t}{\v@lXa}%
    \global\multiply\result@t\f@ctech}#1=\result@t}
\ctr@ld@f\def\n@rmeucTD#1#2{{\Figg@tXY{#2}\getredf@ctTD\f@ctech(\v@lX,\v@lY,\v@lZ)%
    \divide\v@lX\f@ctech\divide\v@lY\f@ctech\divide\v@lZ\f@ctech%
    \v@lXa=\repdecn@mb{\v@lX}\v@lX%
    \v@lYa=\repdecn@mb{\v@lY}\v@lY\v@lZa=\repdecn@mb{\v@lZ}\v@lZ%
    \advance\v@lXa\v@lYa\advance\v@lXa\v@lZa\sqrt@{\global\result@t}{\v@lXa}%
    \global\multiply\result@t\f@ctech}#1=\result@t}
\ctr@ln@m\vecunit@
\ctr@ld@f\def\vecunit@DD#1#2{{\Figg@tXY{#2}\getredf@ctDD\f@ctech(\v@lX,\v@lY)%
    \divide\v@lX\f@ctech\divide\v@lY\f@ctech%
    \Figv@ctCreg#1(\v@lX,\v@lY)\n@rmeucSV{\v@lYa}{#1}%
    \invers@{\v@lXa}{\v@lYa}\edef\v@lv@lXa{\repdecn@mb{\v@lXa}}%
    \v@lX=\v@lv@lXa\v@lX\v@lY=\v@lv@lXa\v@lY%
    \Figv@ctCreg#1(\v@lX,\v@lY)\multiply\v@lYa\f@ctech\global\result@t=\v@lYa}}
\ctr@ld@f\def\vecunit@TD#1#2{{\Figg@tXY{#2}\getredf@ctTD\f@ctech(\v@lX,\v@lY,\v@lZ)%
    \divide\v@lX\f@ctech\divide\v@lY\f@ctech\divide\v@lZ\f@ctech%
    \Figv@ctCreg#1(\v@lX,\v@lY,\v@lZ)\n@rmeucSV{\v@lYa}{#1}%
    \invers@{\v@lXa}{\v@lYa}\edef\v@lv@lXa{\repdecn@mb{\v@lXa}}%
    \v@lX=\v@lv@lXa\v@lX\v@lY=\v@lv@lXa\v@lY\v@lZ=\v@lv@lXa\v@lZ%
    \Figv@ctCreg#1(\v@lX,\v@lY,\v@lZ)\multiply\v@lYa\f@ctech\global\result@t=\v@lYa}}
\ctr@ld@f\def\vecunitC@TD[#1,#2]{\Figg@tXYa{#1}\Figg@tXY{#2}%
    \advance\v@lX-\v@lXa\advance\v@lY-\v@lYa\advance\v@lZ-\v@lZa\c@lvecunitTD}
\ctr@ld@f\def\vecunitCV@TD#1{\Figg@tXY{#1}\c@lvecunitTD}
\ctr@ld@f\def\c@lvecunitTD{\getredf@ctTD\f@ctech(\v@lX,\v@lY,\v@lZ)%
    \divide\v@lX\f@ctech\divide\v@lY\f@ctech\divide\v@lZ\f@ctech%
    \v@lXa=\repdecn@mb{\v@lX}\v@lX%
    \v@lYa=\repdecn@mb{\v@lY}\v@lY\v@lZa=\repdecn@mb{\v@lZ}\v@lZ%
    \advance\v@lXa\v@lYa\advance\v@lXa\v@lZa\sqrt@{\v@lYa}{\v@lXa}%
    \invers@{\v@lXa}{\v@lYa}\edef\v@lv@lXa{\repdecn@mb{\v@lXa}}%
    \v@lX=\v@lv@lXa\v@lX\v@lY=\v@lv@lXa\v@lY\v@lZ=\v@lv@lXa\v@lZ}
\ctr@ln@m\figgetangle
\ctr@ld@f\def\figgetangleDD#1[#2,#3,#4]{\ifGR@cri{\s@uvc@ntr@l\et@tfiggetangleDD\setc@ntr@l{2}%
    \figvectPDD-1[#2,#3]\figvectPDD-2[#2,#4]\vecunit@{-1}{-1}%
    \c@lproscalDD\delt@[-2,-1]\figvectNVDD-1[-1]\c@lproscalDD\v@leur[-2,-1]%
    \arct@n\v@lmax(\delt@,\v@leur)\v@lmax=\rdT@deg\v@lmax%
    \ifdim\v@lmax<\z@\advance\v@lmax\DePI@deg\fi\xdef#1{\repdecn@mb{\v@lmax}}%
    \resetc@ntr@l\et@tfiggetangleDD}\ignorespaces\fi}
\ctr@ld@f\def\figgetangleTD#1[#2,#3,#4,#5]{\ifGR@cri{\s@uvc@ntr@l\et@tfiggetangleTD\setc@ntr@l{2}%
    \figvectPTD-1[#2,#3]\figvectPTD-2[#2,#5]\figvectNVTD-3[-1,-2]%
    \figvectPTD-2[#2,#4]\figvectNVTD-4[-3,-1]%
    \vecunit@{-1}{-1}\c@lproscalTD\delt@[-2,-1]\c@lproscalTD\v@leur[-2,-4]%
    \arct@n\v@lmax(\delt@,\v@leur)\v@lmax=\rdT@deg\v@lmax%
    \ifdim\v@lmax<\z@\advance\v@lmax\DePI@deg\fi\xdef#1{\repdecn@mb{\v@lmax}}%
    \resetc@ntr@l\et@tfiggetangleTD}\ignorespaces\fi}    
\ctr@ld@f\def\figgetdist#1[#2,#3]{\ifGR@cri{\s@uvc@ntr@l\et@tfiggetdist\setc@ntr@l{2}%
    \figvectP-1[#2,#3]\n@rmeuc{\v@lX}{-1}\v@lX=\ptT@unit@\v@lX\xdef#1{\repdecn@mb{\v@lX}}%
    \resetc@ntr@l\et@tfiggetdist}\ignorespaces\fi}
\ctr@ld@f\def\figget#1=#2[#3]{\keln@mun#1|%
    \def\n@mref{a}\ifx\l@debut\n@mref\figgetangle#2[#3]\else
    \def\n@mref{d}\ifx\l@debut\n@mref\figgetdist#2[#3]\else
    \W@rnmeskwd{figget}{#1}\fi\fi\ignorespaces}
\ctr@ld@f\def\Figg@tT#1{\c@ntr@lnum{#1}%
    {\expandafter\expandafter\expandafter\extr@ctT\csname\objc@de\endcsname:%
     \ifnum\B@@ltxt=\z@\ptn@me{#1}\else\csname\objc@de T\endcsname\fi}}
\ctr@ld@f\def\extr@ctT#1,#2,#3/#4:{\def\B@@ltxt{#3}}
\ctr@ld@f\def\Figg@tXY#1{\c@ntr@lnum{#1}%
    \expandafter\expandafter\expandafter\extr@ctC\csname\objc@de\endcsname:}
\ctr@ln@m\extr@ctC
\ctr@ld@f\def\extr@ctCDD#1/#2,#3,#4:{\v@lX=#2\v@lY=#3}
\ctr@ld@f\def\extr@ctCTD#1/#2,#3,#4:{\v@lX=#2\v@lY=#3\v@lZ=#4}
\ctr@ld@f\def\Figg@tXYa#1{\c@ntr@lnum{#1}%
    \expandafter\expandafter\expandafter\extr@ctCa\csname\objc@de\endcsname:}
\ctr@ln@m\extr@ctCa
\ctr@ld@f\def\extr@ctCaDD#1/#2,#3,#4:{\v@lXa=#2\v@lYa=#3}
\ctr@ld@f\def\extr@ctCaTD#1/#2,#3,#4:{\v@lXa=#2\v@lYa=#3\v@lZa=#4}
\ctr@ln@m\t@xt@
\ctr@ld@f\def\figinit#1{\t@stc@tcodech@nge\initpr@lim\Figinit@#1,:\initpss@ttings\ignorespaces}
\ctr@ld@f\def\Figinit@#1,#2:{\setunit@{#1}\def\t@xt@{#2}\ifx\t@xt@\empty\else\Figinit@@#2:\fi}
\ctr@ld@f\def\Figinit@@#1#2:{\if#12 \else\Figs@tproj{#1}\initTD@\fi}
\ctr@ln@w{newif}\ifTr@isDim
\ctr@ld@f\def\UnD@fined{UNDEFINED}
\ctr@ln@m\@utoFN
\ctr@ln@m\@utoFInDone
\ctr@ln@m\disob@unit
\ctr@ld@f\def\initpr@lim{\initb@undb@x\figsetmark{}\figsetptname{$A_{##1}$}\def\Sc@leFact{1}%
    \initDD@\figsetroundcoord{yes}\GR@critrue\expandafter\setupd@te\D@FTupdate:%
    \edef\disob@unit{\UnD@fined}\edef\t@rgetpt{\UnD@fined}\gdef\@utoFInDone{1}\gdef\@utoFN{0}}
\ctr@ld@f\def\initDD@{\Tr@isDimfalse%
    \ifPDFm@ke%
     \let\Ps@rcerc=\Ps@rcercBz%
     \let\Ps@rell=\Ps@rellBz%
    \fi
    \let\c@lDCUn=\c@lDCUnDD%
    \let\c@lDCDeux=\c@lDCDeuxDD%
    \let\c@ldefproj=\relax%
    \let\c@lproscal=\c@lproscalDD%
    \let\c@lprojSP=\relax%
    \let\extr@ctC=\extr@ctCDD%
    \let\extr@ctCa=\extr@ctCaDD%
    \let\extr@ctCF=\extr@ctCFDD%
    \let\Figp@intreg=\Figp@intregDD%
    \let\Figpts@xes=\Figpts@xesDD%
    \let\getredf@ctB=\getredf@ctBDD%
    \let\n@rmeucSV=\n@rmeucSVDD\let\n@rmeuc=\n@rmeucDD\let\n@rmeucC\n@rmeucCDD\let\n@rminf=\n@rminfDD%
    \let\pr@dMatV=\pr@dMatVDD%
    \let\Q@@xes=\Q@@xesDD%
    \let\vecunit@=\vecunit@DD%
    \let\figcoord=\figcoordDD%
    \let\figgetangle=\figgetangleDD%
    \let\figpt=\figptDD%
    \let\figptBezier=\figptBezierDD%
    \let\figptbary=\figptbaryDD%
    \let\figptcirc=\figptcircDD%
    \let\figptcircumcenter=\figptcircumcenterDD%
    \let\figptcopy=\figptcopyDD%
    \let\figptcurvcenter=\figptcurvcenterDD%
    \let\figptell=\figptellDD%
    \let\figptendnormal=\figptendnormalDD%
    \let\figptinterlineplane=\figptinterlineplaneDD%
    \let\figptinterlines=\inters@cDD%
    \let\figptorthocenter=\figptorthocenterDD%
    \let\figptorthoprojline=\figptorthoprojlineDD%
    \let\figptorthoprojplane=\figptorthoprojplaneDD%
    \let\figptrot=\figptrotDD%
    \let\figptscontrol=\figptscontrolDD%
    \let\figptsintercirc=\figptsintercircDD%
    \let\figptsinterlinell=\figptsinterlinellDD%
    \let\figptsorthoprojline=\figptsorthoprojlineDD%
    \let\figptorthoprojplane=\figptorthoprojplaneDD%
    \let\figptsrot=\figptsrotDD%
    \let\figptssym=\figptssymDD%
    \let\figptstra=\figptstraDD%
    \let\figptsym=\figptsymDD%
    \let\figpttraC=\figpttraCDD%
    \let\figpttra=\figpttraDD%
    \let\figptvisilimSL=\figptvisilimSLDD%
    \let\figsetobdist=\figsetobdistDD%
    \let\figsettarget=\figsettargetDD%
    \let\figsetview=\figsetviewDD%
    \let\figvectDBezier=\figvectDBezierDD%
    \let\figvectN=\figvectNDD%
    \let\figvectNV=\figvectNVDD%
    \let\figvectP=\figvectPDD%
    \let\figvectU=\figvectUDD%
    \let\figdrawarccircP=\Q@arccircPDD%
    \let\figdrawarccirc=\Q@arccircDD%
    \let\figdrawarcell=\Q@arcellDD%
    \let\figdrawarcellPA=\Q@arcellPADD%
    \let\figdrawarrowBezier=\Q@arrowBezierDD%
    \let\figdrawarrowcircP=\Q@arrowcircPDD%
    \let\figdrawarrowcirc=\Q@arrowcircDD%
    \let\figdrawarrowhead=\Q@arrowheadDD%
    \let\figdrawarrow=\Q@arrowDD%
    \let\figdrawBezier=\Q@BezierDD%
    \let\figdrawcirc=\Q@circDD%
    \let\figdrawcurve=\Q@curveDD%
    \let\figdrawnormal=\Q@normalDD%
    }
\ctr@ld@f\def\initTD@{\Tr@isDimtrue\initb@undb@xTD\newt@rgetptfalse\newdis@bfalse%
    \let\c@lDCUn=\c@lDCUnTD%
    \let\c@lDCDeux=\c@lDCDeuxTD%
    \let\c@ldefproj=\c@ldefprojTD%
    \let\c@lproscal=\c@lproscalTD%
    \let\extr@ctC=\extr@ctCTD%
    \let\extr@ctCa=\extr@ctCaTD%
    \let\extr@ctCF=\extr@ctCFTD%
    \let\Figp@intreg=\Figp@intregTD%
    \let\Figpts@xes=\Figpts@xesTD%
    \let\getredf@ctB=\getredf@ctBTD%
    \let\n@rmeucSV=\n@rmeucSVTD\let\n@rmeuc=\n@rmeucTD\let\n@rmeucC\n@rmeucCTD\let\n@rminf=\n@rminfTD%
    \let\pr@dMatV=\pr@dMatVTD%
    \let\Q@@xes=\Q@@xesTD%
    \let\vecunit@=\vecunit@TD%
    \let\figcoord=\figcoordTD%
    \let\figgetangle=\figgetangleTD%
    \let\figpt=\figptTD%
    \let\figptBezier=\figptBezierTD%
    \let\figptbary=\figptbaryTD%
    \let\figptcirc=\figptcircTD%
    \let\figptcircumcenter=\figptcircumcenterTD%
    \let\figptcopy=\figptcopyTD%
    \let\figptcurvcenter=\figptcurvcenterTD%
    \let\figptinterlineplane=\figptinterlineplaneTD%
    \let\figptinterlines=\inters@cTD%
    \let\figptorthocenter=\figptorthocenterTD%
    \let\figptorthoprojline=\figptorthoprojlineTD%
    \let\figptorthoprojplane=\figptorthoprojplaneTD%
    \let\figptrot=\figptrotTD%
    \let\figptscontrol=\figptscontrolTD%
    \let\figptsintercirc=\figptsintercircTD%
    \let\figptsorthoprojline=\figptsorthoprojlineTD%
    \let\figptsorthoprojplane=\figptsorthoprojplaneTD%
    \let\figptsrot=\figptsrotTD%
    \let\figptssym=\figptssymTD%
    \let\figptstra=\figptstraTD%
    \let\figptsym=\figptsymTD%
    \let\figpttraC=\figpttraCTD%
    \let\figpttra=\figpttraTD%
    \let\figptvisilimSL=\figptvisilimSLTD%
    \let\figsetobdist=\figsetobdistTD%
    \let\figsettarget=\figsettargetTD%
    \let\figsetview=\figsetviewTD%
    \let\figvectDBezier=\figvectDBezierTD%
    \let\figvectN=\figvectNTD%
    \let\figvectNV=\figvectNVTD%
    \let\figvectP=\figvectPTD%
    \let\figvectU=\figvectUTD%
    \let\figdrawarccircP=\Q@arccircPTD%
    \let\figdrawarccirc=\Q@arccircTD%
    \let\figdrawarcell=\Q@arcellTD%
    \let\figdrawarcellPA=\Q@arcellPATD%
    \let\figdrawarrowBezier=\Q@arrowBezierTD%
    \let\figdrawarrowcircP=\Q@arrowcircPTD%
    \let\figdrawarrowcirc=\Q@arrowcircTD%
    \let\figdrawarrowhead=\Q@arrowheadTD%
    \let\figdrawarrow=\Q@arrowTD%
    \let\figdrawBezier=\Q@BezierTD%
    \let\figdrawcirc=\Q@circTD%
    \let\figdrawcurve=\Q@curveTD%
    }
\ctr@ld@f\def\un@v@ilable#1{\immediate\write16{*** The macro #1 is not available in the current context.}}
\ctr@ld@f\def\figinsert#1{{\def\t@xt@{#1}\relax%
    \ifx\t@xt@\empty\ifnum\@utoFInDone>\z@\Figinsert@\DefGIfilen@me,:\fi%
    \else\expandafter\FiginsertNu@#1 :\fi}\ignorespaces}
\ctr@ld@f\def\FiginsertNu@#1 #2:{\def\t@xt@{#1}\relax\ifx\t@xt@\empty\def\t@xt@{#2}%
    \ifx\t@xt@\empty\ifnum\@utoFInDone>\z@\Figinsert@\DefGIfilen@me,:\fi%
    \else\FiginsertNu@#2:\fi\else\expandafter\FiginsertNd@#1 #2:\fi}
\ctr@ld@f\def\FiginsertNd@#1#2:{\ifcat#1a\Figinsert@#1#2,:\else%
    \ifnum\@utoFInDone>\z@\Figinsert@\DefGIfilen@me,#1#2,:\fi\fi}
\ctr@ln@m\Sc@leFact
\ctr@ld@f\def\Figinsert@#1,#2:{\def\t@xt@{#2}\ifx\t@xt@\empty\xdef\Sc@leFact{1}\else%
    \X@rgdeux@#2\xdef\Sc@leFact{\@rgdeux}\fi%
    \Figdisc@rdLTS{#1}{\t@xt@}\@psfgetbb{\t@xt@}%
    \v@lX=\@psfllx\p@\v@lX=\ptpsT@pt\v@lX\v@lX=\Sc@leFact\v@lX%
    \v@lY=\@psflly\p@\v@lY=\ptpsT@pt\v@lY\v@lY=\Sc@leFact\v@lY%
    \b@undb@x{\v@lX}{\v@lY}%
    \v@lX=\@psfurx\p@\v@lX=\ptpsT@pt\v@lX\v@lX=\Sc@leFact\v@lX%
    \v@lY=\@psfury\p@\v@lY=\ptpsT@pt\v@lY\v@lY=\Sc@leFact\v@lY%
    \b@undb@x{\v@lX}{\v@lY}%
    \ifPDFm@ke\Figinclud@PDF{\t@xt@}{\Sc@leFact}\else%
    \v@lX=\c@nt pt\v@lX=\Sc@leFact\v@lX\edef\F@ct{\repdecn@mb{\v@lX}}%
    \ifx\TeXturesonMacOSltX\special{postscriptfile #1 vscale=\F@ct\space hscale=\F@ct}%
    \else\includegraphics{#1}\fi\fi%
    \message{[\t@xt@]}\ignorespaces}
\ctr@ld@f\def\Figdisc@rdLTS#1#2{\expandafter\Figdisc@rdLTS@#1 :#2}
\ctr@ld@f\def\Figdisc@rdLTS@#1 #2:#3{\def#3{#1}\relax\ifx#3\empty\expandafter\Figdisc@rdLTS@#2:#3\fi}
\ctr@ld@f\def\figinsertE#1{\FiginsertE@#1,:\ignorespaces}
\ctr@ld@f\def\FiginsertE@#1,#2:{{\def\t@xt@{#2}\ifx\t@xt@\empty\xdef\Sc@leFact{1}\else%
    \X@rgdeux@#2\xdef\Sc@leFact{\@rgdeux}\fi%
    \Figdisc@rdLTS{#1}{\t@xt@}\pdfximage{\t@xt@}%
    \setbox\Gb@x=\hbox{\pdfrefximage\pdflastximage}%
    \v@lX=\z@\v@lY=-\Sc@leFact\dp\Gb@x\b@undb@x{\v@lX}{\v@lY}%
    \advance\v@lX\Sc@leFact\wd\Gb@x\advance\v@lY\Sc@leFact\dp\Gb@x%
    \advance\v@lY\Sc@leFact\ht\Gb@x\b@undb@x{\v@lX}{\v@lY}%
    \v@lX=\Sc@leFact\wd\Gb@x\pdfximage width \v@lX {\t@xt@}%
    \rlap{\pdfrefximage\pdflastximage}\message{[\t@xt@]}}\ignorespaces}
\ctr@ld@f\def\X@rgdeux@#1,{\edef\@rgdeux{#1}}
\ctr@ln@m\figpt
\ctr@ld@f\def\figptDD#1:#2(#3,#4){\ifGR@cri\c@ntr@lnum{#1}%
    {\v@lX=#3\unit@\v@lY=#4\unit@\Fig@dmpt{#2}{\z@}}\ignorespaces\fi}
\ctr@ld@f\def\Fig@dmpt#1#2{\def\t@xt@{#1}\ifx\t@xt@\empty\def\B@@ltxt{\z@}%
    \else\expandafter\gdef\csname\objc@de T\endcsname{#1}\def\B@@ltxt{\@ne}\fi%
    \expandafter\xdef\csname\objc@de\endcsname{\ifitis@vect@r\C@dCl@svect%
    \else\C@dCl@spt\fi,\z@,\B@@ltxt/\the\v@lX,\the\v@lY,#2}}
\ctr@ld@f\def\C@dCl@spt{P}
\ctr@ld@f\def\C@dCl@svect{V}
\ctr@ln@m\c@@rdYZ
\ctr@ln@m\c@@rdY
\ctr@ld@f\def\figptTD#1:#2(#3,#4){\ifGR@cri\c@ntr@lnum{#1}%
    \def\c@@rdYZ{#4,0,0}\extrairelepremi@r\c@@rdY\de\c@@rdYZ%
    \extrairelepremi@r\c@@rdZ\de\c@@rdYZ%
    {\v@lX=#3\unit@\v@lY=\c@@rdY\unit@\v@lZ=\c@@rdZ\unit@\Fig@dmpt{#2}{\the\v@lZ}%
    \b@undb@xTD{\v@lX}{\v@lY}{\v@lZ}}\ignorespaces\fi}
\ctr@ln@m\Figp@intreg
\ctr@ld@f\def\Figp@intregDD#1:#2(#3,#4){\c@ntr@lnum{#1}%
    {\result@t=#4\v@lX=#3\v@lY=\result@t\Fig@dmpt{#2}{\z@}}\ignorespaces}
\ctr@ld@f\def\Figp@intregTD#1:#2(#3,#4){\c@ntr@lnum{#1}%
    \def\c@@rdYZ{#4,\z@,\z@}\extrairelepremi@r\c@@rdY\de\c@@rdYZ%
    \extrairelepremi@r\c@@rdZ\de\c@@rdYZ%
    {\v@lX=#3\v@lY=\c@@rdY\v@lZ=\c@@rdZ\Fig@dmpt{#2}{\the\v@lZ}%
    \b@undb@xTD{\v@lX}{\v@lY}{\v@lZ}}\ignorespaces}
\ctr@ln@m\figptBezier
\ctr@ld@f\def\figptBezierDD#1:#2:#3[#4,#5,#6,#7]{\ifGR@cri{\s@uvc@ntr@l\et@tfigptBezierDD%
    \FigptBezier@#3[#4,#5,#6,#7]\Figp@intregDD#1:{#2}(\v@lX,\v@lY)%
    \resetc@ntr@l\et@tfigptBezierDD}\ignorespaces\fi}
\ctr@ld@f\def\figptBezierTD#1:#2:#3[#4,#5,#6,#7]{\ifGR@cri{\s@uvc@ntr@l\et@tfigptBezierTD%
    \FigptBezier@#3[#4,#5,#6,#7]\Figp@intregTD#1:{#2}(\v@lX,\v@lY,\v@lZ)%
    \resetc@ntr@l\et@tfigptBezierTD}\ignorespaces\fi}
\ctr@ld@f\def\FigptBezier@#1[#2,#3,#4,#5]{\setc@ntr@l{2}%
    \edef\T@{#1}\v@leur=\p@\advance\v@leur-#1pt\edef\UNmT@{\repdecn@mb{\v@leur}}%
    \figptcopy-4:/#2/\figptcopy-3:/#3/\figptcopy-2:/#4/\figptcopy-1:/#5/%
    \l@mbd@un=-4 \l@mbd@de=-\thr@@\p@rtent=\m@ne\c@lDecast%
    \l@mbd@un=-4 \l@mbd@de=-\thr@@\p@rtent=-\tw@\c@lDecast%
    \l@mbd@un=-4 \l@mbd@de=-\thr@@\p@rtent=-\thr@@\c@lDecast\Figg@tXY{-4}}
\ctr@ln@m\c@lDCUn
\ctr@ld@f\def\c@lDCUnDD#1#2{\Figg@tXY{#1}\v@lX=\UNmT@\v@lX\v@lY=\UNmT@\v@lY%
    \Figg@tXYa{#2}\advance\v@lX\T@\v@lXa\advance\v@lY\T@\v@lYa%
    \Figp@intregDD#1:(\v@lX,\v@lY)}
\ctr@ld@f\def\c@lDCUnTD#1#2{\Figg@tXY{#1}\v@lX=\UNmT@\v@lX\v@lY=\UNmT@\v@lY\v@lZ=\UNmT@\v@lZ%
    \Figg@tXYa{#2}\advance\v@lX\T@\v@lXa\advance\v@lY\T@\v@lYa\advance\v@lZ\T@\v@lZa%
    \Figp@intregTD#1:(\v@lX,\v@lY,\v@lZ)}
\ctr@ld@f\def\c@lDecast{\relax\ifnum\l@mbd@un<\p@rtent\c@lDCUn{\l@mbd@un}{\l@mbd@de}%
    \advance\l@mbd@un\@ne\advance\l@mbd@de\@ne\c@lDecast\fi}
\ctr@ld@f\def\figptmap#1:#2=#3/#4/#5/{\ifGR@cri{\s@uvc@ntr@l\et@tfigptmap%
    \setc@ntr@l{2}\figvectP-1[#4,#3]\Figg@tXY{-1}%
    \pr@dMatV/#5/\figpttra#1:{#2}=#4/1,-1/%
    \resetc@ntr@l\et@tfigptmap}\ignorespaces\fi}
\ctr@ln@m\pr@dMatV
\ctr@ld@f\def\pr@dMatVDD/#1,#2;#3,#4/{\v@lXa=#1\v@lX\advance\v@lXa#2\v@lY%
    \v@lYa=#3\v@lX\advance\v@lYa#4\v@lY\Figv@ctCreg-1(\v@lXa,\v@lYa)}
\ctr@ld@f\def\pr@dMatVTD/#1,#2,#3;#4,#5,#6;#7,#8,#9/{%
    \v@lXa=#1\v@lX\advance\v@lXa#2\v@lY\advance\v@lXa#3\v@lZ%
    \v@lYa=#4\v@lX\advance\v@lYa#5\v@lY\advance\v@lYa#6\v@lZ%
    \v@lZa=#7\v@lX\advance\v@lZa#8\v@lY\advance\v@lZa#9\v@lZ%
    \Figv@ctCreg-1(\v@lXa,\v@lYa,\v@lZa)}
\ctr@ln@m\figptbary
\ctr@ld@f\def\figptbaryDD#1:#2[#3;#4]{\ifGR@cri{\edef\list@num{#3}\extrairelepremi@r\p@int\de\list@num%
    \s@mme=\z@\@ecfor\c@ef:=#4\do{\advance\s@mme\c@ef}%
    \edef\listec@ef{#4,0}\extrairelepremi@r\c@ef\de\listec@ef%
    \Figg@tXY{\p@int}\divide\v@lX\s@mme\divide\v@lY\s@mme%
    \multiply\v@lX\c@ef\multiply\v@lY\c@ef%
    \@ecfor\p@int:=\list@num\do{\extrairelepremi@r\c@ef\de\listec@ef%
           \Figg@tXYa{\p@int}\divide\v@lXa\s@mme\divide\v@lYa\s@mme%
           \multiply\v@lXa\c@ef\multiply\v@lYa\c@ef%
           \advance\v@lX\v@lXa\advance\v@lY\v@lYa}%
    \Figp@intregDD#1:{#2}(\v@lX,\v@lY)}\ignorespaces\fi}
\ctr@ld@f\def\figptbaryTD#1:#2[#3;#4]{\ifGR@cri{\edef\list@num{#3}\extrairelepremi@r\p@int\de\list@num%
    \s@mme=\z@\@ecfor\c@ef:=#4\do{\advance\s@mme\c@ef}%
    \edef\listec@ef{#4,0}\extrairelepremi@r\c@ef\de\listec@ef%
    \Figg@tXY{\p@int}\divide\v@lX\s@mme\divide\v@lY\s@mme\divide\v@lZ\s@mme%
    \multiply\v@lX\c@ef\multiply\v@lY\c@ef\multiply\v@lZ\c@ef%
    \@ecfor\p@int:=\list@num\do{\extrairelepremi@r\c@ef\de\listec@ef%
           \Figg@tXYa{\p@int}\divide\v@lXa\s@mme\divide\v@lYa\s@mme\divide\v@lZa\s@mme%
           \multiply\v@lXa\c@ef\multiply\v@lYa\c@ef\multiply\v@lZa\c@ef%
           \advance\v@lX\v@lXa\advance\v@lY\v@lYa\advance\v@lZ\v@lZa}%
    \Figp@intregTD#1:{#2}(\v@lX,\v@lY,\v@lZ)}\ignorespaces\fi}
\ctr@ld@f\def\figptbaryR#1:#2[#3;#4]{\ifGR@cri{%
    \v@leur=\z@\@ecfor\c@ef:=#4\do{\maxim@m{\v@lmax}{\c@ef pt}{-\c@ef pt}%
    \ifdim\v@lmax>\v@leur\v@leur=\v@lmax\fi}%
    \ifdim\v@leur<\p@\f@ctech=\@M\else\ifdim\v@leur<\t@n\p@\f@ctech=\@m\else%
    \ifdim\v@leur<\c@nt\p@\f@ctech=\c@nt\else\ifdim\v@leur<\@m\p@\f@ctech=\t@n\else%
    \f@ctech=\@ne\fi\fi\fi\fi%
    \def\listec@ef{0}%
    \@ecfor\c@ef:=#4\do{\sc@lec@nvRI{\c@ef pt}\edef\listec@ef{\listec@ef,\the\s@mme}}%
    \extrairelepremi@r\c@ef\de\listec@ef\figptbary#1:#2[#3;\listec@ef]}\ignorespaces\fi}
\ctr@ld@f\def\sc@lec@nvRI#1{\v@leur=#1\p@rtentiere{\s@mme}{\v@leur}\advance\v@leur-\s@mme\p@%
    \multiply\v@leur\f@ctech\p@rtentiere{\p@rtent}{\v@leur}%
    \multiply\s@mme\f@ctech\advance\s@mme\p@rtent}
\ctr@ln@m\figptcirc
\ctr@ld@f\def\figptcircDD#1:#2:#3;#4(#5){\ifGR@cri{\s@uvc@ntr@l\et@tfigptcircDD%
    \c@lptellDD#1:{#2}:#3;#4,#4(#5)\resetc@ntr@l\et@tfigptcircDD}\ignorespaces\fi}
\ctr@ld@f\def\figptcircTD#1:#2:#3,#4,#5;#6(#7){\ifGR@cri{\s@uvc@ntr@l\et@tfigptcircTD%
    \setc@ntr@l{2}\c@lExtAxes#3,#4,#5(#6)\figptellP#1:{#2}:#3,-4,-5(#7)%
    \resetc@ntr@l\et@tfigptcircTD}\ignorespaces\fi}
\ctr@ln@m\figptcircumcenter
\ctr@ld@f\def\figptcircumcenterDD#1:#2[#3,#4,#5]{\ifGR@cri{\s@uvc@ntr@l\et@tfigptcircumcenterDD%
    \setc@ntr@l{2}\figvectNDD-5[#3,#4]\figptbaryDD-3:[#3,#4;1,1]%
                  \figvectNDD-6[#4,#5]\figptbaryDD-4:[#4,#5;1,1]%
    \resetc@ntr@l{2}\inters@cDD#1:{#2}[-3,-5;-4,-6]%
    \resetc@ntr@l\et@tfigptcircumcenterDD}\ignorespaces\fi}
\ctr@ld@f\def\figptcircumcenterTD#1:#2[#3,#4,#5]{\ifGR@cri{\s@uvc@ntr@l\et@tfigptcircumcenterTD%
    \setc@ntr@l{2}\figvectNTD-1[#3,#4,#5]%
    \figvectPTD-3[#3,#4]\figvectNVTD-5[-1,-3]\figptbaryTD-3:[#3,#4;1,1]%
    \figvectPTD-4[#4,#5]\figvectNVTD-6[-1,-4]\figptbaryTD-4:[#4,#5;1,1]%
    \resetc@ntr@l{2}\inters@cTD#1:{#2}[-3,-5;-4,-6]%
    \resetc@ntr@l\et@tfigptcircumcenterTD}\ignorespaces\fi}
\ctr@ln@m\figptcopy
\ctr@ld@f\def\figptcopyDD#1:#2/#3/{\ifGR@cri{\Figg@tXY{#3}%
    \Figp@intregDD#1:{#2}(\v@lX,\v@lY)}\ignorespaces\fi}
\ctr@ld@f\def\figptcopyTD#1:#2/#3/{\ifGR@cri{\Figg@tXY{#3}%
    \Figp@intregTD#1:{#2}(\v@lX,\v@lY,\v@lZ)}\ignorespaces\fi}
\ctr@ln@m\figptcurvcenter
\ctr@ld@f\def\figptcurvcenterDD#1:#2:#3[#4,#5,#6,#7]{\ifGR@cri{\s@uvc@ntr@l\et@tfigptcurvcenterDD%
    \setc@ntr@l{2}\c@lcurvradDD#3[#4,#5,#6,#7]\edef\Sprim@{\repdecn@mb{\result@t}}%
    \figptBezierDD-1::#3[#4,#5,#6,#7]\figpttraDD#1:{#2}=-1/\Sprim@,-5/%
    \resetc@ntr@l\et@tfigptcurvcenterDD}\ignorespaces\fi}
\ctr@ld@f\def\figptcurvcenterTD#1:#2:#3[#4,#5,#6,#7]{\ifGR@cri{\s@uvc@ntr@l\et@tfigptcurvcenterTD%
    \setc@ntr@l{2}\figvectDBezierTD -5:1,#3[#4,#5,#6,#7]%
    \figvectDBezierTD -6:2,#3[#4,#5,#6,#7]\vecunit@TD{-5}{-5}%
    \edef\Sprim@{\repdecn@mb{\result@t}}\figvectNVTD-1[-6,-5]%
    \figvectNVTD-5[-5,-1]\c@lproscalTD\v@leur[-6,-5]%
    \invers@{\v@leur}{\v@leur}\v@leur=\Sprim@\v@leur\v@leur=\Sprim@\v@leur%
    \figptBezierTD-1::#3[#4,#5,#6,#7]\edef\Sprim@{\repdecn@mb{\v@leur}}%
    \figpttraTD#1:{#2}=-1/\Sprim@,-5/\resetc@ntr@l\et@tfigptcurvcenterTD}\ignorespaces\fi}
\ctr@ld@f\def\c@lcurvradDD#1[#2,#3,#4,#5]{{\figvectDBezierDD -5:1,#1[#2,#3,#4,#5]%
    \figvectDBezierDD -6:2,#1[#2,#3,#4,#5]\vecunit@DD{-5}{-5}%
    \edef\Sprim@{\repdecn@mb{\result@t}}\figvectNVDD-5[-5]\c@lproscalDD\v@leur[-6,-5]%
    \invers@{\v@leur}{\v@leur}\v@leur=\Sprim@\v@leur\v@leur=\Sprim@\v@leur%
    \global\result@t=\v@leur}}
\ctr@ln@m\figptell
\ctr@ld@f\def\figptellDD#1:#2:#3;#4,#5(#6,#7){\ifGR@cri{\s@uvc@ntr@l\et@tfigptell%
    \c@lptellDD#1::#3;#4,#5(#6)\figptrotDD#1:{#2}=#1/#3,#7/%
    \resetc@ntr@l\et@tfigptell}\ignorespaces\fi}
\ctr@ld@f\def\c@lptellDD#1:#2:#3;#4,#5(#6){\c@ssin{\C@}{\S@}{#6}\v@lmin=\C@ pt\v@lmax=\S@ pt%
    \v@lmin=#4\v@lmin\v@lmax=#5\v@lmax%
    \edef\Xc@mp{\repdecn@mb{\v@lmin}}\edef\Yc@mp{\repdecn@mb{\v@lmax}}%
    \setc@ntr@l{2}\figvectC-1(\Xc@mp,\Yc@mp)\figpttraDD#1:{#2}=#3/1,-1/}
\ctr@ld@f\def\figptellP#1:#2:#3,#4,#5(#6){\ifGR@cri{\s@uvc@ntr@l\et@tfigptellP%
    \setc@ntr@l{2}\figvectP-1[#3,#4]\figvectP-2[#3,#5]%
    \v@leur=#6pt\c@lptellP{#3}{-1}{-2}\figptcopy#1:{#2}/-3/%
    \resetc@ntr@l\et@tfigptellP}\ignorespaces\fi}
\ctr@ln@m\@ngle
\ctr@ld@f\def\c@lptellP#1#2#3{\edef\@ngle{\repdecn@mb\v@leur}\c@ssin{\C@}{\S@}{\@ngle}%
    \figpttra-3:=#1/\C@,#2/\figpttra-3:=-3/\S@,#3/}
\ctr@ln@m\figptendnormal
\ctr@ld@f\def\figptendnormalDD#1:#2:#3,#4[#5,#6]{\ifGR@cri{\s@uvc@ntr@l\et@tfigptendnormal%
    \Figg@tXYa{#5}\Figg@tXY{#6}%
    \advance\v@lX-\v@lXa\advance\v@lY-\v@lYa%
    \setc@ntr@l{2}\Figv@ctCreg-1(\v@lX,\v@lY)\vecunit@{-1}{-1}\Figg@tXY{-1}%
    \delt@=#3\unit@\maxim@m{\delt@}{\delt@}{-\delt@}\edef\l@ngueur{\repdecn@mb{\delt@}}%
    \v@lX=\l@ngueur\v@lX\v@lY=\l@ngueur\v@lY%
    \delt@=\p@\advance\delt@-#4pt\edef\l@ngueur{\repdecn@mb{\delt@}}%
    \figptbaryR-1:[#5,#6;#4,\l@ngueur]\Figg@tXYa{-1}%
    \advance\v@lXa\v@lY\advance\v@lYa-\v@lX%
    \setc@ntr@l{1}\Figp@intregDD#1:{#2}(\v@lXa,\v@lYa)\resetc@ntr@l\et@tfigptendnormal}%
    \ignorespaces\fi}
\ctr@ld@f\def\figptexcenter#1:#2[#3,#4,#5]{\ifGR@cri{\let@xte={-}
    \Figptexinsc@nter#1:#2[#3,#4,#5]}\ignorespaces\fi}
\ctr@ld@f\def\figptincenter#1:#2[#3,#4,#5]{\ifGR@cri{\let@xte={}
    \Figptexinsc@nter#1:#2[#3,#4,#5]}\ignorespaces\fi}
\ctr@ld@f
\ctr@ld@f\def\Figptexinsc@nter#1:#2[#3,#4,#5]{%
    \figgetdist\LA@[#4,#5]\figgetdist\LB@[#3,#5]\figgetdist\LC@[#3,#4]%
    \figptbaryR#1:{#2}[#3,#4,#5;\the\let@xte\LA@,\LB@,\LC@]}
\ctr@ln@m\figptinterlineplane
\ctr@ld@f\def\figptinterlineplaneDD{\un@v@ilable{figptinterlineplane}}
\ctr@ld@f\def\figptinterlineplaneTD#1:#2[#3,#4;#5,#6]{\ifGR@cri{\s@uvc@ntr@l\et@tfigptinterlineplane%
    \setc@ntr@l{2}\figvectPTD-1[#3,#5]\vecunit@TD{-2}{#6}%
    \r@pPSTD\v@leur[-2,-1,#4]\edef\v@lcoef{\repdecn@mb{\v@leur}}%
    \figpttraTD#1:{#2}=#3/\v@lcoef,#4/\resetc@ntr@l\et@tfigptinterlineplane}\ignorespaces\fi}
\ctr@ln@m\figptorthocenter
\ctr@ld@f\def\figptorthocenterDD#1:#2[#3,#4,#5]{\ifGR@cri{\s@uvc@ntr@l\et@tfigptorthocenterDD%
    \setc@ntr@l{2}\figvectNDD-3[#3,#4]\figvectNDD-4[#4,#5]%
    \resetc@ntr@l{2}\inters@cDD#1:{#2}[#5,-3;#3,-4]%
    \resetc@ntr@l\et@tfigptorthocenterDD}\ignorespaces\fi}
\ctr@ld@f\def\figptorthocenterTD#1:#2[#3,#4,#5]{\ifGR@cri{\s@uvc@ntr@l\et@tfigptorthocenterTD%
    \setc@ntr@l{2}\figvectNTD-1[#3,#4,#5]%
    \figvectPTD-2[#3,#4]\figvectNVTD-3[-1,-2]%
    \figvectPTD-2[#4,#5]\figvectNVTD-4[-1,-2]%
    \resetc@ntr@l{2}\inters@cTD#1:{#2}[#5,-3;#3,-4]%
    \resetc@ntr@l\et@tfigptorthocenterTD}\ignorespaces\fi}
\ctr@ln@m\figptorthoprojline
\ctr@ld@f\def\figptorthoprojlineDD#1:#2=#3/#4,#5/{\ifGR@cri{\s@uvc@ntr@l\et@tfigptorthoprojlineDD%
    \setc@ntr@l{2}\figvectPDD-3[#4,#5]\figvectNVDD-4[-3]\resetc@ntr@l{2}%
    \inters@cDD#1:{#2}[#3,-4;#4,-3]\resetc@ntr@l\et@tfigptorthoprojlineDD}\ignorespaces\fi}
\ctr@ld@f\def\figptorthoprojlineTD#1:#2=#3/#4,#5/{\ifGR@cri{\s@uvc@ntr@l\et@tfigptorthoprojlineTD%
    \setc@ntr@l{2}\figvectPTD-1[#4,#3]\figvectPTD-2[#4,#5]\vecunit@TD{-2}{-2}%
    \c@lproscalTD\v@leur[-1,-2]\edef\v@lcoef{\repdecn@mb{\v@leur}}%
    \figpttraTD#1:{#2}=#4/\v@lcoef,-2/\resetc@ntr@l\et@tfigptorthoprojlineTD}\ignorespaces\fi}
\ctr@ln@m\figptorthoprojplane
\ctr@ld@f\def\figptorthoprojplaneDD{\un@v@ilable{figptorthoprojplane}}
\ctr@ld@f\def\figptorthoprojplaneTD#1:#2=#3/#4,#5/{\ifGR@cri{\s@uvc@ntr@l\et@tfigptorthoprojplane%
    \setc@ntr@l{2}\figvectPTD-1[#3,#4]\vecunit@TD{-2}{#5}%
    \c@lproscalTD\v@leur[-1,-2]\edef\v@lcoef{\repdecn@mb{\v@leur}}%
    \figpttraTD#1:{#2}=#3/\v@lcoef,-2/\resetc@ntr@l\et@tfigptorthoprojplane}\ignorespaces\fi}
\ctr@ld@f\def\figpthom#1:#2=#3/#4,#5/{\ifGR@cri{\s@uvc@ntr@l\et@tfigpthom%
    \setc@ntr@l{2}\figvectP-1[#4,#3]\figpttra#1:{#2}=#4/#5,-1/%
    \resetc@ntr@l\et@tfigpthom}\ignorespaces\fi}
\ctr@ld@f\def\figptinv#1:#2=#3/#4,#5/{\ifGR@cri{\s@uvc@ntr@l\et@tfigptinv%
    \setc@ntr@l{2}\figvectP-1[#4,#3]\Figg@tXY{-1}%
    \getredf@ctB\f@ctech\n@rmeucC{\delt@}{-1}%
    \delt@=\ptT@unit@\delt@\delt@=\ptT@unit@\delt@%
    \invers@{\delt@}{\delt@}\multiply\f@ctech\f@ctech\divide\delt@\f@ctech%
    \delt@=#5\delt@\edef\v@lcoef{\repdecn@mb{\delt@}}\figpttra#1:{#2}=#4/\v@lcoef,-1/%
    \resetc@ntr@l\et@tfigptinv}\ignorespaces\fi}
\ctr@ln@m\figptrot
\ctr@ld@f\def\figptrotDD#1:#2=#3/#4,#5/{\ifGR@cri{\s@uvc@ntr@l\et@tfigptrotDD%
    \c@ssin{\C@}{\S@}{#5}\setc@ntr@l{2}\figvectPDD-1[#4,#3]\Figg@tXY{-1}%
    \v@lXa=\C@\v@lX\advance\v@lXa-\S@\v@lY%
    \v@lYa=\S@\v@lX\advance\v@lYa\C@\v@lY%
    \Figv@ctCreg-1(\v@lXa,\v@lYa)\figpttraDD#1:{#2}=#4/1,-1/%
    \resetc@ntr@l\et@tfigptrotDD}\ignorespaces\fi}
\ctr@ld@f\def\figptrotTD#1:#2=#3/#4,#5,#6/{\ifGR@cri{\s@uvc@ntr@l\et@tfigptrotTD%
    \c@ssin{\C@}{\S@}{#5}%
    \setc@ntr@l{2}\figptorthoprojplaneTD-3:=#4/#3,#6/\figvectPTD-2[-3,#3]%
    \n@rmeucTD\v@leur{-2}\ifdim\v@leur<\Cepsil@n\Figg@tXYa{#3}\else%
    \edef\v@lcoef{\repdecn@mb{\v@leur}}\figvectNVTD-1[#6,-2]%
    \Figg@tXYa{-1}\v@lXa=\v@lcoef\v@lXa\v@lYa=\v@lcoef\v@lYa\v@lZa=\v@lcoef\v@lZa%
    \v@lXa=\S@\v@lXa\v@lYa=\S@\v@lYa\v@lZa=\S@\v@lZa\Figg@tXY{-2}%
    \advance\v@lXa\C@\v@lX\advance\v@lYa\C@\v@lY\advance\v@lZa\C@\v@lZ%
    \Figg@tXY{-3}\advance\v@lXa\v@lX\advance\v@lYa\v@lY\advance\v@lZa\v@lZ\fi%
    \Figp@intregTD#1:{#2}(\v@lXa,\v@lYa,\v@lZa)\resetc@ntr@l\et@tfigptrotTD}\ignorespaces\fi}
\ctr@ln@m\figptsym
\ctr@ld@f\def\figptsymDD#1:#2=#3/#4,#5/{\ifGR@cri{\s@uvc@ntr@l\et@tfigptsymDD%
    \resetc@ntr@l{2}\figptorthoprojlineDD-5:=#3/#4,#5/\figvectPDD-2[#3,-5]%
    \figpttraDD#1:{#2}=#3/2,-2/\resetc@ntr@l\et@tfigptsymDD}\ignorespaces\fi}
\ctr@ld@f\def\figptsymTD#1:#2=#3/#4,#5/{\ifGR@cri{\s@uvc@ntr@l\et@tfigptsymTD%
    \resetc@ntr@l{2}\figptorthoprojplaneTD-3:=#3/#4,#5/\figvectPTD-2[#3,-3]%
    \figpttraTD#1:{#2}=#3/2,-2/\resetc@ntr@l\et@tfigptsymTD}\ignorespaces\fi}
\ctr@ln@m\figpttra
\ctr@ld@f\def\figpttraDD#1:#2=#3/#4,#5/{\ifGR@cri{\Figg@tXYa{#5}\v@lXa=#4\v@lXa\v@lYa=#4\v@lYa%
    \Figg@tXY{#3}\advance\v@lX\v@lXa\advance\v@lY\v@lYa%
    \Figp@intregDD#1:{#2}(\v@lX,\v@lY)}\ignorespaces\fi}
\ctr@ld@f\def\figpttraTD#1:#2=#3/#4,#5/{\ifGR@cri{\Figg@tXYa{#5}\v@lXa=#4\v@lXa\v@lYa=#4\v@lYa%
    \v@lZa=#4\v@lZa\Figg@tXY{#3}\advance\v@lX\v@lXa\advance\v@lY\v@lYa%
    \advance\v@lZ\v@lZa\Figp@intregTD#1:{#2}(\v@lX,\v@lY,\v@lZ)}\ignorespaces\fi}
\ctr@ln@m\figpttraC
\ctr@ld@f\def\figpttraCDD#1:#2=#3/#4,#5/{\ifGR@cri{\v@lXa=#4\unit@\v@lYa=#5\unit@%
    \Figg@tXY{#3}\advance\v@lX\v@lXa\advance\v@lY\v@lYa%
    \Figp@intregDD#1:{#2}(\v@lX,\v@lY)}\ignorespaces\fi}
\ctr@ld@f\def\figpttraCTD#1:#2=#3/#4,#5,#6/{\ifGR@cri{\v@lXa=#4\unit@\v@lYa=#5\unit@\v@lZa=#6\unit@%
    \Figg@tXY{#3}\advance\v@lX\v@lXa\advance\v@lY\v@lYa\advance\v@lZ\v@lZa%
    \Figp@intregTD#1:{#2}(\v@lX,\v@lY,\v@lZ)}\ignorespaces\fi}
\ctr@ld@f\def\figptsaxes#1:#2(#3){\ifGR@cri{\an@lys@xes#3,:\ifx\t@xt@\empty%
    \ifTr@isDim\Figpts@xes#1:#2(0,#3,0,#3,0,#3)\else\Figpts@xes#1:#2(0,#3,0,#3)\fi%
    \else\Figpts@xes#1:#2(#3)\fi}\ignorespaces\fi}
\ctr@ln@m\Figpts@xes
\ctr@ld@f\def\Figpts@xesDD#1:#2(#3,#4,#5,#6){%
    \s@mme=#1\figpttraC\the\s@mme:$x$=#2/#4,0/%
    \advance\s@mme\@ne\figpttraC\the\s@mme:$y$=#2/0,#6/}
\ctr@ld@f\def\Figpts@xesTD#1:#2(#3,#4,#5,#6,#7,#8){%
    \s@mme=#1\figpttraC\the\s@mme:$x$=#2/#4,0,0/%
    \advance\s@mme\@ne\figpttraC\the\s@mme:$y$=#2/0,#6,0/%
    \advance\s@mme\@ne\figpttraC\the\s@mme:$z$=#2/0,0,#8/}
\ctr@ld@f\def\figptsmap#1=#2/#3/#4/{\ifGR@cri{\s@uvc@ntr@l\et@tfigptsmap%
    \setc@ntr@l{2}\def\list@num{#2}\s@mme=#1%
    \@ecfor\p@int:=\list@num\do{\figvectP-1[#3,\p@int]\Figg@tXY{-1}%
    \pr@dMatV/#4/\figpttra\the\s@mme:=#3/1,-1/\advance\s@mme\@ne}%
    \resetc@ntr@l\et@tfigptsmap}\ignorespaces\fi}
\ctr@ln@m\figptscontrol
\ctr@ld@f\def\figptscontrolDD#1[#2,#3,#4,#5]{\ifGR@cri{\s@uvc@ntr@l\et@tfigptscontrolDD\setc@ntr@l{2}%
    \v@lX=\z@\v@lY=\z@\Figtr@nptDD{-5}{#2}\Figtr@nptDD{2}{#5}%
    \divide\v@lX\@vi\divide\v@lY\@vi%
    \Figtr@nptDD{3}{#3}\Figtr@nptDD{-1.5}{#4}\Figp@intregDD-1:(\v@lX,\v@lY)%
    \v@lX=\z@\v@lY=\z@\Figtr@nptDD{2}{#2}\Figtr@nptDD{-5}{#5}%
    \divide\v@lX\@vi\divide\v@lY\@vi\Figtr@nptDD{-1.5}{#3}\Figtr@nptDD{3}{#4}%
    \s@mme=#1\advance\s@mme\@ne\Figp@intregDD\the\s@mme:(\v@lX,\v@lY)%
    \figptcopyDD#1:/-1/\resetc@ntr@l\et@tfigptscontrolDD}\ignorespaces\fi}
\ctr@ld@f\def\figptscontrolTD#1[#2,#3,#4,#5]{\ifGR@cri{\s@uvc@ntr@l\et@tfigptscontrolTD\setc@ntr@l{2}%
    \v@lX=\z@\v@lY=\z@\v@lZ=\z@\Figtr@nptTD{-5}{#2}\Figtr@nptTD{2}{#5}%
    \divide\v@lX\@vi\divide\v@lY\@vi\divide\v@lZ\@vi%
    \Figtr@nptTD{3}{#3}\Figtr@nptTD{-1.5}{#4}\Figp@intregTD-1:(\v@lX,\v@lY,\v@lZ)%
    \v@lX=\z@\v@lY=\z@\v@lZ=\z@\Figtr@nptTD{2}{#2}\Figtr@nptTD{-5}{#5}%
    \divide\v@lX\@vi\divide\v@lY\@vi\divide\v@lZ\@vi\Figtr@nptTD{-1.5}{#3}\Figtr@nptTD{3}{#4}%
    \s@mme=#1\advance\s@mme\@ne\Figp@intregTD\the\s@mme:(\v@lX,\v@lY,\v@lZ)%
    \figptcopyTD#1:/-1/\resetc@ntr@l\et@tfigptscontrolTD}\ignorespaces\fi}
\ctr@ld@f\def\Figtr@nptDD#1#2{\Figg@tXYa{#2}\v@lXa=#1\v@lXa\v@lYa=#1\v@lYa%
    \advance\v@lX\v@lXa\advance\v@lY\v@lYa}
\ctr@ld@f\def\Figtr@nptTD#1#2{\Figg@tXYa{#2}\v@lXa=#1\v@lXa\v@lYa=#1\v@lYa\v@lZa=#1\v@lZa%
    \advance\v@lX\v@lXa\advance\v@lY\v@lYa\advance\v@lZ\v@lZa}
\ctr@ld@f\def\figptscontrolcurve#1,#2[#3]{\ifGR@cri{\s@uvc@ntr@l\et@tfigptscontrolcurve%
    \def\list@num{#3}\extrairelepremi@r\Ak@\de\list@num%
    \extrairelepremi@r\Ai@\de\list@num\extrairelepremi@r\Aj@\de\list@num%
    \s@mme=#1\figptcopy\the\s@mme:/\Ai@/%
    \setc@ntr@l{2}\figvectP -1[\Ak@,\Aj@]%
    \@ecfor\Ak@:=\list@num\do{\advance\s@mme\@ne\figpttra\the\s@mme:=\Ai@/\curv@roundness,-1/%
       \figvectP -1[\Ai@,\Ak@]\advance\s@mme\@ne\figpttra\the\s@mme:=\Aj@/-\curv@roundness,-1/%
       \advance\s@mme\@ne\figptcopy\the\s@mme:/\Aj@/%
       \edef\Ai@{\Aj@}\edef\Aj@{\Ak@}}\advance\s@mme-#1\divide\s@mme\thr@@%
       \xdef#2{\the\s@mme}%
    \resetc@ntr@l\et@tfigptscontrolcurve}\ignorespaces\fi}
\ctr@ln@m\figptsintercirc
\ctr@ld@f\def\figptsintercircDD#1[#2,#3;#4,#5]{\ifGR@cri{\s@uvc@ntr@l\et@tfigptsintercircDD%
    \setc@ntr@l{2}\let\c@lNVintc=\c@lNVintcDD\Figptsintercirc@#1[#2,#3;#4,#5]%
    \resetc@ntr@l\et@tfigptsintercircDD}\ignorespaces\fi}
\ctr@ld@f\def\figptsintercircTD#1[#2,#3;#4,#5;#6]{\ifGR@cri{\s@uvc@ntr@l\et@tfigptsintercircTD%
    \setc@ntr@l{2}\let\c@lNVintc=\c@lNVintcTD\vecunitC@TD[#2,#6]%
    \Figv@ctCreg-3(\v@lX,\v@lY,\v@lZ)\Figptsintercirc@#1[#2,#3;#4,#5]%
    \resetc@ntr@l\et@tfigptsintercircTD}\ignorespaces\fi}
\ctr@ld@f\def\Figptsintercirc@#1[#2,#3;#4,#5]{\figvectP-1[#2,#4]%
    \vecunit@{-1}{-1}\delt@=\result@t\f@ctech=\result@tent%
    \s@mme=#1\advance\s@mme\@ne\figptcopy#1:/#2/\figptcopy\the\s@mme:/#4/%
    \ifdim\delt@=\z@\else%
    \v@lmin=#3\unit@\v@lmax=#5\unit@\v@leur=\v@lmin\advance\v@leur\v@lmax%
    \ifdim\v@leur>\delt@%
    \v@leur=\v@lmin\advance\v@leur-\v@lmax\maxim@m{\v@leur}{\v@leur}{-\v@leur}%
    \ifdim\v@leur<\delt@%
    \divide\v@lmin\f@ctech\divide\v@lmax\f@ctech\divide\delt@\f@ctech%
    \v@lmin=\repdecn@mb{\v@lmin}\v@lmin\v@lmax=\repdecn@mb{\v@lmax}\v@lmax%
    \invers@{\v@leur}{\delt@}\advance\v@lmax-\v@lmin%
    \v@lmax=-\repdecn@mb{\v@leur}\v@lmax\advance\delt@\v@lmax\delt@=.5\delt@%
    \v@lmax=\delt@\multiply\v@lmax\f@ctech%
    \edef\t@ille{\repdecn@mb{\v@lmax}}\figpttra-2:=#2/\t@ille,-1/%
    \delt@=\repdecn@mb{\delt@}\delt@\advance\v@lmin-\delt@%
    \sqrt@{\v@leur}{\v@lmin}\multiply\v@leur\f@ctech\edef\t@ille{\repdecn@mb{\v@leur}}%
    \c@lNVintc\figpttra#1:=-2/-\t@ille,-1/\figpttra\the\s@mme:=-2/\t@ille,-1/\fi\fi\fi}
\ctr@ld@f\def\c@lNVintcDD{\Figg@tXY{-1}\Figv@ctCreg-1(-\v@lY,\v@lX)} 
\ctr@ld@f\def\c@lNVintcTD{{\Figg@tXY{-3}\v@lmin=\v@lX\v@lmax=\v@lY\v@leur=\v@lZ%
    \Figg@tXY{-1}\c@lprovec{-3}\vecunit@{-3}{-3}
    \Figg@tXY{-1}\v@lmin=\v@lX\v@lmax=\v@lY%
    \v@leur=\v@lZ\Figg@tXY{-3}\c@lprovec{-1}}} 
\ctr@ln@m\figptsinterlinell
\ctr@ld@f\def\figptsinterlinellDD#1[#2,#3,#4,#5;#6,#7]{\ifGR@cri{\s@uvc@ntr@l\et@tfigptsinterlinellDD%
    \figptcopy#1:/#6/\s@mme=#1\advance\s@mme\@ne\figptcopy\the\s@mme:/#7/%
    \v@lmin=#3\unit@\v@lmax=#4\unit@
    \setc@ntr@l{2}\figptbaryDD-4:[#6,#7;1,1]\figptsrotDD-3=-4,#7/#2,-#5/
    \Figg@tXY{-3}\Figg@tXYa{#2}\advance\v@lX-\v@lXa\advance\v@lY-\v@lYa
    \figvectP-1[-3,-2]\Figg@tXYa{-1}\figvectP-3[-4,#7]\Figptsint@rLE{#1}
    \resetc@ntr@l\et@tfigptsinterlinellDD}\ignorespaces\fi}
\ctr@ld@f\def\figptsinterlinellP#1[#2,#3,#4;#5,#6]{\ifGR@cri{\s@uvc@ntr@l\et@tfigptsinterlinellP%
    \figptcopy#1:/#5/\s@mme=#1\advance\s@mme\@ne\figptcopy\the\s@mme:/#6/\setc@ntr@l{2}%
    \figvectP-1[#2,#3]\vecunit@{-1}{-1}\v@lmin=\result@t
    \figvectP-2[#2,#4]\vecunit@{-2}{-2}\v@lmax=\result@t
    \figptbary-4:[#5,#6;1,1]
    \figvectP-3[#2,-4]\c@lproscal\v@lX[-3,-1]\c@lproscal\v@lY[-3,-2]
    \figvectP-3[-4,#6]\c@lproscal\v@lXa[-3,-1]\c@lproscal\v@lYa[-3,-2]
    \Figptsint@rLE{#1}\resetc@ntr@l\et@tfigptsinterlinellP}\ignorespaces\fi}
\ctr@ld@f\def\Figptsint@rLE#1{%
    \getredf@ctDD\f@ctech(\v@lmin,\v@lmax)%
    \getredf@ctDD\p@rtent(\v@lX,\v@lY)\ifnum\p@rtent>\f@ctech\f@ctech=\p@rtent\fi%
    \getredf@ctDD\p@rtent(\v@lXa,\v@lYa)\ifnum\p@rtent>\f@ctech\f@ctech=\p@rtent\fi%
    \divide\v@lmin\f@ctech\divide\v@lmax\f@ctech\divide\v@lX\f@ctech\divide\v@lY\f@ctech%
    \divide\v@lXa\f@ctech\divide\v@lYa\f@ctech%
    \c@rre=\repdecn@mb\v@lXa\v@lmax\mili@u=\repdecn@mb\v@lYa\v@lmin%
    \getredf@ctDD\f@ctech(\c@rre,\mili@u)%
    \c@rre=\repdecn@mb\v@lX\v@lmax\mili@u=\repdecn@mb\v@lY\v@lmin%
    \getredf@ctDD\p@rtent(\c@rre,\mili@u)\ifnum\p@rtent>\f@ctech\f@ctech=\p@rtent\fi%
    \divide\v@lmin\f@ctech\divide\v@lmax\f@ctech\divide\v@lX\f@ctech\divide\v@lY\f@ctech%
    \divide\v@lXa\f@ctech\divide\v@lYa\f@ctech%
    \v@lmin=\repdecn@mb{\v@lmin}\v@lmin\v@lmax=\repdecn@mb{\v@lmax}\v@lmax%
    \edef\G@xde{\repdecn@mb\v@lmin}\edef\P@xde{\repdecn@mb\v@lmax}%
    \c@rre=-\v@lmax\v@leur=\repdecn@mb\v@lY\v@lY\advance\c@rre\v@leur\c@rre=\G@xde\c@rre%
    \v@leur=\repdecn@mb\v@lX\v@lX\v@leur=\P@xde\v@leur\advance\c@rre\v@leur
    \v@lmin=\repdecn@mb\v@lYa\v@lmin\v@lmax=\repdecn@mb\v@lXa\v@lmax%
    \mili@u=\repdecn@mb\v@lX\v@lmax\advance\mili@u\repdecn@mb\v@lY\v@lmin
    \v@lmax=\repdecn@mb\v@lXa\v@lmax\advance\v@lmax\repdecn@mb\v@lYa\v@lmin
    \ifdim\v@lmax>\epsil@n%
    \maxim@m{\v@leur}{\c@rre}{-\c@rre}\maxim@m{\v@lmin}{\mili@u}{-\mili@u}%
    \maxim@m{\v@leur}{\v@leur}{\v@lmin}\maxim@m{\v@lmin}{\v@lmax}{-\v@lmax}%
    \maxim@m{\v@leur}{\v@leur}{\v@lmin}\p@rtentiere{\p@rtent}{\v@leur}\advance\p@rtent\@ne%
    \divide\c@rre\p@rtent\divide\mili@u\p@rtent\divide\v@lmax\p@rtent%
    \delt@=\repdecn@mb{\mili@u}\mili@u\v@leur=\repdecn@mb{\v@lmax}\c@rre%
    \advance\delt@-\v@leur\ifdim\delt@<\z@\else\sqrt@\delt@\delt@%
    \invers@\v@lmax\v@lmax\edef\Uns@rAp{\repdecn@mb\v@lmax}%
    \v@leur=-\mili@u\advance\v@leur-\delt@\v@leur=\Uns@rAp\v@leur%
    \edef\t@ille{\repdecn@mb\v@leur}\figpttra#1:=-4/\t@ille,-3/\s@mme=#1\advance\s@mme\@ne%
    \v@leur=-\mili@u\advance\v@leur\delt@\v@leur=\Uns@rAp\v@leur%
    \edef\t@ille{\repdecn@mb\v@leur}\figpttra\the\s@mme:=-4/\t@ille,-3/\fi\fi}
\ctr@ln@m\figptsorthoprojline
\ctr@ld@f\def\figptsorthoprojlineDD#1=#2/#3,#4/{\ifGR@cri{\s@uvc@ntr@l\et@tfigptsorthoprojlineDD%
    \setc@ntr@l{2}\figvectPDD-3[#3,#4]\figvectNVDD-4[-3]\resetc@ntr@l{2}%
    \def\list@num{#2}\s@mme=#1\@ecfor\p@int:=\list@num\do{%
    \inters@cDD\the\s@mme:[\p@int,-4;#3,-3]\advance\s@mme\@ne}%
    \resetc@ntr@l\et@tfigptsorthoprojlineDD}\ignorespaces\fi}
\ctr@ld@f\def\figptsorthoprojlineTD#1=#2/#3,#4/{\ifGR@cri{\s@uvc@ntr@l\et@tfigptsorthoprojlineTD%
    \setc@ntr@l{2}\figvectPTD-2[#3,#4]\vecunit@TD{-2}{-2}%
    \def\list@num{#2}\s@mme=#1\@ecfor\p@int:=\list@num\do{%
    \figvectPTD-1[#3,\p@int]\c@lproscalTD\v@leur[-1,-2]%
    \edef\v@lcoef{\repdecn@mb{\v@leur}}\figpttraTD\the\s@mme:=#3/\v@lcoef,-2/%
    \advance\s@mme\@ne}\resetc@ntr@l\et@tfigptsorthoprojlineTD}\ignorespaces\fi}
\ctr@ln@m\figptsorthoprojplane
\ctr@ld@f\def\figptsorthoprojplaneDD{\un@v@ilable{figptsorthoprojplane}}
\ctr@ld@f\def\figptsorthoprojplaneTD#1=#2/#3,#4/{\ifGR@cri{\s@uvc@ntr@l\et@tfigptsorthoprojplane%
    \setc@ntr@l{2}\vecunit@TD{-2}{#4}%
    \def\list@num{#2}\s@mme=#1\@ecfor\p@int:=\list@num\do{\figvectPTD-1[\p@int,#3]%
    \c@lproscalTD\v@leur[-1,-2]\edef\v@lcoef{\repdecn@mb{\v@leur}}%
    \figpttraTD\the\s@mme:=\p@int/\v@lcoef,-2/\advance\s@mme\@ne}%
    \resetc@ntr@l\et@tfigptsorthoprojplane}\ignorespaces\fi}
\ctr@ld@f\def\figptshom#1=#2/#3,#4/{\ifGR@cri{\s@uvc@ntr@l\et@tfigptshom%
    \setc@ntr@l{2}\def\list@num{#2}\s@mme=#1%
    \@ecfor\p@int:=\list@num\do{\figvectP-1[#3,\p@int]%
    \figpttra\the\s@mme:=#3/#4,-1/\advance\s@mme\@ne}%
    \resetc@ntr@l\et@tfigptshom}\ignorespaces\fi}
\ctr@ld@f\def\figptsinv#1=#2/#3,#4/{\ifGR@cri{\s@uvc@ntr@l\et@tfigptsinv%
    \setc@ntr@l{2}\def\list@num{#2}\s@mme=#1%
    \@ecfor\p@int:=\list@num\do{\figvectP-1[#3,\p@int]\Figg@tXY{-1}%
    \getredf@ctB\f@ctech\n@rmeucC{\delt@}{-1}%
    \delt@=\ptT@unit@\delt@\delt@=\ptT@unit@\delt@%
    \invers@{\delt@}{\delt@}\multiply\f@ctech\f@ctech\divide\delt@\f@ctech%
    \delt@=#4\delt@\edef\v@lcoef{\repdecn@mb{\delt@}}\figpttra\the\s@mme:=#3/\v@lcoef,-1/%
    \advance\s@mme\@ne}\resetc@ntr@l\et@tfigptsinv}\ignorespaces\fi}
\ctr@ln@m\figptsrot
\ctr@ld@f\def\figptsrotDD#1=#2/#3,#4/{\ifGR@cri{\s@uvc@ntr@l\et@tfigptsrotDD%
    \c@ssin{\C@}{\S@}{#4}\setc@ntr@l{2}\def\list@num{#2}\s@mme=#1%
    \@ecfor\p@int:=\list@num\do{\figvectPDD-1[#3,\p@int]\Figg@tXY{-1}%
    \v@lXa=\C@\v@lX\advance\v@lXa-\S@\v@lY%
    \v@lYa=\S@\v@lX\advance\v@lYa\C@\v@lY%
    \Figv@ctCreg-1(\v@lXa,\v@lYa)\figpttraDD\the\s@mme:=#3/1,-1/\advance\s@mme\@ne}%
    \resetc@ntr@l\et@tfigptsrotDD}\ignorespaces\fi}
\ctr@ld@f\def\figptsrotTD#1=#2/#3,#4,#5/{\ifGR@cri{\s@uvc@ntr@l\et@tfigptsrotTD%
    \c@ssin{\C@}{\S@}{#4}%
    \setc@ntr@l{2}\def\list@num{#2}\s@mme=#1%
    \@ecfor\p@int:=\list@num\do{\figptorthoprojplaneTD-3:=#3/\p@int,#5/%
    \figvectPTD-2[-3,\p@int]%
    \figvectNVTD-1[#5,-2]\n@rmeucTD\v@leur{-2}\edef\v@lcoef{\repdecn@mb{\v@leur}}%
    \Figg@tXYa{-1}\v@lXa=\v@lcoef\v@lXa\v@lYa=\v@lcoef\v@lYa\v@lZa=\v@lcoef\v@lZa%
    \v@lXa=\S@\v@lXa\v@lYa=\S@\v@lYa\v@lZa=\S@\v@lZa\Figg@tXY{-2}%
    \advance\v@lXa\C@\v@lX\advance\v@lYa\C@\v@lY\advance\v@lZa\C@\v@lZ%
    \Figg@tXY{-3}\advance\v@lXa\v@lX\advance\v@lYa\v@lY\advance\v@lZa\v@lZ%
    \Figp@intregTD\the\s@mme:(\v@lXa,\v@lYa,\v@lZa)\advance\s@mme\@ne}%
    \resetc@ntr@l\et@tfigptsrotTD}\ignorespaces\fi}
\ctr@ln@m\figptssym
\ctr@ld@f\def\figptssymDD#1=#2/#3,#4/{\ifGR@cri{\s@uvc@ntr@l\et@tfigptssymDD%
    \setc@ntr@l{2}\figvectPDD-3[#3,#4]\Figg@tXY{-3}\Figv@ctCreg-4(-\v@lY,\v@lX)%
    \resetc@ntr@l{2}\def\list@num{#2}\s@mme=#1%
    \@ecfor\p@int:=\list@num\do{\inters@cDD-5:[#3,-3;\p@int,-4]\figvectPDD-2[\p@int,-5]%
    \figpttraDD\the\s@mme:=\p@int/2,-2/\advance\s@mme\@ne}%
    \resetc@ntr@l\et@tfigptssymDD}\ignorespaces\fi}
\ctr@ld@f\def\figptssymTD#1=#2/#3,#4/{\ifGR@cri{\s@uvc@ntr@l\et@tfigptssymTD%
    \setc@ntr@l{2}\vecunit@TD{-2}{#4}\def\list@num{#2}\s@mme=#1%
    \@ecfor\p@int:=\list@num\do{\figvectPTD-1[\p@int,#3]%
    \c@lproscalTD\v@leur[-1,-2]\v@leur=2\v@leur\edef\v@lcoef{\repdecn@mb{\v@leur}}%
    \figpttraTD\the\s@mme:=\p@int/\v@lcoef,-2/\advance\s@mme\@ne}%
    \resetc@ntr@l\et@tfigptssymTD}\ignorespaces\fi}
\ctr@ln@m\figptstra
\ctr@ld@f\def\figptstraDD#1=#2/#3,#4/{\ifGR@cri{\Figg@tXYa{#4}\v@lXa=#3\v@lXa\v@lYa=#3\v@lYa%
    \def\list@num{#2}\s@mme=#1\@ecfor\p@int:=\list@num\do{\Figg@tXY{\p@int}%
    \advance\v@lX\v@lXa\advance\v@lY\v@lYa%
    \Figp@intregDD\the\s@mme:(\v@lX,\v@lY)\advance\s@mme\@ne}}\ignorespaces\fi}
\ctr@ld@f\def\figptstraTD#1=#2/#3,#4/{\ifGR@cri{\Figg@tXYa{#4}\v@lXa=#3\v@lXa\v@lYa=#3\v@lYa%
    \v@lZa=#3\v@lZa\def\list@num{#2}\s@mme=#1\@ecfor\p@int:=\list@num\do{\Figg@tXY{\p@int}%
    \advance\v@lX\v@lXa\advance\v@lY\v@lYa\advance\v@lZ\v@lZa%
    \Figp@intregTD\the\s@mme:(\v@lX,\v@lY,\v@lZ)\advance\s@mme\@ne}}\ignorespaces\fi}
\ctr@ln@m\figptvisilimSL
\ctr@ld@f\def\figptvisilimSLDD{\un@v@ilable{figptvisilimSL}}
\ctr@ld@f\def\figptvisilimSLTD#1:#2[#3,#4;#5,#6]{\ifGR@cri{\s@uvc@ntr@l\et@tfigptvisilimSLTD%
    \setc@ntr@l{2}\figvectP-1[#3,#4]\n@rminf{\delt@}{-1}%
    \ifcase\CUR@proj\v@lX=\cxa@\p@\v@lY=-\p@\v@lZ=\cxb@\p@
    \Figv@ctCreg-2(\v@lX,\v@lY,\v@lZ)\figvectP-3[#5,#6]\figvectNV-1[-2,-3]%
    \or\figvectP-1[#5,#6]\vecunitCV@TD{-1}\v@lmin=\v@lX\v@lmax=\v@lY
    \v@leur=\v@lZ\v@lX=\cza@\p@\v@lY=\czb@\p@\v@lZ=\czc@\p@\c@lprovec{-1}%
    \or\c@ley@pt{-2}\figvectN-1[#5,#6,-2]\fi
    \edef\Ai@{#3}\edef\Aj@{#4}\figvectP-2[#5,\Ai@]\c@lproscal\v@leur[-1,-2]%
    \ifdim\v@leur>\z@\p@rtent=\@ne\else\p@rtent=\m@ne\fi%
    \figvectP-2[#5,\Aj@]\c@lproscal\v@leur[-1,-2]%
    \ifdim\p@rtent\v@leur>\z@\figptcopy#1:#2/#3/%
    \message{*** \BS@ figptvisilimSL: points are on the same side.}\else%
    \figptcopy-3:/#3/\figptcopy-4:/#4/%
    \loop\figptbary-5:[-3,-4;1,1]\figvectP-2[#5,-5]\c@lproscal\v@leur[-1,-2]%
    \ifdim\p@rtent\v@leur>\z@\figptcopy-3:/-5/\else\figptcopy-4:/-5/\fi%
    \divide\delt@\tw@\ifdim\delt@>\epsil@n\repeat%
    \figptbary#1:#2[-3,-4;1,1]\fi\resetc@ntr@l\et@tfigptvisilimSLTD}\ignorespaces\fi}
\ctr@ld@f\def\c@ley@pt#1{\t@stp@r\ifitis@K\v@lX=\cza@\p@\v@lY=\czb@\p@\v@lZ=\czc@\p@%
    \Figv@ctCreg-1(\v@lX,\v@lY,\v@lZ)\Figp@intreg-2:(\wd\Bt@rget,\ht\Bt@rget,\dp\Bt@rget)%
    \figpttra#1:=-2/-\disob@intern,-1/\else\end\fi}
\ctr@ld@f\def\t@stp@r{\itis@Ktrue\ifnewt@rgetpt\else\itis@Kfalse%
    \message{*** \BS@ figptvisilimXX: target point undefined.}\fi\ifnewdis@b\else%
    \itis@Kfalse\message{*** \BS@ figptvisilimXX: observation distance undefined.}\fi%
    \ifitis@K\else\message{*** This macro must be called after \BS@ figdrawbegin or after
    having set the missing parameter(s) with \BS@ figset proj()}\fi}
\ctr@ld@f\def\figscan#1(#2,#3){{\s@uvc@ntr@l\et@tfigscan\@psfgetbb{#1}\if@psfbbfound\else%
    \def\@psfllx{0}\def\@psflly{20}\def\@psfurx{540}\def\@psfury{640}\fi\figscan@{#2}{#3}%
    \resetc@ntr@l\et@tfigscan}\ignorespaces}
\ctr@ld@f\def\figscan@#1#2{%
    \unit@=\@ne bp\setc@ntr@l{2}\figsetmark{}%
    \def\minst@p{20pt}%
    \v@lX=\@psfllx\p@\v@lX=\Sc@leFact\v@lX\r@undint\v@lX\v@lX%
    \v@lY=\@psflly\p@\v@lY=\Sc@leFact\v@lY\ifdim\v@lY>\z@\r@undint\v@lY\v@lY\fi%
    \delt@=\@psfury\p@\delt@=\Sc@leFact\delt@%
    \advance\delt@-\v@lY\v@lXa=\@psfurx\p@\v@lXa=\Sc@leFact\v@lXa\v@leur=\minst@p%
    \edef\valv@lY{\repdecn@mb{\v@lY}}\edef\LgTr@it{\the\delt@}%
    \loop\ifdim\v@lX<\v@lXa\edef\valv@lX{\repdecn@mb{\v@lX}}%
    \figptDD -1:(\valv@lX,\valv@lY)\figwriten -1:\hbox{\vrule height\LgTr@it}(0)%
    \ifdim\v@leur<\minst@p\else\figsetmark{\raise-8bp\hbox{$\scriptscriptstyle\triangle$}}%
    \figwrites -1:\@ffichnb{0}{\valv@lX}(6)\v@leur=\z@\figsetmark{}\fi%
    \advance\v@leur#1pt\advance\v@lX#1pt\repeat%
    \def\minst@p{10pt}%
    \v@lX=\@psfllx\p@\v@lX=\Sc@leFact\v@lX\ifdim\v@lX>\z@\r@undint\v@lX\v@lX\fi%
    \v@lY=\@psflly\p@\v@lY=\Sc@leFact\v@lY\r@undint\v@lY\v@lY%
    \delt@=\@psfurx\p@\delt@=\Sc@leFact\delt@%
    \advance\delt@-\v@lX\v@lYa=\@psfury\p@\v@lYa=\Sc@leFact\v@lYa\v@leur=\minst@p%
    \edef\valv@lX{\repdecn@mb{\v@lX}}\edef\LgTr@it{\the\delt@}%
    \loop\ifdim\v@lY<\v@lYa\edef\valv@lY{\repdecn@mb{\v@lY}}%
    \figptDD -1:(\valv@lX,\valv@lY)\figwritee -1:\vbox{\hrule width\LgTr@it}(0)%
    \ifdim\v@leur<\minst@p\else\figsetmark{$\triangleright$\kern4bp}%
    \figwritew -1:\@ffichnb{0}{\valv@lY}(6)\v@leur=\z@\figsetmark{}\fi%
    \advance\v@leur#2pt\advance\v@lY#2pt\repeat}
\ctr@ld@f
\ctr@ld@f\def\figscan@E#1(#2,#3){{\s@uvc@ntr@l\et@tfigscan@E%
    \Figdisc@rdLTS{#1}{\t@xt@}\pdfximage{\t@xt@}%
    \setbox\Gb@x=\hbox{\pdfrefximage\pdflastximage}%
    \edef\@psfllx{0}\v@lY=-\dp\Gb@x\edef\@psflly{\repdecn@mb{\v@lY}}%
    \edef\@psfurx{\repdecn@mb{\wd\Gb@x}}%
    \v@lY=\dp\Gb@x\advance\v@lY\ht\Gb@x\edef\@psfury{\repdecn@mb{\v@lY}}%
    \figscan@{#2}{#3}\resetc@ntr@l\et@tfigscan@E}\ignorespaces}
\ctr@ld@f\def\figshowpts[#1,#2]{{\figsetmark{$\bullet$}\figsetptname{\bf ##1}%
    \p@rtent=#2\relax\ifnum\p@rtent<\z@\p@rtent=\z@\fi%
    \s@mme=#1\relax\ifnum\s@mme<\z@\s@mme=\z@\fi%
    \loop\ifnum\s@mme<\p@rtent\pt@rvect{\s@mme}%
    \ifitis@K\figwriten{\the\s@mme}:(4pt)\fi\advance\s@mme\@ne\repeat%
    \pt@rvect{\s@mme}\ifitis@K\figwriten{\the\s@mme}:(4pt)\fi}\ignorespaces}
\ctr@ld@f\def\pt@rvect#1{\set@bjc@de{#1}%
    \expandafter\expandafter\expandafter\inqpt@rvec\csname\objc@de\endcsname:}
\ctr@ld@f\def\inqpt@rvec#1#2:{\if#1\C@dCl@spt\itis@Ktrue\else\itis@Kfalse\fi}
\ctr@ld@f\def\figshowsettings{{%
    \immediate\write16{====================================================================}%
    \immediate\write16{ Current settings are (DDV means "with dynamic default value"):}%
    \immediate\write16{ --- GENERAL ---}%
    \immediate\write16{Scale factor and Unit = \unit@util\space (\the\unit@)
     \space -> \BS@ figinit{ScaleFactorUnit}}%
    \immediate\write16{Update mode = \ifGRupdatem@de yes\else no\fi
     \space-> \BS@ figset(update=yes/no) or \BS@ figsetdefault(update=yes/no)}%
    \immediate\write16{ --- WRITING ---}%
    \immediate\write16{Implicit point name = \ptn@me{i} \space-> \BS@ figset write(ptname={Name})}%
    \immediate\write16{Point marker = \the\c@nsymb \space -> \BS@ figset write(mark=Mark)}%
    \immediate\write16{Print rounded coordinates = \ifr@undcoord yes\else no\fi
     \space-> \BS@ figset write(roundcoord=yes/no)}%
    \immediate\write16{ --- GRAPHICAL (general) ---}%
    \immediate\write16{Color = \CUR@color \space-> \BS@ figset(color=ColorDefinition)}%
    \immediate\write16{Filling mode = \iffillm@de yes\else no\fi
     \space-> \BS@ figset(fillmode=yes/no)}%
    \immediate\write16{Line join = \CUR@join \space-> \BS@ figset(join=miter/round/bevel)}%
    \immediate\write16{Line style = \CUR@dash \space-> \BS@ figset(dash=Index/Pattern)}%
    \immediate\write16{Line width = \CUR@width
     \space-> \BS@ figset(width=real in PostScript units)}%
    \immediate\write16{ --- GRAPHICAL (specific) ---}%
    \immediate\write16{Altitude (all the following attributes are DDV):}%
    \immediate\write16{ Base line color =
     \ifx\DDV@blcolor\D@FTref general color\else\DDV@blcolor\fi
     \space-> \BS@ figset altitude(blcolor=ColorDefinition)}%
    \immediate\write16{ Base line style =
     \ifx\DDV@bldash\D@FTref general style\else\DDV@bldash\fi
     \space-> \BS@ figset altitude(bldash=Index/Pattern)}%
    \immediate\write16{ Base line width =
     \ifx\DDV@blwidth\D@FTref general width\else\DDV@blwidth\fi
     \space-> \BS@ figset altitude(blwidth=real in PostScript units)}%
    \immediate\write16{ Square line color =
     \ifx\DDV@sqcolor\D@FTref general color\else\DDV@sqcolor\fi
     \space-> \BS@ figset altitude(sqcolor=ColorDefinition)}%
    \immediate\write16{ Square line style =
     \ifx\DDV@sqdash\D@FTref general style\else\DDV@sqdash\fi
     \space-> \BS@ figset altitude(sqdash=Index/Pattern)}%
    \immediate\write16{ Square line width =
     \ifx\DDV@sqwidth\D@FTref general width\else\DDV@sqwidth\fi
     \space-> \BS@ figset altitude(sqwidth=real in PostScript units)}%
    \immediate\write16{Arrowhead:}%
    \immediate\write16{ (half-)Angle = \@rrowheadangle
     \space-> \BS@ figset arrowhead(angle=real in degrees)}%
    \immediate\write16{ Filling mode = \if@rrowhfill yes\else no\fi
     \space-> \BS@ figset arrowhead(fillmode=yes/no)}%
    \immediate\write16{ "Outside" = \if@rrowhout yes\else no\fi
     \space-> \BS@ figset arrowhead(out=yes/no)}%
    \immediate\write16{ Length = \@rrowheadlength
     \if@rrowratio\space(not active)\else\space(active)\fi
     \space-> \BS@ figset arrowhead(length=real in user coord.)}%
    \immediate\write16{ Ratio = \@rrowheadratio
     \if@rrowratio\space(active)\else\space(not active)\fi
     \space-> \BS@ figset arrowhead(ratio=real in [0,1])}%
    \immediate\write16{Curve:}%
    \immediate\write16{ Roundness = \curv@roundness
     \space-> \BS@ figset curve(roundness=real in [0,0.5])}%
    \immediate\write16{Flow chart:}%
    \immediate\write16{ Arrow position = \@rrowp@s
     \space-> \BS@ figset flowchart(arrowposition=real in [0,1])}%
    \immediate\write16{ Arrow reference point = \ifcase\@rrowr@fpt start\else end\fi
     \space-> \BS@ figset flowchart(arrowrefpt = start/end)}%
    \immediate\write16{ Background color = \fcbgc@lor
     \space-> \BS@ figset flowchart(bgcolor=ColorDefinition)}%
    \immediate\write16{ Line type = \ifcase\fclin@typ@ curve\else polygon\fi
     \space-> \BS@ figset flowchart(line=polygon/curve)}%
    \immediate\write16{ Padding = (\Xp@dd, \Yp@dd)
     \space-> \BS@ figset flowchart(padding = real in user coord.)}%
    \immediate\write16{\space\space\space\space(or
     \BS@ figset flowchart(xpadding=real, ypadding=real) )}%
    \immediate\write16{ Radius = \fclin@r@d
     \space-> \BS@ figset flowchart(radius=positive real in user coord.)}%
    \immediate\write16{ Shape = \fcsh@pe
     \space-> \BS@ figset flowchart(shape = rectangle, ellipse or lozenge)}%
    \immediate\write16{ Thickness color (DDV) = 
     \ifx\DDV@thickcolor\D@FTref general color\else\DDV@thickcolor\fi
     \space-> \BS@ figset flowchart(thickcolor=ColorDefinition)}%
    \immediate\write16{ Thickness = \thickn@ss
     \space-> \BS@ figset flowchart(thickness = real in user coord.)}%
    \immediate\write16{Mesh:}%
    \immediate\write16{ Diagonal = \c@ntrolmesh
     \space-> \BS@ figset mesh(diag=integer in {-1,0,1})}%
    \immediate\write16{ Lines color (DDV) =
     \ifx\DDV@meshcolor\D@FTref general color\else\DDV@meshcolor\fi
     \space-> \BS@ figset mesh(color=ColorDefinition)}%
    \immediate\write16{ Lines style (DDV) =
     \ifx\DDV@meshdash\D@FTref general style\else\DDV@meshdash\fi
     \space-> \BS@ figset mesh(dash=Index/Pattern)}%
    \immediate\write16{ Lines width (DDV) =
     \ifx\DDV@meshwidth\D@FTref general width\else\DDV@meshwidth\fi
     \space-> \BS@ figset mesh(width=real in PostScript units)}%
    \immediate\write16{Trimesh:}%
    \immediate\write16{ Lines color (DDV) =
     \ifx\DDV@tmeshcolor\D@FTref general color\else\DDV@tmeshcolor\fi
     \space-> \BS@ figset trimesh(color=ColorDefinition)}%
    \immediate\write16{ Lines style (DDV) =
     \ifx\DDV@tmeshdash\D@FTref general style\else\DDV@tmeshdash\fi
     \space-> \BS@ figset trimesh(dash=Index/Pattern)}%
    \immediate\write16{ Lines width (DDV) =
     \ifx\DDV@tmeshwidth\D@FTref general width\else\DDV@tmeshwidth\fi
     \space-> \BS@ figset trimesh(width=real in PostScript units)}%
    \ifTr@isDim%
    \immediate\write16{ --- 3D to 2D PROJECTION ---}%
    \immediate\write16{Projection : \typ@proj \space-> \BS@ figinit{ScaleFactorUnit, ProjType}}%
    \immediate\write16{Longitude (psi) = \v@lPsi \space-> \BS@ figset proj(psi=real in degrees)}%
    \ifcase\CUR@proj\immediate\write16{Depth coeff. (Lambda)
     \space = \v@lTheta \space-> \BS@ figset proj(lambda=real in [0,1])}%
    \else\immediate\write16{Latitude (theta)
     \space = \v@lTheta \space-> \BS@ figset proj(theta=real in degrees)}%
    \fi%
    \ifnum\CUR@proj=\tw@%
    \immediate\write16{Observation distance = \disob@unit
     \space-> \BS@ figset proj(dist=real in user coord.)}%
    \immediate\write16{Target point = \t@rgetpt \space-> \BS@ figset proj(targetpt=pt number)}%
     \v@lX=\ptT@unit@\wd\Bt@rget\v@lY=\ptT@unit@\ht\Bt@rget\v@lZ=\ptT@unit@\dp\Bt@rget%
    \immediate\write16{ Its coordinates are
     (\repdecn@mb{\v@lX}, \repdecn@mb{\v@lY}, \repdecn@mb{\v@lZ})}%
    \fi%
    \fi%
    \immediate\write16{====================================================================}%
    \ignorespaces}}
\ctr@ln@w{newif}\ifitis@vect@r
\ctr@ld@f\def\figvectC#1(#2,#3){{\itis@vect@rtrue\figpt#1:(#2,#3)}\ignorespaces}
\ctr@ld@f\def\Figv@ctCreg#1(#2,#3){{\itis@vect@rtrue\Figp@intreg#1:(#2,#3)}\ignorespaces}
\ctr@ln@m\figvectDBezier
\ctr@ld@f\def\figvectDBezierDD#1:#2,#3[#4,#5,#6,#7]{\ifGR@cri{\s@uvc@ntr@l\et@tfigvectDBezierDD%
    \FigvectDBezier@#2,#3[#4,#5,#6,#7]\v@lX=\c@ef\v@lX\v@lY=\c@ef\v@lY%
    \Figv@ctCreg#1(\v@lX,\v@lY)\resetc@ntr@l\et@tfigvectDBezierDD}\ignorespaces\fi}
\ctr@ld@f\def\figvectDBezierTD#1:#2,#3[#4,#5,#6,#7]{\ifGR@cri{\s@uvc@ntr@l\et@tfigvectDBezierTD%
    \FigvectDBezier@#2,#3[#4,#5,#6,#7]\v@lX=\c@ef\v@lX\v@lY=\c@ef\v@lY\v@lZ=\c@ef\v@lZ%
    \Figv@ctCreg#1(\v@lX,\v@lY,\v@lZ)\resetc@ntr@l\et@tfigvectDBezierTD}\ignorespaces\fi}
\ctr@ld@f\def\FigvectDBezier@#1,#2[#3,#4,#5,#6]{\setc@ntr@l{2}%
    \edef\T@{#2}\v@leur=\p@\advance\v@leur-#2pt\edef\UNmT@{\repdecn@mb{\v@leur}}%
    \ifnum#1=\tw@\def\c@ef{6}\else\def\c@ef{3}\fi%
    \figptcopy-4:/#3/\figptcopy-3:/#4/\figptcopy-2:/#5/\figptcopy-1:/#6/%
    \l@mbd@un=-4 \l@mbd@de=-\thr@@\p@rtent=\m@ne\c@lDecast%
    \ifnum#1=\tw@\c@lDCDeux{-4}{-3}\c@lDCDeux{-3}{-2}\c@lDCDeux{-4}{-3}\else%
    \l@mbd@un=-4 \l@mbd@de=-\thr@@\p@rtent=-\tw@\c@lDecast%
    \c@lDCDeux{-4}{-3}\fi\Figg@tXY{-4}}
\ctr@ln@m\c@lDCDeux
\ctr@ld@f\def\c@lDCDeuxDD#1#2{\Figg@tXY{#2}\Figg@tXYa{#1}%
    \advance\v@lX-\v@lXa\advance\v@lY-\v@lYa\Figp@intregDD#1:(\v@lX,\v@lY)}
\ctr@ld@f\def\c@lDCDeuxTD#1#2{\Figg@tXY{#2}\Figg@tXYa{#1}\advance\v@lX-\v@lXa%
    \advance\v@lY-\v@lYa\advance\v@lZ-\v@lZa\Figp@intregTD#1:(\v@lX,\v@lY,\v@lZ)}
\ctr@ln@m\figvectN
\ctr@ld@f\def\figvectNDD#1[#2,#3]{\ifGR@cri{\Figg@tXYa{#2}\Figg@tXY{#3}%
    \advance\v@lX-\v@lXa\advance\v@lY-\v@lYa%
    \Figv@ctCreg#1(-\v@lY,\v@lX)}\ignorespaces\fi}
\ctr@ld@f\def\figvectNTD#1[#2,#3,#4]{\ifGR@cri{\vecunitC@TD[#2,#4]\v@lmin=\v@lX\v@lmax=\v@lY%
    \v@leur=\v@lZ\vecunitC@TD[#2,#3]\c@lprovec{#1}}\ignorespaces\fi}
\ctr@ln@m\figvectNV
\ctr@ld@f\def\figvectNVDD#1[#2]{\ifGR@cri{\Figg@tXY{#2}\Figv@ctCreg#1(-\v@lY,\v@lX)}\ignorespaces\fi}
\ctr@ld@f\def\figvectNVTD#1[#2,#3]{\ifGR@cri{\vecunitCV@TD{#3}\v@lmin=\v@lX\v@lmax=\v@lY%
    \v@leur=\v@lZ\vecunitCV@TD{#2}\c@lprovec{#1}}\ignorespaces\fi}
\ctr@ln@m\figvectP
\ctr@ld@f\def\figvectPDD#1[#2,#3]{\ifGR@cri{\Figg@tXYa{#2}\Figg@tXY{#3}%
    \advance\v@lX-\v@lXa\advance\v@lY-\v@lYa%
    \Figv@ctCreg#1(\v@lX,\v@lY)}\ignorespaces\fi}
\ctr@ld@f\def\figvectPTD#1[#2,#3]{\ifGR@cri{\Figg@tXYa{#2}\Figg@tXY{#3}%
    \advance\v@lX-\v@lXa\advance\v@lY-\v@lYa\advance\v@lZ-\v@lZa%
    \Figv@ctCreg#1(\v@lX,\v@lY,\v@lZ)}\ignorespaces\fi}
\ctr@ln@m\figvectU
\ctr@ld@f\def\figvectUDD#1[#2]{\ifGR@cri{\n@rmeuc\v@leur{#2}\invers@\v@leur\v@leur%
    \delt@=\repdecn@mb{\v@leur}\unit@\edef\v@ldelt@{\repdecn@mb{\delt@}}%
    \Figg@tXY{#2}\v@lX=\v@ldelt@\v@lX\v@lY=\v@ldelt@\v@lY%
    \Figv@ctCreg#1(\v@lX,\v@lY)}\ignorespaces\fi}
\ctr@ld@f\def\figvectUTD#1[#2]{\ifGR@cri{\n@rmeuc\v@leur{#2}\invers@\v@leur\v@leur%
    \delt@=\repdecn@mb{\v@leur}\unit@\edef\v@ldelt@{\repdecn@mb{\delt@}}%
    \Figg@tXY{#2}\v@lX=\v@ldelt@\v@lX\v@lY=\v@ldelt@\v@lY\v@lZ=\v@ldelt@\v@lZ%
    \Figv@ctCreg#1(\v@lX,\v@lY,\v@lZ)}\ignorespaces\fi}
\ctr@ld@f\def\figvisu#1#2#3{\c@ldefproj\initb@undb@x\xdef\figforTeXFigno{\figforTeXnextFigno}%
    \s@mme=\figforTeXnextFigno\advance\s@mme\@ne\xdef\figforTeXnextFigno{\number\s@mme}%
    \setbox\b@xvisu=\hbox{\ifnum\@utoFN>\z@\figinsert{}\gdef\@utoFInDone{0}\fi\ignorespaces#3}%
    \gdef\@utoFInDone{1}\gdef\@utoFN{0}%
    \v@lXa=-\c@@rdYmin\v@lYa=\c@@rdYmax\advance\v@lYa-\c@@rdYmin%
    \v@lX=\c@@rdXmax\advance\v@lX-\c@@rdXmin%
    \setbox#1=\hbox{#2}\v@lY=-\v@lX\maxim@m{\v@lX}{\v@lX}{\wd#1}%
    \advance\v@lY\v@lX\divide\v@lY\tw@\advance\v@lY-\c@@rdXmin%
    \setbox#1=\vbox{\parindent\z@\hsize=\v@lX\vskip\v@lYa%
    \rlap{\hskip\v@lY\smash{\raise\v@lXa\box\b@xvisu}}%
    \def\t@xt@{#2}\ifx\t@xt@\empty\else\medskip\centerline{#2}\fi}\wd#1=\v@lX}
\ctr@ld@f\def\figDecrementFigno{{\xdef\figforTeXnextFigno{\figforTeXFigno}%
    \s@mme=\figforTeXFigno\advance\s@mme\m@ne\xdef\figforTeXFigno{\number\s@mme}}}
\ctr@ln@w{newbox}\Bt@rget\setbox\Bt@rget=\null
\ctr@ln@w{newbox}\BminTD@\setbox\BminTD@=\null
\ctr@ln@w{newbox}\BmaxTD@\setbox\BmaxTD@=\null
\ctr@ln@w{newif}\ifnewt@rgetpt\ctr@ln@w{newif}\ifnewdis@b
\ctr@ld@f\def\b@undb@xTD#1#2#3{%
    \relax\ifdim#1<\wd\BminTD@\global\wd\BminTD@=#1\fi%
    \relax\ifdim#2<\ht\BminTD@\global\ht\BminTD@=#2\fi%
    \relax\ifdim#3<\dp\BminTD@\global\dp\BminTD@=#3\fi%
    \relax\ifdim#1>\wd\BmaxTD@\global\wd\BmaxTD@=#1\fi%
    \relax\ifdim#2>\ht\BmaxTD@\global\ht\BmaxTD@=#2\fi%
    \relax\ifdim#3>\dp\BmaxTD@\global\dp\BmaxTD@=#3\fi}
\ctr@ld@f\def\c@ldefdisob{{\ifdim\wd\BminTD@<\maxdimen\v@leur=\wd\BmaxTD@\advance\v@leur-\wd\BminTD@%
    \delt@=\ht\BmaxTD@\advance\delt@-\ht\BminTD@\maxim@m{\v@leur}{\v@leur}{\delt@}%
    \delt@=\dp\BmaxTD@\advance\delt@-\dp\BminTD@\maxim@m{\v@leur}{\v@leur}{\delt@}%
    \v@leur=5\v@leur\else\v@leur=800pt\fi\c@ldefdisob@{\v@leur}}}
\ctr@ln@m\disob@intern
\ctr@ln@m\disob@
\ctr@ln@m\divf@ctproj
\ctr@ld@f\def\c@ldefdisob@#1{{\v@leur=#1\ifdim\v@leur<\p@\v@leur=800pt\fi%
    \xdef\disob@intern{\repdecn@mb{\v@leur}}%
    \delt@=\ptT@unit@\v@leur\xdef\disob@unit{\repdecn@mb{\delt@}}%
    \f@ctech=\@ne\loop\ifdim\v@leur>\t@n pt\divide\v@leur\t@n\multiply\f@ctech\t@n\repeat%
    \xdef\disob@{\repdecn@mb{\v@leur}}\xdef\divf@ctproj{\the\f@ctech}}%
    \global\newdis@btrue}
\ctr@ln@m\t@rgetpt
\ctr@ld@f\def\c@ldeft@rgetpt{\newt@rgetpttrue\def\t@rgetpt{CenterBoundBox}{%
    \delt@=\wd\BmaxTD@\advance\delt@-\wd\BminTD@\divide\delt@\tw@%
    \v@leur=\wd\BminTD@\advance\v@leur\delt@\global\wd\Bt@rget=\v@leur%
    \delt@=\ht\BmaxTD@\advance\delt@-\ht\BminTD@\divide\delt@\tw@%
    \v@leur=\ht\BminTD@\advance\v@leur\delt@\global\ht\Bt@rget=\v@leur%
    \delt@=\dp\BmaxTD@\advance\delt@-\dp\BminTD@\divide\delt@\tw@%
    \v@leur=\dp\BminTD@\advance\v@leur\delt@\global\dp\Bt@rget=\v@leur}}
\ctr@ln@m\c@ldefproj
\ctr@ld@f\def\c@ldefprojTD{\ifnewt@rgetpt\else\c@ldeft@rgetpt\fi\ifnewdis@b\else\c@ldefdisob\fi}
\ctr@ld@f\def\c@lprojcav{
    \v@lZa=\cxa@\v@lY\advance\v@lX\v@lZa%
    \v@lZa=\cxb@\v@lY\v@lY=\v@lZ\advance\v@lY\v@lZa\ignorespaces}
\ctr@ln@m\v@lcoef
\ctr@ld@f\def\c@lprojrea{
    \advance\v@lX-\wd\Bt@rget\advance\v@lY-\ht\Bt@rget\advance\v@lZ-\dp\Bt@rget%
    \v@lZa=\cza@\v@lX\advance\v@lZa\czb@\v@lY\advance\v@lZa\czc@\v@lZ%
    \divide\v@lZa\divf@ctproj\advance\v@lZa\disob@ pt\invers@{\v@lZa}{\v@lZa}%
    \v@lZa=\disob@\v@lZa\edef\v@lcoef{\repdecn@mb{\v@lZa}}%
    \v@lXa=\cxa@\v@lX\advance\v@lXa\cxb@\v@lY\v@lXa=\v@lcoef\v@lXa%
    \v@lY=\cyb@\v@lY\advance\v@lY\cya@\v@lX\advance\v@lY\cyc@\v@lZ%
    \v@lY=\v@lcoef\v@lY\v@lX=\v@lXa\ignorespaces}
\ctr@ld@f\def\c@lprojort{
    \v@lXa=\cxa@\v@lX\advance\v@lXa\cxb@\v@lY%
    \v@lY=\cyb@\v@lY\advance\v@lY\cya@\v@lX\advance\v@lY\cyc@\v@lZ%
    \v@lX=\v@lXa\ignorespaces}
\ctr@ld@f\def\Figptpr@j#1:#2/#3/{{\Figg@tXY{#3}\superc@lprojSP%
    \Figp@intregDD#1:{#2}(\v@lX,\v@lY)}\ignorespaces}
\ctr@ln@m\figsetobdist
\ctr@ld@f\def\figsetobdistDD{\un@v@ilable{figsetobdist}}
\ctr@ld@f\def\figsetobdistTD(#1){{\ifCUR@PS\W@rnmesIgn{figset proj(dist=...)}%
    \else\v@leur=#1\unit@\c@ldefdisob@{\v@leur}\fi}\ignorespaces}
\ctr@ln@m\c@lprojSP
\ctr@ln@m\CUR@proj
\ctr@ln@m\typ@proj
\ctr@ln@m\superc@lprojSP
\ctr@ld@f\def\Figs@tproj#1{%
    \if#13 \def@ultproj\else\if#1c\def@ultproj%
    \else\if#1o\xdef\CUR@proj{1}\xdef\typ@proj{orthogonal}%
         \figsetviewTD(\def@ultpsi,\def@ulttheta)%
         \global\let\c@lprojSP=\c@lprojort\global\let\superc@lprojSP=\c@lprojort%
    \else\if#1r\xdef\CUR@proj{2}\xdef\typ@proj{realistic}%
         \figsetviewTD(\def@ultpsi,\def@ulttheta)%
         \global\let\c@lprojSP=\c@lprojrea\global\let\superc@lprojSP=\c@lprojrea%
    \else\def@ultproj\message{*** Unknown projection. Cavalier projection assumed.}%
    \fi\fi\fi\fi}
\ctr@ld@f\def\def@ultproj{\xdef\CUR@proj{0}\xdef\typ@proj{cavalier}\figsetviewTD(\def@ultpsi,0.5)%
         \global\let\c@lprojSP=\c@lprojcav\global\let\superc@lprojSP=\c@lprojcav}
\ctr@ln@m\figsettarget
\ctr@ld@f\def\figsettargetDD{\un@v@ilable{figsettarget}}
\ctr@ld@f\def\figsettargetTD[#1]{{\ifCUR@PS\W@rnmesIgn{figset proj(targetpt=...)}%
    \else\global\newt@rgetpttrue\xdef\t@rgetpt{#1}\Figg@tXY{#1}\global\wd\Bt@rget=\v@lX%
    \global\ht\Bt@rget=\v@lY\global\dp\Bt@rget=\v@lZ\fi}\ignorespaces}
\ctr@ln@m\figsetview
\ctr@ld@f\def\figsetviewDD{\un@v@ilable{figsetview}}
\ctr@ld@f\def\figsetviewTD(#1){\ifCUR@PS\W@rnmesIgn{figset proj(Psi|Theta|Lambda=...)}%
     \else\Figsetview@#1,:\fi\ignorespaces}
\ctr@ld@f\def\Figsetview@#1,#2:{{\xdef\v@lPsi{#1}\def\t@xt@{#2}%
    \ifx\t@xt@\empty\def\@rgdeux{\v@lTheta}\else\X@rgdeux@#2\fi%
    \c@ssin{\costhet@}{\sinthet@}{#1}\v@lmin=\costhet@ pt\v@lmax=\sinthet@ pt%
    \ifcase\CUR@proj%
    \v@leur=\@rgdeux\v@lmin\xdef\cxa@{\repdecn@mb{\v@leur}}%
    \v@leur=\@rgdeux\v@lmax\xdef\cxb@{\repdecn@mb{\v@leur}}\v@leur=\@rgdeux pt%
    \relax\ifdim\v@leur>\p@\message{*** Lambda too large ! See \BS@ figset proj() !}\fi%
    \else%
    \v@lmax=-\v@lmax\xdef\cxa@{\repdecn@mb{\v@lmax}}\xdef\cxb@{\costhet@}%
    \ifx\t@xt@\empty\edef\@rgdeux{\def@ulttheta}\fi\c@ssin{\C@}{\S@}{\@rgdeux}%
    \v@lmax=-\S@ pt%
    \v@leur=\v@lmax\v@leur=\costhet@\v@leur\xdef\cya@{\repdecn@mb{\v@leur}}%
    \v@leur=\v@lmax\v@leur=\sinthet@\v@leur\xdef\cyb@{\repdecn@mb{\v@leur}}%
    \xdef\cyc@{\C@}\v@lmin=-\C@ pt%
    \v@leur=\v@lmin\v@leur=\costhet@\v@leur\xdef\cza@{\repdecn@mb{\v@leur}}%
    \v@leur=\v@lmin\v@leur=\sinthet@\v@leur\xdef\czb@{\repdecn@mb{\v@leur}}%
    \xdef\czc@{\repdecn@mb{\v@lmax}}\fi%
    \xdef\v@lTheta{\@rgdeux}}}
\ctr@ld@f\def\def@ultpsi{40}
\ctr@ld@f\def\def@ulttheta{25}
\ctr@ln@m\l@debut
\ctr@ln@m\n@mref
\ctr@ld@f\def\Figsetpr@j#1=#2|{\keln@mtr#1|%
    \def\n@mref{dep}\ifx\l@debut\n@mref\Figsetd@p{#2}\else
    \def\n@mref{dis}\ifx\l@debut\n@mref%
     \ifnum\CUR@proj=\tw@\figsetobdist(#2)\else\Figset@rr\fi\else
    \def\n@mref{lam}\ifx\l@debut\n@mref\Figsetd@p{#2}\else
    \def\n@mref{lat}\ifx\l@debut\n@mref\Figsetth@{#2}\else
    \def\n@mref{lon}\ifx\l@debut\n@mref\figsetview(#2)\else
    \def\n@mref{psi}\ifx\l@debut\n@mref\figsetview(#2)\else
    \def\n@mref{tar}\ifx\l@debut\n@mref%
     \ifnum\CUR@proj=\tw@\figsettarget[#2]\else\Figset@rr\fi\else
    \def\n@mref{the}\ifx\l@debut\n@mref\Figsetth@{#2}\else
    \W@rnmesAttr{figset proj}{#1}\fi\fi\fi\fi\fi\fi\fi\fi}
\ctr@ld@f\def\Figsetd@p#1{\ifnum\CUR@proj=\z@\figsetview(\v@lPsi,#1)\else\Figset@rr\fi}
\ctr@ld@f\def\Figsetth@#1{\ifnum\CUR@proj=\z@\Figset@rr\else\figsetview(\v@lPsi,#1)\fi}
\ctr@ld@f\def\Figset@rr{\message{*** \BS@ figset proj(): Attribute "\n@mref" ignored, incompatible
    with current projection}}
\ctr@ld@f\def\initb@undb@xTD{\wd\BminTD@=\maxdimen\ht\BminTD@=\maxdimen\dp\BminTD@=\maxdimen%
    \wd\BmaxTD@=-\maxdimen\ht\BmaxTD@=-\maxdimen\dp\BmaxTD@=-\maxdimen}
\ctr@ln@w{newbox}\Gb@x      
\ctr@ln@w{newbox}\Gb@xSC    
\ctr@ln@w{newtoks}\c@nsymb  
\ctr@ln@w{newif}\ifr@undcoord\ctr@ln@w{newif}\ifunitpr@sent
\ctr@ld@f\def\unssqrttw@{0.707106 }
\ctr@ld@f\def\figAst{\raise-1.15ex\hbox{$\ast$}}
\ctr@ld@f\def\figBullet{\raise-1.15ex\hbox{$\bullet$}}
\ctr@ld@f\def\figCirc{\raise-1.15ex\hbox{$\circ$}}
\ctr@ld@f\def\figDiamond{\raise-1.15ex\hbox{$\diamond$}}%
\ctr@ld@f\def\boxit#1#2{\leavevmode\hbox{\vrule\vbox{\hrule\vglue#1%
    \vtop{\hbox{\kern#1{#2}\kern#1}\vglue#1\hrule}}\vrule}}
\ctr@ld@f
\ctr@ld@f
\ctr@ld@f\def\c@nterpt{\ignorespaces%
    \kern-.5\wd\Gb@xSC%
    \raise-.5\ht\Gb@xSC\rlap{\hbox{\raise.5\dp\Gb@xSC\hbox{\copy\Gb@xSC}}}%
    \kern .5\wd\Gb@xSC\ignorespaces}
\ctr@ld@f\def\b@undb@xSC#1#2{{\v@lXa=#1\v@lYa=#2%
    \v@leur=\ht\Gb@xSC\advance\v@leur\dp\Gb@xSC%
    \advance\v@lXa-.5\wd\Gb@xSC\advance\v@lYa-.5\v@leur\b@undb@x{\v@lXa}{\v@lYa}%
    \advance\v@lXa\wd\Gb@xSC\advance\v@lYa\v@leur\b@undb@x{\v@lXa}{\v@lYa}}}
\ctr@ln@m\Dist@n
\ctr@ln@m\l@suite
\ctr@ld@f\def\@keldist#1#2{\edef\Dist@n{#2}\y@tiunit{\Dist@n}%
    \ifunitpr@sent#1=\Dist@n\else#1=\Dist@n\unit@\fi}
\ctr@ld@f\def\y@tiunit#1{\unitpr@sentfalse\expandafter\y@tiunit@#1:}
\ctr@ld@f\def\y@tiunit@#1#2:{\ifcat#1a\unitpr@senttrue\else\def\l@suite{#2}%
    \ifx\l@suite\empty\else\y@tiunit@#2:\fi\fi}
\ctr@ln@m\figcoord
\ctr@ld@f\def\figcoordDD#1{{\v@lX=\ptT@unit@\v@lX\v@lY=\ptT@unit@\v@lY%
    \ifr@undcoord\ifcase#1\v@leur=0.5pt\or\v@leur=0.05pt\or\v@leur=0.005pt%
    \or\v@leur=0.0005pt\else\v@leur=\z@\fi%
    \ifdim\v@lX<\z@\advance\v@lX-\v@leur\else\advance\v@lX\v@leur\fi%
    \ifdim\v@lY<\z@\advance\v@lY-\v@leur\else\advance\v@lY\v@leur\fi\fi%
    (\@ffichnb{#1}{\repdecn@mb{\v@lX}},\ifmmode\else\thinspace\fi%
    \@ffichnb{#1}{\repdecn@mb{\v@lY}})}}
\ctr@ld@f\def\@ffichnb#1#2{{\def\@@ffich{\@ffich#1(}\edef\n@mbre{#2}%
    \expandafter\@@ffich\n@mbre)}}
\ctr@ld@f\def\@ffich#1(#2.#3){{#2\ifnum#1>\z@.\fi\def\dig@ts{#3}\s@mme=\z@%
    \loop\ifnum\s@mme<#1\expandafter\@ffichdec\dig@ts:\advance\s@mme\@ne\repeat}}
\ctr@ld@f\def\@ffichdec#1#2:{\relax#1\def\dig@ts{#20}}
\ctr@ld@f\def\figcoordTD#1{{\v@lX=\ptT@unit@\v@lX\v@lY=\ptT@unit@\v@lY\v@lZ=\ptT@unit@\v@lZ%
    \ifr@undcoord\ifcase#1\v@leur=0.5pt\or\v@leur=0.05pt\or\v@leur=0.005pt%
    \or\v@leur=0.0005pt\else\v@leur=\z@\fi%
    \ifdim\v@lX<\z@\advance\v@lX-\v@leur\else\advance\v@lX\v@leur\fi%
    \ifdim\v@lY<\z@\advance\v@lY-\v@leur\else\advance\v@lY\v@leur\fi%
    \ifdim\v@lZ<\z@\advance\v@lZ-\v@leur\else\advance\v@lZ\v@leur\fi\fi%
    (\@ffichnb{#1}{\repdecn@mb{\v@lX}},\ifmmode\else\thinspace\fi%
     \@ffichnb{#1}{\repdecn@mb{\v@lY}},\ifmmode\else\thinspace\fi%
     \@ffichnb{#1}{\repdecn@mb{\v@lZ}})}}
\ctr@ld@f\def\figsetroundcoord#1{\expandafter\Figsetr@undcoord#1:\ignorespaces}
\ctr@ld@f\def\Figsetr@undcoord#1#2:{\if#1n\r@undcoordfalse\else\r@undcoordtrue\fi}
\ctr@ld@f\def\Figsetwr@te#1=#2|{\keln@mun#1|%
    \def\n@mref{m}\ifx\l@debut\n@mref\figsetmark{#2}\else
    \def\n@mref{p}\ifx\l@debut\n@mref\figsetptname{#2}\else
    \def\n@mref{r}\ifx\l@debut\n@mref\figsetroundcoord{#2}\else
    \W@rnmesAttr{figset write}{#1}\fi\fi\fi}
\ctr@ld@f\def\figsetmark#1{\c@nsymb={#1}\setbox\Gb@xSC=\hbox{\the\c@nsymb}\ignorespaces}
\ctr@ln@m\ptn@me
\ctr@ld@f\def\figsetptname#1{\def\ptn@me##1{#1}\ignorespaces}
\ctr@ld@f\def\FigWrit@L#1:#2(#3,#4){\ignorespaces\@keldist\v@leur{#3}\@keldist\delt@{#4}%
    \C@rp@r@m\def\list@num{#1}\@ecfor\p@int:=\list@num\do{\FigWrit@pt{\p@int}{#2}}}
\ctr@ld@f\def\FigWrit@pt#1#2{\FigWp@r@m{#1}{#2}\Vc@rrect\figWp@si%
    \ifdim\wd\Gb@xSC>\z@\b@undb@xSC{\v@lX}{\v@lY}\fi\figWBB@x}
\ctr@ld@f\def\FigWp@r@m#1#2{\Figg@tXY{#1}%
    \setbox\Gb@x=\hbox{\def\t@xt@{#2}\ifx\t@xt@\empty\Figg@tT{#1}\else#2\fi}\c@lprojSP}
\ctr@ld@f\let\Vc@rrect=\relax
\ctr@ld@f\let\C@rp@r@m=\relax
\ctr@ld@f\def\figwrite[#1]#2{{\ignorespaces\def\list@num{#1}\@ecfor\p@int:=\list@num\do{%
    \setbox\Gb@x=\hbox{\def\t@xt@{#2}\ifx\t@xt@\empty\Figg@tT{\p@int}\else#2\fi}%
    \Figwrit@{\p@int}}}\ignorespaces}
\ctr@ld@f\def\Figwrit@#1{\Figg@tXY{#1}\c@lprojSP%
    \rlap{\kern\v@lX\raise\v@lY\hbox{\unhcopy\Gb@x}}\v@leur=\v@lY%
    \advance\v@lY\ht\Gb@x\b@undb@x{\v@lX}{\v@lY}\advance\v@lX\wd\Gb@x%
    \v@lY=\v@leur\advance\v@lY-\dp\Gb@x\b@undb@x{\v@lX}{\v@lY}}
\ctr@ld@f\def\figwritec[#1]#2{{\ignorespaces\def\list@num{#1}%
    \@ecfor\p@int:=\list@num\do{\Figwrit@c{\p@int}{#2}}}\ignorespaces}
\ctr@ld@f\def\Figwrit@c#1#2{\FigWp@r@m{#1}{#2}%
    \rlap{\kern\v@lX\raise\v@lY\hbox{\rlap{\kern-.5\wd\Gb@x%
    \raise-.5\ht\Gb@x\hbox{\raise.5\dp\Gb@x\hbox{\unhcopy\Gb@x}}}}}%
    \v@leur=\ht\Gb@x\advance\v@leur\dp\Gb@x%
    \advance\v@lX-.5\wd\Gb@x\advance\v@lY-.5\v@leur\b@undb@x{\v@lX}{\v@lY}%
    \advance\v@lX\wd\Gb@x\advance\v@lY\v@leur\b@undb@x{\v@lX}{\v@lY}}
\ctr@ld@f\def\figwritep[#1]{{\ignorespaces\def\list@num{#1}\setbox\Gb@x=\hbox{\c@nterpt}%
    \@ecfor\p@int:=\list@num\do{\Figwrit@{\p@int}}}\ignorespaces}
\ctr@ld@f\def\figwritew#1:#2(#3){\figwritegcw#1:{#2}(#3,0pt)}
\ctr@ld@f\def\figwritee#1:#2(#3){\figwritegce#1:{#2}(#3,0pt)}
\ctr@ld@f\def\figwriten#1:#2(#3){{\def\Vc@rrect{\v@lZ=\v@leur\advance\v@lZ\dp\Gb@x}%
    \Figwrit@NS#1:{#2}(#3)}\ignorespaces}
\ctr@ld@f\def\figwrites#1:#2(#3){{\def\Vc@rrect{\v@lZ=-\v@leur\advance\v@lZ-\ht\Gb@x}%
    \Figwrit@NS#1:{#2}(#3)}\ignorespaces}
\ctr@ld@f\def\Figwrit@NS#1:#2(#3){\let\figWp@si=\FigWp@siNS\let\figWBB@x=\FigWBB@xNS%
    \FigWrit@L#1:{#2}(#3,0pt)}
\ctr@ld@f\def\FigWp@siNS{\rlap{\kern\v@lX\raise\v@lY\hbox{\rlap{\kern-.5\wd\Gb@x%
    \raise\v@lZ\hbox{\unhcopy\Gb@x}}\c@nterpt}}}
\ctr@ld@f\def\FigWBB@xNS{\advance\v@lY\v@lZ%
    \advance\v@lY-\dp\Gb@x\advance\v@lX-.5\wd\Gb@x\b@undb@x{\v@lX}{\v@lY}%
    \advance\v@lY\ht\Gb@x\advance\v@lY\dp\Gb@x%
    \advance\v@lX\wd\Gb@x\b@undb@x{\v@lX}{\v@lY}}
\ctr@ld@f\def\figwritenw#1:#2(#3){{\let\figWp@si=\FigWp@sigW\let\figWBB@x=\FigWBB@xgWE%
    \def\C@rp@r@m{\v@leur=\unssqrttw@\v@leur\delt@=\v@leur%
    \ifdim\delt@=\z@\delt@=\epsil@n\fi}\let@xte={-}\FigWrit@L#1:{#2}(#3,0pt)}\ignorespaces}
\ctr@ld@f\def\figwritesw#1:#2(#3){{\let\figWp@si=\FigWp@sigW\let\figWBB@x=\FigWBB@xgWE%
    \def\C@rp@r@m{\v@leur=\unssqrttw@\v@leur\delt@=-\v@leur%
    \ifdim\delt@=\z@\delt@=-\epsil@n\fi}\let@xte={-}\FigWrit@L#1:{#2}(#3,0pt)}\ignorespaces}
\ctr@ld@f\def\figwritene#1:#2(#3){{\let\figWp@si=\FigWp@sigE\let\figWBB@x=\FigWBB@xgWE%
    \def\C@rp@r@m{\v@leur=\unssqrttw@\v@leur\delt@=\v@leur%
    \ifdim\delt@=\z@\delt@=\epsil@n\fi}\let@xte={}\FigWrit@L#1:{#2}(#3,0pt)}\ignorespaces}
\ctr@ld@f\def\figwritese#1:#2(#3){{\let\figWp@si=\FigWp@sigE\let\figWBB@x=\FigWBB@xgWE%
    \def\C@rp@r@m{\v@leur=\unssqrttw@\v@leur\delt@=-\v@leur%
    \ifdim\delt@=\z@\delt@=-\epsil@n\fi}\let@xte={}\FigWrit@L#1:{#2}(#3,0pt)}\ignorespaces}
\ctr@ld@f\def\figwritegw#1:#2(#3,#4){{\let\figWp@si=\FigWp@sigW\let\figWBB@x=\FigWBB@xgWE%
    \let@xte={-}\FigWrit@L#1:{#2}(#3,#4)}\ignorespaces}
\ctr@ld@f\def\figwritege#1:#2(#3,#4){{\let\figWp@si=\FigWp@sigE\let\figWBB@x=\FigWBB@xgWE%
    \let@xte={}\FigWrit@L#1:{#2}(#3,#4)}\ignorespaces}
\ctr@ld@f\def\FigWp@sigW{\v@lXa=\z@\v@lYa=\ht\Gb@x\advance\v@lYa\dp\Gb@x%
    \ifdim\delt@>\z@\relax%
    \rlap{\kern\v@lX\raise\v@lY\hbox{\rlap{\kern-\wd\Gb@x\kern-\v@leur%
          \raise\delt@\hbox{\raise\dp\Gb@x\hbox{\unhcopy\Gb@x}}}\c@nterpt}}%
    \else\ifdim\delt@<\z@\relax\v@lYa=-\v@lYa%
    \rlap{\kern\v@lX\raise\v@lY\hbox{\rlap{\kern-\wd\Gb@x\kern-\v@leur%
          \raise\delt@\hbox{\raise-\ht\Gb@x\hbox{\unhcopy\Gb@x}}}\c@nterpt}}%
    \else\v@lXa=-.5\v@lYa%
    \rlap{\kern\v@lX\raise\v@lY\hbox{\rlap{\kern-\wd\Gb@x\kern-\v@leur%
          \raise-.5\ht\Gb@x\hbox{\raise.5\dp\Gb@x\hbox{\unhcopy\Gb@x}}}\c@nterpt}}%
    \fi\fi}
\ctr@ld@f\def\FigWp@sigE{\v@lXa=\z@\v@lYa=\ht\Gb@x\advance\v@lYa\dp\Gb@x%
    \ifdim\delt@>\z@\relax%
    \rlap{\kern\v@lX\raise\v@lY\hbox{\c@nterpt\kern\v@leur%
          \raise\delt@\hbox{\raise\dp\Gb@x\hbox{\unhcopy\Gb@x}}}}%
    \else\ifdim\delt@<\z@\relax\v@lYa=-\v@lYa%
    \rlap{\kern\v@lX\raise\v@lY\hbox{\c@nterpt\kern\v@leur%
          \raise\delt@\hbox{\raise-\ht\Gb@x\hbox{\unhcopy\Gb@x}}}}%
    \else\v@lXa=-.5\v@lYa%
    \rlap{\kern\v@lX\raise\v@lY\hbox{\c@nterpt\kern\v@leur%
          \raise-.5\ht\Gb@x\hbox{\raise.5\dp\Gb@x\hbox{\unhcopy\Gb@x}}}}%
    \fi\fi}
\ctr@ld@f\def\FigWBB@xgWE{\advance\v@lY\delt@%
    \advance\v@lX\the\let@xte\v@leur\advance\v@lY\v@lXa\b@undb@x{\v@lX}{\v@lY}%
    \advance\v@lX\the\let@xte\wd\Gb@x\advance\v@lY\v@lYa\b@undb@x{\v@lX}{\v@lY}}
\ctr@ld@f\def\figwritegcw#1:#2(#3,#4){{\let\figWp@si=\FigWp@sigcW\let\figWBB@x=\FigWBB@xgcWE%
    \let@xte={-}\FigWrit@L#1:{#2}(#3,#4)}\ignorespaces}
\ctr@ld@f\def\figwritegce#1:#2(#3,#4){{\let\figWp@si=\FigWp@sigcE\let\figWBB@x=\FigWBB@xgcWE%
    \let@xte={}\FigWrit@L#1:{#2}(#3,#4)}\ignorespaces}
\ctr@ld@f\def\FigWp@sigcW{\rlap{\kern\v@lX\raise\v@lY\hbox{\rlap{\kern-\wd\Gb@x\kern-\v@leur%
     \raise-.5\ht\Gb@x\hbox{\raise\delt@\hbox{\raise.5\dp\Gb@x\hbox{\unhcopy\Gb@x}}}}%
     \c@nterpt}}}
\ctr@ld@f\def\FigWp@sigcE{\rlap{\kern\v@lX\raise\v@lY\hbox{\c@nterpt\kern\v@leur%
    \raise-.5\ht\Gb@x\hbox{\raise\delt@\hbox{\raise.5\dp\Gb@x\hbox{\unhcopy\Gb@x}}}}}}
\ctr@ld@f\def\FigWBB@xgcWE{\v@lZ=\ht\Gb@x\advance\v@lZ\dp\Gb@x%
    \advance\v@lX\the\let@xte\v@leur\advance\v@lY\delt@\advance\v@lY.5\v@lZ%
    \b@undb@x{\v@lX}{\v@lY}%
    \advance\v@lX\the\let@xte\wd\Gb@x\advance\v@lY-\v@lZ\b@undb@x{\v@lX}{\v@lY}}
\ctr@ld@f\def\figwritebn#1:#2(#3){{\def\Vc@rrect{\v@lZ=\v@leur}\Figwrit@NS#1:{#2}(#3)}\ignorespaces}
\ctr@ld@f\def\figwritebs#1:#2(#3){{\def\Vc@rrect{\v@lZ=-\v@leur}\Figwrit@NS#1:{#2}(#3)}\ignorespaces}
\ctr@ld@f\def\figwritebw#1:#2(#3){{\let\figWp@si=\FigWp@sibW\let\figWBB@x=\FigWBB@xbWE%
    \let@xte={-}\FigWrit@L#1:{#2}(#3,0pt)}\ignorespaces}
\ctr@ld@f\def\figwritebe#1:#2(#3){{\let\figWp@si=\FigWp@sibE\let\figWBB@x=\FigWBB@xbWE%
    \let@xte={}\FigWrit@L#1:{#2}(#3,0pt)}\ignorespaces}
\ctr@ld@f\def\FigWp@sibW{\rlap{\kern\v@lX\raise\v@lY\hbox{\rlap{\kern-\wd\Gb@x\kern-\v@leur%
          \hbox{\unhcopy\Gb@x}}\c@nterpt}}}
\ctr@ld@f\def\FigWp@sibE{\rlap{\kern\v@lX\raise\v@lY\hbox{\c@nterpt\kern\v@leur%
          \hbox{\unhcopy\Gb@x}}}}
\ctr@ld@f\def\FigWBB@xbWE{\v@lZ=\ht\Gb@x\advance\v@lZ\dp\Gb@x%
    \advance\v@lX\the\let@xte\v@leur\advance\v@lY\ht\Gb@x\b@undb@x{\v@lX}{\v@lY}%
    \advance\v@lX\the\let@xte\wd\Gb@x\advance\v@lY-\v@lZ\b@undb@x{\v@lX}{\v@lY}}
\ctr@ln@w{newread}\frf@g  \ctr@ln@w{newwrite}\fwf@g
\ctr@ln@w{newif}\ifCUR@PS
\ctr@ln@w{newif}\ifGR@cri
\ctr@ln@w{newif}\ifUse@llipse
\ctr@ln@w{newif}\ifGRdebugm@de \GRdebugm@defalse 
\ctr@ln@w{newif}\ifPDFm@ke
\ifx\pdfliteral\undefined\else\ifnum\pdfoutput>\z@\PDFm@ketrue\fi\fi
\ctr@ld@f\def\initPDF@rDVI{%
\ifPDFm@ke
 \let\figscan=\figscan@E
 \let\newGr@FN=\newGr@FNPDF
 \ctr@ld@f\def\c@mcurveto{c}
 \ctr@ld@f\def\c@mfill{f}
 \ctr@ld@f\def\c@mgsave{q}
 \ctr@ld@f\def\c@mgrestore{Q}
 \ctr@ld@f\def\c@mlineto{l}
 \ctr@ld@f\def\c@mmoveto{m}
 \ctr@ld@f\def\c@msetgray{g}     \ctr@ld@f\def\c@msetgrayStroke{G}
 \ctr@ld@f\def\c@msetcmykcolor{k}\ctr@ld@f\def\c@msetcmykcolorStroke{K}
 \ctr@ld@f\def\c@msetrgbcolor{rg}\ctr@ld@f\def\c@msetrgbcolorStroke{RG}
 \ctr@ld@f\def\d@fprimarC@lor{\CUR@color\space\CUR@colorc@md%
               \space\CUR@color\space\CUR@colorc@mdStroke}
 \ctr@ld@f\def\c@msetdash{d}
 \ctr@ld@f\def\c@msetlinejoin{j}
 \ctr@ld@f\def\c@msetlinewidth{w}
 \ctr@ld@f\def\f@gclosestroke{\immediate\write\fwf@g{s}}
 \ctr@ld@f\def\f@gfill{\immediate\write\fwf@g{\fillc@md}}
 \ctr@ld@f\def\f@gnewpath{}
 \ctr@ld@f\def\f@gstroke{\immediate\write\fwf@g{S}}
\else
 \let\figinsertE=\figinsert
 \let\newGr@FN=\newGr@FNDVI
 \ctr@ld@f\def\c@mcurveto{curveto}
 \ctr@ld@f\def\c@mfill{fill}
 \ctr@ld@f\def\c@mgsave{gsave}
 \ctr@ld@f\def\c@mgrestore{grestore}
 \ctr@ld@f\def\c@mlineto{lineto}
 \ctr@ld@f\def\c@mmoveto{moveto}
 \ctr@ld@f\def\c@msetgray{setgray}          \ctr@ld@f\def\c@msetgrayStroke{}
 \ctr@ld@f\def\c@msetcmykcolor{setcmykcolor}\ctr@ld@f\def\c@msetcmykcolorStroke{}
 \ctr@ld@f\def\c@msetrgbcolor{setrgbcolor}  \ctr@ld@f\def\c@msetrgbcolorStroke{}
 \ctr@ld@f\def\d@fprimarC@lor{\CUR@color\space\CUR@colorc@md}
 \ctr@ld@f\def\c@msetdash{setdash}
 \ctr@ld@f\def\c@msetlinejoin{setlinejoin}
 \ctr@ld@f\def\c@msetlinewidth{setlinewidth}
 \ctr@ld@f\def\f@gclosestroke{\immediate\write\fwf@g{closepath\space stroke}}
 \ctr@ld@f\def\f@gfill{\immediate\write\fwf@g{\fillc@md}}
 \ctr@ld@f\def\f@gnewpath{\immediate\write\fwf@g{newpath}}
 \ctr@ld@f\def\f@gstroke{\immediate\write\fwf@g{stroke}}
\fi}
\ctr@ld@f\def\c@pypsfile#1#2{\c@pyfil@{\immediate\write#1}{#2}}
\ctr@ld@f\def\Figinclud@PDF#1#2{\openin\frf@g=#1\pdfliteral{q #2 0 0 #2 0 0 cm}%
    \c@pyfil@{\pdfliteral}{\frf@g}\pdfliteral{Q}\closein\frf@g}
\ctr@ln@w{newif}\ifmored@ta
\ctr@ln@m\bl@nkline
\ctr@ld@f\def\c@pyfil@#1#2{\def\bl@nkline{\par}{\catcode`\%=12
    \loop\ifeof#2\mored@tafalse\else\mored@tatrue\immediate\read#2 to\tr@c
    \ifx\tr@c\bl@nkline\else#1{\tr@c}\fi\fi\ifmored@ta\repeat}}
\ctr@ld@f\def\keln@mun#1#2|{\def\l@debut{#1}\def\l@suite{#2}}
\ctr@ld@f\def\keln@mde#1#2#3|{\def\l@debut{#1#2}\def\l@suite{#3}}
\ctr@ld@f\def\keln@mtr#1#2#3#4|{\def\l@debut{#1#2#3}\def\l@suite{#4}}
\ctr@ld@f\def\keln@mqu#1#2#3#4#5|{\def\l@debut{#1#2#3#4}\def\l@suite{#5}}
\ctr@ld@f\let\@psffilein=\frf@g 
\ctr@ln@w{newif}\if@psffileok    
\ctr@ln@w{newif}\if@psfbbfound   
\ctr@ln@w{newif}\if@psfverbose   
\@psfverbosetrue
\ctr@ln@m\@psfllx \ctr@ln@m\@psflly
\ctr@ln@m\@psfurx \ctr@ln@m\@psfury
\ctr@ln@m\resetcolonc@tcode
\ctr@ld@f\def\@psfgetbb#1{\global\@psfbbfoundfalse%
\global\def\@psfllx{0}\global\def\@psflly{0}%
\global\def\@psfurx{30}\global\def\@psfury{30}%
\openin\@psffilein=#1\relax
\ifeof\@psffilein\errmessage{I couldn't open #1, will ignore it}\else
   \edef\resetcolonc@tcode{\catcode`\noexpand\:\the\catcode`\:\relax}%
   {\@psffileoktrue \chardef\other=12
    \def\do##1{\catcode`##1=\other}\dospecials \catcode`\ =10 \resetcolonc@tcode
    \loop
       \read\@psffilein to \@psffileline
       \ifeof\@psffilein\@psffileokfalse\else
          \expandafter\@psfaux\@psffileline:. \\%
       \fi
   \if@psffileok\repeat
   \if@psfbbfound\else
    \if@psfverbose\message{No bounding box comment in #1; using defaults}\fi\fi
   }\closein\@psffilein\fi}%
\ctr@ln@m\@psfbblit
\ctr@ln@m\@psfpercent
{\catcode`\%=12 \global\let\@psfpercent=
\ctr@ln@m\@psfaux
\long\def\@psfaux#1#2:#3\\{\ifx#1\@psfpercent
   \def\testit{#2}\ifx\testit\@psfbblit
      \@psfgrab #3 . . . \\%
      \@psffileokfalse
      \global\@psfbbfoundtrue
   \fi\else\ifx#1\par\else\@psffileokfalse\fi\fi}%
\ctr@ld@f\def\@psfempty{}%
\ctr@ld@f\def\@psfgrab #1 #2 #3 #4 #5\\{%
\global\def\@psfllx{#1}\ifx\@psfllx\@psfempty
      \@psfgrab #2 #3 #4 #5 .\\\else
   \global\def\@psflly{#2}%
   \global\def\@psfurx{#3}\global\def\@psfury{#4}\fi}%
\ctr@ld@f\def\PSwrit@cmd#1#2#3{{\Figg@tXY{#1}\c@lprojSP\b@undb@x{\v@lX}{\v@lY}%
    \v@lX=\ptT@ptps\v@lX\v@lY=\ptT@ptps\v@lY%
    \immediate\write#3{\repdecn@mb{\v@lX}\space\repdecn@mb{\v@lY}\space#2}}}
\ctr@ld@f\def\PSwrit@cmdS#1#2#3#4#5{{\Figg@tXY{#1}\c@lprojSP\b@undb@x{\v@lX}{\v@lY}%
    \global\result@t=\v@lX\global\result@@t=\v@lY%
    \v@lX=\ptT@ptps\v@lX\v@lY=\ptT@ptps\v@lY%
    \immediate\write#3{\repdecn@mb{\v@lX}\space\repdecn@mb{\v@lY}\space#2}}%
    \edef#4{\the\result@t}\edef#5{\the\result@@t}}
\ctr@ld@f\def\update@ttr#1#2#3{\Figdisc@rdLTS{#3}{\n@mref}%
    \ifx\n@mref\D@FTref#2{#1}\else#2{#3}\fi}
\ctr@ld@f\def\D@FTref{default}
\ctr@ld@f\def\W@rnmesAttr#1#2{%
    \immediate\write16{*** Unknown attribute: \BS@ #1(..., #2=...)}}
\ctr@ld@f\def\W@rnmeskwd#1#2{%
    \immediate\write16{*** Unknown keyword #2 in \BS@ #1}}
\ctr@ld@f\def\W@rnmesIgn#1{\immediate\write16{*** \BS@ #1 is ignored inside a
     \BS@ figdrawbegin-\BS@ figdrawend block.}}
\ctr@ld@f\def\Psset@lti#1=#2|{\keln@mtr#1|%
    \def\n@mref{blc}\ifx\l@debut\n@mref\update@ttr\D@FTref\P@setblcolor{#2}\else
    \def\n@mref{bld}\ifx\l@debut\n@mref\update@ttr\D@FTref\P@setbldash{#2}\else
    \def\n@mref{blw}\ifx\l@debut\n@mref\update@ttr\D@FTref\P@setblwidth{#2}\else
    \def\n@mref{sqc}\ifx\l@debut\n@mref\update@ttr\D@FTref\P@setsqcolor{#2}\else
    \def\n@mref{sqd}\ifx\l@debut\n@mref\update@ttr\D@FTref\P@setsqdash{#2}\else
    \def\n@mref{sqw}\ifx\l@debut\n@mref\update@ttr\D@FTref\P@setsqwidth{#2}\else
    \W@rnmesAttr{figset altitude}{#1}\fi\fi\fi\fi\fi\fi}
\ctr@ln@m\DDV@blcolor
\ctr@ld@f\def\P@setblcolor#1{\edef\DDV@blcolor{#1}}
\ctr@ln@m\DDV@bldash
\ctr@ld@f\def\P@setbldash#1{\edef\DDV@bldash{#1}}
\ctr@ln@m\DDV@blwidth
\ctr@ld@f\def\P@setblwidth#1{\edef\DDV@blwidth{#1}}
\ctr@ln@m\DDV@sqcolor
\ctr@ld@f\def\P@setsqcolor#1{\edef\DDV@sqcolor{#1}}
\ctr@ln@m\DDV@sqdash
\ctr@ld@f\def\P@setsqdash#1{\edef\DDV@sqdash{#1}}
\ctr@ln@m\DDV@sqwidth
\ctr@ld@f\def\P@setsqwidth#1{\edef\DDV@sqwidth{#1}}
\ctr@ld@f\def\figdrawaltitude#1[#2,#3,#4]{{\ifCUR@PS\ifGR@cri%
    \PSc@mment{altitude Square Dim=#1, Triangle=[#2 / #3,#4]}%
    \s@uvc@ntr@l\et@tpsaltitude\resetc@ntr@l{2}\figptorthoprojline-5:=#2/#3,#4/%
    \figvectP -1[#3,#4]\n@rminf{\v@leur}{-1}\vecunit@{-3}{-1}%
    \figvectP -1[-5,#3]\n@rminf{\v@lmin}{-1}\figvectP -2[-5,#4]\n@rminf{\v@lmax}{-2}%
    \ifdim\v@lmin<\v@lmax\s@mme=#3\else\v@lmax=\v@lmin\s@mme=#4\fi%
    \figvectP -4[-5,#2]\vecunit@{-4}{-4}\delt@=#1\unit@%
    \edef\t@ille{\repdecn@mb{\delt@}}\figpttra-1:=-5/\t@ille,-3/%
    \figptstra-3=-5,-1/\t@ille,-4/\figdrawline[#2,-5]%
    \Pss@tspecifSt{color=\DDV@sqcolor,dash=\DDV@sqdash,width=\DDV@sqwidth}%
    \figdrawline[-1,-2,-3]%
    \Psrest@reSt{color=\DDV@sqcolor,dash=\DDV@sqdash,width=\DDV@sqwidth}%
    \ifdim\v@leur<\v@lmax%
    \Pss@tspecifSt{color=\DDV@blcolor,dash=\DDV@bldash,width=\DDV@blwidth}%
    \figdrawline[-5,\the\s@mme]%
    \Psrest@reSt{color=\DDV@blcolor,dash=\DDV@bldash,width=\DDV@blwidth}%
    \fi\PSc@mment{End altitude}\resetc@ntr@l\et@tpsaltitude\fi\fi}}
\ctr@ld@f\def\Ps@rcerc#1;#2(#3,#4){\ellBB@x#1;#2,#2(#3,#4,0)%
    \f@gnewpath{\delt@=#2\unit@\delt@=\ptT@ptps\delt@%
    \BdingB@xfalse%
    \PSwrit@cmd{#1}{\repdecn@mb{\delt@}\space #3\space #4\space arc}{\fwf@g}}}
\ctr@ln@m\figdrawarccirc
\ctr@ld@f\def\Q@arccircDD#1;#2(#3,#4){\ifCUR@PS\ifGR@cri%
    \PSc@mment{arccircDD Center=#1 ; Radius=#2 (Ang1=#3, Ang2=#4)}%
    \iffillm@de\Ps@rcerc#1;#2(#3,#4)%
    \f@gfill%
    \else\Ps@rcerc#1;#2(#3,#4)\f@gstroke\fi%
    \PSc@mment{End arccircDD}\fi\fi}
\ctr@ld@f\def\Q@arccircTD#1,#2,#3;#4(#5,#6){{\ifCUR@PS\ifGR@cri\s@uvc@ntr@l\et@tpsarccircTD%
    \PSc@mment{arccircTD Center=#1,P1=#2,P2=#3 ; Radius=#4 (Ang1=#5, Ang2=#6)}%
    \setc@ntr@l{2}\c@lExtAxes#1,#2,#3(#4)\Q@arcellPATD#1,-4,-5(#5,#6)%
    \PSc@mment{End arccircTD}\resetc@ntr@l\et@tpsarccircTD\fi\fi}}
\ctr@ld@f\def\c@lExtAxes#1,#2,#3(#4){%
    \figvectPTD-5[#1,#2]\vecunit@{-5}{-5}\figvectNTD-4[#1,#2,#3]\vecunit@{-4}{-4}%
    \figvectNVTD-3[-4,-5]\delt@=#4\unit@\edef\r@yon{\repdecn@mb{\delt@}}%
    \figpttra-4:=#1/\r@yon,-5/\figpttra-5:=#1/\r@yon,-3/}
\ctr@ln@m\figdrawarccircP
\ctr@ld@f\def\Q@arccircPDD#1;#2[#3,#4]{{\ifCUR@PS\ifGR@cri\s@uvc@ntr@l\et@tpsarccircPDD%
    \PSc@mment{arccircPDD Center=#1; Radius=#2, [P1=#3, P2=#4]}%
    \Ps@ngleparam#1;#2[#3,#4]\ifdim\v@lmin>\v@lmax\advance\v@lmax\DePI@deg\fi%
    \edef\@ngdeb{\repdecn@mb{\v@lmin}}\edef\@ngfin{\repdecn@mb{\v@lmax}}%
    \figdrawarccirc#1;\r@dius(\@ngdeb,\@ngfin)%
    \PSc@mment{End arccircPDD}\resetc@ntr@l\et@tpsarccircPDD\fi\fi}}
\ctr@ld@f\def\Q@arccircPTD#1;#2[#3,#4,#5]{{\ifCUR@PS\ifGR@cri\s@uvc@ntr@l\et@tpsarccircPTD%
    \PSc@mment{arccircPTD Center=#1; Radius=#2, [P1=#3, P2=#4, P3=#5]}%
    \setc@ntr@l{2}\c@lExtAxes#1,#3,#5(#2)\figdrawarcellPP#1,-4,-5[#3,#4]%
    \PSc@mment{End arccircPTD}\resetc@ntr@l\et@tpsarccircPTD\fi\fi}}
\ctr@ld@f\def\Ps@ngleparam#1;#2[#3,#4]{\setc@ntr@l{2}%
    \figvectPDD-1[#1,#3]\vecunit@{-1}{-1}\Figg@tXY{-1}\arct@n\v@lmin(\v@lX,\v@lY)%
    \figvectPDD-2[#1,#4]\vecunit@{-2}{-2}\Figg@tXY{-2}\arct@n\v@lmax(\v@lX,\v@lY)%
    \v@lmin=\rdT@deg\v@lmin\v@lmax=\rdT@deg\v@lmax%
    \v@leur=#2pt\maxim@m{\mili@u}{-\v@leur}{\v@leur}%
    \edef\r@dius{\repdecn@mb{\mili@u}}}
\ctr@ld@f\def\Ps@rcercBz#1;#2(#3,#4){\Ps@rellBz#1;#2,#2(#3,#4,0)}
\ctr@ld@f\def\Ps@rellBz#1;#2,#3(#4,#5,#6){%
    \ellBB@x#1;#2,#3(#4,#5,#6)\BdingB@xfalse%
    \c@lNbarcs{#4}{#5}\v@leur=#4pt\setc@ntr@l{2}\figptell-13::#1;#2,#3(#4,#6)%
    \f@gnewpath\PSwrit@cmd{-13}{\c@mmoveto}{\fwf@g}%
    \s@mme=\z@\bcl@rellBz#1;#2,#3(#6)\BdingB@xtrue}
\ctr@ld@f\def\bcl@rellBz#1;#2,#3(#4){\relax%
    \ifnum\s@mme<\p@rtent\advance\s@mme\@ne%
    \advance\v@leur\delt@\edef\@ngle{\repdecn@mb\v@leur}\figptell-14::#1;#2,#3(\@ngle,#4)%
    \advance\v@leur\delt@\edef\@ngle{\repdecn@mb\v@leur}\figptell-15::#1;#2,#3(\@ngle,#4)%
    \advance\v@leur\delt@\edef\@ngle{\repdecn@mb\v@leur}\figptell-16::#1;#2,#3(\@ngle,#4)%
    \figptscontrolDD-18[-13,-14,-15,-16]%
    \PSwrit@cmd{-18}{}{\fwf@g}\PSwrit@cmd{-17}{}{\fwf@g}%
    \PSwrit@cmd{-16}{\c@mcurveto}{\fwf@g}%
    \figptcopyDD-13:/-16/\bcl@rellBz#1;#2,#3(#4)\fi}
\ctr@ld@f\def\Ps@rell#1;#2,#3(#4,#5,#6){\ellBB@x#1;#2,#3(#4,#5,#6)%
    \f@gnewpath{\v@lmin=#2\unit@\v@lmin=\ptT@ptps\v@lmin%
    \v@lmax=#3\unit@\v@lmax=\ptT@ptps\v@lmax\BdingB@xfalse%
    \PSwrit@cmd{#1}%
    {#6\space\repdecn@mb{\v@lmin}\space\repdecn@mb{\v@lmax}\space #4\space #5\space ellipse}{\fwf@g}}%
    \global\Use@llipsetrue}
\ctr@ln@m\figdrawarcell
\ctr@ld@f\def\Q@arcellDD#1;#2,#3(#4,#5,#6){{\ifCUR@PS\ifGR@cri%
    \PSc@mment{arcellDD Center=#1 ; XRad=#2, YRad=#3 (Ang1=#4, Ang2=#5, Inclination=#6)}%
    \iffillm@de\Ps@rell#1;#2,#3(#4,#5,#6)%
    \f@gfill%
    \else\Ps@rell#1;#2,#3(#4,#5,#6)\f@gstroke\fi%
    \PSc@mment{End arcellDD}\fi\fi}}
\ctr@ld@f\def\Q@arcellTD#1;#2,#3(#4,#5,#6){{\ifCUR@PS\ifGR@cri\s@uvc@ntr@l\et@tpsarcellTD%
    \PSc@mment{arcellTD Center=#1 ; XRad=#2, YRad=#3 (Ang1=#4, Ang2=#5, Inclination=#6)}%
    \setc@ntr@l{2}\figpttraC -8:=#1/#2,0,0/\figpttraC -7:=#1/0,#3,0/%
    \figvectC -4(0,0,1)\figptsrot -8=-8,-7/#1,#6,-4/\Q@arcellPATD#1,-8,-7(#4,#5)%
    \PSc@mment{End arcellTD}\resetc@ntr@l\et@tpsarcellTD\fi\fi}}
\ctr@ln@m\figdrawarcellPA
\ctr@ld@f\def\Q@arcellPADD#1,#2,#3(#4,#5){{\ifCUR@PS\ifGR@cri\s@uvc@ntr@l\et@tpsarcellPADD%
    \PSc@mment{arcellPADD Center=#1,PtAxis1=#2,PtAxis2=#3 (Ang1=#4, Ang2=#5)}%
    \setc@ntr@l{2}\figvectPDD-1[#1,#2]\vecunit@DD{-1}{-1}\v@lX=\ptT@unit@\result@t%
    \edef\XR@d{\repdecn@mb{\v@lX}}\Figg@tXY{-1}\arct@n\v@lmin(\v@lX,\v@lY)%
    \v@lmin=\rdT@deg\v@lmin\edef\Inclin@{\repdecn@mb{\v@lmin}}%
    \figgetdist\YR@d[#1,#3]\Q@arcellDD#1;\XR@d,\YR@d(#4,#5,\Inclin@)%
    \PSc@mment{End arcellPADD}\resetc@ntr@l\et@tpsarcellPADD\fi\fi}}
\ctr@ld@f\def\Q@arcellPATD#1,#2,#3(#4,#5){{\ifCUR@PS\ifGR@cri\s@uvc@ntr@l\et@tpsarcellPATD%
    \PSc@mment{arcellPATD Center=#1,PtAxis1=#2,PtAxis2=#3 (Ang1=#4, Ang2=#5)}%
    \iffillm@de\Ps@rellPATD#1,#2,#3(#4,#5)%
    \f@gfill%
    \else\Ps@rellPATD#1,#2,#3(#4,#5)\f@gstroke\fi%
    \PSc@mment{End arcellPATD}\resetc@ntr@l\et@tpsarcellPATD\fi\fi}}
\ctr@ld@f\def\Ps@rellPATD#1,#2,#3(#4,#5){\let\c@lprojSP=\relax%
    \setc@ntr@l{2}\figvectPTD-1[#1,#2]\figvectPTD-2[#1,#3]\c@lNbarcs{#4}{#5}%
    \v@leur=#4pt\c@lptellP{#1}{-1}{-2}\Figptpr@j-5:/-3/%
    \f@gnewpath\PSwrit@cmdS{-5}{\c@mmoveto}{\fwf@g}{\X@un}{\Y@un}%
    \edef\C@nt@r{#1}\s@mme=\z@\bcl@rellPATD}
\ctr@ld@f\def\bcl@rellPATD{\relax%
    \ifnum\s@mme<\p@rtent\advance\s@mme\@ne%
    \advance\v@leur\delt@\c@lptellP{\C@nt@r}{-1}{-2}\Figptpr@j-4:/-3/%
    \advance\v@leur\delt@\c@lptellP{\C@nt@r}{-1}{-2}\Figptpr@j-6:/-3/%
    \advance\v@leur\delt@\c@lptellP{\C@nt@r}{-1}{-2}\Figptpr@j-3:/-3/%
    \v@lX=\z@\v@lY=\z@\Figtr@nptDD{-5}{-5}\Figtr@nptDD{2}{-3}%
    \divide\v@lX\@vi\divide\v@lY\@vi%
    \Figtr@nptDD{3}{-4}\Figtr@nptDD{-1.5}{-6}\v@lmin=\v@lX\v@lmax=\v@lY%
    \v@lX=\z@\v@lY=\z@\Figtr@nptDD{2}{-5}\Figtr@nptDD{-5}{-3}%
    \divide\v@lX\@vi\divide\v@lY\@vi\Figtr@nptDD{-1.5}{-4}\Figtr@nptDD{3}{-6}%
    \BdingB@xfalse%
    \Figp@intregDD-4:(\v@lmin,\v@lmax)\PSwrit@cmdS{-4}{}{\fwf@g}{\X@de}{\Y@de}%
    \Figp@intregDD-4:(\v@lX,\v@lY)\PSwrit@cmdS{-4}{}{\fwf@g}{\X@tr}{\Y@tr}%
    \BdingB@xtrue\PSwrit@cmdS{-3}{\c@mcurveto}{\fwf@g}{\X@qu}{\Y@qu}%
    \B@zierBB@x{1}{\Y@un}(\X@un,\X@de,\X@tr,\X@qu)%
    \B@zierBB@x{2}{\X@un}(\Y@un,\Y@de,\Y@tr,\Y@qu)%
    \edef\X@un{\X@qu}\edef\Y@un{\Y@qu}\figptcopyDD-5:/-3/\bcl@rellPATD\fi}
\ctr@ld@f\def\c@lNbarcs#1#2{%
    \delt@=#2pt\advance\delt@-#1pt\maxim@m{\v@lmax}{\delt@}{-\delt@}%
    \v@leur=\v@lmax\divide\v@leur45 \p@rtentiere{\p@rtent}{\v@leur}\advance\p@rtent\@ne%
    \s@mme=\p@rtent\multiply\s@mme\thr@@\divide\delt@\s@mme}
\ctr@ld@f\def\figdrawarcellPP#1,#2,#3[#4,#5]{{\ifCUR@PS\ifGR@cri\s@uvc@ntr@l\et@tpsarcellPP%
    \PSc@mment{arcellPP Center=#1,PtAxis1=#2,PtAxis2=#3 [Point1=#4, Point2=#5]}%
    \setc@ntr@l{2}\figvectP-2[#1,#3]\vecunit@{-2}{-2}\v@lmin=\result@t%
    \invers@{\v@lmax}{\v@lmin}%
    \figvectP-1[#1,#2]\vecunit@{-1}{-1}\v@leur=\result@t%
    \v@leur=\repdecn@mb{\v@lmax}\v@leur\edef\AsB@{\repdecn@mb{\v@leur}}
    \c@lAngle{#1}{#4}{\v@lmin}\edef\@ngdeb{\repdecn@mb{\v@lmin}}%
    \c@lAngle{#1}{#5}{\v@lmax}\ifdim\v@lmin>\v@lmax\advance\v@lmax\DePI@deg\fi%
    \edef\@ngfin{\repdecn@mb{\v@lmax}}\figdrawarcellPA#1,#2,#3(\@ngdeb,\@ngfin)%
    \PSc@mment{End arcellPP}\resetc@ntr@l\et@tpsarcellPP\fi\fi}}
\ctr@ld@f\def\c@lAngle#1#2#3{\figvectP-3[#1,#2]%
    \c@lproscal\delt@[-3,-1]\c@lproscal\v@leur[-3,-2]%
    \v@leur=\AsB@\v@leur\arct@n#3(\delt@,\v@leur)#3=\rdT@deg#3}
\ctr@ln@w{newif}\if@rrowratio\@rrowratiotrue
\ctr@ln@w{newif}\if@rrowhfill
\ctr@ln@w{newif}\if@rrowhout
\ctr@ld@f\def\Psset@rrowhe@d#1=#2|{\keln@mun#1|%
    \def\n@mref{a}\ifx\l@debut\n@mref\update@ttr\D@FTarrowheadangle\Q@s@tarrowheadangle{#2}\else
    \def\n@mref{f}\ifx\l@debut\n@mref\update@ttr\D@FTarrowheadfill\Q@s@tarrowheadfill{#2}\else
    \def\n@mref{l}\ifx\l@debut\n@mref\update@ttr\D@FTarrowheadlength\Q@s@tarrowheadlength{#2}\else
    \def\n@mref{o}\ifx\l@debut\n@mref\update@ttr\D@FTarrowheadout\Q@s@tarrowheadout{#2}\else
    \def\n@mref{r}\ifx\l@debut\n@mref\update@ttr\D@FTarrowheadratio\Q@s@tarrowheadratio{#2}\else
    \W@rnmesAttr{figset arrowhead}{#1}\fi\fi\fi\fi\fi}
\ctr@ln@m\@rrowheadangle
\ctr@ln@m\C@AHANG \ctr@ln@m\S@AHANG \ctr@ln@m\UNSS@N
\ctr@ld@f\def\Q@s@tarrowheadangle#1{\edef\@rrowheadangle{#1}{\c@ssin{\C@}{\S@}{#1}%
    \xdef\C@AHANG{\C@}\xdef\S@AHANG{\S@}\v@lmax=\S@ pt%
    \invers@{\v@leur}{\v@lmax}\maxim@m{\v@leur}{\v@leur}{-\v@leur}%
    \xdef\UNSS@N{\the\v@leur}}}
\ctr@ld@f\def\Q@s@tarrowheadfill#1{\expandafter\set@rrowhfill#1:}
\ctr@ld@f\def\set@rrowhfill#1#2:{\if#1n\@rrowhfillfalse\else\@rrowhfilltrue\fi}
\ctr@ld@f\def\Q@s@tarrowheadout#1{\expandafter\set@rrowhout#1:}
\ctr@ld@f\def\set@rrowhout#1#2:{\if#1n\@rrowhoutfalse\else\@rrowhouttrue\fi}
\ctr@ln@m\@rrowheadlength
\ctr@ld@f\def\Q@s@tarrowheadlength#1{\edef\@rrowheadlength{#1}\@rrowratiofalse}
\ctr@ln@m\@rrowheadratio
\ctr@ld@f\def\Q@s@tarrowheadratio#1{\edef\@rrowheadratio{#1}\@rrowratiotrue}
\ctr@ln@m\D@FTarrowheadlength
\ctr@ld@f\def\figresetarrowhead{%
    \Q@s@tarrowheadangle{\D@FTarrowheadangle}%
    \Q@s@tarrowheadfill{\D@FTarrowheadfill}%
    \Q@s@tarrowheadout{\D@FTarrowheadout}%
    \Q@s@tarrowheadratio{\D@FTarrowheadratio}%
    \d@fm@cdim\D@FTarrowheadlength{\D@FTh@rdahlength}
    \Q@s@tarrowheadlength{\D@FTarrowheadlength}}
\ctr@ld@f\def\D@FTarrowheadratio{0.1}
\ctr@ld@f\def\D@FTarrowheadangle{20}
\ctr@ld@f\def\D@FTarrowheadfill{no}
\ctr@ld@f\def\D@FTarrowheadout{no}
\ctr@ld@f\def\D@FTh@rdahlength{8pt}
\ctr@ln@m\figdrawarrow
\ctr@ld@f\def\Q@arrowDD[#1,#2]{{\ifCUR@PS\ifGR@cri\s@uvc@ntr@l\et@tpsarrow%
    \PSc@mment{arrowDD [Pt1,Pt2]=[#1,#2]}\Q@s@tfillmode{no}%
    \Q@arrowheadDD[#1,#2]\setc@ntr@l{2}\figdrawline[#1,-3]%
    \PSc@mment{End arrowDD}\resetc@ntr@l\et@tpsarrow\fi\fi}}
\ctr@ld@f\def\Q@arrowTD[#1,#2]{{\ifCUR@PS\ifGR@cri\s@uvc@ntr@l\et@tpsarrowTD%
    \PSc@mment{arrowTD [Pt1,Pt2]=[#1,#2]}\resetc@ntr@l{2}%
    \Figptpr@j-5:/#1/\Figptpr@j-6:/#2/\let\c@lprojSP=\relax\Q@arrowDD[-5,-6]%
    \PSc@mment{End arrowTD}\resetc@ntr@l\et@tpsarrowTD\fi\fi}}
\ctr@ln@m\figdrawarrowhead
\ctr@ld@f\def\Q@arrowheadDD[#1,#2]{{\ifCUR@PS\ifGR@cri\s@uvc@ntr@l\et@tpsarrowheadDD%
    \if@rrowhfill\def\@hangle{-\@rrowheadangle}\else\def\@hangle{\@rrowheadangle}\fi%
    \if@rrowratio%
    \if@rrowhout\def\@hratio{-\@rrowheadratio}\else\def\@hratio{\@rrowheadratio}\fi%
    \PSc@mment{arrowheadDD Ratio=\@hratio, Angle=\@hangle, [Pt1,Pt2]=[#1,#2]}%
    \Ps@rrowhead\@hratio,\@hangle[#1,#2]%
    \else%
    \if@rrowhout\def\@hlength{-\@rrowheadlength}\else\def\@hlength{\@rrowheadlength}\fi%
    \PSc@mment{arrowheadDD Length=\@hlength, Angle=\@hangle, [Pt1,Pt2]=[#1,#2]}%
    \Ps@rrowheadfd\@hlength,\@hangle[#1,#2]%
    \fi%
    \PSc@mment{End arrowheadDD}\resetc@ntr@l\et@tpsarrowheadDD\fi\fi}}
\ctr@ld@f\def\Q@arrowheadTD[#1,#2]{{\ifCUR@PS\ifGR@cri\s@uvc@ntr@l\et@tpsarrowheadTD%
    \PSc@mment{arrowheadTD [Pt1,Pt2]=[#1,#2]}\resetc@ntr@l{2}%
    \Figptpr@j-5:/#1/\Figptpr@j-6:/#2/\let\c@lprojSP=\relax\Q@arrowheadDD[-5,-6]%
    \PSc@mment{End arrowheadTD}\resetc@ntr@l\et@tpsarrowheadTD\fi\fi}}
\ctr@ld@f\def\Ps@rrowhead#1,#2[#3,#4]{\v@leur=#1\p@\maxim@m{\v@leur}{\v@leur}{-\v@leur}%
    \ifdim\v@leur>\Cepsil@n{
    \PSc@mment{@rrowhead Ratio=#1, Angle=#2, [Pt1,Pt2]=[#3,#4]}\v@leur=\UNSS@N%
    \v@leur=\CUR@width\v@leur\v@leur=\ptpsT@pt\v@leur\delt@=.5\v@leur
    \setc@ntr@l{2}\figvectPDD-3[#4,#3]%
    \Figg@tXY{-3}\v@lX=#1\v@lX\v@lY=#1\v@lY\Figv@ctCreg-3(\v@lX,\v@lY)%
    \vecunit@{-4}{-3}\mili@u=\result@t%
    \ifdim#2pt>\z@\v@lXa=-\C@AHANG\delt@%
     \edef\c@ef{\repdecn@mb{\v@lXa}}\figpttraDD-3:=-3/\c@ef,-4/\fi%
    \edef\c@ef{\repdecn@mb{\delt@}}%
    \v@lXa=\mili@u\v@lXa=\C@AHANG\v@lXa%
    \v@lYa=\ptpsT@pt\p@\v@lYa=\CUR@width\v@lYa\v@lYa=\sDcc@ngle\v@lYa%
    \advance\v@lXa-\v@lYa\gdef\sDcc@ngle{0}%
    \ifdim\v@lXa>\v@leur\edef\c@efendpt{\repdecn@mb{\v@leur}}%
    \else\edef\c@efendpt{\repdecn@mb{\v@lXa}}\fi%
    \Figg@tXY{-3}\v@lmin=\v@lX\v@lmax=\v@lY%
    \v@lXa=\C@AHANG\v@lmin\v@lYa=\S@AHANG\v@lmax\advance\v@lXa\v@lYa%
    \v@lYa=-\S@AHANG\v@lmin\v@lX=\C@AHANG\v@lmax\advance\v@lYa\v@lX%
    \setc@ntr@l{1}\Figg@tXY{#4}\advance\v@lX\v@lXa\advance\v@lY\v@lYa%
    \setc@ntr@l{2}\Figp@intregDD-2:(\v@lX,\v@lY)%
    \v@lXa=\C@AHANG\v@lmin\v@lYa=-\S@AHANG\v@lmax\advance\v@lXa\v@lYa%
    \v@lYa=\S@AHANG\v@lmin\v@lX=\C@AHANG\v@lmax\advance\v@lYa\v@lX%
    \setc@ntr@l{1}\Figg@tXY{#4}\advance\v@lX\v@lXa\advance\v@lY\v@lYa%
    \setc@ntr@l{2}\Figp@intregDD-1:(\v@lX,\v@lY)%
    \ifdim#2pt<\z@\fillm@detrue\figdrawline[-2,#4,-1]
    \else\figptstraDD-3=#4,-2,-1/\c@ef,-4/\s@uvdash{\typ@dash}\Q@s@tdash{\D@FTdash}%
    \figdrawline[-2,-3,-1]\Q@s@tdash{\typ@dash}\fi
    \ifdim#1pt>\z@\figpttraDD-3:=#4/\c@efendpt,-4/\else\figptcopyDD-3:/#4/\fi%
    \PSc@mment{End @rrowhead}}\fi}
\ctr@ld@f\def\sDcc@ngle{0}
\ctr@ld@f\def\Ps@rrowheadfd#1,#2[#3,#4]{{%
    \PSc@mment{@rrowheadfd Length=#1, Angle=#2, [Pt1,Pt2]=[#3,#4]}%
    \setc@ntr@l{2}\figvectPDD-1[#3,#4]\n@rmeucDD{\v@leur}{-1}\v@leur=\ptT@unit@\v@leur%
    \invers@{\v@leur}{\v@leur}\v@leur=#1\v@leur\edef\R@tio{\repdecn@mb{\v@leur}}%
    \Ps@rrowhead\R@tio,#2[#3,#4]\PSc@mment{End @rrowheadfd}}}
\ctr@ln@m\figdrawarrowBezier
\ctr@ld@f\def\Q@arrowBezierDD[#1,#2,#3,#4]{{\ifCUR@PS\ifGR@cri\s@uvc@ntr@l\et@tpsarrowBezierDD%
    \PSc@mment{arrowBezierDD Control points=#1,#2,#3,#4}\setc@ntr@l{2}%
    \if@rrowratio\c@larclengthDD\v@leur,10[#1,#2,#3,#4]\else\v@leur=\z@\fi%
    \Ps@rrowB@zDD\v@leur[#1,#2,#3,#4]%
    \PSc@mment{End arrowBezierDD}\resetc@ntr@l\et@tpsarrowBezierDD\fi\fi}}
\ctr@ld@f\def\Q@arrowBezierTD[#1,#2,#3,#4]{{\ifCUR@PS\ifGR@cri\s@uvc@ntr@l\et@tpsarrowBezierTD%
    \PSc@mment{arrowBezierTD Control points=#1,#2,#3,#4}\resetc@ntr@l{2}%
    \Figptpr@j-7:/#1/\Figptpr@j-8:/#2/\Figptpr@j-9:/#3/\Figptpr@j-10:/#4/%
    \let\c@lprojSP=\relax\ifnum\CUR@proj<\tw@\Q@arrowBezierDD[-7,-8,-9,-10]%
    \else\f@gnewpath\PSwrit@cmd{-7}{\c@mmoveto}{\fwf@g}%
    \if@rrowratio\c@larclengthDD\mili@u,10[-7,-8,-9,-10]\else\mili@u=\z@\fi%
    \p@rtent=\NBz@rcs\advance\p@rtent\m@ne\subB@zierTD\p@rtent[#1,#2,#3,#4]%
    \f@gstroke%
    \advance\v@lmin\p@rtent\delt@
    \v@leur=\v@lmin\advance\v@leur0.33333 \delt@\edef\unti@rs{\repdecn@mb{\v@leur}}%
    \v@leur=\v@lmin\advance\v@leur0.66666 \delt@\edef\deti@rs{\repdecn@mb{\v@leur}}%
    \figptcopyDD-8:/-10/\c@lsubBzarc\unti@rs,\deti@rs[#1,#2,#3,#4]%
    \figptcopyDD-8:/-4/\figptcopyDD-9:/-3/\Ps@rrowB@zDD\mili@u[-7,-8,-9,-10]\fi%
    \PSc@mment{End arrowBezierTD}\resetc@ntr@l\et@tpsarrowBezierTD\fi\fi}}
\ctr@ld@f\def\c@larclengthDD#1,#2[#3,#4,#5,#6]{{\p@rtent=#2\figptcopyDD-5:/#3/%
    \delt@=\p@\divide\delt@\p@rtent\c@rre=\z@\v@leur=\z@\s@mme=\z@%
    \loop\ifnum\s@mme<\p@rtent\advance\s@mme\@ne\advance\v@leur\delt@%
    \edef\T@{\repdecn@mb{\v@leur}}\figptBezierDD-6::\T@[#3,#4,#5,#6]%
    \figvectPDD-1[-5,-6]\n@rmeucDD{\mili@u}{-1}\advance\c@rre\mili@u%
    \figptcopyDD-5:/-6/\repeat\global\result@t=\ptT@unit@\c@rre}#1=\result@t}
\ctr@ld@f\def\Ps@rrowB@zDD#1[#2,#3,#4,#5]{{\Q@s@tfillmode{no}%
    \if@rrowratio\delt@=\@rrowheadratio#1\else\delt@=\@rrowheadlength pt\fi%
    \v@leur=\C@AHANG\delt@\edef\R@dius{\repdecn@mb{\v@leur}}%
    \FigptintercircB@zDD-5::0,\R@dius[#5,#4,#3,#2]%
    \Q@s@tarrowheadlength{\repdecn@mb{\delt@}}\Q@arrowheadDD[-5,#5]%
    \let\n@rmeuc=\n@rmeucDD\figgetdist\R@dius[#5,-3]%
    \FigptintercircB@zDD-6::0,\R@dius[#5,#4,#3,#2]%
    \figptBezierDD-5::0.33333[#5,#4,#3,#2]\figptBezierDD-3::0.66666[#5,#4,#3,#2]%
    \figptscontrolDD-5[-6,-5,-3,#2]\Q@BezierDD1[-6,-5,-4,#2]}}
\ctr@ln@m\figdrawarrowcirc
\ctr@ld@f\def\Q@arrowcircDD#1;#2(#3,#4){{\ifCUR@PS\ifGR@cri\s@uvc@ntr@l\et@tpsarrowcircDD%
    \PSc@mment{arrowcircDD Center=#1 ; Radius=#2 (Ang1=#3,Ang2=#4)}%
    \Q@s@tfillmode{no}\Pscirc@rrowhead#1;#2(#3,#4)%
    \setc@ntr@l{2}\figvectPDD -4[#1,-3]\vecunit@{-4}{-4}%
    \Figg@tXY{-4}\arct@n\v@lmin(\v@lX,\v@lY)%
    \v@lmin=\rdT@deg\v@lmin\v@leur=#4pt\advance\v@leur-\v@lmin%
    \maxim@m{\v@leur}{\v@leur}{-\v@leur}%
    \ifdim\v@leur>\DemiPI@deg\relax\ifdim\v@lmin<#4pt\advance\v@lmin\DePI@deg%
    \else\advance\v@lmin-\DePI@deg\fi\fi\edef\ar@ngle{\repdecn@mb{\v@lmin}}%
    \ifdim#3pt<#4pt\figdrawarccirc#1;#2(#3,\ar@ngle)\else\figdrawarccirc#1;#2(\ar@ngle,#3)\fi%
    \PSc@mment{End arrowcircDD}\resetc@ntr@l\et@tpsarrowcircDD\fi\fi}}
\ctr@ld@f\def\Q@arrowcircTD#1,#2,#3;#4(#5,#6){{\ifCUR@PS\ifGR@cri\s@uvc@ntr@l\et@tpsarrowcircTD%
    \PSc@mment{arrowcircTD Center=#1,P1=#2,P2=#3 ; Radius=#4 (Ang1=#5, Ang2=#6)}%
    \resetc@ntr@l{2}\c@lExtAxes#1,#2,#3(#4)\let\c@lprojSP=\relax%
    \figvectPTD-11[#1,-4]\figvectPTD-12[#1,-5]\c@lNbarcs{#5}{#6}%
    \if@rrowratio\v@lmax=\degT@rd\v@lmax\edef\D@lpha{\repdecn@mb{\v@lmax}}\fi%
    \advance\p@rtent\m@ne\mili@u=\z@%
    \v@leur=#5pt\c@lptellP{#1}{-11}{-12}\Figptpr@j-9:/-3/%
    \f@gnewpath\PSwrit@cmdS{-9}{\c@mmoveto}{\fwf@g}{\X@un}{\Y@un}%
    \edef\C@nt@r{#1}\s@mme=\z@\bcl@rcircTD\f@gstroke%
    \advance\v@leur\delt@\c@lptellP{#1}{-11}{-12}\Figptpr@j-5:/-3/%
    \advance\v@leur\delt@\c@lptellP{#1}{-11}{-12}\Figptpr@j-6:/-3/%
    \advance\v@leur\delt@\c@lptellP{#1}{-11}{-12}\Figptpr@j-10:/-3/%
    \figptscontrolDD-8[-9,-5,-6,-10]%
    \if@rrowratio\c@lcurvradDD0.5[-9,-8,-7,-10]\advance\mili@u\result@t%
    \maxim@m{\mili@u}{\mili@u}{-\mili@u}\mili@u=\ptT@unit@\mili@u%
    \mili@u=\D@lpha\mili@u\advance\p@rtent\@ne\divide\mili@u\p@rtent\fi%
    \Ps@rrowB@zDD\mili@u[-9,-8,-7,-10]%
    \PSc@mment{End arrowcircTD}\resetc@ntr@l\et@tpsarrowcircTD\fi\fi}}
\ctr@ld@f\def\bcl@rcircTD{\relax%
    \ifnum\s@mme<\p@rtent\advance\s@mme\@ne%
    \advance\v@leur\delt@\c@lptellP{\C@nt@r}{-11}{-12}\Figptpr@j-5:/-3/%
    \advance\v@leur\delt@\c@lptellP{\C@nt@r}{-11}{-12}\Figptpr@j-6:/-3/%
    \advance\v@leur\delt@\c@lptellP{\C@nt@r}{-11}{-12}\Figptpr@j-10:/-3/%
    \figptscontrolDD-8[-9,-5,-6,-10]\BdingB@xfalse%
    \PSwrit@cmdS{-8}{}{\fwf@g}{\X@de}{\Y@de}\PSwrit@cmdS{-7}{}{\fwf@g}{\X@tr}{\Y@tr}%
    \BdingB@xtrue\PSwrit@cmdS{-10}{\c@mcurveto}{\fwf@g}{\X@qu}{\Y@qu}%
    \if@rrowratio\c@lcurvradDD0.5[-9,-8,-7,-10]\advance\mili@u\result@t\fi%
    \B@zierBB@x{1}{\Y@un}(\X@un,\X@de,\X@tr,\X@qu)%
    \B@zierBB@x{2}{\X@un}(\Y@un,\Y@de,\Y@tr,\Y@qu)%
    \edef\X@un{\X@qu}\edef\Y@un{\Y@qu}\figptcopyDD-9:/-10/\bcl@rcircTD\fi}
\ctr@ld@f\def\Pscirc@rrowhead#1;#2(#3,#4){{%
    \PSc@mment{circ@rrowhead Center=#1 ; Radius=#2 (Ang1=#3,Ang2=#4)}%
    \v@leur=#2\unit@\edef\s@glen{\repdecn@mb{\v@leur}}\v@lY=\z@\v@lX=\v@leur%
    \resetc@ntr@l{2}\Figv@ctCreg-3(\v@lX,\v@lY)\figpttraDD-5:=#1/1,-3/%
    \figptrotDD-5:=-5/#1,#4/%
    \figvectPDD-3[#1,-5]\Figg@tXY{-3}\v@leur=\v@lX%
    \ifdim#3pt<#4pt\v@lX=\v@lY\v@lY=-\v@leur\else\v@lX=-\v@lY\v@lY=\v@leur\fi%
    \Figv@ctCreg-3(\v@lX,\v@lY)\vecunit@{-3}{-3}%
    \if@rrowratio\v@leur=#4pt\advance\v@leur-#3pt\maxim@m{\mili@u}{-\v@leur}{\v@leur}%
    \mili@u=\degT@rd\mili@u\v@leur=\s@glen\mili@u\edef\s@glen{\repdecn@mb{\v@leur}}%
    \mili@u=#2\mili@u\mili@u=\@rrowheadratio\mili@u\else\mili@u=\@rrowheadlength pt\fi%
    \figpttraDD-6:=-5/\s@glen,-3/\v@leur=#2pt\v@leur=2\v@leur%
    \invers@{\v@leur}{\v@leur}\c@rre=\repdecn@mb{\v@leur}\mili@u
    \mili@u=\c@rre\mili@u=\repdecn@mb{\c@rre}\mili@u%
    \v@leur=\p@\advance\v@leur-\mili@u
    \invers@{\mili@u}{2\v@leur}\delt@=\c@rre\delt@=\repdecn@mb{\mili@u}\delt@%
    \xdef\sDcc@ngle{\repdecn@mb{\delt@}}
    \sqrt@{\mili@u}{\v@leur}\arct@n\v@leur(\mili@u,\c@rre)%
    \v@leur=\rdT@deg\v@leur
    \ifdim#3pt<#4pt\v@leur=-\v@leur\fi%
    \if@rrowhout\v@leur=-\v@leur\fi\edef\cor@ngle{\repdecn@mb{\v@leur}}%
    \figptrotDD-6:=-6/-5,\cor@ngle/\Q@arrowheadDD[-6,-5]%
    \PSc@mment{End circ@rrowhead}}}
\ctr@ln@m\figdrawarrowcircP
\ctr@ld@f\def\Q@arrowcircPDD#1;#2[#3,#4]{{\ifCUR@PS\ifGR@cri%
    \PSc@mment{arrowcircPDD Center=#1; Radius=#2, [P1=#3,P2=#4]}%
    \s@uvc@ntr@l\et@tpsarrowcircPDD\Ps@ngleparam#1;#2[#3,#4]%
    \ifdim\v@leur>\z@\ifdim\v@lmin>\v@lmax\advance\v@lmax\DePI@deg\fi%
    \else\ifdim\v@lmin<\v@lmax\advance\v@lmin\DePI@deg\fi\fi%
    \edef\@ngdeb{\repdecn@mb{\v@lmin}}\edef\@ngfin{\repdecn@mb{\v@lmax}}%
    \figdrawarrowcirc#1;\r@dius(\@ngdeb,\@ngfin)%
    \PSc@mment{End arrowcircPDD}\resetc@ntr@l\et@tpsarrowcircPDD\fi\fi}}
\ctr@ld@f\def\Q@arrowcircPTD#1;#2[#3,#4,#5]{{\ifCUR@PS\ifGR@cri\s@uvc@ntr@l\et@tpsarrowcircPTD%
    \PSc@mment{arrowcircPTD Center=#1; Radius=#2, [P1=#3,P2=#4,P3=#5]}%
    \figgetangleTD\@ngfin[#1,#3,#4,#5]\v@leur=#2pt%
    \maxim@m{\mili@u}{-\v@leur}{\v@leur}\edef\r@dius{\repdecn@mb{\mili@u}}%
    \ifdim\v@leur<\z@\v@lmax=\@ngfin pt\advance\v@lmax-\DePI@deg%
    \edef\@ngfin{\repdecn@mb{\v@lmax}}\fi\Q@arrowcircTD#1,#3,#5;\r@dius(0,\@ngfin)%
    \PSc@mment{End arrowcircPTD}\resetc@ntr@l\et@tpsarrowcircPTD\fi\fi}}
\ctr@ld@f\def\figdrawaxes#1(#2){{\ifCUR@PS\ifGR@cri\s@uvc@ntr@l\et@tpsaxes%
    \PSc@mment{axes Origin=#1 Range=(#2)}\an@lys@xes#2,:\resetc@ntr@l{2}%
    \ifx\t@xt@\empty\ifTr@isDim\Q@@xes#1(0,#2,0,#2,0,#2)\else\Q@@xes#1(0,#2,0,#2)\fi%
    \else\Q@@xes#1(#2)\fi\PSc@mment{End axes}\resetc@ntr@l\et@tpsaxes\fi\fi}}
\ctr@ld@f\def\an@lys@xes#1,#2:{\def\t@xt@{#2}}
\ctr@ln@m\Q@@xes
\ctr@ld@f\def\Q@@xesDD#1(#2,#3,#4,#5){%
    \figpttraC-5:=#1/#2,0/\figpttraC-6:=#1/#3,0/\Q@arrowDD[-5,-6]%
    \figpttraC-5:=#1/0,#4/\figpttraC-6:=#1/0,#5/\Q@arrowDD[-5,-6]}
\ctr@ld@f\def\Q@@xesTD#1(#2,#3,#4,#5,#6,#7){%
    \figpttraC-7:=#1/#2,0,0/\figpttraC-8:=#1/#3,0,0/\Q@arrowTD[-7,-8]%
    \figpttraC-7:=#1/0,#4,0/\figpttraC-8:=#1/0,#5,0/\Q@arrowTD[-7,-8]%
    \figpttraC-7:=#1/0,0,#6/\figpttraC-8:=#1/0,0,#7/\Q@arrowTD[-7,-8]}
\ctr@ln@m\newGr@FN
\ctr@ld@f\def\newGr@FNPDF#1{\s@mme=\Gr@FNb\advance\s@mme\@ne\xdef\Gr@FNb{\number\s@mme}}
\ctr@ld@f\def\newGr@FNDVI#1{\newGr@FNPDF{}\xdef#1{\jobname GI\Gr@FNb.anx}}
\ctr@ld@f\def\figdrawbegin#1{\newGr@FN\DefGIfilen@me\gdef\@utoFN{0}%
    \def\t@xt@{#1}\relax\ifx\t@xt@\empty\GRupdatem@detrue%
    \gdef\@utoFN{1}\Psb@ginfig\DefGIfilen@me\else\expandafter\Psb@ginfigNu@#1 :\fi}
\ctr@ld@f\def\Psb@ginfigNu@#1 #2:{\def\t@xt@{#1}\relax\ifx\t@xt@\empty\def\t@xt@{#2}%
    \ifx\t@xt@\empty\GRupdatem@detrue\gdef\@utoFN{1}\Psb@ginfig\DefGIfilen@me%
    \else\Psb@ginfigNu@#2:\fi\else\Psb@ginfig{#1}\fi}
\ctr@ln@m\PSfilen@me \ctr@ln@m\auxfilen@me
\ctr@ld@f\def\Psb@ginfig#1{\ifCUR@PS\else%
    \edef\PSfilen@me{#1}\edef\auxfilen@me{\jobname.anx}%
    \ifGRupdatem@de\GR@critrue\else\openin\frf@g=\PSfilen@me\relax%
    \ifeof\frf@g\GR@critrue\else\GR@crifalse\fi\closein\frf@g\fi%
    \CUR@PStrue\c@ldefproj\expandafter\setupd@te\D@FTupdate:%
    \ifGR@cri\initb@undb@x%
    \immediate\openout\fwf@g=\auxfilen@me\initpss@ttings\fi%
    \fi}
\ctr@ld@f\def\Gr@FNb{0}
\ctr@ld@f\def\figforTeXFileno{\Gr@FNb}
\ctr@ld@f\def\figforTeXFigno{0 }
\ctr@ld@f\def\figforTeXnextFigno{1 }
\ctr@ld@f\edef\DefGIfilen@me{\jobname GI.anx}
\ctr@ld@f\def\initpss@ttings{\figreset{altitude,arrowhead,curve,general,flowchart,mesh,trimesh}%
    \Use@llipsefalse}
\ctr@ld@f\def\B@zierBB@x#1#2(#3,#4,#5,#6){{\c@rre=\t@n\epsil@n
    \v@lmax=#4\advance\v@lmax-#5\v@lmax=\thr@@\v@lmax\advance\v@lmax#6\advance\v@lmax-#3%
    \mili@u=#4\mili@u=-\tw@\mili@u\advance\mili@u#3\advance\mili@u#5%
    \v@lmin=#4\advance\v@lmin-#3\maxim@m{\v@leur}{-\v@lmax}{\v@lmax}%
    \maxim@m{\delt@}{-\mili@u}{\mili@u}\maxim@m{\v@leur}{\v@leur}{\delt@}%
    \maxim@m{\delt@}{-\v@lmin}{\v@lmin}\maxim@m{\v@leur}{\v@leur}{\delt@}%
    \ifdim\v@leur>\c@rre\invers@{\v@leur}{\v@leur}\edef\Uns@rM@x{\repdecn@mb{\v@leur}}%
    \v@lmax=\Uns@rM@x\v@lmax\mili@u=\Uns@rM@x\mili@u\v@lmin=\Uns@rM@x\v@lmin%
    \maxim@m{\v@leur}{-\v@lmax}{\v@lmax}\ifdim\v@leur<\c@rre%
    \maxim@m{\v@leur}{-\mili@u}{\mili@u}\ifdim\v@leur<\c@rre\else%
    \invers@{\mili@u}{\mili@u}\v@leur=-0.5\v@lmin%
    \v@leur=\repdecn@mb{\mili@u}\v@leur\m@jBBB@x{\v@leur}{#1}{#2}(#3,#4,#5,#6)\fi%
    \else\delt@=\repdecn@mb{\mili@u}\mili@u\v@leur=\repdecn@mb{\v@lmax}\v@lmin%
    \advance\delt@-\v@leur\ifdim\delt@<\z@\else\invers@{\v@lmax}{\v@lmax}%
    \edef\Uns@rAp{\repdecn@mb{\v@lmax}}\sqrt@{\delt@}{\delt@}%
    \v@leur=-\mili@u\advance\v@leur\delt@\v@leur=\Uns@rAp\v@leur%
    \m@jBBB@x{\v@leur}{#1}{#2}(#3,#4,#5,#6)%
    \v@leur=-\mili@u\advance\v@leur-\delt@\v@leur=\Uns@rAp\v@leur%
    \m@jBBB@x{\v@leur}{#1}{#2}(#3,#4,#5,#6)\fi\fi\fi}}
\ctr@ld@f\def\m@jBBB@x#1#2#3(#4,#5,#6,#7){{\relax\ifdim#1>\z@\ifdim#1<\p@%
    \edef\T@{\repdecn@mb{#1}}\v@lX=\p@\advance\v@lX-#1\edef\UNmT@{\repdecn@mb{\v@lX}}%
    \v@lX=#4\v@lY=#5\v@lZ=#6\v@lXa=#7\v@lX=\UNmT@\v@lX\advance\v@lX\T@\v@lY%
    \v@lY=\UNmT@\v@lY\advance\v@lY\T@\v@lZ\v@lZ=\UNmT@\v@lZ\advance\v@lZ\T@\v@lXa%
    \v@lX=\UNmT@\v@lX\advance\v@lX\T@\v@lY\v@lY=\UNmT@\v@lY\advance\v@lY\T@\v@lZ%
    \v@lX=\UNmT@\v@lX\advance\v@lX\T@\v@lY%
    \ifcase#2\or\v@lY=#3\or\v@lY=\v@lX\v@lX=#3\fi\b@undb@x{\v@lX}{\v@lY}\fi\fi}}
\ctr@ld@f\def\PsB@zier#1[#2]{{\f@gnewpath%
    \s@mme=\z@\def\list@num{#2,0}\extrairelepremi@r\p@int\de\list@num%
    \PSwrit@cmdS{\p@int}{\c@mmoveto}{\fwf@g}{\X@un}{\Y@un}\p@rtent=#1\bclB@zier}}
\ctr@ld@f\def\bclB@zier{\relax%
    \ifnum\s@mme<\p@rtent\advance\s@mme\@ne\BdingB@xfalse%
    \extrairelepremi@r\p@int\de\list@num\PSwrit@cmdS{\p@int}{}{\fwf@g}{\X@de}{\Y@de}%
    \extrairelepremi@r\p@int\de\list@num\PSwrit@cmdS{\p@int}{}{\fwf@g}{\X@tr}{\Y@tr}%
    \BdingB@xtrue%
    \extrairelepremi@r\p@int\de\list@num\PSwrit@cmdS{\p@int}{\c@mcurveto}{\fwf@g}{\X@qu}{\Y@qu}%
    \B@zierBB@x{1}{\Y@un}(\X@un,\X@de,\X@tr,\X@qu)%
    \B@zierBB@x{2}{\X@un}(\Y@un,\Y@de,\Y@tr,\Y@qu)%
    \edef\X@un{\X@qu}\edef\Y@un{\Y@qu}\bclB@zier\fi}
\ctr@ln@m\figdrawBezier
\ctr@ld@f\def\Q@BezierDD#1[#2]{\ifCUR@PS\ifGR@cri%
    \PSc@mment{BezierDD N arcs=#1, Control points=#2}%
    \iffillm@de\PsB@zier#1[#2]%
    \f@gfill%
    \else\PsB@zier#1[#2]\f@gstroke\fi%
    \PSc@mment{End BezierDD}\fi\fi}
\ctr@ln@m\et@tpsBezierTD
\ctr@ld@f\def\Q@BezierTD#1[#2]{\ifCUR@PS\ifGR@cri\s@uvc@ntr@l\et@tpsBezierTD%
    \PSc@mment{BezierTD N arcs=#1, Control points=#2}%
    \iffillm@de\PsB@zierTD#1[#2]%
    \f@gfill%
    \else\PsB@zierTD#1[#2]\f@gstroke\fi%
    \PSc@mment{End BezierTD}\resetc@ntr@l\et@tpsBezierTD\fi\fi}
\ctr@ld@f\def\PsB@zierTD#1[#2]{\ifnum\CUR@proj<\tw@\PsB@zier#1[#2]\else\PsB@zier@TD#1[#2]\fi}
\ctr@ld@f\def\PsB@zier@TD#1[#2]{{\f@gnewpath%
    \s@mme=\z@\def\list@num{#2,0}\extrairelepremi@r\p@int\de\list@num%
    \let\c@lprojSP=\relax\setc@ntr@l{2}\Figptpr@j-7:/\p@int/%
    \PSwrit@cmd{-7}{\c@mmoveto}{\fwf@g}%
    \loop\ifnum\s@mme<#1\advance\s@mme\@ne\extrairelepremi@r\p@intun\de\list@num%
    \extrairelepremi@r\p@intde\de\list@num\extrairelepremi@r\p@inttr\de\list@num%
    \subB@zierTD\NBz@rcs[\p@int,\p@intun,\p@intde,\p@inttr]\edef\p@int{\p@inttr}\repeat}}
\ctr@ld@f\def\subB@zierTD#1[#2,#3,#4,#5]{\delt@=\p@\divide\delt@\NBz@rcs\v@lmin=\z@%
    {\Figg@tXY{-7}\edef\X@un{\the\v@lX}\edef\Y@un{\the\v@lY}%
    \s@mme=\z@\loop\ifnum\s@mme<#1\advance\s@mme\@ne%
    \v@leur=\v@lmin\advance\v@leur0.33333 \delt@\edef\unti@rs{\repdecn@mb{\v@leur}}%
    \v@leur=\v@lmin\advance\v@leur0.66666 \delt@\edef\deti@rs{\repdecn@mb{\v@leur}}%
    \advance\v@lmin\delt@\edef\trti@rs{\repdecn@mb{\v@lmin}}%
    \figptBezierTD-8::\trti@rs[#2,#3,#4,#5]\Figptpr@j-8:/-8/%
    \c@lsubBzarc\unti@rs,\deti@rs[#2,#3,#4,#5]\BdingB@xfalse%
    \PSwrit@cmdS{-4}{}{\fwf@g}{\X@de}{\Y@de}\PSwrit@cmdS{-3}{}{\fwf@g}{\X@tr}{\Y@tr}%
    \BdingB@xtrue\PSwrit@cmdS{-8}{\c@mcurveto}{\fwf@g}{\X@qu}{\Y@qu}%
    \B@zierBB@x{1}{\Y@un}(\X@un,\X@de,\X@tr,\X@qu)%
    \B@zierBB@x{2}{\X@un}(\Y@un,\Y@de,\Y@tr,\Y@qu)%
    \edef\X@un{\X@qu}\edef\Y@un{\Y@qu}\figptcopyDD-7:/-8/\repeat}}
\ctr@ld@f\def\NBz@rcs{2}
\ctr@ld@f\def\c@lsubBzarc#1,#2[#3,#4,#5,#6]{\figptBezierTD-5::#1[#3,#4,#5,#6]%
    \figptBezierTD-6::#2[#3,#4,#5,#6]\Figptpr@j-4:/-5/\Figptpr@j-5:/-6/%
    \figptscontrolDD-4[-7,-4,-5,-8]}
\ctr@ln@m\figdrawcirc
\ctr@ld@f\def\Q@circDD#1(#2){\ifCUR@PS\ifGR@cri\PSc@mment{circDD Center=#1 (Radius=#2)}%
    \Q@arccircDD#1;#2(0,360)\PSc@mment{End circDD}\fi\fi}
\ctr@ld@f\def\Q@circTD#1,#2,#3(#4){\ifCUR@PS\ifGR@cri%
    \PSc@mment{circTD Center=#1,P1=#2,P2=#3 (Radius=#4)}%
    \Q@arccircTD#1,#2,#3;#4(0,360)\PSc@mment{End circTD}\fi\fi}
\ctr@ln@m\p@urcent
{\catcode`\%=12\gdef\p@urcent{
\ctr@ld@f\def\PSc@mment#1{\ifGRdebugm@de\immediate\write\fwf@g{\p@urcent\space#1}\fi}
\ctr@ln@m\acc@louv \ctr@ln@m\acc@lfer
{\catcode`\[=1\catcode`\{=12\gdef\acc@louv[{}}
{\catcode`\]=2\catcode`\}=12\gdef\acc@lfer{}]]
\ctr@ld@f\def\PSdict@{\ifUse@llipse%
    \immediate\write\fwf@g{/ellipsedict 9 dict def ellipsedict /mtrx matrix put}%
    \immediate\write\fwf@g{/ellipse \acc@louv ellipsedict begin}%
    \immediate\write\fwf@g{ /endangle exch def /startangle exch def}%
    \immediate\write\fwf@g{ /yrad exch def /xrad exch def}%
    \immediate\write\fwf@g{ /rotangle exch def /y exch def /x exch def}%
    \immediate\write\fwf@g{ /savematrix mtrx currentmatrix def}%
    \immediate\write\fwf@g{ x y translate rotangle rotate xrad yrad scale}%
    \immediate\write\fwf@g{ 0 0 1 startangle endangle arc}%
    \immediate\write\fwf@g{ savematrix setmatrix end\acc@lfer def}%
    \fi\PShe@der{EndProlog}}
\ctr@ld@f\def\Pssetc@rve#1=#2|{\keln@mun#1|%
    \def\n@mref{r}\ifx\l@debut\n@mref\update@ttr\D@FTroundness\Q@s@troundness{#2}\else
    \W@rnmesAttr{figset curve}{#1}\fi}
\ctr@ln@m\curv@roundness
\ctr@ld@f\def\Q@s@troundness#1{\edef\curv@roundness{#1}}
\ctr@ld@f\def\D@FTroundness{0.2} 
\ctr@ln@m\figdrawcurve
\ctr@ld@f\def\Q@curveDD[#1]{{\ifCUR@PS\ifGR@cri\PSc@mment{curveDD Points=#1}%
    \s@uvc@ntr@l\et@tpscurveDD%
    \iffillm@de\Psc@rveDD\curv@roundness[#1]%
    \f@gfill%
    \else\Psc@rveDD\curv@roundness[#1]\f@gstroke\fi%
    \PSc@mment{End curveDD}\resetc@ntr@l\et@tpscurveDD\fi\fi}}
\ctr@ld@f\def\Q@curveTD[#1]{{\ifCUR@PS\ifGR@cri%
    \PSc@mment{curveTD Points=#1}\s@uvc@ntr@l\et@tpscurveTD\let\c@lprojSP=\relax%
    \iffillm@de\Psc@rveTD\curv@roundness[#1]%
    \f@gfill%
    \else\Psc@rveTD\curv@roundness[#1]\f@gstroke\fi%
    \PSc@mment{End curveTD}\resetc@ntr@l\et@tpscurveTD\fi\fi}}
\ctr@ld@f\def\Psc@rveDD#1[#2]{%
    \def\list@num{#2}\extrairelepremi@r\Ak@\de\list@num%
    \extrairelepremi@r\Ai@\de\list@num\extrairelepremi@r\Aj@\de\list@num%
    \f@gnewpath\PSwrit@cmdS{\Ai@}{\c@mmoveto}{\fwf@g}{\X@un}{\Y@un}%
    \setc@ntr@l{2}\figvectPDD -1[\Ak@,\Aj@]%
    \@ecfor\Ak@:=\list@num\do{\figpttraDD-2:=\Ai@/#1,-1/\BdingB@xfalse%
       \PSwrit@cmdS{-2}{}{\fwf@g}{\X@de}{\Y@de}%
       \figvectPDD -1[\Ai@,\Ak@]\figpttraDD-2:=\Aj@/-#1,-1/%
       \PSwrit@cmdS{-2}{}{\fwf@g}{\X@tr}{\Y@tr}\BdingB@xtrue%
       \PSwrit@cmdS{\Aj@}{\c@mcurveto}{\fwf@g}{\X@qu}{\Y@qu}%
       \B@zierBB@x{1}{\Y@un}(\X@un,\X@de,\X@tr,\X@qu)%
       \B@zierBB@x{2}{\X@un}(\Y@un,\Y@de,\Y@tr,\Y@qu)%
       \edef\X@un{\X@qu}\edef\Y@un{\Y@qu}\edef\Ai@{\Aj@}\edef\Aj@{\Ak@}}}
\ctr@ld@f\def\Psc@rveTD#1[#2]{\ifnum\CUR@proj<\tw@\Psc@rvePPTD#1[#2]\else\Psc@rveCPTD#1[#2]\fi}
\ctr@ld@f\def\Psc@rvePPTD#1[#2]{\setc@ntr@l{2}%
    \def\list@num{#2}\extrairelepremi@r\Ak@\de\list@num\Figptpr@j-5:/\Ak@/%
    \extrairelepremi@r\Ai@\de\list@num\Figptpr@j-3:/\Ai@/%
    \extrairelepremi@r\Aj@\de\list@num\Figptpr@j-4:/\Aj@/%
    \f@gnewpath\PSwrit@cmdS{-3}{\c@mmoveto}{\fwf@g}{\X@un}{\Y@un}%
    \figvectPDD -1[-5,-4]%
    \@ecfor\Ak@:=\list@num\do{\Figptpr@j-5:/\Ak@/\figpttraDD-2:=-3/#1,-1/%
       \BdingB@xfalse\PSwrit@cmdS{-2}{}{\fwf@g}{\X@de}{\Y@de}%
       \figvectPDD -1[-3,-5]\figpttraDD-2:=-4/-#1,-1/%
       \PSwrit@cmdS{-2}{}{\fwf@g}{\X@tr}{\Y@tr}\BdingB@xtrue%
       \PSwrit@cmdS{-4}{\c@mcurveto}{\fwf@g}{\X@qu}{\Y@qu}%
       \B@zierBB@x{1}{\Y@un}(\X@un,\X@de,\X@tr,\X@qu)%
       \B@zierBB@x{2}{\X@un}(\Y@un,\Y@de,\Y@tr,\Y@qu)%
       \edef\X@un{\X@qu}\edef\Y@un{\Y@qu}\figptcopyDD-3:/-4/\figptcopyDD-4:/-5/}}
\ctr@ld@f\def\Psc@rveCPTD#1[#2]{\setc@ntr@l{2}%
    \def\list@num{#2}\extrairelepremi@r\Ak@\de\list@num%
    \extrairelepremi@r\Ai@\de\list@num\extrairelepremi@r\Aj@\de\list@num%
    \Figptpr@j-7:/\Ai@/%
    \f@gnewpath\PSwrit@cmd{-7}{\c@mmoveto}{\fwf@g}%
    \figvectPTD -9[\Ak@,\Aj@]%
    \@ecfor\Ak@:=\list@num\do{\figpttraTD-10:=\Ai@/#1,-9/%
       \figvectPTD -9[\Ai@,\Ak@]\figpttraTD-11:=\Aj@/-#1,-9/%
       \subB@zierTD\NBz@rcs[\Ai@,-10,-11,\Aj@]\edef\Ai@{\Aj@}\edef\Aj@{\Ak@}}}
\ctr@ld@f\def\figdrawend{\ifCUR@PS\ifGR@cri\immediate\closeout\fwf@g%
    \immediate\openout\fwf@g=\PSfilen@me\relax%
    \ifPDFm@ke\PSBdingB@x\else%
    \immediate\write\fwf@g{\p@urcent\string!PS-Adobe-2.0 EPSF-2.0}%
    \PShe@der{Creator\string: TeX (fig4tex.tex)}%
    \PShe@der{Title\string: \PSfilen@me}%
    \PShe@der{CreationDate\string: \the\day/\the\month/\the\year}%
    \PSBdingB@x%
    \PShe@der{EndComments}\PSdict@\fi%
    \immediate\write\fwf@g{\c@mgsave}%
    \openin\frf@g=\auxfilen@me\c@pypsfile\fwf@g\frf@g\closein\frf@g%
    \immediate\write\fwf@g{\c@mgrestore}%
    \PSc@mment{End of file.}\immediate\closeout\fwf@g%
    \immediate\openout\fwf@g=\auxfilen@me\immediate\closeout\fwf@g%
    \immediate\write16{File \PSfilen@me\space created.}\fi\fi\CUR@PSfalse\GR@critrue}
\ctr@ld@f\def\PShe@der#1{\immediate\write\fwf@g{\p@urcent\p@urcent#1}}
\ctr@ld@f\def\PSBdingB@x{{\v@lX=\ptT@ptps\c@@rdXmin\v@lY=\ptT@ptps\c@@rdYmin%
     \v@lXa=\ptT@ptps\c@@rdXmax\v@lYa=\ptT@ptps\c@@rdYmax%
     \PShe@der{BoundingBox\string: \repdecn@mb{\v@lX}\space\repdecn@mb{\v@lY}%
     \space\repdecn@mb{\v@lXa}\space\repdecn@mb{\v@lYa}}}}
\ctr@ld@f\def\figdrawfcconnect[#1]{{\ifCUR@PS\ifGR@cri\PSc@mment{fcconnect Points=#1}%
    \Q@s@tfillmode{no}\s@uvc@ntr@l\et@tpsfcconnect\resetc@ntr@l{2}%
    \fcc@nnect@[#1]\resetc@ntr@l\et@tpsfcconnect\PSc@mment{End fcconnect}\fi\fi}}
\ctr@ld@f\def\fcc@nnect@[#1]{\let\N@rm=\n@rmeucDD\def\list@num{#1}%
    \extrairelepremi@r\Ai@\de\list@num\edef\pr@m{\Ai@}\v@leur=\z@\p@rtent=\@ne\c@llgtot%
    \ifcase\fclin@typ@\edef\list@num{[\pr@m,#1,\Ai@}\expandafter\figdrawcurve\list@num]%
    \else\ifdim\fclin@r@d\p@>\z@\Pslin@conge[#1]\else\figdrawline[#1]\fi\fi%
    \v@leur=\@rrowp@s\v@leur\edef\list@num{#1,\Ai@,0}%
    \extrairelepremi@r\Ai@\de\list@num\mili@u=\epsil@n\c@llgpart%
    \advance\mili@u-\epsil@n\advance\mili@u-\delt@\advance\v@leur-\mili@u%
    \ifcase\fclin@typ@\invers@\mili@u\delt@%
    \ifnum\@rrowr@fpt>\z@\advance\delt@-\v@leur\v@leur=\delt@\fi%
    \v@leur=\repdecn@mb\v@leur\mili@u\edef\v@lt{\repdecn@mb\v@leur}%
    \extrairelepremi@r\Ak@\de\list@num%
    \figvectPDD-1[\pr@m,\Aj@]\figpttraDD-6:=\Ai@/\curv@roundness,-1/%
    \figvectPDD-1[\Ak@,\Ai@]\figpttraDD-7:=\Aj@/\curv@roundness,-1/%
    \delt@=\@rrowheadlength\p@\delt@=\C@AHANG\delt@\edef\R@dius{\repdecn@mb{\delt@}}%
    \ifcase\@rrowr@fpt%
    \FigptintercircB@zDD-8::\v@lt,\R@dius[\Ai@,-6,-7,\Aj@]\Q@arrowheadDD[-5,-8]\else%
    \FigptintercircB@zDD-8::\v@lt,\R@dius[\Aj@,-7,-6,\Ai@]\Q@arrowheadDD[-8,-5]\fi%
    \else\advance\delt@-\v@leur%
    \p@rtentiere{\p@rtent}{\delt@}\edef\C@efun{\the\p@rtent}%
    \p@rtentiere{\p@rtent}{\v@leur}\edef\C@efde{\the\p@rtent}%
    \figptbaryDD-5:[\Ai@,\Aj@;\C@efun,\C@efde]\ifcase\@rrowr@fpt%
    \delt@=\@rrowheadlength\unit@\delt@=\C@AHANG\delt@\edef\t@ille{\repdecn@mb{\delt@}}%
    \figvectPDD-2[\Ai@,\Aj@]\vecunit@{-2}{-2}\figpttraDD-5:=-5/\t@ille,-2/\fi%
    \Q@arrowheadDD[\Ai@,-5]\fi}
\ctr@ld@f\def\c@llgtot{\@ecfor\Aj@:=\list@num\do{\figvectP-1[\Ai@,\Aj@]\N@rm\delt@{-1}%
    \advance\v@leur\delt@\advance\p@rtent\@ne\edef\Ai@{\Aj@}}}
\ctr@ld@f\def\c@llgpart{\extrairelepremi@r\Aj@\de\list@num\figvectP-1[\Ai@,\Aj@]\N@rm\delt@{-1}%
    \advance\mili@u\delt@\ifdim\mili@u<\v@leur\edef\pr@m{\Ai@}\edef\Ai@{\Aj@}\c@llgpart\fi}
\ctr@ld@f\def\Pslin@conge[#1]{\ifnum\p@rtent>\tw@{\def\list@num{#1}%
    \extrairelepremi@r\Ai@\de\list@num\extrairelepremi@r\Aj@\de\list@num%
    \figptcopy-6:/\Ai@/\figvectP-3[\Ai@,\Aj@]\vecunit@{-3}{-3}\v@lmax=\result@t%
    \@ecfor\Ak@:=\list@num\do{\figvectP-4[\Aj@,\Ak@]\vecunit@{-4}{-4}%
    \minim@m\v@lmin\v@lmax\result@t\v@lmax=\result@t%
    \det@rm\delt@[-3,-4]\maxim@m\mili@u{\delt@}{-\delt@}\ifdim\mili@u>\Cepsil@n%
    \ifdim\delt@>\z@\figgetangleDD\Angl@[\Aj@,\Ak@,\Ai@]\else%
    \figgetangleDD\Angl@[\Aj@,\Ai@,\Ak@]\fi%
    \v@leur=\PI@deg\advance\v@leur-\Angl@\p@\divide\v@leur\tw@%
    \edef\Angl@{\repdecn@mb\v@leur}\c@ssin{\C@}{\S@}{\Angl@}\v@leur=\fclin@r@d\unit@%
    \v@leur=\S@\v@leur\mili@u=\C@\p@\invers@\mili@u\mili@u%
    \v@leur=\repdecn@mb{\mili@u}\v@leur%
    \minim@m\v@leur\v@leur\v@lmin\edef\t@ille{\repdecn@mb{\v@leur}}%
    \figpttra-5:=\Aj@/-\t@ille,-3/\figdrawline[-6,-5]\figpttra-6:=\Aj@/\t@ille,-4/%
    \figvectNVDD-3[-3]\figvectNVDD-8[-4]\inters@cDD-7:[-5,-3;-6,-8]%
    \ifdim\delt@>\z@\figdrawarccircP-7;\fclin@r@d[-5,-6]\else\figdrawarccircP-7;\fclin@r@d[-6,-5]\fi%
    \else\figdrawline[-6,\Aj@]\figptcopy-6:/\Aj@/\fi
    \edef\Ai@{\Aj@}\edef\Aj@{\Ak@}\figptcopy-3:/-4/}\figdrawline[-6,\Aj@]}\else\figdrawline[#1]\fi}
\ctr@ld@f\def\figdrawfcnode[#1]#2{{\ifCUR@PS\ifGR@cri\PSc@mment{fcnode Points=#1}%
    \s@uvc@ntr@l\et@tpsfcnode\resetc@ntr@l{2}%
    \def\t@xt@{#2}\ifx\t@xt@\empty\def\g@tt@xt{\setbox\Gb@x=\hbox{\Figg@tT{\p@int}}}%
    \else\def\g@tt@xt{\setbox\Gb@x=\hbox{#2}}\fi%
    \v@lmin=\h@rdfcXp@dd\advance\v@lmin\Xp@dd\unit@\multiply\v@lmin\tw@%
    \v@lmax=\h@rdfcYp@dd\advance\v@lmax\Yp@dd\unit@\multiply\v@lmax\tw@%
    \Figv@ctCreg-8(\unit@,-\unit@)\def\list@num{#1}%
    \delt@=\CUR@width bp\divide\delt@\tw@%
    \fcn@de\PSc@mment{End fcnode}\resetc@ntr@l\et@tpsfcnode\fi\fi}}
\ctr@ld@f\def\d@butn@de{\g@tt@xt\v@lX=\wd\Gb@x%
    \v@lY=\ht\Gb@x\advance\v@lY\dp\Gb@x\advance\v@lX\v@lmin\advance\v@lY\v@lmax}
\ctr@ld@f\def\fcn@deE{%
    \@ecfor\p@int:=\list@num\do{\d@butn@de\v@lX=\unssqrttw@\v@lX\v@lY=\unssqrttw@\v@lY%
    \ifdim\thickn@ss\p@>\z@
    \v@lXa=\v@lX\advance\v@lXa\delt@\v@lXa=\ptT@unit@\v@lXa\edef\XR@d{\repdecn@mb\v@lXa}%
    \v@lYa=\v@lY\advance\v@lYa\delt@\v@lYa=\ptT@unit@\v@lYa\edef\YR@d{\repdecn@mb\v@lYa}%
    \arct@n\v@leur(\v@lXa,\v@lYa)\v@leur=\rdT@deg\v@leur\edef\@nglde{\repdecn@mb\v@leur}%
    {\c@lptellDD-2::\p@int;\XR@d,\YR@d(\@nglde)}
    \advance\v@leur-\PI@deg\edef\@nglun{\repdecn@mb\v@leur}%
    {\c@lptellDD-3::\p@int;\XR@d,\YR@d(\@nglun)}%
    \figptstra-6=-3,-2,\p@int/\thickn@ss,-8/\Q@s@tfillmode{yes}%
    \Pss@tspecifSt{color=\DDV@thickcolor}%
    \figdrawline[-2,-3,-6,-5]\figdrawarcell-4;\XR@d,\YR@d(\@nglun,\@nglde,0)%
    \Psrest@reSt{color=\DDV@thickcolor}\fi
    \v@lX=\ptT@unit@\v@lX\v@lY=\ptT@unit@\v@lY%
    \edef\XR@d{\repdecn@mb\v@lX}\edef\YR@d{\repdecn@mb\v@lY}%
    \Q@s@tfillmode{yes}\Pss@tspecifSt{color=\fcbgc@lor}%
    \figdrawarcell\p@int;\XR@d,\YR@d(0,360,0)%
    \Q@s@tfillmode{no}\Psrest@reSt{color=\fcbgc@lor}\figdrawarcell\p@int;\XR@d,\YR@d(0,360,0)}}
\ctr@ld@f\def\fcn@deL{\delt@=\ptT@unit@\delt@\edef\t@ille{\repdecn@mb\delt@}%
    \@ecfor\p@int:=\list@num\do{\Figg@tXYa{\p@int}\d@butn@de%
    \ifdim\v@lX>\v@lY\itis@Ktrue\else\itis@Kfalse\fi%
    \advance\v@lXa-\v@lX\Figp@intreg-1:(\v@lXa,\v@lYa)%
    \advance\v@lXa\v@lX\advance\v@lYa-\v@lY\Figp@intreg-2:(\v@lXa,\v@lYa)%
    \advance\v@lXa\v@lX\advance\v@lYa\v@lY\Figp@intreg-3:(\v@lXa,\v@lYa)%
    \advance\v@lXa-\v@lX\advance\v@lYa\v@lY\Figp@intreg-4:(\v@lXa,\v@lYa)%
    \ifdim\thickn@ss\p@>\z@
    \Figg@tXYa{\p@int}\Q@s@tfillmode{yes}\Pss@tspecifSt{color=\DDV@thickcolor}%
    \c@lpt@xt{-1}{-4}\c@lpt@xt@\v@lXa\v@lYa\v@lX\v@lY\c@rre\delt@%
    \Figp@intregDD-9:(\v@lZ,\v@lYa)\Figp@intregDD-11:(\v@lZa,\v@lYa)%
    \c@lpt@xt{-4}{-3}\c@lpt@xt@\v@lYa\v@lXa\v@lY\v@lX\delt@\c@rre%
    \Figp@intregDD-12:(\v@lXa,\v@lZ)\Figp@intregDD-10:(\v@lXa,\v@lZa)%
    \ifitis@K\figptstra-7=-9,-10,-11/\thickn@ss,-8/\figdrawline[-9,-11,-5,-6,-7]\else%
    \figptstra-7=-10,-11,-12/\thickn@ss,-8/\figdrawline[-10,-12,-5,-6,-7]\fi%
    \Psrest@reSt{color=\DDV@thickcolor}\fi
    \Q@s@tfillmode{yes}\Pss@tspecifSt{color=\fcbgc@lor}\figdrawline[-1,-2,-3,-4]%
    \Q@s@tfillmode{no}\Psrest@reSt{color=\fcbgc@lor}\figdrawline[-1,-2,-3,-4,-1]}}
\ctr@ld@f\def\c@lpt@xt#1#2{\figvectN-7[#1,#2]\vecunit@{-7}{-7}\figpttra-5:=#1/\t@ille,-7/%
    \figvectP-7[#1,#2]\Figg@tXY{-7}\c@rre=\v@lX\delt@=\v@lY\Figg@tXY{-5}}
\ctr@ld@f\def\c@lpt@xt@#1#2#3#4#5#6{\v@lZ=#6\invers@{\v@lZ}{\v@lZ}\v@leur=\repdecn@mb{#5}\v@lZ%
    \v@lZ=#2\advance\v@lZ-#4\mili@u=\repdecn@mb{\v@leur}\v@lZ%
    \v@lZ=#3\advance\v@lZ\mili@u\v@lZa=-\v@lZ\advance\v@lZa\tw@#1}
\ctr@ld@f\def\fcn@deR{\@ecfor\p@int:=\list@num\do{\Figg@tXYa{\p@int}\d@butn@de%
    \advance\v@lXa-0.5\v@lX\advance\v@lYa-0.5\v@lY\Figp@intreg-1:(\v@lXa,\v@lYa)%
    \advance\v@lXa\v@lX\Figp@intreg-2:(\v@lXa,\v@lYa)%
    \advance\v@lYa\v@lY\Figp@intreg-3:(\v@lXa,\v@lYa)%
    \advance\v@lXa-\v@lX\Figp@intreg-4:(\v@lXa,\v@lYa)%
    \ifdim\thickn@ss\p@>\z@
    \Q@s@tfillmode{yes}\Pss@tspecifSt{color=\DDV@thickcolor}%
    \Figv@ctCreg-5(-\delt@,-\delt@)\figpttra-9:=-1/1,-5/%
    \Figv@ctCreg-5(\delt@,-\delt@)\figpttra-10:=-2/1,-5/%
    \Figv@ctCreg-5(\delt@,\delt@)\figpttra-11:=-3/1,-5/%
    \figptstra-7=-9,-10,-11/\thickn@ss,-8/\figdrawline[-9,-11,-5,-6,-7]%
    \Psrest@reSt{color=\DDV@thickcolor}\fi
    \Q@s@tfillmode{yes}\Pss@tspecifSt{color=\fcbgc@lor}\figdrawline[-1,-2,-3,-4]%
    \Q@s@tfillmode{no}\Psrest@reSt{color=\fcbgc@lor}\figdrawline[-1,-2,-3,-4,-1]}}
\ctr@ld@f\def\Pssetfl@wchart#1=#2|{\keln@mtr#1|%
    \def\n@mref{arr}\ifx\l@debut\n@mref\expandafter\keln@mtr\l@suite|%
     \def\n@mref{owp}\ifx\l@debut\n@mref\update@ttr\D@FTfcarrowposition\P@setfcarrowposition{#2}\else
     \def\n@mref{owr}\ifx\l@debut\n@mref\update@ttr\D@FTfcarrowrefpt\P@setfcarrowrefpt{#2}\else
     \W@rnmesAttr{figset flowchart}{#1}\fi\fi\else%
    \def\n@mref{bgc}\ifx\l@debut\n@mref\update@ttr\D@FTfcbgcolor\P@setfcbgcolor{#2}\else
    \def\n@mref{lin}\ifx\l@debut\n@mref\update@ttr\D@FTfcline\P@setfcline{#2}\else
    \def\n@mref{pad}\ifx\l@debut\n@mref\update@ttr\D@FTfcxpadding\P@setfcxpadding{#2}%
                                       \update@ttr\D@FTfcypadding\P@setfcypadding{#2}\else
    \def\n@mref{rad}\ifx\l@debut\n@mref\update@ttr\D@FTfcradius\P@setfcradius{#2}\else
    \def\n@mref{sha}\ifx\l@debut\n@mref\update@ttr\D@FTfcshape\P@setfcshape{#2}\else
    \def\n@mref{thi}\ifx\l@debut\n@mref\expandafter\keln@mtr\l@suite|%
     \def\n@mref{ckc}\ifx\l@debut\n@mref\update@ttr\D@FTref\P@setfcthickcolor{#2}\else
     \def\n@mref{ckn}\ifx\l@debut\n@mref\update@ttr\D@FTfcthickness\P@setfcthickness{#2}\else
     \W@rnmesAttr{figset flowchart}{#1}\fi\fi\else%
    \def\n@mref{xpa}\ifx\l@debut\n@mref\update@ttr\D@FTfcxpadding\P@setfcxpadding{#2}\else
    \def\n@mref{ypa}\ifx\l@debut\n@mref\update@ttr\D@FTfcypadding\P@setfcypadding{#2}\else
    \W@rnmesAttr{figset flowchart}{#1}\fi\fi\fi\fi\fi\fi\fi\fi\fi}
\ctr@ln@m\@rrowp@s
\ctr@ld@f\def\P@setfcarrowposition#1{\edef\@rrowp@s{#1}}
\ctr@ln@m\@rrowr@fpt
\ctr@ld@f\def\P@setfcarrowrefpt#1{\setfcr@fpt#1|}
\ctr@ld@f\def\setfcr@fpt#1#2|{\if#1e\def\@rrowr@fpt{1}\else\def\@rrowr@fpt{0}\fi}
\ctr@ln@m\fcbgc@lor
\ctr@ld@f\def\P@setfcbgcolor#1{\edef\fcbgc@lor{#1}}
\ctr@ln@m\fclin@typ@
\ctr@ld@f\def\P@setfcline#1{\setfccurv@#1|}
\ctr@ld@f\def\setfccurv@#1#2|{\if#1c\def\fclin@typ@{0}\else\def\fclin@typ@{1}\fi}
\ctr@ln@m\fclin@r@d
\ctr@ld@f\def\P@setfcradius#1{\edef\fclin@r@d{#1}}
\ctr@ln@m\fcn@de \ctr@ln@m\fcsh@pe
\ctr@ln@m\h@rdfcXp@dd \ctr@ln@m\h@rdfcYp@dd
\ctr@ld@f\def\P@setfcshape#1{\setfcshap@#1|}
\ctr@ld@f\def\setfcshap@#1#2|{%
    \if#1e\let\fcn@de=\fcn@deE\def\h@rdfcXp@dd{4pt}\def\h@rdfcYp@dd{4pt}%
     \edef\fcsh@pe{ellipse}\else%
    \if#1l\let\fcn@de=\fcn@deL\def\h@rdfcXp@dd{4pt}\def\h@rdfcYp@dd{4pt}%
     \edef\fcsh@pe{lozenge}\else%
          \let\fcn@de=\fcn@deR\def\h@rdfcXp@dd{6pt}\def\h@rdfcYp@dd{6pt}%
     \edef\fcsh@pe{rectangle}\fi\fi}
\ctr@ln@m\DDV@thickcolor
\ctr@ld@f\def\P@setfcthickcolor#1{\edef\DDV@thickcolor{#1}}
\ctr@ln@m\thickn@ss
\ctr@ld@f\def\P@setfcthickness#1{\edef\thickn@ss{#1}}
\ctr@ln@m\Xp@dd
\ctr@ld@f\def\P@setfcxpadding#1{\edef\Xp@dd{#1}}
\ctr@ln@m\Yp@dd
\ctr@ld@f\def\P@setfcypadding#1{\edef\Yp@dd{#1}}
\ctr@ld@f\def\figdrawline[#1]{{\ifCUR@PS\ifGR@cri\PSc@mment{line Points=#1}%
    \let\figdrawlign@=\Pslign@P\Pslin@{#1}\PSc@mment{End line}\fi\fi}}
\ctr@ld@f\def\figdrawlineF#1{{\ifCUR@PS\ifGR@cri\PSc@mment{lineF Filename=#1}%
    \let\figdrawlign@=\Pslign@F\Pslin@{#1}\PSc@mment{End lineF}\fi\fi}}
\ctr@ld@f\def\figdrawlineC(#1){{\ifCUR@PS\ifGR@cri\PSc@mment{lineC}%
    \let\figdrawlign@=\Pslign@C\Pslin@{#1}\PSc@mment{End lineC}\fi\fi}}
\ctr@ld@f\def\Pslin@#1{\iffillm@de\figdrawlign@{#1}%
    \f@gfill%
    \else\figdrawlign@{#1}\ifx\derp@int\premp@int%
    \f@gclosestroke%
    \else\f@gstroke\fi\fi}
\ctr@ld@f\def\Pslign@P#1{\def\list@num{#1}\extrairelepremi@r\p@int\de\list@num%
    \edef\premp@int{\p@int}\f@gnewpath%
    \PSwrit@cmd{\p@int}{\c@mmoveto}{\fwf@g}%
    \@ecfor\p@int:=\list@num\do{\PSwrit@cmd{\p@int}{\c@mlineto}{\fwf@g}%
    \edef\derp@int{\p@int}}}
\ctr@ld@f\def\Pslign@F#1{\s@uvc@ntr@l\et@tPslign@F\setc@ntr@l{2}\openin\frf@g=#1\relax%
    \ifeof\frf@g\message{*** File #1 not found !}\end\else%
    \read\frf@g to\tr@c\edef\premp@int{\tr@c}\expandafter\extr@ctCF\tr@c:%
    \f@gnewpath\PSwrit@cmd{-1}{\c@mmoveto}{\fwf@g}%
    \loop\read\frf@g to\tr@c\ifeof\frf@g\mored@tafalse\else\mored@tatrue\fi%
    \ifmored@ta\expandafter\extr@ctCF\tr@c:\PSwrit@cmd{-1}{\c@mlineto}{\fwf@g}%
    \edef\derp@int{\tr@c}\repeat\fi\closein\frf@g\resetc@ntr@l\et@tPslign@F}
\ctr@ln@m\extr@ctCF
\ctr@ld@f\def\extr@ctCFDD#1 #2:{\v@lX=#1\unit@\v@lY=#2\unit@\Figp@intregDD-1:(\v@lX,\v@lY)}
\ctr@ld@f\def\extr@ctCFTD#1 #2 #3:{\v@lX=#1\unit@\v@lY=#2\unit@\v@lZ=#3\unit@%
    \Figp@intregTD-1:(\v@lX,\v@lY,\v@lZ)}
\ctr@ld@f\def\Pslign@C#1{\s@uvc@ntr@l\et@tPslign@C\setc@ntr@l{2}%
    \def\list@num{#1}\extrairelepremi@r\p@int\de\list@num%
    \edef\premp@int{\p@int}\f@gnewpath%
    \expandafter\Pslign@C@\p@int:\PSwrit@cmd{-1}{\c@mmoveto}{\fwf@g}%
    \@ecfor\p@int:=\list@num\do{\expandafter\Pslign@C@\p@int:%
    \PSwrit@cmd{-1}{\c@mlineto}{\fwf@g}\edef\derp@int{\p@int}}%
    \resetc@ntr@l\et@tPslign@C}
\ctr@ld@f\def\Pslign@C@#1 #2:{{\def\t@xt@{#1}\ifx\t@xt@\empty\Pslign@C@#2:
    \else\extr@ctCF#1 #2:\fi}}
\ctr@ld@f\def\Pssetm@sh#1=#2|{\keln@mde#1|%
    \def\n@mref{co}\ifx\l@debut\n@mref\update@ttr\D@FTref\P@setmeshcolor{#2}\else
    \def\n@mref{da}\ifx\l@debut\n@mref\update@ttr\D@FTref\P@setmeshdash{#2}\else
    \def\n@mref{di}\ifx\l@debut\n@mref\update@ttr\D@FTmeshdiag\Q@s@tmeshdiag{#2}\else
    \def\n@mref{wi}\ifx\l@debut\n@mref\update@ttr\D@FTref\P@setmeshwidth{#2}\else
    \W@rnmesAttr{figset mesh}{#1}\fi\fi\fi\fi}
\ctr@ln@m\c@ntrolmesh
\ctr@ld@f\def\Q@s@tmeshdiag#1{\edef\c@ntrolmesh{#1}}
\ctr@ld@f\def\D@FTmeshdiag{0}    
\ctr@ln@m\DDV@meshcolor
\ctr@ld@f\def\P@setmeshcolor#1{\edef\DDV@meshcolor{#1}}
\ctr@ln@m\DDV@meshdash
\ctr@ld@f\def\P@setmeshdash#1{\edef\DDV@meshdash{#1}}
\ctr@ln@m\DDV@meshwidth
\ctr@ld@f\def\P@setmeshwidth#1{\edef\DDV@meshwidth{#1}}
\ctr@ld@f\def\figdrawmesh#1,#2[#3,#4,#5,#6]{{\ifCUR@PS\ifGR@cri%
    \PSc@mment{mesh N1=#1, N2=#2, Quadrangle=[#3,#4,#5,#6]}\s@uvc@ntr@l\et@tpsmesh%
    \Pss@tspecifSt{color=\DDV@meshcolor,dash=\DDV@meshdash,width=\DDV@meshwidth}%
    \setc@ntr@l{2}%
    \ifnum#1>\@ne\Psmeshp@rt#1[#3,#4,#5,#6]\fi%
    \ifnum#2>\@ne\Psmeshp@rt#2[#4,#5,#6,#3]\fi%
    \ifnum\c@ntrolmesh>\z@\Psmeshdi@g#1,#2[#3,#4,#5,#6]\fi%
    \ifnum\c@ntrolmesh<\z@\Psmeshdi@g#2,#1[#4,#5,#6,#3]\fi%
    \Psrest@reSt{color=\DDV@meshcolor,dash=\DDV@meshdash,width=\DDV@meshwidth}%
    \figdrawline[#3,#4,#5,#6,#3]\PSc@mment{End mesh}\resetc@ntr@l\et@tpsmesh\fi\fi}}
\ctr@ld@f\def\Psmeshp@rt#1[#2,#3,#4,#5]{{\l@mbd@un=\@ne\l@mbd@de=#1\loop%
    \ifnum\l@mbd@un<#1\advance\l@mbd@de\m@ne\figptbary-1:[#2,#3;\l@mbd@de,\l@mbd@un]%
    \figptbary-2:[#5,#4;\l@mbd@de,\l@mbd@un]\figdrawline[-1,-2]\advance\l@mbd@un\@ne\repeat}}
\ctr@ld@f\def\Psmeshdi@g#1,#2[#3,#4,#5,#6]{\figptcopy-2:/#3/\figptcopy-3:/#6/%
    \l@mbd@un=\z@\l@mbd@de=#1\loop\ifnum\l@mbd@un<#1%
    \advance\l@mbd@un\@ne\advance\l@mbd@de\m@ne\figptcopy-1:/-2/\figptcopy-4:/-3/%
    \figptbary-2:[#3,#4;\l@mbd@de,\l@mbd@un]%
    \figptbary-3:[#6,#5;\l@mbd@de,\l@mbd@un]\Psmeshdi@gp@rt#2[-1,-2,-3,-4]\repeat}
\ctr@ld@f\def\Psmeshdi@gp@rt#1[#2,#3,#4,#5]{{\l@mbd@un=\z@\l@mbd@de=#1\loop%
    \ifnum\l@mbd@un<#1\figptbary-5:[#2,#5;\l@mbd@de,\l@mbd@un]%
    \advance\l@mbd@de\m@ne\advance\l@mbd@un\@ne%
    \figptbary-6:[#3,#4;\l@mbd@de,\l@mbd@un]\figdrawline[-5,-6]\repeat}}
\ctr@ln@m\figdrawnormal
\ctr@ld@f\def\Q@normalDD#1,#2[#3,#4]{{\ifCUR@PS\ifGR@cri%
    \PSc@mment{normal Length=#1, Lambda=#2 [Pt1,Pt2]=[#3,#4]}%
    \s@uvc@ntr@l\et@tpsnormal\resetc@ntr@l{2}\figptendnormal-6::#1,#2[#3,#4]%
    \figptcopyDD-5:/-1/\figdrawarrow[-5,-6]%
    \PSc@mment{End normal}\resetc@ntr@l\et@tpsnormal\fi\fi}}
\ctr@ld@f\def\figreset#1{\trtlis@rg{#1}{\Psreset@}}
\ctr@ld@f\def\Psreset@#1|{\def\t@xt@{#1}\ifx\t@xt@\empty\P@resetg@n
    \else\keln@mde#1|%
    \def\n@mref{al}\ifx\l@debut\n@mref%
        \figset altitude(blcolor=default,bldash=default,blwidth=default,%
        sqcolor=default,sqdash=default,sqwidth=default)\else
    \def\n@mref{ar}\ifx\l@debut\n@mref\figresetarrowhead\else
    \def\n@mref{cu}\ifx\l@debut\n@mref\figset curve(roundness=\D@FTroundness)\else
    \def\n@mref{ge}\ifx\l@debut\n@mref\P@resetg@n\else
    \def\n@mref{fl}\ifx\l@debut\n@mref%
        \figset flowchart(arrowp=\D@FTfcarrowposition,arrowr=\D@FTfcarrowrefpt,%
	bgcolor=\D@FTfcbgcolor,line=\D@FTfcline,radius=\D@FTfcradius,%
	shape=\D@FTfcshape,thickcolor=default,thickness=\D@FTfcthickness,%
	xpadd=\D@FTfcxpadding,ypadd=\D@FTfcypadding)\else
    \def\n@mref{me}\ifx\l@debut\n@mref\figset mesh(diag=\D@FTmeshdiag,%
        color=default,dash=default,width=default)\else
    \def\n@mref{tr}\ifx\l@debut\n@mref%
        \figset trimesh(color=default,dash=default,width=default)\else
    \W@rnmeskwd{figreset}{#1}\fi\fi\fi\fi\fi\fi\fi\fi}
\ctr@ld@f\def\P@resetg@n{\figset (color=\D@FTcolor,dash=\D@FTdash,fill=\D@FTfill,%
    join=\D@FTjoin,width=\D@FTwidth)}
\ctr@ld@f\def\figset#1(#2){\def\t@xt@{#1}\ifx\t@xt@\empty\trtlis@rg{#2}{\Pssetg@n}
    \else\keln@mde#1|%
    \def\n@mref{al}\ifx\l@debut\n@mref\trtlis@rg{#2}{\Psset@lti}\else
    \def\n@mref{ar}\ifx\l@debut\n@mref\trtlis@rg{#2}{\Psset@rrowhe@d}\else
    \def\n@mref{cu}\ifx\l@debut\n@mref\trtlis@rg{#2}{\Pssetc@rve}\else
    \def\n@mref{fl}\ifx\l@debut\n@mref\trtlis@rg{#2}{\Pssetfl@wchart}\else
    \def\n@mref{ge}\ifx\l@debut\n@mref\trtlis@rg{#2}{\Pssetg@n}\else
    \def\n@mref{me}\ifx\l@debut\n@mref\trtlis@rg{#2}{\Pssetm@sh}\else
    \def\n@mref{pr}\ifx\l@debut\n@mref\ifCUR@PS\W@rnmesIgn{figset proj(...)}%
     \else\trtlis@rg{#2}{\Figsetpr@j}\fi\else
    \def\n@mref{tr}\ifx\l@debut\n@mref\trtlis@rg{#2}{\Pssettrim@sh}\else
    \def\n@mref{wr}\ifx\l@debut\n@mref\let\M@cro=\Figsetwr@te\trtlis@rgtok{#2,|}\else
    \W@rnmeskwd{figset}{#1}\fi\fi\fi\fi\fi\fi\fi\fi\fi\fi\ignorespaces}
\ctr@ld@f\def\figsetdefault#1(#2){\ifCUR@PS\W@rnmesIgn{figsetdefault}\else%
    \def\t@xt@{#1}\ifx\t@xt@\empty\trtlis@rg{#2}{\Pssd@g@n}\else\keln@mun#1|
    \def\n@mref{a}\ifx\l@debut\n@mref\trtlis@rg{#2}{\Pssd@@rrowhe@d}\else
    \def\n@mref{c}\ifx\l@debut\n@mref\trtlis@rg{#2}{\Pssd@c@rve}\else
    \def\n@mref{g}\ifx\l@debut\n@mref\trtlis@rg{#2}{\Pssd@g@n}\else
    \def\n@mref{f}\ifx\l@debut\n@mref\trtlis@rg{#2}{\Pssd@fl@wchart}\else
    \def\n@mref{m}\ifx\l@debut\n@mref\trtlis@rg{#2}{\Pssd@m@sh}\else
    \W@rnmeskwd{figsetdefault}{#1}\fi\fi\fi\fi\fi\fi\initpss@ttings\fi}
\ctr@ld@f\def\Pssd@g@n#1=#2|{\keln@mun#1|%
    \def\n@mref{c}\ifx\l@debut\n@mref\edef\D@FTcolor{#2}\else
    \def\n@mref{d}\ifx\l@debut\n@mref\edef\D@FTdash{#2}\else
    \def\n@mref{f}\ifx\l@debut\n@mref\edef\D@FTfill{#2}\else
    \def\n@mref{j}\ifx\l@debut\n@mref\edef\D@FTjoin{#2}\else
    \def\n@mref{u}\ifx\l@debut\n@mref\edef\D@FTupdate{#2}\Q@s@tupdate{#2}\else
    \def\n@mref{w}\ifx\l@debut\n@mref\edef\D@FTwidth{#2}\else
    \W@rnmesAttr{figsetdefault}{#1}\fi\fi\fi\fi\fi\fi}
\ctr@ld@f\def\Pssd@@rrowhe@d#1=#2|{\keln@mun#1|%
    \def\n@mref{a}\ifx\l@debut\n@mref\edef\D@FTarrowheadangle{#2}\else
    \def\n@mref{f}\ifx\l@debut\n@mref\edef\D@FTarrowheadfill{#2}\else
    \def\n@mref{l}\ifx\l@debut\n@mref\y@tiunit{#2}\ifunitpr@sent%
     \edef\D@FTh@rdahlength{#2}\else\edef\D@FTh@rdahlength{#2pt}%
     \message{*** \BS@ figsetdefault (..., #1=#2, ...) : unit is missing, pt is assumed.}%
     \fi\else
    \def\n@mref{o}\ifx\l@debut\n@mref\edef\D@FTarrowheadout{#2}\else
    \def\n@mref{r}\ifx\l@debut\n@mref\edef\D@FTarrowheadratio{#2}\else
    \W@rnmesAttr{figsetdefault arrowhead}{#1}\fi\fi\fi\fi\fi}
\ctr@ld@f\def\Pssd@c@rve#1=#2|{\keln@mun#1|%
    \def\n@mref{r}\ifx\l@debut\n@mref\edef\D@FTroundness{#2}\else%
    \W@rnmesAttr{figsetdefault curve}{#1}\fi}
\ctr@ld@f\def\Pssd@fl@wchart#1=#2|{\keln@mtr#1|%
    \def\n@mref{arr}\ifx\l@debut\n@mref\expandafter\keln@mtr\l@suite|%
     \def\n@mref{owp}\ifx\l@debut\n@mref\edef\D@FTfcarrowposition{#2}\else
     \def\n@mref{owr}\ifx\l@debut\n@mref\edef\D@FTfcarrowrefpt{#2}\else
                     \W@rnmesAttr{figsetdefault flowchart}{#1}\fi\fi\else%
    \def\n@mref{bgc}\ifx\l@debut\n@mref\edef\D@FTfcbgcolor{#2}\else
    \def\n@mref{lin}\ifx\l@debut\n@mref\edef\D@FTfcline{#2}\else
    \def\n@mref{pad}\ifx\l@debut\n@mref\edef\D@FTfcxpadding{#2}%
                    \edef\D@FTfcypadding{#2}\else
    \def\n@mref{rad}\ifx\l@debut\n@mref\edef\D@FTfcradius{#2}\else
    \def\n@mref{sha}\ifx\l@debut\n@mref\edef\D@FTfcshape{#2}\else
    \def\n@mref{thi}\ifx\l@debut\n@mref\expandafter\keln@mtr\l@suite|%
     \def\n@mref{ckn}\ifx\l@debut\n@mref\edef\D@FTfcthickness{#2}\else
                     \W@rnmesAttr{figsetdefault flowchart}{#1}\fi\else%
    \def\n@mref{xpa}\ifx\l@debut\n@mref\edef\D@FTfcxpadding{#2}\else
    \def\n@mref{ypa}\ifx\l@debut\n@mref\edef\D@FTfcypadding{#2}\else
    \W@rnmesAttr{figsetdefault flowchart}{#1}\fi\fi\fi\fi\fi\fi\fi\fi\fi}
\ctr@ld@f\def\D@FTfcarrowposition{0.5}
\ctr@ld@f\def\D@FTfcarrowrefpt{start}
\ctr@ld@f\def\D@FTfcbgcolor{1}
\ctr@ld@f\def\D@FTfcline{polygon}
\ctr@ld@f\def\D@FTfcradius{0}
\ctr@ld@f\def\D@FTfcshape{rectangle}
\ctr@ld@f\def\D@FTfcthickness{0}
\ctr@ld@f\def\D@FTfcxpadding{0}
\ctr@ld@f\def\D@FTfcypadding{0}
\ctr@ld@f\def\Pssd@m@sh#1=#2|{\keln@mun#1|%
    \def\n@mref{d}\ifx\l@debut\n@mref\edef\D@FTmeshdiag{#2}\else%
    \W@rnmesAttr{figsetdefault mesh}{#1}\fi}
\ctr@ln@w{newif}\iffillm@de
\ctr@ld@f\def\Q@s@tfillmode#1{\expandafter\setfillm@de#1:}
\ctr@ld@f\def\setfillm@de#1#2:{\if#1n\fillm@defalse\else\fillm@detrue\fi}
\ctr@ld@f\def\D@FTfill{no}     
\ctr@ln@w{newif}\ifGRupdatem@de
\ctr@ld@f\def\Q@s@tupdate#1{\ifCUR@PS\W@rnmesIgn{figset (update=...)}%
    \else\expandafter\setupd@te#1:\fi}
\ctr@ld@f\def\setupd@te#1#2:{\if#1n\GRupdatem@defalse\else\GRupdatem@detrue\fi}
\ctr@ld@f\def\D@FTupdate{no}     
\ctr@ln@m\CUR@color \ctr@ln@m\CUR@colorc@md
\ctr@ld@f\def\s@uvcolor#1{\edef#1{\CUR@color}}
\ctr@ld@f\def\D@FTcolor{0}       
\ctr@ld@f\def\Pssetc@lor#1{\ifGR@cri\result@tent=\@ne\expandafter\c@lnbV@l#1 :%
    \def\CUR@color{}\def\CUR@colorc@md{}%
    \ifcase\result@tent\or\Q@s@tgray{#1}\or\or\Q@s@trgb{#1}\or\Q@s@tcmyk{#1}\fi\fi}
\ctr@ln@m\CUR@colorc@mdStroke
\ctr@ld@f\def\Q@s@tcmyk#1{\ifGR@cri\def\CUR@color{#1}\def\CUR@colorc@md{\c@msetcmykcolor}%
    \def\CUR@colorc@mdStroke{\c@msetcmykcolorStroke}%
    \ifCUR@PS\PSc@mment{setcmyk Color=#1}\us@primarC@lor\fi\fi}
\ctr@ld@f\def\Q@s@trgb#1{\ifGR@cri\def\CUR@color{#1}\def\CUR@colorc@md{\c@msetrgbcolor}%
    \def\CUR@colorc@mdStroke{\c@msetrgbcolorStroke}%
    \ifCUR@PS\PSc@mment{setrgb Color=#1}\us@primarC@lor\fi\fi}
\ctr@ld@f\def\Q@s@tgray#1{\ifGR@cri\def\CUR@color{#1}\def\CUR@colorc@md{\c@msetgray}%
    \def\CUR@colorc@mdStroke{\c@msetgrayStroke}%
    \ifCUR@PS\PSc@mment{setgray Gray level=#1}\us@primarC@lor\fi\fi}
\ctr@ln@m\fillc@md
\ctr@ld@f\def\us@primarC@lor{\immediate\write\fwf@g{\d@fprimarC@lor}%
    \let\fillc@md=\prfillc@md}
\ctr@ld@f\def\prfillc@md{\d@fprimarC@lor\space\c@mfill}
\ctr@ld@f\def\c@lnbV@l#1 #2:{\def\t@xt@{#1}\relax\ifx\t@xt@\empty\c@lnbV@l#2:
    \else\c@lnbV@l@#1 #2:\fi}
\ctr@ld@f\def\c@lnbV@l@#1 #2:{\def\t@xt@{#2}\ifx\t@xt@\empty%
    \def\t@xt@{#1}\ifx\t@xt@\empty\advance\result@tent\m@ne\fi
    \else\advance\result@tent\@ne\c@lnbV@l@#2:\fi}
\ctr@ld@f\def\Blackcmyk{0 0 0 1}
\ctr@ld@f\def\Whitecmyk{0 0 0 0}
\ctr@ld@f\def\Cyancmyk{1 0 0 0}
\ctr@ld@f\def\Magentacmyk{0 1 0 0}
\ctr@ld@f\def\Yellowcmyk{0 0 1 0}
\ctr@ld@f\def\Redcmyk{0 1 1 0}
\ctr@ld@f\def\Greencmyk{1 0 1 0}
\ctr@ld@f\def\Bluecmyk{1 1 0 0}
\ctr@ld@f\def\Graycmyk{0 0 0 0.50}
\ctr@ld@f\def\BrickRedcmyk{0 0.89 0.94 0.28} 
\ctr@ld@f\def\Browncmyk{0 0.81 1 0.60} 
\ctr@ld@f\def\ForestGreencmyk{0.91 0 0.88 0.12} 
\ctr@ld@f\def\Goldenrodcmyk{ 0 0.10 0.84 0} 
\ctr@ld@f\def\Marooncmyk{0 0.87 0.68 0.32} 
\ctr@ld@f\def\Orangecmyk{0 0.61 0.87 0} 
\ctr@ld@f\def\Purplecmyk{0.45 0.86 0 0} 
\ctr@ld@f\def\RoyalBluecmyk{1. 0.50 0 0} 
\ctr@ld@f\def\Violetcmyk{0.79 0.88 0 0} 
\ctr@ld@f\def\Blackrgb{0 0 0}
\ctr@ld@f\def\Whitergb{1 1 1}
\ctr@ld@f\def\Redrgb{1 0 0}
\ctr@ld@f\def\Greenrgb{0 1 0}
\ctr@ld@f\def\Bluergb{0 0 1}
\ctr@ld@f\def\Cyanrgb{0 1 1}
\ctr@ld@f\def\Magentargb{1 0 1}
\ctr@ld@f\def\Yellowrgb{1 1 0}
\ctr@ld@f\def\Grayrgb{0.5 0.5 0.5}
\ctr@ld@f\def\Chocolatergb{0.824 0.412 0.118}
\ctr@ld@f\def\DarkGoldenrodrgb{0.722 0.525 0.043}
\ctr@ld@f\def\DarkOrangergb{1 0.549 0}
\ctr@ld@f\def\Firebrickrgb{0.698 0.133 0.133}
\ctr@ld@f\def\ForestGreenrgb{0.133 0.545 0.133}
\ctr@ld@f\def\Goldrgb{1 0.843 0}
\ctr@ld@f\def\HotPinkrgb{1 0.412 0.706}
\ctr@ld@f\def\Maroonrgb{0.690 0.188 0.376}
\ctr@ld@f\def\Pinkrgb{1 0.753 0.796}
\ctr@ld@f\def\RoyalBluergb{0.255 0.412 0.882}
\ctr@ld@f\def\Pssetg@n#1=#2|{\keln@mun#1|%
    \def\n@mref{c}\ifx\l@debut\n@mref\update@ttr\D@FTcolor\Pssetc@lor{#2}\else
    \def\n@mref{d}\ifx\l@debut\n@mref\update@ttr\D@FTdash\Q@s@tdash{#2}\else
    \def\n@mref{f}\ifx\l@debut\n@mref\update@ttr\D@FTfill\Q@s@tfillmode{#2}\else
    \def\n@mref{j}\ifx\l@debut\n@mref\update@ttr\D@FTjoin\Q@s@tjoin{#2}\else
    \def\n@mref{u}\ifx\l@debut\n@mref\update@ttr\D@FTupdate\Q@s@tupdate{#2}\else
    \def\n@mref{w}\ifx\l@debut\n@mref\update@ttr\D@FTwidth\Q@s@twidth{#2}\else
    \W@rnmesAttr{figset}{#1}\fi\fi\fi\fi\fi\fi}
\ctr@ln@m\CUR@dash
\ctr@ld@f\def\s@uvdash#1{\edef#1{\CUR@dash}}
\ctr@ld@f\def\D@FTdash{1}        
\ctr@ld@f\def\Q@s@tdash#1{\ifGR@cri\edef\CUR@dash{#1}\ifCUR@PS\expandafter\Pssetd@sh#1 :\fi\fi}
\ctr@ld@f\def\Pssetd@shI#1{\PSc@mment{setdash Index=#1}\ifcase#1%
    \or\immediate\write\fwf@g{[] 0 \c@msetdash}
    \or\immediate\write\fwf@g{[6 2] 0 \c@msetdash}
    \or\immediate\write\fwf@g{[4 2] 0 \c@msetdash}
    \or\immediate\write\fwf@g{[2 2] 0 \c@msetdash}
    \or\immediate\write\fwf@g{[1 2] 0 \c@msetdash}
    \or\immediate\write\fwf@g{[2 4] 0 \c@msetdash}
    \or\immediate\write\fwf@g{[3 5] 0 \c@msetdash}
    \or\immediate\write\fwf@g{[3 3] 0 \c@msetdash}
    \or\immediate\write\fwf@g{[3 5 1 5] 0 \c@msetdash}
    \or\immediate\write\fwf@g{[6 4 2 4] 0 \c@msetdash}
    \fi}
\ctr@ld@f\def\Pssetd@sh#1 #2:{{\def\t@xt@{#1}\ifx\t@xt@\empty\Pssetd@sh#2:
    \else\def\t@xt@{#2}\ifx\t@xt@\empty\Pssetd@shI{#1}\else\s@mme=\@ne\def\debutp@t{#1}%
    \an@lysd@sh#2:\ifodd\s@mme\edef\debutp@t{\debutp@t\space\finp@t}\def\finp@t{0}\fi%
    \PSc@mment{setdash Pattern=#1 #2}%
    \immediate\write\fwf@g{[\debutp@t] \finp@t\space\c@msetdash}\fi\fi}}
\ctr@ld@f\def\an@lysd@sh#1 #2:{\def\t@xt@{#2}\ifx\t@xt@\empty\def\finp@t{#1}\else%
    \edef\debutp@t{\debutp@t\space#1}\advance\s@mme\@ne\an@lysd@sh#2:\fi}
\ctr@ln@m\CUR@width
\ctr@ld@f\def\s@uvwidth#1{\edef#1{\CUR@width}}
\ctr@ld@f\def\D@FTwidth{0.4}     
\ctr@ld@f\def\Q@s@twidth#1{\ifGR@cri\edef\CUR@width{#1}\ifCUR@PS%
    \PSc@mment{setwidth Width=#1}\immediate\write\fwf@g{#1 \c@msetlinewidth}\fi\fi}
\ctr@ln@m\CUR@join
\ctr@ld@f\def\s@uvjoin#1{\edef#1{\CUR@join}}
\ctr@ld@f\def\D@FTjoin{miter}   
\ctr@ld@f\def\Q@s@tjoin#1{\ifGR@cri\edef\CUR@join{#1}\ifCUR@PS\expandafter\Pssetj@in#1:\fi\fi}
\ctr@ld@f\def\Pssetj@in#1#2:{\PSc@mment{setjoin join=#1}%
    \if#1r\def\t@xt@{1}\else\if#1b\def\t@xt@{2}\else\def\t@xt@{0}\fi\fi%
    \immediate\write\fwf@g{\t@xt@\space\c@msetlinejoin}}
\ctr@ld@f\def\Pss@tspecifSt#1{\trtlis@rg{#1}{\Pss@tspecifSt@}}
\ctr@ld@f\def\Pss@tspecifSt@#1=#2|{\keln@mun#1|%
    \def\n@mref{c}\ifx\l@debut\n@mref\def\n@mref{#2}\ifx\n@mref\D@FTref\else%
     \s@uvcolor{\typ@color}\Pssetc@lor{#2}\fi\else
    \def\n@mref{d}\ifx\l@debut\n@mref\def\n@mref{#2}\ifx\n@mref\D@FTref\else%
     \s@uvdash{\typ@dash}\Q@s@tdash{#2}\fi\else
    \def\n@mref{j}\ifx\l@debut\n@mref\def\n@mref{#2}\ifx\n@mref\D@FTref\else%
     \s@uvjoin{\typ@join}\Q@s@tjoin{#2}\fi\else
    \def\n@mref{w}\ifx\l@debut\n@mref\def\n@mref{#2}\ifx\n@mref\D@FTref\else%
     \s@uvwidth{\typ@width}\Q@s@twidth{#2}\fi\else
    \W@rnmeskwd{Pss@tspecifSt}{#1}\fi\fi\fi\fi}
\ctr@ld@f\def\Psrest@reSt#1{\trtlis@rg{#1}{\Psrest@reSt@}}
\ctr@ld@f\def\Psrest@reSt@#1=#2|{\keln@mun#1|%
    \def\n@mref{c}\ifx\l@debut\n@mref\def\n@mref{#2}\ifx\n@mref\D@FTref\else%
     \Pssetc@lor{\typ@color}\fi\else
    \def\n@mref{d}\ifx\l@debut\n@mref\def\n@mref{#2}\ifx\n@mref\D@FTref\else%
     \Q@s@tdash{\typ@dash}\fi\else
    \def\n@mref{j}\ifx\l@debut\n@mref\def\n@mref{#2}\ifx\n@mref\D@FTref\else%
     \Q@s@tjoin{\typ@join}\fi\else
    \def\n@mref{w}\ifx\l@debut\n@mref\def\n@mref{#2}\ifx\n@mref\D@FTref\else%
     \Q@s@twidth{\typ@width}\fi\else
    \W@rnmeskwd{Psrest@reSt}{#1}\fi\fi\fi\fi}
\ctr@ld@f\def\Pssettrim@sh#1=#2|{\keln@mde#1|%
    \def\n@mref{co}\ifx\l@debut\n@mref\update@ttr\D@FTref\P@settmeshcolor{#2}\else
    \def\n@mref{da}\ifx\l@debut\n@mref\update@ttr\D@FTref\P@settmeshdash{#2}\else
    \def\n@mref{wi}\ifx\l@debut\n@mref\update@ttr\D@FTref\P@settmeshwidth{#2}\else
    \W@rnmesAttr{figset trimesh}{#1}\fi\fi\fi}
\ctr@ln@m\DDV@tmeshcolor
\ctr@ld@f\def\P@settmeshcolor#1{\edef\DDV@tmeshcolor{#1}}
\ctr@ln@m\DDV@tmeshdash
\ctr@ld@f\def\P@settmeshdash#1{\edef\DDV@tmeshdash{#1}}
\ctr@ln@m\DDV@tmeshwidth
\ctr@ld@f\def\P@settmeshwidth#1{\edef\DDV@tmeshwidth{#1}}
\ctr@ld@f\def\figdrawtrimesh#1[#2,#3,#4]{{\ifCUR@PS\ifGR@cri%
    \PSc@mment{trimesh Type=#1, Triangle=[#2,#3,#4]}%
    \s@uvc@ntr@l\et@tpstrimesh\ifnum#1>\@ne%
    \Pss@tspecifSt{color=\DDV@tmeshcolor,dash=\DDV@tmeshdash,width=\DDV@tmeshwidth}%
    \setc@ntr@l{2}%
    \Pstrimeshp@rt#1[#2,#3,#4]\Pstrimeshp@rt#1[#3,#4,#2]\Pstrimeshp@rt#1[#4,#2,#3]%
    \Psrest@reSt{color=\DDV@tmeshcolor,dash=\DDV@tmeshdash,width=\DDV@tmeshwidth}%
    \fi\figdrawline[#2,#3,#4,#2]%
    \PSc@mment{End trimesh}\resetc@ntr@l\et@tpstrimesh\fi\fi}}
\ctr@ld@f\def\Pstrimeshp@rt#1[#2,#3,#4]{{\l@mbd@un=\@ne\l@mbd@de=#1\loop\ifnum\l@mbd@de>\@ne%
    \advance\l@mbd@de\m@ne\figptbary-1:[#2,#3;\l@mbd@de,\l@mbd@un]%
    \figptbary-2:[#2,#4;\l@mbd@de,\l@mbd@un]\figdrawline[-1,-2]%
    \advance\l@mbd@un\@ne\repeat}}
\initpr@lim\initpss@ttings\initPDF@rDVI
\ctr@ln@w{newbox}\figBoxA
\ctr@ln@w{newbox}\figBoxB
\ctr@ln@w{newbox}\figBoxC
\catcode`\@=12

\pssetdefault(update=yes)
\newbox\figbox

\usepackage{color}
\definecolor{gr}{rgb}   {0.,   0.6,   0.25 }
\definecolor{mg}{rgb}   {0.85,  0.,    0.85}
\definecolor{marin}{rgb}   {0.,   0.65,   0.25}
\definecolor{rouge}{rgb}   {0.8,   0.,   0.}
\definecolor{orange}{rgb}   {0.8,   0.4,   0.}

\newcommand{\Gr}{\color{gr}}
\newcommand{\Mg}{\color{mg}}
\newcommand{\Bk}{\color{black}}
\newcommand{\Rd}{\color{red}}
\newcommand{\Bl}{\color{blue}}
\newcommand{\Or}{\color{orange}}

\newtheorem{theorem}{Theorem}[section]

\theoremstyle{definition}
\newtheorem{definition}[theorem]{Definition}

\newtheorem{assumption}[theorem]{Assumption}
\newtheorem{postulate}[theorem]{Postulate}

\theoremstyle{remark}
\newtheorem{remark}[theorem]{Remark}
\newtheorem{example}[theorem]{Example}

\numberwithin{equation}{section}

\newcommand{\ee}{\hskip0.15ex}
\newcommand{\me}{\hskip-0.15ex}
\newcommand{\dd}[1]{_{\raise-0.3ex\hbox{$\scriptstyle #1$}}}
\newcommand{\di}{\displaystyle}

\renewcommand{\Re}{\operatorname{\rm Re}}
\renewcommand{\Im}{\operatorname{\rm Im}}

\newcommand{\vpha}{\left.\vphantom{T^{j_0}_{j_0}}\!\!\right.}

\newcommand {\Norm}[2]{ \mathchoice
    {|\ee #1\ee|\dd{#2}}
    {| #1 |_{#2}}
    {| #1 |_{#2}}
    {| #1 |_{#2}} }
\newcommand {\DNorm}[2]{ \mathchoice
    {\|\ee #1\ee\|\dd{#2}}
    {\| #1 \|_{#2}}
    {\| #1 \|_{#2}}
    {\| #1 \|_{#2}} }
\newcommand {\Normc}[2]{ \mathchoice
    {|\ee #1\ee|\dd{#2}^2}
    {| #1 |_{#2}^2}
    {| #1 |_{#2}^2}
    {| #1 |_{#2}^2} }
\newcommand {\DNormc}[2]{ \mathchoice
    {\|\ee #1\ee\|\dd{#2}^2}
    {\| #1 \|_{#2}^2}
    {\| #1 \|_{#2}^2}
    {\| #1 \|_{#2}^2} }

\newcommand\rd{{\mathrm d}}
\newcommand\re{{\mathrm e}}

\renewcommand{\div}{\operatorname{\rm div}}
\newcommand{\grad}{\operatorname{\textbf{grad}}}
\newcommand{\curl}{\operatorname{\textbf{curl}}}
\newcommand{\curls}{\operatorname{\rm curl}}
\newcommand{\Cinf}{\mathscr{C}^\infty}

\renewcommand\P{{\mathbb P}}
\newcommand\Sbb{{\mathbb S}}
\newcommand\C{{\mathbb C}}
\newcommand\R{{\mathbb R}}
\newcommand\N{{\mathbb N}}
\newcommand\Q{{\mathbb Q}}
\newcommand\T{\mathbb{T}}
\newcommand\Z{{\mathbb Z}}

\newcommand\A{{\mathscr A}}
\newcommand\B{\mathsf{B}}
\newcommand\BB{{\mathscr B}}
\newcommand\sC{{\mathscr C}}
\newcommand\cC{{\mathcal C}}
\newcommand\D{\mathrm{D}}
\newcommand\DD{{\mathscr D}}
\newcommand\E{{\mathcal E}}
\newcommand\F{{\mathscr F}}
\newcommand\G{{\mathscr G}}
\newcommand\GG{{\mathcal G}}
\newcommand\cI{{\mathcal I}}
\newcommand\cH{{\mathcal H}}
\newcommand\J{{\mathcal J}}
\newcommand\cK{\boldsymbol{\mathcal K}}
\renewcommand\L{{\mathscr L}}
\newcommand\cO{{\mathcal O}}
\newcommand\PP{{\mathscr P}}
\newcommand\QQ{{\mathcal Q}}
\newcommand\RR{{\mathcal R}}
\newcommand\cS{{\mathcal S}}
\newcommand\sT{{\mathscr T}}
\newcommand\cV{{\mathcal V}}
\newcommand\cZ{{\mathcal Z}}

\newcommand\VV{{\boldsymbol V}}
\newcommand\WW{{\boldsymbol W}}
\newcommand\XX{{\boldsymbol X}}

\newcommand{\ba}{\boldsymbol{a}}
\newcommand{\bA}{\boldsymbol{\mathsf{A}}}
\newcommand{\bB}{\boldsymbol{\mathsf{B}}}
\newcommand{\bb}{\boldsymbol{b}}
\newcommand{\bc}{{\boldsymbol{c}}}
\newcommand{\be}{{\boldsymbol{\mathsf{E}}}}
\newcommand{\bF}{{\boldsymbol{\mathsf{F}}}}
\newcommand{\bff}{\boldsymbol{f}}
\newcommand{\bK}{\boldsymbol{\mathsf{K}}}
\newcommand{\bM}{\boldsymbol{\mathsf{M}}}
\newcommand{\bn}{\boldsymbol{\mathsf{N}}}
\newcommand{\bp}{\boldsymbol{\mathsf{P}}}
\newcommand{\bu}{\boldsymbol{\mathsf{u}}}
\newcommand{\bv}{\boldsymbol{\mathsf{v}}}
\newcommand{\bX}{\boldsymbol{\mathsf{X}}}
\newcommand{\bV}{\boldsymbol{\mathsf{V}}}
\newcommand{\bW}{\boldsymbol{\mathsf{V}}^\varepsilon}
\newcommand{\bt}{\boldsymbol{t}}
\newcommand{\bw}{\boldsymbol{w}}
\newcommand{\bx}{{\boldsymbol{x}}}
\newcommand{\by}{{\boldsymbol{y}}}
\newcommand{\bz}{{\boldsymbol{z}}}

\newcommand{\rb}{{\mathsf{b}}}
\newcommand{\rc}{{\mathsf{c}}}
\newcommand{\rB}{{\mathsf{B}}}
\newcommand{\rD}{{\mathsf{D}}}
\newcommand{\sh}{{g}}
\newcommand{\rh}{{\mathsf{g}}}
\newcommand{\rH}{{\mathsf{H}}}
\newcommand{\rK}{{\mathsf{K}}}
\newcommand{\rM}{{\mathsf{M}}}
\newcommand{\rL}{{\mathsf{L}}}
\newcommand{\ru}{{\mathsf{u}}}
\newcommand{\rv}{{\mathsf{v}}}
\newcommand{\rR}{{\mathsf{R}}}
\newcommand{\rT}{{\mathsf{T}}}
\newcommand{\rU}{{\mathsf{U}}}
\newcommand{\rV}{{\mathsf{V}}}
\newcommand{\rW}{{\mathsf{W}}}

\newcommand{\tbu}{\bu^*}
\newcommand{\tru}{\tilde\ru}
\newcommand{\kk}{[k]}

\renewcommand\H{{\boldsymbol H}}

\newcommand\ttau{\boldsymbol{\tau}}
\newcommand\zetaz{\boldsymbol{\zeta}}
\newcommand\gammag{\boldsymbol{\gamma}}
\newcommand\rhor{\boldsymbol{\rho}}
\newcommand\zetazt{\zetaz^*}
\newcommand\etat{\eta^*}
\newcommand\zetat{\zeta^*}

\newcommand\sep{\hskip-0.2ex,\hskip0.15ex}

\newcommand {\gA}{\mathfrak{A}}
\newcommand {\gB}{\boldsymbol{\mathfrak{B}}}
\newcommand {\gC}{\mathfrak{C}}
\newcommand {\gE}{\mathfrak{E}}
\newcommand {\gF}{\mathfrak{F}}
\newcommand {\gK}{{\boldsymbol{\mathfrak K}}}
\newcommand {\gL}{\boldsymbol{\mathfrak L}}
\newcommand {\gM}{{\boldsymbol{\mathfrak M}}}
\newcommand {\gP}{{\mathfrak P}}
\newcommand {\gR}{{\mathfrak R}}
\newcommand {\gS}{\sigma}
\newcommand {\gU}{{\mathfrak U}}
\newcommand {\gW}{{\mathfrak W}}
\newcommand {\gZ}{{\mathfrak Z}}

\newcommand {\hm}{\mathsf{h}}
\newcommand {\am}{\mathsf{a}}
\newcommand {\mm}{\mathsf{m}}
\newcommand {\hB}{\breve{\mathfrak B}}
\newcommand {\hM}{\mathfrak{H}}
\newcommand {\hK}{\mathfrak{A}}

\newcommand {\aB}{\rB}
\newcommand {\aL}{\rL}
\newcommand {\aK}{\rK}
\newcommand {\aM}{\rM}

\newcommand{\comp}{\mathrm{comp}}
\newcommand{\diam}{\mathrm{diam}}
\newcommand{\dist}{\mathrm{dist}}
\newcommand{\meas}{\mathrm{meas}}
\newcommand{\cor}{\mathsf{cor}}

\newcommand{\ocirc}{\overset{\circ}}

\newcommand{\lamone}[1]{\lambda_1\big[#1\big]}
\newcommand{\zetone}[1]{\eta_0\big[#1\big]}
\newcommand{\near}[1]{\lfloor#1\rceil}

\usepackage{soul}

\begin{document}

\title[Free vibrations of axisymmetric shells\,:
parabolic and elliptic cases]{\bf Free vibrations of axisymmetric shells\,:\\[0.5em]
parabolic and elliptic cases }  

\author{Marie Chaussade-Beaudouin}

\author{Monique Dauge}

\author{Erwan Faou}

\author{Zohar Yosibash}

\address{MC, MD \& EF: Irmar, (Cnrs, Inria) Universit\'{e} de Rennes 1, \hfill\break
Campus de Beaulieu,
35042 Rennes Cedex, France}

\address{ZY: Ben-Gurion University of the Negev, Dept. of Mechanical Engineering, \hfill\break
             POBox 653, Beer-Sheva 84105, Israel}

\email{\href{mailto:monique.dauge@univ-rennes1.fr}{monique.dauge@univ-rennes1.fr}}
\urladdr{\href{http://perso.univ-rennes1.fr/monique.dauge/}{http://perso.univ-rennes1.fr/monique.dauge/}}

\email{\href{mailto:Erwan.Faou@inria.fr}{Erwan.Faou@inria.fr}}
\urladdr{\href{http://www.irisa.fr/ipso/perso/faou/}{http://www.irisa.fr/ipso/perso/faou/}}

\email{\href{mailto:zohary@bgu.ac.il}{zohary@bgu.ac.il}}
\urladdr{\href{http://www.bgu.ac.il/~zohary/}{http://www.bgu.ac.il/~zohary/}}

\date{9 December 2016}



\begin{abstract}
Approximate eigenpairs (quasimodes) of axisymmetric thin elastic domains with laterally clamped boundary conditions (Lam\'e system) are determined by an asymptotic analysis as the thickness ($2\varepsilon$) tends to zero. The departing point is the Koiter shell model that we reduce by asymptotic analysis to a scalar model
that depends on two parameters: the angular frequency $k$ and the half-thickness $\varepsilon$. Optimizing $k$ for each chosen $\varepsilon$, we find power laws for $k$ in function of $\varepsilon$ that provide the smallest eigenvalues of the scalar reductions.
Corresponding eigenpairs generate quasimodes for the 3D Lam\'e system by means of several reconstruction operators, including boundary layer terms. Numerical experiments demonstrate that in many cases the constructed eigenpair corresponds to the first eigenpair of the Lam\'e system.

Geometrical conditions are necessary to this approach: The Gaussian curvature has to be nonnegative and the azimuthal curvature has to dominate the meridian curvature in any point of the midsurface. In this case, the first eigenvector admits progressively larger oscillation in the angular variable as $\varepsilon$ tends to $0$. Its angular frequency exhibits a power law relation
of the form $k=\gamma \varepsilon^{-\beta}$ with $\beta=\frac14$ in the parabolic case (cylinders and trimmed cones), and the various $\beta$s  $\frac25$, $\frac37$, and $\frac13$ in the elliptic case.
For these cases where the mathematical analysis is applicable, numerical examples that illustrate the theoretical results are presented.
\end{abstract}

\maketitle


\section{Introduction}
Shells are three-dimensional thin objects widely addressed in the literature in mechanics, engineering as well as in mathematics. According to any classical definition, a shell is determined by its midsurface $\cS$ and a thickness parameter $\varepsilon$:  The shell denoted by $\Omega^\varepsilon$ is obtained by thickening $\cS$ on either side by $\varepsilon$ along unit normals to $\cS$. Like most of references, we assume that $\Omega^\varepsilon$ is made of a linear homogeneous isotropic material and we furthermore consider clamped boundary conditions along its lateral boundary.

In this paper, the behavior of the fundamental vibration mode of such a shell is investigated as $\varepsilon$ tends to $0$. We consider {\em free vibration modes}, that is, eigenpairs $(\lambda,\bu)$ of the 3D Lam\'e system $\gL$ in $\Omega^\varepsilon$ complemented by suitable boundary conditions. Here $\lambda$ is the square of the eigenfrequency and $\bu$ the eigen-displacement. The thin domain limit $\varepsilon\to0$ pertains to ``shell theory''. 

Shell theory consists of finding surface models, i.e., systems of equations posed on $\cS$, approximating the 3D Lam\'e system $\gL$ on $\Omega^\varepsilon$ when $\varepsilon$ tends to $0$. This approach was started for plates (the case when $\cS$ is flat) by Kirchhoff, Reissner and Mindlin see for instance \cite{Kirchhoff1,Reissner1,Mindlin} respectively. When the structure is a genuine shell for which the midsurface has nonzero curvature, the problem is even more difficult and was first tackled in the seminal works of Koiter, John, Naghdi and Novozhilov in the sixties \cite{Koiter5,Koiter2,Koiter3,John2,Novozhilov,Naghdi1}. A large literature developed afterwards aimed at laying more rigorous mathematical bases to shell theory see for instance the works of Sanchez-Palencia, Sanchez-Hubert \cite{SanchezPalencia3,SanchezPalencia4,SanchezPalencia5,SanchezHubertSanchezPalencia2}, Ciarlet, Lods, Mardare, Miara \cite{CiarletLods7,CiarletLods8,CiarletLods9,LodsMardare4} and the book \cite{Ciarlet10}, and more recently Dauge, Faou \cite{Faou1,Faou2,DaFa1}. Most of  these  works apply to the static problem,  and the results strongly depend on the geometrical nature of the shell (namely parabolic, elliptic or hyperbolic according to the Gaussian curvature $K$ of $\cS$ being zero, positive or negative).

Much fewer works were devoted to free vibrations of thin shells. Plates were addressed beforehand, see \cite{CiarletKesavan,DDFR99}. For shells and more general thin structures, let us quote Soedel \cite{SoedelEncyclopedia,SoedelBook}. To the best of our knowledge, theoretical works devoted to the asymptotic analysis of eigenmodes in thin elastic shells were associated with a surface model, such as the Koiter model.  

Recall that the Koiter model \cite{Koiter5,Koiter2} takes the form: 
\begin{equation}
\label{chopin}
    \gK(\varepsilon) = \gM + \varepsilon^2 \gB,
\end{equation}
where $\gM$ is the membrane operator, $\gB$ the bending operator, and $\varepsilon$ the half-thickness of the shell.
These two operators are $3\times3$ systems posed on $\cS$, 
acting on 3-component vector fields $\zetaz$. When these fields are represented in surface fitted components $\zeta_\alpha$ and $\zeta_3$ (the tangential and normal components), these two operators display special structures. For plates, they uncouple: $\gM$ amounts to a $2\times2$ Lam\'e system acting on tangential components $\zeta_\alpha$ and $\gB$ is a multiple of the biharmonic operator $\Delta^2$ acting on the sole normal component $\zeta_3$. For general shells, the membrane operator $\gM$ is of order 2 on tangential components $\zeta_\alpha$, but of order $0$ on the normal component $\zeta_3$. The bending operator $\gB$ has a complementary role: It is order $4$ on $\zeta_3$.

In \cite{SanchezHubertSanchezPalencia2}, the essential spectrum of the membrane operator $\gM$ (the set of $\lambda$'s such that $\gM-\lambda$ is not Fredholm)  was characterized in the elliptic, parabolic, and hyperbolic cases.  The series of papers by Artioli, Beir{\~a}o Da Veiga, Hakula and Lovadina \cite{BeiraoLovadina2008,ArtioliBeiraoHakulaLovadina2008,ArtioliBeiraoHakulaLovadina2009} investigated the first eigenvalue of models like $\gK(\varepsilon)$. Effective results hold for {\em axisymmetric}  shells with clamped lateral boundary: Defining the order $\alpha$ of a positive function $\varepsilon\mapsto\lambda(\varepsilon)$, continuous on $(0,\varepsilon_0]$, by the conditions
\begin{equation}
\label{order}
   \forall\eta>0,\quad
   \lim_{\varepsilon\to0^+} \lambda(\varepsilon)\,\varepsilon^{-\alpha+\eta} = 0
   \quad\mbox{and}\quad
   \lim_{\varepsilon\to0^+} \lambda(\varepsilon)\,\varepsilon^{-\alpha-\eta} = \infty
\end{equation}
they proved that $\alpha=0$ in the elliptic case, $\alpha=1$ for parabolic case, and $\alpha=\frac23$ in the hyperbolic case.

\subsection{Axisymmetric shells}
Besides their natural interest in structural mechanics, isotropic axisymmetric shells have the nice property that all 3D Lam\'e eigenpairs $(\lambda,\bu)$ can be classified by their {\em azimuthal frequency} $k$ (aka {\em angular frequency}). Indeed, the 3D Lam\'e system $\gL$ as well as the membrane and bending operators $\gM$ and $\gB$ can be diagonalized by Fourier decomposition with respect to the azimuthal angle $\varphi$, see \cite{ChDaFaYo16b} for example. So, in particular, the azimuthal frequency $k(\varepsilon)$ of the first eigenvector makes sense. Based on some analytical calculations it was known that $k(\varepsilon)$ may have a non trivial behavior: Quoting W.\ Soedel \cite{SoedelEncyclopedia}
{\em ``[We observe] a phenomenon which is particular to many deep shells, namely that the lowest natural frequency does not correspond to the simplest natural mode, as is typically the case for rods, beams, and plates.''} In other words, $k(\varepsilon)$ is not zero as it would be for a simpler operator like the Laplacian, see also \cite{ChDaFaYo16b}.

For axisymmetric shells Beirao et al.\ and Artioli et al. \cite{BeiraoLovadina2008,ArtioliBeiraoHakulaLovadina2008,ArtioliBeiraoHakulaLovadina2009} investigated by numerical simulations the azimuthal frequency $k(\varepsilon)$ of the first eigenvector of $\gK(\varepsilon)$: Like in the phenomenon of {\em sensitivity} \cite{SanchezHubertSanchezPalencia2}, the lowest eigenvalues are associated with eigenvectors with growing angular frequencies and $k(\varepsilon)$ exhibits a negative power law of type $\varepsilon^{-\beta}$, for which \cite{ArtioliBeiraoHakulaLovadina2009} identifies the exponents $\beta=\frac14$ for cylinders (see also \cite{BeiraoHakulaPitkaranta2008} for some theoretical arguments), $\beta=\frac25$ for a particular family of elliptic shells, and $\beta=\frac13$ for another particular family of hyperbolic shells.

Similarly to the aforementioned publications, we consider here axisymmetric shells whose mid-surface $\cS$ is parametrized by a smooth positive function $f$ representing the radius as a function of the axial variable:
\begin{equation}
\label{eq:surfpara}
\begin{array}{ccccc}
    F : & \cI\times\T &\longrightarrow &\cS \\
    & (z,\varphi)  &\longmapsto& (f(z)\cos\varphi, \,f(z)\sin\varphi,\,z).
\end{array}
\end{equation}
Here $\cI$ is the parametrization interval and $\T$ is the torus $\R/2\pi\Z$.
\begin{figure}[ht]
\includegraphics[scale=0.5]{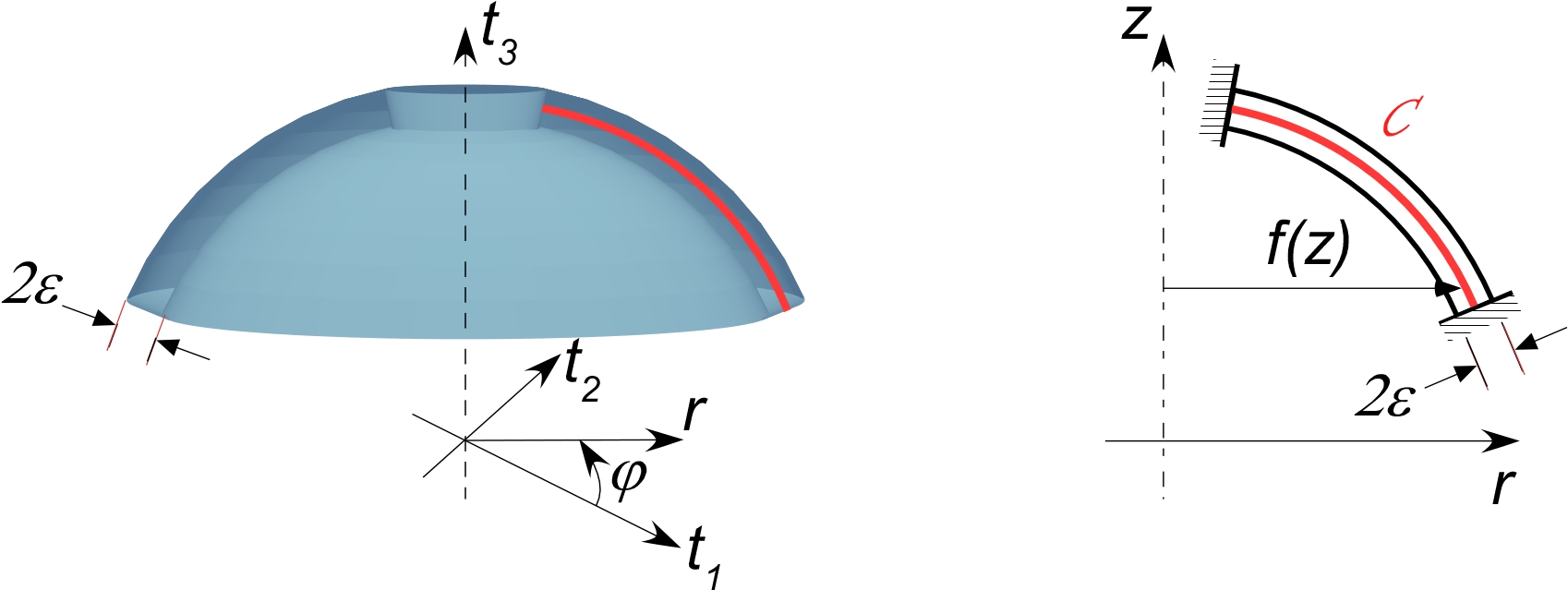}
   \caption{Axisymmetric shell $\Omega^\varepsilon$ with Cartesian and cylindrical coordinates (left) and 
   the meridian domain 
   $\omega^\varepsilon$ with its midcurve $\cC$ parametrized by the equation $r=f(z)$ (left).}
\label{f.axiShellNotations}
\end{figure}
In  Figure \ref{f.axiShellNotations} are represented an instance of 3D shell $\Omega^\varepsilon$, together with its meridian domain $\omega^\varepsilon$. The 2D domain $\omega^\varepsilon$ has the meridian set $\cC$ of the midsurface $\cS$ as meridian curve.

We focus on cases when sensitivity may show up, i.e., when the azimuthal frequency $k(\varepsilon)$ of the first eigenvector is likely to tend to infinity as the thickness tends to $0$. As will be shown, the rules driving this phenomenon are far to be  straightforward, and depend in a non trivial manner on the geometry of the shell: In the sole elliptic case, we show that there exist at least three distinct power laws for $k(\varepsilon)$. This is the expression of some {\em bending} effects and  may sound as a paradox since for elliptic shells the membrane is an elliptic system in the sense of Agmon, Douglis and Nirenberg \cite{ADN2}, see \cite{Genevey}.
However, there also exist elliptic shells for which $k(\varepsilon)$ remains constant, see the computations for a spherical cap in \cite{DaFaYo}.

\subsection{High frequency analysis}
Our departing point is a high frequency analysis (in $k$) of the membrane operator $\gM$ on surfaces $\cS$ with a parametrization of type \eqref{eq:surfpara}. By the Fourier decomposition naturally induced by the cylindrical symmetry, we define $\gM^k$  as  the membrane operator acting at the frequency $k\in\N$ and we perform a {\em scalar reduction} of the eigenproblem by a special factorization in a formal series algebra in powers of the small parameter $\frac{1}{k}$.  This mathematical tool, developed  for cylindrical shells  in the PhD thesis \cite{beaudouin:tel} of the first author, reduces the original eigenproblem (which is a $3\times3$ system) to a {\em scalar} eigenproblem posed on the transverse component of the displacement. This way, we can construct in a variety of parabolic and elliptic cases a new explicit scalar differential operator $\hM^k$  whose first eigenvalue $\lamone{\hM^k}$ has a computable asymptotics as $k\to\infty$
\begin{equation}
\label{eq:mu1hMk}
   \lamone{\hM^k} = \hm_0 + \hm_1 k^{-\eta_1} + \cO(k^{-\eta_2}),
   \quad 0<\eta_1<\eta_2.
\end{equation}
In \eqref{eq:mu1hMk} all coefficients and exponents depend on shell's geometry, i.e.\ on the function $f$ in \eqref{eq:surfpara}. This leads to a {\em quasimode construction} for $\gM^k$ that is valid for all parabolic shells of type \eqref{eq:surfpara} and all elliptic shells with azimuthal curvature dominating.
The operator $\hM^k$ strongly depends on the nature of the shell:
\begin{equation}
\label{eq:hMk}
   \begin{cases}
   \hM^k = k^{-4}\rH_4 & \mbox{in the parabolic case (i.e., when $f''=0$)}, \\
   \hM^k = \rH_0 + k^{-2}\rH_2 & \mbox{in the elliptic case (i.e., when $f''<0$)}. \\
   \end{cases}
\end{equation}
with explicit operators $\rH_0$, $\rH_2$ and $\rH_4$, cf.\ sect.\ \ref{ss:memred} for formulas. Let us mention at this point that for hyperbolic shells such a suitable scalar reduction $\hM^k$ cannot be found. 

This membrane scalar reduction induces a {\em Koiter-like scalar reduced operator}  $\hK(\varepsilon)$ for the shell that we define at the frequency $k$ by
\begin{equation}
\label{eq:hKk}
   \hK^k(\varepsilon) = \hM^k + \varepsilon^2k^4\rB_0
\end{equation}
where the function $\rB_0$ is positive and explicit ($k^4$ corresponding to the leading order in the Fourier expansion of the bending operator $\gB$).
Then the lowest eigenvalue of $\hK(\varepsilon)$ is the infimum on all angular frequencies of the first eigenvalues of $\hK^k(\varepsilon)$:
\begin{equation}
\label{eq:mu1hH}
   \lamone{\hK(\varepsilon)} = \inf_{k\in\N} \lamone{ \hK^k(\varepsilon)}.
\end{equation}
In all relevant parabolic cases (i.e., cylinders and cones) and a variety of elliptic cases, we prove in this paper:
\begin{itemize}
\item[{\em (i)}] The infimum in \eqref{eq:mu1hH} is reached for $k=\near{k(\varepsilon)}$, the nearest integer from $k(\varepsilon)$, with $k(\varepsilon)$ satisfying a power law of the form
\begin{equation}
\label{eq:keps}
   k(\varepsilon) = \gamma \varepsilon^{-\beta} + \cO(\varepsilon^{-\beta'}),
   \quad 0\le\beta'<\beta,
\end{equation}
with $\beta$ depending only on $f$ and $\gamma$ positive. The exponent $\beta$ is calculated so to equilibrate $k^{-\eta_1}$ (cf.\ \eqref{eq:mu1hMk})  and $\varepsilon^2k^4\equiv k^{4-2/\beta}$, which yields:
\begin{equation}
\label{eq:beta}
   \framebox{$\displaystyle\beta = \frac{2}{4+\eta_1}$}
\end{equation}
\item[{\em (ii)}] The smallest eigenvalue of the reduced scalar model $\hK(\varepsilon)$ has an asymptotic expansion of the form, as $\varepsilon\to0$
\begin{equation}
\label{eq:mu1hK}
   \lamone{\hK(\varepsilon)} = \am_0 + \am_1 \varepsilon^{\alpha_1} + \cO(\varepsilon^{\alpha_2}),
   \quad 0<\alpha_1<\alpha_2,
\end{equation}
where $\am_0$ coincides with the coefficient $\hm_0$ present in \eqref{eq:mu1hMk} and $\alpha_1$ is given by the formula (replace $k$ with $\varepsilon^{-\beta}$ into the term $k^{-\eta_1}$ in \eqref{eq:mu1hMk})
\begin{equation}
\label{eq:alpha}
   \framebox{$\displaystyle\alpha_1 = \eta_1\beta = \frac{2\eta_1}{4+\eta_1}$}
\end{equation}
\item[{\em (iii)}] The corresponding eigenvector $\zetone{\hK(\varepsilon)}$ has a multiscale expansion in variables $z$ and $\varphi$ that involves 1 or 2 scales in $z$ (including or not boundary layers), depending on the parametrization $f$, i.e.\ on the geometry of $\cS$.
\end{itemize}

Once the asymptotic expansions for the Koiter scalar reduced operator $\hK(\varepsilon)$ is resolved we
construct  quasimodes for the full Koiter model $\gK(\varepsilon)$. Then, by energy estimates linking surfacic and 3D models similar to those of \cite{DaFa1}, we find a sort of {\em quasi-eigenvector} $\bu^\varepsilon$ whose Rayley quotient provides an asymptotic upper bound $\mm_1(\varepsilon)$ for the first eigenvalue $\lamone{\gL(\varepsilon)}$ of the 3D Lam\'e system $\gL$ in the shell $\Omega^\varepsilon$. This upper bound is given by the first two terms in \eqref{eq:mu1hK}:
\begin{equation}
\label{eq:quasi}
   \mm_1(\varepsilon) = \am_0 + \am_1 \varepsilon^{\alpha_1}.
\end{equation}
To make the analysis more complete, we perform numerical simulations. They aim at comparing true eigenpairs with quasimodes $(\mm_1(\varepsilon),\bu^\varepsilon)$.
To this end computations are performed at three different levels:
\begin{itemize}
\item[(1D)] We calculate $\am_0$, $\am_1$ of \eqref{eq:mu1hK} and $\gamma$ of \eqref{eq:keps}. We either use explicit analytical formulas when available, or compute numerically the spectrum of the 1D scalar reduced operators $\hK^k(\varepsilon)$ through a 1D finite element method applied to an auxiliary operator.
\item[(2D)] The Fourier decomposition of the 3D Lam\'e system $\gL$ in the shell $\Omega^\varepsilon$ provides a family $\gL^k$, $k\in\N$, of $3\times3$ systems posed on the 2D meridian domain $\omega^\varepsilon$.
We discretize these systems by a 2D finite element method in $\omega^\varepsilon$ for collections of integers $k\in\{0,1,\ldots, K_\varepsilon\}$ depending on the thickness $\varepsilon$, and compute the lowest eigenvalue $\lamone{\gL^k(\varepsilon)}$. This procedure provides an approximation of $\lamone{\gL(\varepsilon)}$ and of $k(\varepsilon)$ through the formula
\[
   \lamone{\gL(\varepsilon)} = \min_{k=0}^{K_\varepsilon} \lamone{\gL^k(\varepsilon)} \quad\mbox{and}\quad
   k(\varepsilon) = \arg\min_{k=0}^{K_\varepsilon}\lamone{\gL^k(\varepsilon)}.
\]
This method is a Fourier spectral discretization of the 3D problem. Note that in \cite{ArtioliBeiraoHakulaLovadina2009} a 1D Fourier spectral method is used for the discretization of the surfacic Koiter and Naghdi models.
\item[(3D)] We compute the first eigenvalue $\lamone{\gL(\varepsilon)}$ of the 3D Lam\'e system $\gL$ in the shell using directly a 3D finite element method in $\Omega^\varepsilon$.
\end{itemize}
This combination of simulations show that, in a number of cases, the theoretical quasimode $(\mm_1(\varepsilon),\bu^\varepsilon)$ is a good approximation of the true first eigenpair of $\gL(\varepsilon)$. 

\subsection{Specification in the parabolic and elliptic cases}
In the Lam\'e system we use the engineering notations of the material parameters: $E$ is the Young modulus and $\nu$ is the Poisson ratio. The shells to which our analysis apply are uniquely defined by the function $f$ and the interval $\cI$ in \eqref{eq:surfpara}. The inverse parametrization (the axial variable function of the radius) would not provide distinct cases where our analysis is applicable. 

$\bullet$ The parabolic cases are those for which $f''=0$ on $\cI$. So $f$ is affine. The midsurface $\cS$ is developable. We classify parabolic cases in two types:
\begin{enumerate}
\item\textbf{\textit{`Cylinder'}} $f$ is constant;
\item\textbf{\textit{`Cone'}}  $f$ is affine and not constant.
\end{enumerate}

$\bullet$ The elliptic cases are those for which $f''<0$ on $\cI$. To conduct our analysis, we assume moreover that the azimuthal curvature dominates the meridian curvature (admissible cases), which amounts to
\begin{equation}
\label{eq:azi}
    1+f'^2+ff'' \geq   0.
\end{equation}
We discriminate admissible elliptic cases by the behavior of the function $\rH_0$ that is the first term of the scalar reduction $\hM^k$, cf.\ \eqref{eq:hMk},
\begin{equation}
\label{eq:H0}
   \rH_0 = E \frac{f''^2}{(1+f'^2)^3}\,,
\end{equation}
classifying them in three generic types:
\begin{enumerate}
\item \textbf{\textit{`Toroidal'}} $\rH_0$  is constant.
\item \textbf{\textit{`Gauss'}} $\rH_0$ is not constant and reaches its minimum at $z_0$ inside $\cI$ and not on its boundary $\partial\cI$, with the exception of cases for which $\rH''_0$ or $1+f'^2+ff''$ are zero at $z_0$.
\item \textbf{\textit{`Airy'}} $\rH_0$ is not constant and reaches its minimum at $z_0$ in the boundary $\partial\cI$,  with the exception of cases for which $\rH'_0$ or $1+f'^2+ff''$ are zero at $z_0$.
\end{enumerate}

\begin{figure}[ht]
\includegraphics[scale=0.15]{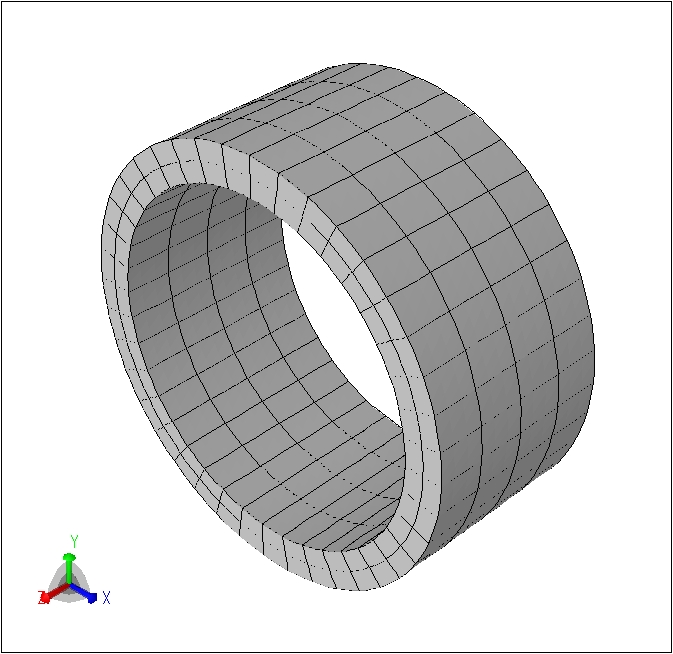}\hskip-1em%
\includegraphics[scale=0.15]{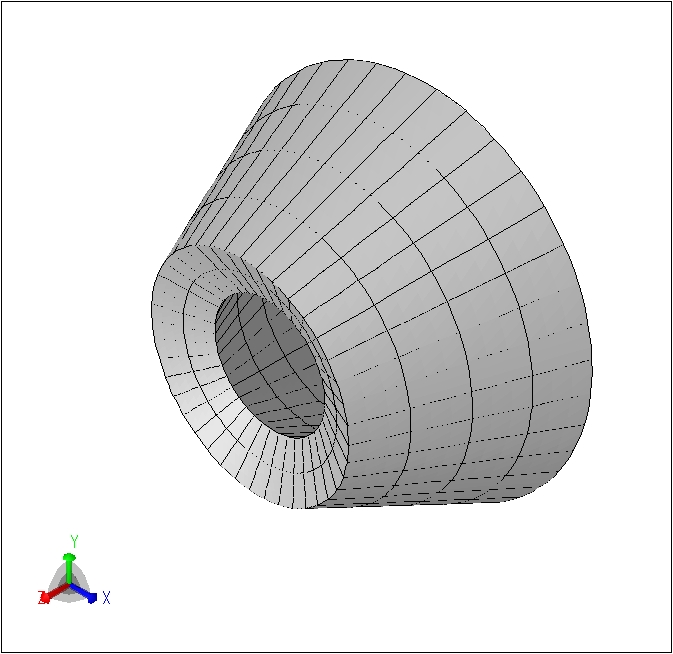}\hskip-1em%
\includegraphics[scale=0.15]{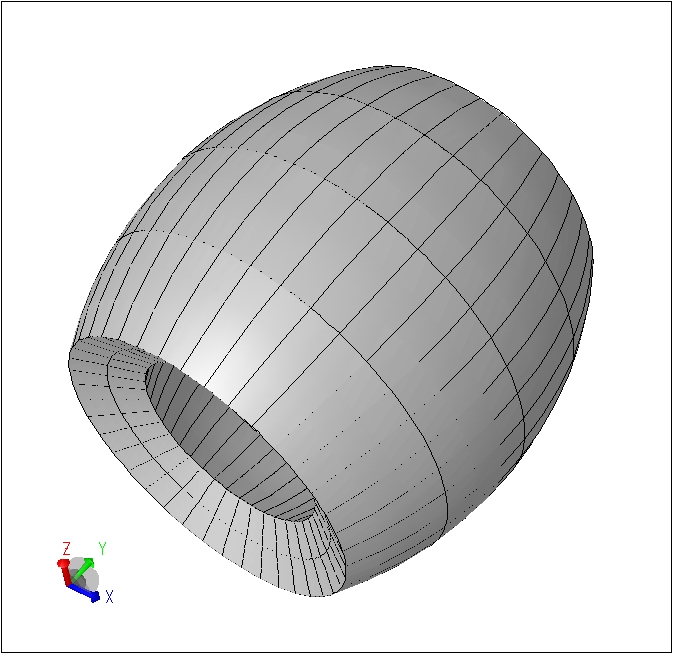}\hskip-1em%
\includegraphics[scale=0.15]{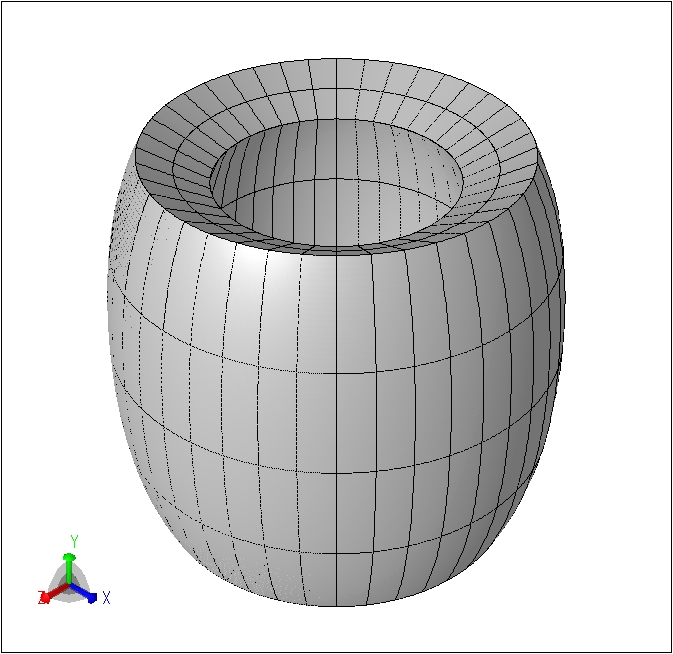}\hskip-1em%
\includegraphics[scale=0.15]{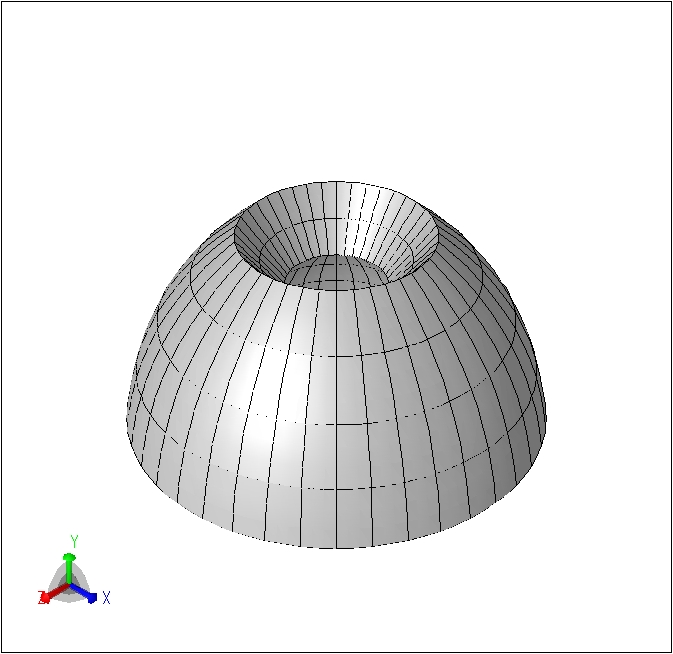}%
   \caption{The five models A, B, D, H, L, used for computations (here $\varepsilon=0.2$).}
\label{f:models}\end{figure}

We summarize  in Table \ref{tab:1} our main theoretical results on the exponents $\eta_1$, $\beta$, $\alpha_1$, on the azimuthal frequency $k(\varepsilon)$, and on the quasi-eigenvalue (qev) $\mm_1(\varepsilon)$. The exponents $\alpha$ of \cite{ArtioliBeiraoHakulaLovadina2008} are confirmed ($1$ in the parabolic cases and $0$ in the elliptic cases). Inspired by \cite{ArtioliBeiraoHakulaLovadina2009}, we mention in the table the factor $\rR$ representing the ratio (Bending Energy)/(Total Energy). This ratio is asymptotically represented by, cf \eqref{eq:hKk}
\begin{equation}
\label{eq:Rratio}
   \rR = \frac{\varepsilon^2k^4\langle \rB_0\eta_0,\eta_0\rangle}
   {\langle\hK^k(\varepsilon) \eta_0,\eta_0\rangle}
   \quad\mbox{for $k=k(\varepsilon)$ and $\eta_0$ the corresponding eigenvector of $\hK^k(\varepsilon)$.}
\end{equation}
The names of models used for numerical simulations are also mentioned in this table, whereas in Figure \ref{f:models} we represent these models in their 3D version for $\varepsilon=0.2$.
\begin{table}
{\small
\begin{tabular}{l@{\quad}ll@{\quad}l@{\quad}llllll}
\hline
Type (Model)& $\eta_1$ & $\beta$   & $\alpha_1$ & $\am_0$ & $\am_1$
            & $k(\varepsilon)$ 
            & $\mm_1(\varepsilon)$ & $\rR$ \\
\hline
\underline{\sc Parabolic} \\
`Cylinder' (A) & $4$        & $\frac14$ & $1$        & $0$   & explicit wrt 1D ev's
            & $\gamma\varepsilon^{-1/4}$ 
            & $\am_1\varepsilon$    & $\frac12$  \\
`Cone' (B)  & $4$        & $\frac14$ & $1$        & $0$   & optimization of 1D ev's
            & $\gamma\varepsilon^{-1/4}$ 
            & $\am_1\varepsilon$    & $\frac12$  \\
\underline{\sc Elliptic} \\
`Toroidal' (D) & $2$        & $\frac13$ & $\frac23$  & $\rH_0$        & optimization of 1D ev's
            & $\gamma\varepsilon^{-1/3}$ 
            & $\am_0 + \am_1\varepsilon^{2/3}$   & $\delta\varepsilon^{2/3}$ \\
`Gauss' (H) & $1$        & $\frac25$ & $\frac25$  & $\rH_0(z_0)$   & explicit
            & $\gamma\varepsilon^{-2/5}$ 
            & $\am_0 + \am_1\varepsilon^{2/5}$   & $\delta\varepsilon^{2/5}$  \\
`Airy' (L)  & $\frac23$  & $\frac37$ & $\frac27$  & $\rH_0(z_0)$   & explicit
            & $\gamma\varepsilon^{-3/7}$ 
            & $\am_0 + \am_1\varepsilon^{2/7}$   & $\delta\varepsilon^{2/7}$ \\
\hline\\
\end{tabular}
}
\caption{Summary of exponents $\eta_1$, $\beta$, $\alpha_1$, frequency $k(\varepsilon)$, qev $\mm_1(\varepsilon)$ and ratio of energies $\rR$ \eqref{eq:Rratio}. Coefficients $\gamma$ and $\delta$ are determined by the 1D reduction.}
\label{tab:1}
\end{table}

\subsection{Overview of main notation. Plan of the paper}
To relieve the complexity of notation, we gather here some definitions relating to coordinate systems, operators, and spectrum, before presenting the plan of the paper.

\subsubsection{Coordinates}
We use three systems of coordinates:

$\bullet$ Cartesian coordinates $\bt=(t_1,t_2,t_3)\in\R^3$ with coordinate vectors $\be_{t_1}$, $\be_{t_2}$, $\be_{t_3}$.

$\bullet$ Cylindrical coordinates $(r,\varphi,\tau)\in\R^+\times\T\times\R$ related to Cartesian coordinates by relations
\begin{equation}
\label{eq:cyl1}
   (t_1,t_2,t_3) = \sT(r,\varphi,\tau) \quad\mbox{with}\quad
   t_1 = r\cos\varphi,\ \ t_2 = r\sin\varphi,
   \ \  t_3=\tau\,.
\end{equation}
The coordinate vectors associated with the transformation $\sT$ are $\be_r=\partial_r\sT$, $\be_\varphi=\partial_\varphi\sT$, and $\be_\tau=\partial_\tau\sT$. We have
\begin{equation}
\label{eq:cyl2}
   \be_r = \be_{t_1}\cos\varphi + \be_{t_2}\sin\varphi,\quad
   \be_\varphi = -r\be_{t_1}\sin\varphi + r\be_{t_2}\cos\varphi,
   \quad\mbox{and}\quad \be_\tau=\be_{t_3}\,.
\end{equation}

$\bullet$ Normal coordinates $(x_1,x_2,x_3)$, specified as $(z,\varphi,x_3)$ in our case. Such coordinates are related to the surface $\cS$ and a chosen unit normal field $\bn$ to $\cS$. The variable $x_3$ is the coordinate along $\bn$. The variables $(x_1,x_2)$, specified as $(z,\varphi)$ in our case, parametrize the surface. The full transformation $\F:(z,\varphi,x_3)\mapsto(t_1,t_2,t_3)$ sends the product $\cI\times\T\times(-\varepsilon,\varepsilon)$ onto the shell $\Omega^\varepsilon$ and is explicitly given by
\begin{equation}
\label{eq:para}
   t_1 = \big(f(z) + x_3\ \tfrac{1}{s(z)} \big) \cos\varphi,\quad
   t_2 = \big(f(z) + x_3\ \tfrac{1}{s(z)} \big) \sin\varphi,\quad
   t_3 = z - x_3 \ \tfrac{f'(z)}{s(z)},
\end{equation}
where $s=\sqrt{1+f'^2}$. The restriction of $\F$ on the surface $\cS$ (corresponding to $x_3 = 0$) gives back $F$ \eqref{eq:surfpara}. The coordinate vectors associated with the transformation $\F$ are $\partial_z\F=:\be_z$, $\partial_\varphi\sT$ that coincides with $\be_\varphi$ above, and $\partial_3\F=:\be_3$. On the surface $\cS$, $x_3=0$ and $\be_3$ coincides with $\bn$, whereas $\be_z$ and $\be_\varphi$ are tangent to $\cS$.

These three systems of coordinates determine the contravariant components of a displacement $\bu$ in each of these  systems by identities
\begin{equation}
\label{eq:contra}
   \bu = \ru^{t_1}\be_{t_1} + \ru^{t_2}\be_{t_2} + \ru^{t_3}\be_{t_3}
   = \ru^{r}\be_r+\ru^{\varphi}\be_\varphi + \ru^{\tau}\be_\tau
   = \ru^{z}\be_z+\ru^{\varphi}\be_\varphi + \ru^{3}\be_3\,.
\end{equation}

The cylindrical and normal systems of coordinates are suitable for angular Fourier decomposition $\T\ni\varphi\mapsto k\in\Z$. The Fourier coefficient of rank $k$ of a function $\ru$ is denoted by $\ru^k$
\begin{equation}
\label{eq:Fou1}
   \ru^k = \frac{1}{2\pi} \int_{0}^{2\pi} \ru(\varphi) \,e^{- i k\varphi} \,\rd \varphi.
\end{equation}
For functions on $\Omega^\varepsilon$, the Fourier coefficients are defined on the meridian domain $\omega^\varepsilon\subset\R^2$ of $\Omega^\varepsilon$.

Concerning 3D displacements $\bu$ defined on $\Omega^\varepsilon$ or surface displacements $\zetaz$ defined on $\cS$, we have first to expand them in a suitable system of coordinates (cylindric or normal) and then calculate Fourier coefficients of their components, see \cite{ChDaFaYo16b}: for instance
\begin{equation}
\label{eq:Fou2}
\begin{gathered}
   \bu^k = (\ru^{r})^k\be_r+(\ru^{\varphi})^k\be_\varphi + (\ru^{\tau})^k\be_\tau
   \quad\!\mbox{with}\!\quad
   (\ru^{r})^k(r,\tau) = \frac{1}{2\pi} \int_{0}^{2\pi} 
   \ru^r(r,\varphi,\tau) \,e^{- i k\varphi} \,\rd \varphi,
   \ \mbox{...}\hskip-1em\\
   \zetaz^k = (\zeta^{z})^k\be_z+(\zeta^{\varphi})^k\be_\varphi + (\zeta^{3})^k\bn
   \quad\mbox{with}\quad
   (\zeta^{z})^k(z) = \frac{1}{2\pi} \int_{0}^{2\pi} \zeta^z(z,\varphi) \,e^{- i k\varphi} \,\rd \varphi,
   \ \mbox{...}
\end{gathered}
\end{equation}

\subsubsection{Operators}
\label{ss:op}
We manipulate a collection of operators and their Fourier symbols. The Lam\'e system $\gL$ acting on 3D displacements $\bu$ defined on the shell $\Omega^\varepsilon$ is particularized as $\gL(\varepsilon)$. After angular Fourier decomposition, we obtain the family of $3\times3$ operators $\gL^k(\varepsilon)$ defined on the meridian domain $\omega^\varepsilon$. On the surface $\cS$ we have the membrane, bending and Koiter operators $\gM$, $\gB$ and $\gK(\varepsilon)$. They act on 3-component surface displacements $\zetaz$. On the meridian curve $\cC$ of $\cS$, we have the corresponding families $\gM^k$, $\gB^k$ and $\gK^k(\varepsilon)$. Finally, on the meridian curve $\cC$, we have our scalar reductions $\hM^k$ and $\hK^k(\varepsilon)=\hM^k + \varepsilon^2k^4\rB_0$ acting on functions $\eta$. We go from a higher model to a lower one by reduction, and the converse way by reconstruction. For instance we go from $\bu$ to $\zetaz$ by restriction to $\cS$. The converse way uses the reconstruction operator $\rU$ \eqref{eq:U}. For any chosen integer $k$, we go from $\zetaz^k$ to $\eta^k$ by selecting the normal component of $\zetaz^k$. The converse way uses the reconstruction operators $\bV[k]$ that we will construct.

\subsubsection{Spectrum}
We denote by $\sigma(\gA)$ and $\sigma_{\sf ess}(\gA)$ the spectrum and the essential spectrum of a selfadjoint operator $\gA$, respectively, which means the set of $\lambda$'s such that $\gA-\lambda$ is not invertible and not Fredholm, respectively. If moreover, $\gA$ is non-negative we denote by $\lambda_1[\gA]$ its lowest eigenvalue.

\subsubsection{Outline} 
After the present introduction, we revisit in sect.~\ref{s:2a} the linear shell theory in general with a brief introduction of 3D (Lam\'e) and surfacic (Koiter, membrane, bending) problems, and in sect.~\ref{s:2} we particularize formulas for axisymmetric shells. In sect.~\ref{s:5} we set the principles of the high frequency analysis, in sect.~\ref{s:6} and \ref{s:7} we address more particularly the parabolic and elliptic cases, respectively. In sect.~\ref{s:8} we present numerical experiments addressing a model for each of the five main types described above. We conclude in sect.~\ref{s:9}. We provide in Appendix \ref{app:A} details on the factorization in formal series leading to the scalar reduction and in Appendix \ref{app:B} variational formulations in the meridian domain $\omega^\varepsilon$ of the Fourier operator coefficients $\gL^k$ of the 3D Lam\'e system.

\section{Essentials on shell theory}
\label{s:2a}

Recall that Cartesian coordinates of a point $\bp\in\R^3$ are denoted by $\bt=(t_1,t_2,t_3)$. A shell $\Omega^\varepsilon$ is a three-dimensional object defined by its midsurface $\cS$ and its thickness parameter $\varepsilon$ in the following way: We assume that $\cS$ is smooth and orientable, so that there exists a smooth unit normal field $\bp\mapsto\bn(\bp)$ on $\cS$ and so that for $\varepsilon>0$ small enough the following map is one to one and smooth
\begin{equation}
\label{1E1}
\begin{array}{cccc}
   \Phi : & \cS\times(-\varepsilon,\varepsilon) &\to &\Omega^\varepsilon \\
   & (\bp,x_3) &\mapsto& \bt=\bp+x_3\,\bn(\bp).
\end{array}
\end{equation}
The boundary of $\Omega^\varepsilon$ has two parts:
\begin{enumerate}
\item Its lateral boundary $\partial_0\Omega^\varepsilon := \Phi\big(\partial\cS\times(-\varepsilon,\varepsilon)\big)$,
\item The rest of its boundary (natural boundary) $\partial_1\Omega^\varepsilon := \partial\Omega^\varepsilon\setminus \partial_0\Omega^\varepsilon$.
\end{enumerate}

\subsection{3D vibration modes}
\label{ss:3}
On the domain $\Omega^\varepsilon$, we consider the Lam\'e operator associated with an isotropic and homogeneous material  with Young coefficient $E$ and Poisson coefficient $\nu$. This means that
the material tensor is given by
\begin{equation}
\label{eq:Amat}
   A^{ijk\ell} = \tfrac{E\nu}{(1+\nu)(1 - 2\nu)} \delta^{ij}\delta^{k\ell}
   + \tfrac{E}{2(1 + \nu)} (\delta^{ik}\delta^{j\ell} + \delta^{i\ell}\delta^{jk}).
\end{equation}
For {\em clamped} boundary conditions the variational space is
\begin{equation}
\label{eq:V}
V(\Omega^\varepsilon) := \{\bu = (\ru_{t_1},\ru_{t_2},\ru_{t_3}) \in H^1(\Omega^\varepsilon)^3 \, , \quad
 \bu = 0 \quad \mbox{on}\quad \partial_0\Omega^\varepsilon\}.
\end{equation}
For a given displacement field $\bu$
let $e_{ij}(\bu) = \frac12 (\partial_i \ru_{t_j} + \partial_j \ru_{t_i})$ be the strain tensor, where $\partial_i$ stands for the partial derivative with respect to $t_i$.
The Lam\'e energy scalar product between two displacements $\bu$ and $\tbu$ is given by
\begin{equation}
\label{eq:en3D}
   a^\varepsilon_{\aL}(\bu,\tbu) =
   \int_{\Omega^\varepsilon} A^{ijk\ell} e_{ij}(\bu) \,e_{k\ell}(\tbu) \, \rd \Omega^\varepsilon \,,
\end{equation}
using the summation convention of repeated indices.
The three-dimensional modal problem can be written in variational form as:
{\sl Find $(\bu, \lambda)$ in $V(\Omega^\varepsilon) \times \R$ with $\bu\neq0$ such that}
\begin{equation}
\label{3E2}
\forall\, \tbu \in V(\Omega^\varepsilon), \quad
   a^\varepsilon_{\aL}(\bu,\tbu)  =
   \lambda \int_{\Omega^\varepsilon} \ru^{t_i} \ru^*_{t_i} \, \rd \Omega^\varepsilon.
\end{equation}
The strong formulation of \eqref{3E2} can be written as
   $\gL(\varepsilon)\bu = \lambda\bu$,
where $\gL(\varepsilon)$ is the Lam\'e system
\begin{equation}
\label{eq:Lame}
   \gL = -\tfrac{E}{2(1+\nu)(1 - 2\nu)} \big( (1-2\nu)\Delta + \nabla\div \big)
\end{equation}
set on $\Omega^\varepsilon$ and associated with Dirichlet BC's on $\partial_0\Omega^\varepsilon$ and natural BC's on the rest of the boundary.
Its spectrum $\gS(\gL(\varepsilon)$ is discrete and positive. Let $\lamone{\gL(\varepsilon)}$ be its first eigenvalue. It is obtained by the minimum Rayleigh quotient
\[
   \lamone{\gL(\varepsilon)} = \min_{\bu\in V(\Omega^\varepsilon)}
   \frac{a^\varepsilon_{\aL}(\bu,\bu)}{\DNormc{\bu}{L^2(\Omega^\varepsilon)}}\,.
\]

\subsection{Surfacic shell models}
\label{ss:2D}
The key operators of the reduction to the midsurface $\cS$, namely the membrane and bending operators, are defined via intrinsic geometrical objects attached to $\cS$. To introduce them, we need generic parametrizations $F:(x_\alpha)_{\alpha\in1,2}\to\bt$ acting from maps neighborhoods $\cV$ into the midsurface $\cS$. Associated tangent coordinate vector fields are
\[
   \be_\alpha = \partial_{\alpha} \F,\quad \alpha=1,2,\quad\mbox{with}\quad
   \partial_\alpha=\frac{\partial}{\partial x_\alpha}.
\]
Completed by the unit normal field $\bn$ they form a basis $\{\be_1,\be_2,\bn\}$ in each point of $\cS$. The metric tensor $(a_{\alpha\beta})$ and the curvature tensor $(b_{\alpha\beta})$ are given by
\[
   a_{\alpha\beta} = \langle \be_\alpha,\be_\beta\rangle
   \quad\mbox{and}\quad
   b_{\alpha\beta}=\langle \partial_{\alpha\beta}\F,\bn\rangle\,.
\]
Denoting by $(a^{\alpha\beta})$ the inverse of $(a_{\alpha\beta})$, the curvature (symmetric) matrix is defined by
\[
   (b^\alpha_\beta) \quad\mbox{with}\quad b^\alpha_\beta = a^{\alpha\gamma}b_{\gamma\beta}.
\]
The eigenvalues $\kappa_1$ and $\kappa_2$ of the matrix $(b^\alpha_\beta)$ are called the {\em principal curvatures} of $\cS$ and their product is the {\em Gaussian curvature} $K$. Here comes the classification of shells: If $K\equiv0$, the shell is {\em parabolic}, if $K>0$, the shell is {\em elliptic}, if $K<0$, the shell is {\em hyperbolic}. Finally let $R$ denote the minimal radius of curvature of $\cS$
\begin{equation}
\label{eq:R}
   R = \inf_{\bp\in\cS} \Big\{ \min\{|\kappa_1(\bp)|^{-1},|\kappa_2(\bp)|^{-1}\} \Big\}.
\end{equation}
The basis $\{\be_\alpha,\bn\}$ determines contravariant components $(\zeta^\alpha,\zeta^3)$ of a vector field $\zetaz$ on $\cS$:
\[
   \zetaz = \zeta^{t_i}\be_{t_i} = \zeta^\alpha\be_\alpha + \zeta^3\bn\,.
\]
The covariant components are $(\zeta_\alpha,\zeta_3)$ with $\zeta_\alpha = a_{\alpha\beta}\zeta^\beta$ and $\zeta_3=\zeta^3$.
The surfacic rigidity tensor on $\cS$ is given by
$$
   M^{\alpha\beta\sigma\delta}= \tfrac{\nu E}{1-\nu^2}
   a^{\alpha\beta}a^{\sigma\delta}+\tfrac{E}{2(1+\nu)}(a^{\alpha\sigma}a^{\beta\delta}+a^{\alpha\delta}a^{\beta\sigma}).
$$
Note that, even if $\cS$ is flat ($a^{\alpha\beta}=\delta^{\alpha\beta}$), $M$ is different than the 3D rigidity tensor $A$.

\subsubsection{Membrane operator}
The variational space associated with the membrane operator is
\begin{equation}
\label{eq:VM}
   V_\aM(\cS) =  H^1_0(\cS) \times H^1_0(\cS) \times L^2(\cS).
\end{equation}
For an element $\zetaz=(\zeta_\alpha,\zeta_3)$ in $V_\aM(\cS)$, the change of metric tensor $\gammag = \gamma_{\alpha\beta}(\zetaz)$ is given by
\[
   \gamma_{\alpha\beta}(\zetaz)=\tfrac12(\D_\alpha\zeta_\beta+\D_\beta\zeta_\alpha)-b_{\alpha\beta}\zeta_3,
\]
where $\D_\alpha$ is the covariant derivative on $\cS$, see \cite{doCarmo1,doCarmo2,Spivak}. The membrane energy scalar product is defined as
\[
   a_\aM(\zetaz,\zetazt) =
   \int_\cS M^{\alpha \beta\sigma\delta} \gamma_{\alpha \beta}(\zetaz) \,\gamma_{\sigma\delta}(\zetazt) \,\rd \cS\,.
\]
Here the volume form $\rd \cS$ is $\sqrt{|\det(a_{\alpha\beta})|}\,\rd x_1\rd x_2$.
The variational formulation of the modal problem associated with the membrane operator $\gM$ is given by

{\sl Find $(\zetaz, \Lambda)$ with  $\zetaz \in V_\aM(\cS)\setminus\{0\}$ and $\Lambda \in \R$ such that for all $\zetazt \in V_\aM(\cS)$, }
\begin{equation}
\label{eq:membev}
   a_\aM(\zetaz,\zetazt) =
   \Lambda \int_\cS ( \zeta^\beta \zetat_\beta + \zeta^3 \zetat_3 ) \,\rd \cS.
\end{equation}

\subsubsection{Bending operator and Koiter model}
The variational space associated with the bending operator is
\begin{equation}
\label{eq:VB}
   V_\aB(\cS) =  H^1_0(\cS) \times H^1_0(\cS) \times H^2_0(\cS).
\end{equation}
For an element $\zetaz=(\zeta_\alpha,\zeta_3)$ in $V_\aB(\cS)$, the change of curvature tensor $\rhor = \rho_{\alpha\beta}(\zetaz)$ is given by
\[
   \rho_{\alpha\beta}(\zetaz) =
   \D_{\alpha}\D_{\beta}\zeta_3+\D_{\alpha}(b_{\beta}^\delta \zeta_{\delta})
   +b_{\alpha}^\delta \D_{\beta}\zeta_{\delta}-b_\alpha^\delta b_{\beta\delta}\zeta_3.
\]
The bending operator $\gB$ acts on the variational space $V_\aB(\cS) $ and its energy scalar product is
$$
   a_\aB(\zetaz,\zetazt) =
   \frac13\int_\cS M^{\alpha \beta\sigma\delta} \rho_{\alpha \beta}(\zetaz) \rho_{\sigma\delta}(\zetazt)
   \,\rd\cS \,.
$$
For any positive $\varepsilon$, the Koiter operator $\gK(\varepsilon)$ is defined as $\gM +\varepsilon^2\gB$.
It can  be shown, see \cite{BernadouCiarlet}, that $\gK(\varepsilon)$ is elliptic with multi-order on $ V_\aB(\cS)$ in the sense of Agmon-Douglis-Nirenberg \cite{ADN2}. The corresponding Koiter energy scalar product is
\begin{equation}
\label{eq:en2D}
   a^\varepsilon_{\aK}(\zetaz,\zetazt) =
   2\varepsilon\, a_\aM(\zetaz,\zetazt) + 2\varepsilon^3a_\aB(\zetaz,\zetazt).
\end{equation}

\subsection{Reconstruction operators from the midsurface to the shell}
The parametrizations of the midsurface induce local system of normal coordinates $(x_\alpha,x_3)$ inside the shell and, correspondingly, the covariant components $\ru_\alpha$ and $\ru_3$ of a displacement $\bu$.
The rationale of the shell theory is to deduce by an explicit procedure a solution $\bu$ of the 3D Lam\'e system posed on the shell from a solution $\zetaz$ of the Koiter model posed on the midsurface. This is done via a {\em reconstruction operator} $\rU$, cf \cite{Koiter2,Koiter3} and \cite{DaFa1}.
With any displacement $\zetaz(x_\alpha)$ defined on the midsurface $\cS$, $\rU$ associates a 3D displacement $\bu$ depending on the three coordinates $(x_\alpha,x_3)$ in $\Omega^\varepsilon$.
The operator $\rU$ is defined by
\begin{equation}
\label{eq:U}
   \rU = \rT \circ \rW
\end{equation}
where $\rW$ is the shifted reconstruction operator
\begin{equation}
\label{eq:W}
   \rW\zetaz = \left\{
   \begin{array}{l}
   \zeta_{\sigma} - x_3 (\D_{\sigma} \zeta_{3} + b_{\sigma}^{\alpha} \zeta_{\alpha}),\\[1ex]
   \zeta_3 - \frac{\nu}{1-\nu}\, x_3\,  \gamma_{\alpha}^{\alpha}(\zetaz)
   + \frac{\nu}{2-2\nu} \,x_3^{2} \,\rho_{\alpha}^{\alpha}(\zetaz),
   \end{array}
   \right.
\end{equation}
and $\rT:\zetaz\mapsto\rT\zetaz$ is the shifter defined as $(\rT\zetaz)_\sigma=\zeta_\sigma-x_3 b^\alpha_\sigma\zeta_\alpha$ and $(\rT\zetaz)_3=\zeta_3$, see \cite{Naghdi1}. The Koiter elastic energy of $\zetaz$ is a good approximation of the 3D elastic energy of $\rU\zetaz$, cf.\ \cite[Theorem A.1]{DaFa1}:
For any $\zetaz\in (H^2\times H^2\times H^3)\cap V_B(\cS)$, there holds, with non-dimensional constant $A$
\begin{equation}
\label{eq:ener}
   \big\vert  a^\varepsilon_{\aK}(\zetaz,\zetaz) -  a^\varepsilon_{\aL}(\rU\zetaz,\rU\zetaz)\big\vert
   \leq
   A \,a^\varepsilon_{\aK}(\zetaz,\zetaz)\, \Big(\frac{\varepsilon} {R} 
   + \frac{\varepsilon^{2}}{L^2} \Big)   ,
\end{equation}
where $R$ is the minimal radius of curvature \eqref{eq:R} of $\cS$, and $L$ is the {\em wave length} for $\zetaz$ defined as the largest constant such that the following ``inverse estimates'' hold
\begin{equation}
\label{eq:L}
    L\,\Norm{\gammag}{H^1(\cS)} \le \DNorm{\gammag}{L^2(\cS)}
    \quad\mbox{and}\quad
    L\,\Norm{\rhor}{H^1(\cS)} \le \DNorm{\rhor}{L^2(\cS)}\, .
\end{equation}
Note that for $\zetaz\in V_\aB(\cS)$, the first two components of $\rU\zetaz$ satisfy the Dirichlet condition on $\partial_0\Omega^\varepsilon$, whereas the third one does not need to satisfy it. In order to remedy that, we add a corrector term $\bu^\cor$ to $\rW\zetaz$ to compensate for the nonzero trace $g=- \frac{\nu}{1-\nu}\, x_3\,  \gamma_{\alpha}^{\alpha}(\zetaz) + \frac{\nu}{2-2\nu} \,x_3^{2} \,\rho_{\alpha}^{\alpha}(\zetaz)\big|_{\partial\cS}$. This corrector term is constructed and its energy estimated in \cite[sect.\,7]{DaFa1}. It has a simple tensor product form and exhibits the typical 3D boundary layer scale $d/\varepsilon$ with $d=\dist(\bp,\partial_0 \Omega^\varepsilon)$:
\[
   \bu^\cor = \Big( 0, 0, g\, \chi\Big( \frac{d}{\varepsilon}\Big)\Big)^\top
   \quad\mbox{with}\quad
   \chi\in C^\infty_0(\R),\ \chi(0)=1.
\]
``True'' boundary layer terms live at the same scale, decay exponentially, but have a non-tensor form in variables $(d,x_3)$, see \cite{DaugeGruais2,DDFR99} for plates and \cite{Faou2} for elliptic shells. Nevertheless this expression for $\bu^\cor$ suffices to obtain good estimates: There holds
\[
   a^\varepsilon_{\aL}(\bu^\cor,\bu^\cor) \le A \,a^\varepsilon_{\aK}(\zetaz,\zetaz) \,
   \, \Big(\frac{\varepsilon}{\ell} + \frac{\varepsilon^3}{\ell^3} \Big),
\]
for $\ell$ the lateral wave length of $\zetaz$ defined as the largest constant such that
\begin{equation}
\label{eq:ell}
    \ell\Normc{\gammag}{L^2(\partial \cS)} + \ell^3\Normc{\gammag}{H^1(\partial \cS)}
    \leq \DNormc{\gammag}{L^2(\cS)}\quad\mbox{and}\quad
    \ell\Normc{\rhor}{L^2(\partial \cS)} + \ell^3\Normc{\rhor}{H^1(\partial \cS)}
    \leq \DNormc{\rhor}{L^2(\cS)}\,.
\end{equation}

\begin{example}
\label{ex:bd}
Let $G$ in $H^2(\R_+)$ be such that $G\equiv0$ for $t\ge1$. Let $k\in\N$ and $0<\tau<\tau_0$ for $\tau_0$ small enough. The function $g(z,\varphi)$ defined on $\cS$ as
\[
   g(z,\varphi) = e^{ik\varphi}\,G\Big( \frac{d}{\tau}\Big)
\]
satisfies the estimates
$L\,\Norm{g}{H^1(\cS)} \le \DNorm{g}{L^2(\cS)}$ and
$\ell\Normc{g}{L^2(\partial \cS)} + \ell^3\Normc{g}{H^1(\partial \cS)} \leq \DNormc{g}{L^2(\cS)}$
for $L$ and $\ell$ larger than $c(G)\min\{\tau,k^{-1}\}$ where the positive constant $c(G)$ is independent of $\tau$ and $k$.
\end{example}

In the present work, we are interested in comparing surfacic and 3D Rayleigh quotients so we introduce the following notations
\[
   Q^\varepsilon_{\aL}(\bu) = \frac{a^\varepsilon_{\aL}(\bu,\bu)}{\DNormc{\bu}{L^2(\Omega^\varepsilon)}}\,,
   \quad\mbox{and}\quad
   Q^\varepsilon_{\aK}(\zetaz) =
   \frac{a^\varepsilon_{\aK}(\zetaz,\zetaz)}{2\varepsilon\DNormc{\zetaz}{L^2(\cS)}}\,.
\]
By similar inequalities as in \cite{DaFa1} we can prove the following relative estimate

\begin{theorem}
\label{th:rayl}
{\em (i)\ } For all $\zetaz\in (H^2\times H^2\times H^3)\cap V_B(\cS)$ and with $\rU$ defined in \eqref{eq:U} we set
\[
   \ocirc{\rU}\zetaz = \rU\zetaz - \bu^\cor.
\]
Then $\ocirc{\rU}\zetaz$ belongs to the 3D variational space $V(\Omega^\varepsilon)$.
With $L$ and $\ell$ the wave lengths  \eqref{eq:L} and \eqref{eq:ell}, let us assume $\varepsilon\le L$ and $\varepsilon\le\ell$. We also assume $Q^\varepsilon_{\aK}(\zetaz)\le EM$ for a chosen constant $M\ge1$ independent of $\varepsilon$. Then we have the relative estimates between Rayleigh quotients for $\varepsilon$ small enough
\begin{equation}
\label{eq:rayl}
   \big\vert  Q^\varepsilon_{\aK}(\zetaz) -
   Q^\varepsilon_{\aL}(\ocirc{\rU}\zetaz)\big\vert
   \leq
   A' \,Q^\varepsilon_{\aK}(\zetaz)\,
   \Big( \frac{\varepsilon}{R} + \frac{\varepsilon^{2}}{L^2}
   + \Big(\frac{\varepsilon}{\ell}\Big)^{1/2} + \varepsilon \sqrt{M} \,\Big)   \,,
\end{equation}
with a constant $A'$ independent of $\varepsilon$ and $\zetaz$.

{\em (ii)\ } If $\zetaz$ belongs to $(H^2_0\times H^2_0\times H^3_0)(\cS)$, the boundary corrector $\bu^\cor$ is zero and the above estimates do not involve the term $\sqrt{\varepsilon/\ell}$ any more.
\end{theorem}

This theorem allows to find upper bounds for the first 3D eigenvalue $\lambda^{\varepsilon}_1$ if we know convenient energy minimizers $\zetaz^\varepsilon$ for the Koiter model $\gK(\varepsilon)$ and if we have the relevant information about their wave lengths.

\section{Axisymmetric shells}
\label{s:2}
An axisymmetric shell is invariant by rotation around an axis that we may choose as $t_3$. Recall that $(r,\varphi,\tau)\in\R^+\times\T\times\R$ denote associated cylindrical coordinates satisfying relations \eqref{eq:cyl1} and coordinate vectors are $\be_r$, $\be_\varphi$, and $\be_\tau$ given by \eqref{eq:cyl2}.
Accordingly, the (contravariant) {\em cylindrical components} of a displacement $\bu = \ru^{t_i}\be_{t_i}$ are $(\ru^r,\ru^\varphi,\ru^\tau)$ so that $\bu =\ru^{r}\be_r+\ru^{\varphi}\be_\varphi + \ru^{\tau}\be_\tau$. In particular
the {\em radial component} of $\bu$ is given by
\begin{equation}
\label{eq:rad}
   \ru^r = \ru^{t_1}\cos\varphi + \ru^{t_2}\sin\varphi.
\end{equation}
The components $\ru^\varphi$ and $\ru^\tau$ are called azimuthal and axial, respectively.

An axisymmetric domain $\Omega\subset\R^3$ is associated with a meridian domain $\omega\subset\R^+\times\R$ so that
\begin{equation}
\label{1E2}
   \Omega = \{\bx\in\R^3,\quad (r,\tau)\in\omega\ \ \mbox{and}\ \ \varphi\in\T\}.
\end{equation}

\subsection{Axisymmetric parametrization}
For a shell $\Omega^\varepsilon$ that is axisymmetric, let $\omega^\varepsilon$ be its meridian domain. The midsurface $\cS$ of $\Omega^\varepsilon$ is axisymmetric too. Let $\cC$ be its meridian domain. We have a relation similar to \eqref{1E1}
\begin{equation}
\label{1E3}
   \Phi : \quad\cC\times(-\varepsilon,\varepsilon)\ni\big((r,\tau),x_3\big) \;\;\longmapsto\;\; (r,\tau)+x_3\bn(r,\tau)\in\omega^\varepsilon.
\end{equation}
The meridian midsurface $\cC$ is a curve in the halfplane $\R^+\times\R$.

\begin{assumption}
\label{as:1}
Let $\cI$ denote any bounded interval and let $z$ be the variable in $\cI$.
\\
(i) The curve $\cC$ can be parametrized by one map defined on $\cI$ by a smooth function $f$:
\begin{equation}
\label{1E4}
\begin{array}{ccc}
   \cI &\longrightarrow &\cC \\
   z &\longmapsto&  (r,\tau) = (f(z),z)
\end{array}
\quad\mbox{with }\quad
f:z\mapsto r=f(z).
\end{equation}
(ii) The shells are disjoint from the rotation axis, i.e., there exists $R_{\min}>0$ such that $f\ge R_{\min}$.
\end{assumption}

\begin{remark}
We impose condition (ii) to avoid technical difficulties due to the singularity at the origin.
We have observed that, if we keep this condition, the inverse parametrization $z = g(r)$ does not bring new examples in the framework that we investigate in this paper. For instance {\em annular plates} pertain to this inverse parametrization, but they fall in \cite{DDFR99} that provides a complete eigenvalue asymptotics.
\end{remark}

The parametrization \eqref{1E4} of the meridian curve $\cC$ provides a parametrization of the meridian domain $\omega^\varepsilon$ by $\cI\times(-\varepsilon,\varepsilon)$: Let us introduce the arc-length
\begin{equation}
\label{eq:sz}
   s(z) = \sqrt{1+f'(z)^2}, \quad z \in \cI\,.
\end{equation}
The unit normal vector $\bn$ to $\cC$ at the point $(r,\tau)=(f(z),z)$ is given by $(\frac{1}{s(z)},-\frac{f'(z)}{s(z)})$ and the parametrization by
\[
  \cI\times(-\varepsilon,\varepsilon)\ni (z,x_3) \;\; \longmapsto \;\;
  \Big(f(z) + x_3\ \tfrac{1}{s(z)},z - x_3 \ \tfrac{f'(z)}{s(z)}\Big) \in\omega^\varepsilon\,.
\]
The parametrization \eqref{1E4} also induces the parametrization $F$ \eqref{eq:surfpara} of
the midsurface $\cS$ by the variables $(z,\varphi)\in \cI\times\T$. The unit normal vector $\bn$ to $\cS$ at the point $F(z,\varphi)$ is given by
\[
   \bn = s(z)^{-1} (\be_r - f'(z)\be_\tau)
\]
while tangent coordinate vectors are $\be_z=\partial_zF$ and $\be_\varphi=\partial_\varphi F$, i.e.
\[
   \be_z = f'(z)\be_r + \be_\tau
\]
while $\be_\varphi$ coincides the coordinate vector of same name corresponding to cylindrical coordinates \eqref{eq:cyl2}.
The metric tensor $(a_{\alpha\beta})$ is given by $\langle \be_\alpha,\be_\beta\rangle$ with $\alpha,\beta\in\{z,\varphi\}$, i.e.
\begin{equation}
\label{metric}
\begin{pmatrix}
a_{zz} & a_{z\varphi} \\ a_{\varphi z} & a_{\varphi\varphi}
\end{pmatrix}(z)
 = \begin{pmatrix}
s(z)^2 & 0 \\ 0 &  f(z)^2
\end{pmatrix}.
\end{equation}
The curvature tensor and Gaussian curvature $K$ are respectively given by
\begin{equation}
\label{2E2}
\begin{pmatrix}
b^z_{z} & b^z_{\varphi} \\ b^\varphi_{z} & b^\varphi_{\varphi}
\end{pmatrix}(z)
 = \begin{pmatrix}
f''(z)s(z)^{-3} & 0 \\ 0 & - f(z)^{-1}s(z)^{-1}
\end{pmatrix}
\quad\mbox{and}\quad
K(z) = - \frac{ f''(z)}{f(z) s(z)^4}
\,.
\end{equation}
So the curvature tensor is in diagonal form, and $K$ is simply the product of its diagonal elements.

\begin{definition}
\label{def:curv}
We call $b^z_z$ the {\em meridian curvature} and $b^\varphi_\varphi$ the {\em azimuthal curvature}.
\end{definition}

Since we have assumed that $f\ge R_0>0$,  all terms are bounded  and we find that
\begin{enumerate}
\item[(1)]  If $f''\equiv0$, i.e.\ $f$ is affine, the shell is (nondegenerate) parabolic. If $f$ is constant, the shell is a cylinder, if not it is a truncated cone (without conical point!).
\item[(2)]  If $f''<0$, the shell is elliptic.
\item[(3)]  If $f''>0$, the shell is hyperbolic.
\end{enumerate}

\subsection{Surfacic axisymmetric models in normal coordinates}
\label{ss:4}
Relations \eqref{eq:surfpara} and \eqref{1E1} define {\em normal coordinates} $(z,\varphi,x_3)$ in the thin shell $\Omega^\varepsilon$. For example when the midsurface $\cS$ is a cylinder ($f$ constant), the normal coordinates are a permutation of standard coordinates: $(z,\varphi,x_3)=(\tau,\varphi,r)$. The associate (contravariant) decomposition of surface displacement fields $\zetaz$ is written as $\zetaz = \zeta^z\be_z + \zeta^\varphi\be_\varphi + \zeta^3\bn$, where $\zeta^3$ is the component of the displacement in the normal direction $\bn$ to the midsurface, $\zeta^z$  and  $\zeta^\varphi$ the meridian and azimuthal components respectively, defined so that there holds
\[
   \zeta^{t_1}\be_{t_1} + \zeta^{t_2}\be_{t_2} + \zeta^{t_3}\be_{t_3}
   =  \zeta^z\be_z + \zeta^\varphi\be_\varphi + \zeta^3\bn\,.
\]
Note that the azimuthal component is the same as defined by cylindrical coordinates. The covariant components are
\[
   \zeta_z = s^2\zeta^z,\quad
   \zeta_\varphi = f^2\zeta^\varphi,
   \quad\mbox{and}\quad\zeta_3=\zeta^3.
\]
The change of metric tensor $\gamma_{\alpha\beta}(\zetaz)$ has the expression in normal coordinates
\begin{equation}
\label{eqgamma}
\begin{array}{rcl}
   \gamma_{zz}(\zetaz)&=&
   \displaystyle\partial_{z}\zeta_{z}-\frac{f'f''}{s^2}\zeta_{z}-\frac{f''}{s}\zeta_3\\[1ex]
   \gamma_{z\varphi}(\zetaz)&=&
   \displaystyle\frac{1}{2}(\partial_{z} \zeta_{\varphi}
   +\partial_\varphi \zeta_{z})-\frac{f'}{f}\zeta_{\varphi}\\[1ex]
   \gamma_{\varphi\varphi}(\zetaz)&=&
   \displaystyle \partial_\varphi\zeta_{\varphi}+\frac{f f'}{s^2}\zeta_{z}+\frac{f}{s}\zeta_3\,,
\end{array}
\end{equation}
while the change of curvature tensor $\rho_{\alpha\beta}(\zetaz)$ is written as
\begin{equation}
\label{rho}
\begin{array}{rcl}
   \rho_{zz}(\zetaz)&=&
   \displaystyle\partial_{z}^2 \zeta_3-\frac{f''^2}{s^4}\zeta_3
   +\frac{2f''}{s^3}\partial_{z}\zeta_{z}+\frac{f'''s^2-5f'f''^2}{s^5}\zeta_{z}\\[1ex]
   \rho_{\varphi\varphi}(\zetaz)&=&
   \displaystyle\partial_{\varphi}^2\zeta_3-\frac{1}{s^2}\zeta_3
   -\frac{2}{fs}\partial_{\varphi}\zeta_{\varphi} -\frac{ 2f'}{s^3}\zeta_{z}\\[1ex]
   \rho_{z\varphi}(\zetaz)&=&
   \displaystyle\partial_{z\varphi}\zeta_3 + \frac{f''}{s^3}\partial_{\varphi}\zeta_{z}
   -\frac{1}{f s}\partial_{z}\zeta_{\varphi}+\frac{2f'}{f^2s}\zeta_{\varphi}\,.
\end{array}
\end{equation}

\section{Principles of construction: High frequency analysis}
\label{s:5}
The construction is based on the following postulate:

\begin{postulate}
\label{4P1}
The eigenmodes associated with the smallest vibrations are strongly oscillating in the angular variable $\varphi$ and this oscillation is dominating.
\end{postulate}

This means that if this postulate happens to be true for certain families of shells, our construction will provide rigorous quasimodes and, moreover, these quasimodes are candidates to be associated with lowest energy eigenpairs.
We may notice that Postulate \ref{4P1} is wrong for planar shells. But it appears to be true for nondegenerate parabolic shells and some subclasses of elliptic shells.

\subsection{Angular Fourier decomposition}
We can perform a discrete Fourier decomposition in the shell $\Omega^\varepsilon \equiv \omega^\varepsilon\times\T$ and in its midsurface $\cS\equiv\cC\times\T\cong\cI\times\T$.
For a displacement $\bu$ defined on $\Omega^\varepsilon$, and its Fourier coefficient of order $k\in\Z$ is denoted by $\bu^k$ and defined on $\omega^\varepsilon$, see \eqref{eq:Fou2}. Likewise, a surface displacement $\zetaz$ defined on $\cS$, and its Fourier coefficient of order $k$ is denoted by $\zetaz^k$ and defined on the curve $\cC$.
This Fourier decomposition diagonalizes the Lam\'e system $\gL$  with respect to the angular modes $\re^{ik\varphi}$, $k\in\Z$, due to the relation:
\[
    (\gL\bu)^k = \gL^k \bu^k.
\]
Similar properties hold with the membrane and bending operators $\gM$ and $\gB$ defined on the spaces $V_\aM(\cS)$ and $V_\aB(\cS)$, composing the Koiter operator $\gK(\varepsilon)$.
Recall from sect.\ref{ss:op} that $\gL^k(\varepsilon)$, $\gM^k$, $\gB^k$ and $\gK^k(\varepsilon)$, are the angular Fourier decomposition of $\gL(\varepsilon)$, $\gM$, $\gB$ and $\gK(\varepsilon)$, respectively.

The (non decreasing) collections of the eigenvalues of $\gL^k(\varepsilon)$ for all $k\in\Z$ gives back all eigenvalues of $\gL(\varepsilon)$. Note that since $\gL$ is real valued, the eigenvalues for $k$ and $-k$ are identical. Thus $\lamone{\gL(\varepsilon)} = \inf_{k\in\N} \lamone{\gL^k(\varepsilon)}$ and we denote by $k(\varepsilon)$ the smallest natural integer $k$ such that 
\[
   \lamone{\gL(\varepsilon)} = \lamone{\gL^{k(\varepsilon)}(\varepsilon)}\,.
\]
Postulate \ref{4P1} means that $k(\varepsilon)\to\infty$ as $\varepsilon\to0$.

\subsection{High frequency analysis of the membrane operator}
The eigenmode membrane equation \eqref{eq:membev} at azimuthal frequency $k$ takes the form
\begin{equation}
\label{eq:membevk}
   \gM^k \zetaz^k = \Lambda^k\bA\zetaz^k
\end{equation}
where $\bA$ is the mass matrix
\begin{equation}
\label{eq:A}
\bA = \begin{pmatrix} a^{zz} & 0 & 0 \\ 0 & a^{\varphi\varphi}& 0 \\ 0 & 0 & 1 \end{pmatrix} = \begin{pmatrix} s^{-2} & 0 & 0 \\ 0 & f^{-2} & 0 \\ 0 & 0 & 1 \end{pmatrix}.
\end{equation}
We construct quasimodes for $\gM^k$ as $k\to\infty$, i.e. pairs $(\breve\Lambda^k,\breve\zetaz{}^k)$ with $\breve\zetaz{}^k$ in the domain of the operator $\gM^k$ and satisfying the estimates
\[
   \DNorm{(\gM^k-\breve\Lambda^k)\,\breve\zetaz{}^k}{L^2(\cS)} \le
   \delta(k) \DNorm{\breve\zetaz{}^k}{L^2(\cS)}
   \quad\mbox{with}\quad
   \delta(k)/\breve\Lambda^k \to 0 \mbox{ as } k\to\infty.
\]

Now we consider the membrane operator as a formal series with respect to $k$
\begin{equation}
\label{bach}
   \gM^k = k^2\bM_0 + k\bM_1 + \bM_2 \equiv \bM\kk, \quad \mbox{with}\quad
   \bM \kk  = k^{2}\sum_{n\in\N} k^{-n}\bM_n \,,
\end{equation}
and try to solve \eqref{eq:membevk} in the formal series algebra:
\begin{equation}
\label{eq:membeqfs}
   \bM\kk\zetaz\kk = \Lambda\kk\bA\zetaz\kk.
\end{equation}
Here the multiplication of formal series is the Cauchy product: For two formal series $a\kk = \sum_n k^{-n} a_n$ and $b\kk = \sum_{n} k^{-n} b_n$, the coefficients of the series $a\kk \,b\kk = \sum_{n} k^{-n} c_n$ are given by $c_n = \sum_{\ell +m=n} a_\ell b_{m}$.

The director $\bM_0$ of the series $\bM\kk$ is given in parametrization $r=f(z)$ by
\begin{equation}
\label{eq:M0}
\bM_0=\frac{E}{1-\nu^2}\begin{pmatrix}
\frac{1-\nu}{2f^2s^2} & 0 & 0\\
0 &\frac{1}{f^4}& 0\\
0 & 0 & 0
\end{pmatrix}.
\end{equation}
Its kernel is given by all triples $\zetaz$ of the form $(0,0,\zeta_3)^\top$. This is the reason why we look for a reduction of the eigenvalue problem for $\bM$ to a scalar eigenvalue problem set on the normal component $\zeta_3$. The key is a factorization process in the formal series algebra proved in \cite[Chap.3]{beaudouin:tel},
\begin{equation}
\label{eq:reduc}
\bM\kk \bV\kk -  \Lambda\kk \bA\bV\kk= \bV_{0} \circ (\rH\kk - \Lambda\kk)\,.
\end{equation}
Here $\bV\kk$ is a (formal series of) reconstruction operators whose first term $\bV_0$ is the embedding $\bV_0 \eta = (0,0,\eta)^\top$  in the kernel of $\bM_0$, and $\rH\kk$ is the scalar reduction.

\begin{theorem}
\label{reducform}
Let be a formal series with real coefficients :
$$
\Lambda[k] = \sum_{n \geq 0} k ^{-n}\Lambda_n.
$$
For $n\geq 1$, there exist operators
$\rV_{n,z},\rV_{n,\varphi}: C^\infty(\overline\cI) \to C^\infty(\overline\cI)$ of order $n-1$, polynomial in $\Lambda_j$, for $j\leq n-3$, and for $n\geq 0$ scalar operators  $\rH_{n}: C^\infty(\overline\cI) \to C^\infty(\overline\cI)$   of order $n$, polynomial in $\Lambda_j$, for $j\leq n-2$  such that if we set :
$$
\bV[k]=\sum_{n\geq 0}k^{-n} \bV_{n}\quad \mbox{with}\quad  \bV_{n}=(\rV_{n,z},\rV_{n,\varphi},0)^\top,
\quad\mbox{and}\quad
\rH[k]=\sum_{n\geq 0} k^{-n}\rH_{n}
$$
we have \eqref{eq:reduc} in the sense of formal series.
\end{theorem}
See Appendix \ref{app:A} for more details on this theorem.

With the scalar reduction $\rH\kk$ is associated the formal series problem
\begin{equation}
\label{scalreducfs}
\rH\kk \,\eta\kk = \Lambda\kk \,\eta\kk
\end{equation}
where $\eta\kk = \sum_{n \geq 0} k^{-n} \eta_n$ is a scalar formal series. The previous theorem shows that any solution to \eqref{scalreducfs} provides a solution $\zetaz\kk = \bV\kk \eta\kk$ to \eqref{eq:membeqfs}.

The cornerstone of our quasimodes construction for $\gM^k$ as $k\to\infty$ is to construct a solution $\eta\kk$ of the problem \eqref{scalreducfs}. This relies on the possibility to extract an elliptic operator $\hM^k$ with compact resolvent from the first terms of the series $\rH\kk$ as we describe in several geometrical situations later on.

\begin{remark}
\label{rem:ess}
The  essential spectrum $\sigma_{\sf ess}(\gM^k)$ of the membrane operator $\gM^k$ at frequency $k$ can be determined explicitly thanks to \cite[Th.4.5]{AtLangMenShka}. It depends only on its principal part, which coincides with the (multi-degree) principal part of $\bM_2$, and is given by the range of $\frac{E}{f(z)^2s(z)^2}$ for $z\in\cI$, see \cite[sect.\,2.7]{beaudouin:tel} for details.
With formula \eqref{2E2}, we note the relation with the azimuthal curvature
\begin{equation}
\label{eq:essMk}
   \sigma_{\sf ess}(\gM^k)=\big\{ E\,b^\varphi_\varphi(z)^2\;,\quad z\in\cI\big\}.
\end{equation}
As a consequence of Assumption \ref{as:1}, the minimum of $\sigma_{\sf ess}(\gM^k)$ is positive.
\end{remark}

\subsection{High frequency analysis of the Koiter operator}
Similar to the membrane operator $\bM\kk$, the bending operator expands as
$$
\gB^k = k^4 \bB_0 + \sum_{n = 1}^4 k^{4 - n}\bB_n \equiv\bB\kk,
$$
with first term
\begin{equation}
\label{eq:B0}
\bB_0=\begin{pmatrix}
0 & 0 & 0\\
0 & 0 & 0\\
0 & 0 & \rB_0
\end{pmatrix}\quad\mbox{with}\quad \rB_0=\frac{E}{1-\nu^2}\,\frac{1}{3 f^4}\,.
\end{equation}
We notice that we have the commutation relation
\[
   \bB_0 \bV\kk = \bV_0\rB_0\,.
\]
Therefore the identity \eqref{eq:reduc} implies for all $\varepsilon$ the identity
\begin{equation}
\label{eq:reduB}
   \big(\bM\kk +\varepsilon^2k^4\bB_0\big)\bV\kk -  \Lambda\kk \bA\bV\kk =
   \bV_{0} \circ \big(\rH\kk +\varepsilon^2k^4\rB_0 - \Lambda\kk\big)\,.
\end{equation}
Thus the same factorization as for the membrane operator will generate the quasimode constructions for the Koiter operator as soon as the higher order terms of $\gB^k$ correspond to perturbation terms. This is related to Postulate \ref{4P1}. The identity \eqref{eq:reduB} motivates the formula \eqref{eq:hKk} defining the reduced Koiter operator $\hK^k(\varepsilon) = \hM^k + \varepsilon^2k^4\rB_0$.
In the following two sections we provide $\hK^k(\varepsilon)$ and its lowest eigenvalues in several well defined cases.

\section{Nondegenerate parabolic case.}
\label{s:6}
We assume in addition to Assumption \ref{as:1}
\begin{equation}
\label{5E1}
   f(z) =  T z + R_0,\quad z\in\cI,\quad\mbox{with}\quad
   R_0>0,\ \ T\in\R\,.
\end{equation}
If $T=0$, the corresponding surface $\cS$ is a cylinder of radius $R_0$ and the minimal radius of curvature $R$ \eqref{eq:R} equals to $R_0$. So we write $f=R$ in the cylinder case. If $T\neq0$, the surface $\cS$ is a truncated cone. The arc length \eqref{eq:sz} is $s=\sqrt{1+T^2}$. In this section, we address successively the membrane scalar reduction, the Koiter scalar reduction, and finally the reconstruction of quasimodes into the shell $\Omega^\varepsilon$, providing an upper bound for $\lamone{\gL(\varepsilon)}$.

\subsection{Membrane scalar reduction in the parabolic case}
The first terms $\rH_n$ of the scalar formal series reduction of the membrane operator have been explicitly calculated in \cite{beaudouin:tel} in the cylindrical case $T=0$ and have the following expression in the general parabolic case:
\begin{align}
\label{5E2}
\rH_0 = \rH_1 = \rH_2 = \rH_3 = 0\quad\mbox{and}\quad
\rH_4(z,\partial_{z})=E\Big(\frac{f^2}{s^6} \partial_{z}^4+\frac{6f' f}{s^6} \partial_{z}^3
+\frac{6f'^2}{s^6} \partial_{z}^2 \Big) \,.
\end{align}
It is relevant to notice that $\rH_4$ is selfadjoint on $H^2_0(\cI)$ with respect to the natural measure $\rd\cI = f(z)s(z)\,\rd z\,$, since there holds
\begin{equation}
\label{5E3}
   \big\langle \rH_4\eta,\etat \big\rangle_{\cI} =
   \frac{E}{(1+T^2)^3} \int_{\cI} f(z)^2
   \,\partial^2_z\eta\,\partial^2_z\etat\,\rd\cI\,.
\end{equation}
This also proves that $\rH_4$ is positive.
The Dirichlet boundary conditions $\eta=\partial_z\eta=0$ on $\partial\cI$ are the right conditions to implement the  membrane boundary condition $\zeta_\alpha=0$ on $\partial\cI$  through the reconstruction operators $\bV_n$, see \eqref{V1}\,--\,\eqref{V2}. The eigenvalue formal series $\Lambda\kk$ starts with $\Lambda_4$ that is the first eigenvalue of $\rH_4$:
\begin{equation}
\label{eq:Lamh}
   \Lambda_0=\Lambda_1=\Lambda_2=\Lambda_3=0\quad\mbox{and}\quad \Lambda_4>0.
\end{equation}
The pair $(\breve\Lambda{}^k,\breve\zetaz{}^k)$
\begin{equation}
\label{eq:qmm}
   \breve\zetaz{}^k = \sum_{0\le n+m\le6} k^{-n-m}\bV_n\eta_m
   \quad \mbox{and}\quad
   \breve\Lambda{}^k = k^{-4}\Lambda_4
\end{equation}
with $(\Lambda_4,\eta_0)$ an eigenpair of $\rH_4$, and $\eta_m$ ($m=1,\ldots,6$)  constructed by induction  so that the membrane boundary conditions $\breve\zeta{}^k_\alpha=0$ are satisfied, is a quasimode for $\gM^k$. For instance, in the cylindrical case $f=R$, the triple $\breve\zetaz{}^k$ takes the form
\begin{equation}\label{uk}\footnotesize
   \breve\zetaz{}^k = \begin{pmatrix}
   0\\ 0 \\ \eta_{0}
   \end{pmatrix}
   +\frac{i}{k}\begin{pmatrix}
   0\\ R\eta_{0} \\ 0
   \end{pmatrix}
   +\frac{1}{k^2}\begin{pmatrix}
   -R \eta'_{0} \\ 0 \\ \eta_{2}
   \end{pmatrix}
   +\frac{i}{k^3}\begin{pmatrix}
   0\\ -\nu R^3\eta''_{0} + R\eta_{2}\\ \eta_{3}
   \end{pmatrix}
   -\frac{1}{k^4}\begin{pmatrix}
   (\nu\!+\!2)R^3\eta'''_0+R\eta'_2\\ R\eta_{3}\\ \eta_{4}
   \end{pmatrix}+\ldots\normalsize
\end{equation}
and the boundary conditions are, for  $z\in \partial \cI$
\begin{equation}
\label{eq:bcz}
   \eta_0(z)=0,\;\; \eta'_0(z)=0,\;\;
   \eta_2(z)=\nu R^2\eta''_0(z),\;\; \eta'_2(z)=(\nu+2)R^2\eta'''_0(z),\;\;
   \eta_3(z)=0,\ldots
\end{equation}

Recall that the minimum of the essential spectrum of $\gM^k$ is positive by Remark \ref{rem:ess}.
For $|k|$ large enough, $\gM^k$ has therefore at least an eigenvalue $\cong \Lambda_4k^{-4}$ under its essential spectrum and
\begin{equation}
\label{eq:memb5.1}
   \dist\big(k^{-4}\Lambda_4,\gS(\gM^k)\big) \lesssim k^{-5},\quad k\to\infty.
\end{equation}

\subsection{Koiter scalar reduction in the parabolic case}
The leading term of the series $\rH(k)$ is $ \hM^k =k^{-4}\rH_4$, as mentioned in the introduction, see \eqref{eq:hMk}. So, the leading term of the scalar reduction of the Koiter operator is, cf.\ \eqref{eq:reduB}
\begin{equation}
\label{5E4}
   \hK^k(\varepsilon) =
   k^{-4}\rH_4 + \varepsilon^2 k^4 \rB_0 =
   k^{-4}\rH_4 + \frac{\varepsilon^2}{3} \,\frac{E}{1-\nu^2} \, \frac{k^4}{f^4} \,.
\end{equation}
The operator $\hK^k(\varepsilon)$ is a priori defined for integers $k$, nevertheless it makes sense for any real number $k$, like all the other operators $\gM^k$, $\gB^k$ and $\gK^k(\varepsilon)$. We keep this extended framework all along this subsection. All functions and vector fields are defined on the parametric interval $\cI$ with variable $z$.

\subsubsection{Optimizing $k$}
The operator $\hK^k(\varepsilon)$ is self-adjoint on $H^2_0(\cI)$ real-valued and positive. Let $\lamone{\hK^k(\varepsilon)}$ denote its smallest eigenvalue. For any chosen $\varepsilon$ we look for $k_{\min}=k(\varepsilon)$ realizing the minimum $\mu_1^{\hK}(\varepsilon)$ of $\lamone{\hK^k(\varepsilon)}$ if it exists:
\[
   \mu_1^{\hK}(\varepsilon) = \lamone{\hK^{k(\varepsilon)}(\varepsilon)} = \min_{k\in\R_+}\lamone{\hK^k(\varepsilon)}.
\]
To ``homogenize'' the terms $k^{-4}$ and $\varepsilon^2k^4$ let us define $\gamma(\varepsilon)$ by setting
\begin{equation}
\label{5E5}
   \gamma(\varepsilon) = \varepsilon^{1/4} k(\varepsilon),
\end{equation}
so that we look equivalently for $\gamma(\varepsilon)$. There holds
\[
   \hK^{k(\varepsilon)}(\varepsilon) = k(\varepsilon)^{-4}\rH_4 + \varepsilon^2 k(\varepsilon)^4 \rB_0 =
   \varepsilon \Big(\frac{1}{\gamma(\varepsilon)^4} \rH_4 + \gamma(\varepsilon)^4 \rB_0\Big).
\]
Therefore $\gamma(\varepsilon)$ does not depend on $\varepsilon$. Let $\mu_1(\gamma)$ be the first eigenvalue of the operator
\begin{equation}
\label{5E6}
   \frac{1}{\gamma^4} \rH_4 + \gamma^4 \rB_0\,.
\end{equation}
The function $\gamma\mapsto\mu_1(\gamma)$ is continuous and, since $\rH_4$ and $\rB_0$ are {\em positive}, it tends to infinity as $\gamma$ tends to $0$ or to $+\infty$. Therefore we can define $\gamma_{\min}$  as the (smallest) positive constant such that $\mu_1(\gamma)$ is minimum
\begin{equation}
\label{5E6b}
   \mu_1(\gamma_{\min}) = \min_{\gamma\in\R_+}\mu_1(\gamma) =: \am_1 \,.
\end{equation}
Thus $k(\varepsilon)$ satisfies a power law that yields a formula for the minimal first eigenvalue $\mu_1^{\hK}(\varepsilon)$:
\begin{equation}
\label{5E7}
   k(\varepsilon) = \varepsilon^{-1/4} \gamma_{\min}\quad\mbox{and}\quad
   \mu_1^{\hK}(\varepsilon) = \am_1\varepsilon.
\end{equation}
Let $\eta_0$ be a corresponding eigenvector.  By definition
\begin{equation}
\label{eq:eta0}
   \eta_0\in H^2_0(\cI) \ \ \mbox{first eigenvector of }\ \ 
   \frac{1}{\gamma^4_{\min}} \rH_4 + \gamma^4_{\min} \rB_0 = \varepsilon^{-1}\hK^{k(\varepsilon)}(\varepsilon)\,.
\end{equation}
Note that $\mu_1(\gamma_{\min})$ coincides with the minimum of the Rayleigh quotients associated with $\eta_0$:
\begin{equation}
\label{eq:Rqeta0}
   \mu_1(\gamma_{\min}) = \min_{\gamma\in\R_+} 
   \frac{\langle \gamma^{-4} \rH_4\eta_0 + \gamma^4 \rB_0\eta_0,\eta_0\rangle}
   {\langle\eta_0,\eta_0\rangle}
\end{equation}
Therefore $\gamma_{\min}$ equilibrates the two terms in the numerator, which proves that
the ratio $\rR$ \eqref{eq:Rratio} between bending energy and total energy is equal to $\frac12$:
\begin{equation}
\label{eq:Rratiopara}
   \rR = \frac{\langle  \gamma^4_{\min} \rB_0\eta_0,\eta_0\rangle}
   {\langle \gamma^{-4}_{\min} \rH_4\eta_0 + \gamma^4_{\min} \rB_0\eta_0,\eta_0\rangle} = \frac12\,.
\end{equation}

\subsubsection{Case of cylinders}
In the cylindrical case $T = 0$, formulas are more explicit because $f$ is constant. So everything can be written as a function of the first Dirichlet eigenvalue $\mu_1^{\sf bilap}$ of the bilaplacian operator $\Delta^2$ on $H^2_0(\cI)$ as we explain now. We have
\[
   \rH_4 = ER^2\,\Delta^2 \quad\mbox{and}\quad
   \rB_0 = \frac{E}{1-\nu^2} \, \frac{1}{3R^4} \, .
\]
So the eigenvalue of $\frac{1}{\gamma^4} \rH_4 + \gamma^4 \rB_0$ is
\begin{equation}
\label{eq:twoterm}
   \mu_1(\gamma)=\frac1{\gamma^4} ER^2 \mu_1^{\sf bilap} + \gamma^4\, \frac{E}{1-\nu^2} \,\frac{1}{3R^4}.
\end{equation}
It is minimum for $\gamma_{\min}$ such that
\begin{equation}
\label{eq:gammin}
   \gamma^4_{\min} = R^3 \sqrt{3(1-\nu^2)\, \mu_1^{\sf bilap}}
\end{equation}
and we find that the minimum eigenvalue \eqref{5E6b} is
\begin{equation}
\label{5E9}
   \mu_1(\gamma_{\min}) =  \frac{2E}{R} \sqrt{\frac{ \mu_1^{\sf bilap}}{3(1-\nu^2)}} =: \am_1\,.
\end{equation}
Thus
\begin{equation}
\label{5E10}
   k(\varepsilon) = \varepsilon^{-1/4} R^{3/4} 
\big(3(1-\nu^2)\,
   \mu_1^{\sf bilap}
\big)^{1/8} .
\end{equation}

\begin{remark}
Denote by $\mu^{\sf bilap}$ the first eigenvalue of $\Delta^2$ on the unit interval $(0,1)$. We have the relation $\mu_1^{\sf bilap} =\mu^{\sf bilap}\, L^{-4}$ with the length $L$ of the interval $\cI$. 
\end{remark}

\subsubsection{Reconstruction of vectors from scalars. Membrane boundary conditions} In order to reconstruct fields $\zetaz^k$ from the scalar eigenvector $\eta_0$ \eqref{eq:eta0}, we 
convert the law \eqref{5E7} giving $k$ as a function of $\varepsilon$ into a law giving $\varepsilon$ as a function of $k$
\begin{equation}
\label{eq:epsk}
   \varepsilon = k^{-4} \gamma^4_{\min}
\end{equation}
and insert it into the identity \eqref{eq:reduB}. We obtain
\begin{equation}
\label{scriabine}
   \big(\bM\kk +\gamma^{8}_{\min}k^{-4}\bB_0\big)\bV\kk -  \Lambda\kk \bA\bV\kk =
   \bV_{0} \circ \big(\rH\kk +\gamma^{8}_{\min}k^{-4}\rB_0 - \Lambda\kk\big)\,.
\end{equation}
So the series $\Lambda\kk$ starts with the first eigenvalue $\Lambda_4=\gamma^{4}_{\min}\am_1$ of the operator $\rH_4+\gamma^{8}_{\min}\rB_0$. Then $\eta_0$ \eqref{eq:eta0} is an associated eigenvector. Like before, but now with this new $\eta_0$, and $k=k(\varepsilon)$, there exist further terms $\eta_1$, \ldots, $\eta_6$ such that the pair $(\breve\Lambda{}^k,\breve\zetaz{}^k)$ defined by \eqref{eq:qmm} is a quasimode for $\gM^k+\gamma^{8}_{\min}k^{-4}\bB_0 = \gM^k+\varepsilon^2k^{4}\bB_0$ with membrane boundary conditions. Since with law \eqref{eq:epsk} the terms $\varepsilon^2(\gB^{k(\varepsilon)}-{k(\varepsilon)}^4\bB_0)$ are of order ${k(\varepsilon)}^{-5}$ or higher, the same pair
\begin{equation}
\label{eq:qmb0}
   \breve\Lambda{}^{k(\varepsilon)} = k(\varepsilon)^{-4}\gamma^{4}_{\min}\am_1 = \varepsilon\am_1
   \quad\mbox{and}\quad
   \breve\zetaz{}^{k(\varepsilon)} = \big( 0,0,\eta_0\big)^{\!\top} 
   + \ \mbox{higher order terms in $k(\varepsilon)^{-1}$}
\end{equation}
is a quasimode for the full Koiter operator $\gK^{k(\varepsilon)}(\varepsilon)$, but still with the sole {\em membrane boundary conditions}.

\subsubsection{Quasimodes for the Koiter model at angular frequency $k(\varepsilon)$. Bending boundary layers}
\label{sss:bbl}
The full bending boundary conditions $\zeta_3=0$ and $\zeta'_3=0$ on $\partial\cI$ cannot be implemented in general for the quasimodes $(\breve\Lambda{}^{k(\varepsilon)},\breve\zetaz{}^{k(\varepsilon)})$. The singularly perturbed nature of the Koiter operator causes the loss of these boundary conditions between the bending and membrane operator. Solutions of the Koiter model, just as eigenvectors, incorporate boundary layer terms. In all cases investigated in this paper, these terms exist at the scale $d/\sqrt{\varepsilon}$ with $d=\dist(z,\partial\cI)$. Such a scaling appears in \cite{PitkarantaMatacheSchwab} in a variety of nondegenerate cases (the boundary of $\partial\cS$ is noncharacteristic for the curvature). It is rigorously analyzed in \cite{Faou2} in the case of static clamped elliptic shells.

More precisely, the scaled variable is (for $\cI=(z_-,z_+)$)
\begin{equation}
\label{eq:Z}
   Z = \frac{d}{\sqrt{\varepsilon}}
   \quad\mbox{with}\quad
   d = z_+-z \;\;\mbox{or}\;\; z-z_-,
\end{equation}
according as we consider the localization at the end $z_0=z_+$ or $z_0=z_-$ of the interval $\cI$.
In view of law \eqref{5E7}, we can write the operator $\gK^{k(\varepsilon)}(\varepsilon)$ as a series in powers of $\varepsilon^{1/4}$. In the rapid variable $Z$, there holds $\partial_z G(Z) = \varepsilon^{-1/2}G'$ for any profile $G(Z)$, which provides a new formal series $\cK[\varepsilon^{1/4}]$. Its leading term $\cK_0$ is compatible with the full bending boundary conditions at $Z=0$. It has the following form in the cylindrical case $f=R$
$$
   \cK_0 = \frac{E}{1-\nu^2} \begin{pmatrix}
   -\partial_{Z}^2 & 0 & \frac{\nu}{R}\partial_{Z}\\
   0 & -\frac{1-\nu}{2R^2}\partial_{Z}^2 & 0\\
   -\frac{\nu}{R}\partial_{Z} & 0 &\frac{1}{R^2}+\frac{1}{3}\partial_{Z}^4
\end{pmatrix}.
$$
It allows to construct a series of exponentially decreasing vector functions $\mathbf{G}[\varepsilon^{1/4}]$  satisfying a formal series relation of the type $\cK[\varepsilon^{1/4}] \mathbf{G}[\varepsilon^{1/4}] = \Lambda[\varepsilon^{1/4}]\mathbf{G}[\varepsilon^{1/4}] $,  that compensate for the missing traces of $\breve\zetaz{}^{k(\varepsilon)}$, see \cite[Section 5.6]{beaudouin:tel}. Our ``true'' quasimode has now the form $(\Lambda{}^{k(\varepsilon)},\zetaz{}^{k(\varepsilon)}(\varepsilon))$ with
\begin{equation}
\label{eq:qmb}
   \Lambda^{k(\varepsilon)} =\breve\Lambda{}^{k(\varepsilon)}=\am_1\varepsilon
   \ \quad\mbox{and}\quad\ 
   \zetaz^{k(\varepsilon)}(\varepsilon)(z) = \breve\zetaz{}^{k(\varepsilon)}(z)
   + \chi(d)\sum_{n=2}^6 \varepsilon^{n/4} \mathbf{G}_n(Z)\,.
\end{equation}
Here $\chi$ is a smooth cut-off that localizes near the boundary $\partial\cI$.
The outcome is the spectral estimate
\begin{equation}
\label{eq:koi5.2}
   \dist\big(\am_1\varepsilon\,,\,\gS(\gK^{k(\varepsilon)}(\varepsilon))\big) \lesssim \varepsilon^{5/4}
   \;\;\mbox{with}\;\; k(\varepsilon)=\varepsilon^{-1/4}\gamma_{\min},
   \quad \mbox{as}\quad\varepsilon\to0.
\end{equation}

\subsection{3D reconstruction and Rayleigh quotients}
We construct a three-component vector field on the surface $\cS$ by setting in normal coordinates
\[
   \zetaz^\varepsilon(z,\varphi) = e^{ik\varphi} \zetaz{}^k(z)
   \quad\mbox{with $\;\zetaz^k= \zetaz^{k(\varepsilon)}(\varepsilon)\;$ \eqref{eq:qmb}, \eqref{eq:qmb0} \ and
   $\;k=\near{k(\varepsilon)}=\near{\varepsilon^{-1/4}\gamma_{\min}}$.}
\]
By construction, $\zetaz^\varepsilon$ belongs to the variational space $V_\aB(\cS)$, and by the elliptic regularity of the Koiter problem, it also belongs to $(H^2\times H^2\times H^3)(\cS)$. So we may apply the reconstruction operator introduced in Theorem \ref{th:rayl}: Set
\[
   \bu^\varepsilon = \ocirc\rU\zetaz^\varepsilon.
\]
To take advantage of the comparison \eqref{eq:rayl} between the Rayleigh quotients of $\zetaz^\varepsilon$ and $\bu^\varepsilon$, we have to exhibit the behavior of the wave lengths  $ L = L^\varepsilon$ \eqref{eq:L} and $\ell = \ell^\varepsilon$ \eqref{eq:ell}  of $\zetaz^\varepsilon$ as $\varepsilon\to0$. Following the construction of the fields $\zetaz^\varepsilon$, we see that they all originate from an eigenfunction $\eta_0$ that does not depend on $\varepsilon$. The nontrivial behavior of $L^\varepsilon$ and $\ell^\varepsilon$ arises from, cf.\ Example \ref{ex:bd}:
\begin{itemize}
\item The Koiter boundary layer terms $\mathbf{G}_n(Z)=\mathbf{G}_n(d/\varepsilon^{1/2})$ that contribute a term in $\varepsilon^{1/2}$,
\item The azimuthal oscillation $e^{ik\varphi}$ that contributes a term in $k^{-1}\simeq\varepsilon^{1/4}$.
\end{itemize}
As a result we find in the nondegenerate parabolic case
$
   L^\varepsilon,\;\ell^\varepsilon \gtrsim \varepsilon^{1/2}.
$
So the assumptions of Theorem \ref{th:rayl} are uniformly satisfied for the family $(\zetaz^\varepsilon)_\varepsilon$ and the estimate \eqref{eq:rayl} reads now
\[
   \big\vert  Q^\varepsilon_{\aK}(\zetaz^\varepsilon) -
   Q^\varepsilon_{\aL}(\bu^\varepsilon)\big\vert
   \lesssim
   \varepsilon^{1/4} \,Q^\varepsilon_{\aK}(\zetaz^\varepsilon) \lesssim \varepsilon^{5/4}\, .
\]
Stricto sensu, we have at hand a family of 3D displacements $\bu^\varepsilon$ with azimuthal frequency $\near{k(\varepsilon)}\equiv \varepsilon^{-1/4}\gamma_{\min}$ such that
\[
  \big\vert  Q^\varepsilon_{\aL}(\bu^\varepsilon) - \am_1\varepsilon\big\vert \lesssim \varepsilon^{5/4}.
\]
So, with $\am_1$ and $\gamma=\gamma_{\min}$ defined in \eqref{5E6b}, we have proved the results summarized in the first two lines of Table \ref{tab:1}. 
By construction, in normal coordinates:
\begin{equation}
\label{eq:qmcyl}
   \bu^\varepsilon\big|_{\cS}(z,\varphi) = e^{i\near{k(\varepsilon)}\varphi} \begin{pmatrix}
   0,\ 0,\  \eta_0(z)
   \end{pmatrix}^\top \ \ \mbox{modulo higher order terms as $\varepsilon\to0$,}
\end{equation}
with $\eta_0$ the generating scalar eigenvector \eqref{eq:eta0}. 
Our numerical experiments (Model A, sect.\,\ref{modelA}, and Model B, sect.\,\ref{modelB}) suggest that, in fact, $(\am_1\varepsilon,\bu^\varepsilon)$ is an approximation of the first 3D eigenpair.

\section{Elliptic case (small meridian curvature)}
\label{s:7}
The elliptic case in parametrization $r=f(z)$, $z\in\cI$, corresponds to the situation $f''<0$ on $\overline\cI$.
After an exposition of the general principles of scalar reduction in the elliptic case, we address separately three different families of axisymmetric shells: Gaussian, Airy and toroidal.

\subsection{Membrane scalar reduction in the general case}
\label{ss:memred}
When the parametrizing function $f$ is not affine, i.e., when $f''\not\equiv0$, the scalar reduction of the membrane operator has non-vanishing first terms as follows:
\begin{equation}
\label{6E1}
\rH_0(z,\partial_z) = E \frac{f''^2}{s^6},\quad\; \rH_1(z,\partial_z) = 0,
\quad\;
\rH_{2}(z,\partial_{z}) = \rH_2^{(2)}\!(z)\,\partial^2_z + \rH_2^{(1)}\!(z)\,\partial_z + \rH_2^{(0)}\!(z)
\end{equation}
with
\begin{equation}
\label{6E2}
\left\{
\begin{aligned}
\rH_2^{(2)}\!(z)&=2E\Big(\frac{ff''}{s^6}+\frac{f^2f''^2}{s^8}\Big)\\[1ex]
\rH_2^{(1)}\!(z)&=2E\Big(\frac{2f'f''}{s^6}+\frac{ff'''}{s^6} -\frac{2ff'f''^2}{s^8}+\frac{2f^2f''f'''}{s^8}-\frac{7f^2f'f''^3}{s^{10}}\Big)\\[1ex]
\rH_2^{(0)}\!(z)&=E\Big(
- \frac{10f'^2f''^2}{s^8} 
+\frac{4f'f'''}{s^6}+\frac{2f'^2f''}{fs^6}-\frac{(\nu-2)ff'^2f''^3}{s^{10}} -\frac{5ff'f''f'''}{s^8}\\[1ex]
&\hskip 4ex
+\frac{ff^{(4)} }{s^6} + \frac{2f^2f''f^{(4)}}{s^8} 
+\frac{36f^2f'^2f''^4}{s^{12}}+\frac{(\nu-2)ff''^3}{s^8}-\frac{6f^2f''^4}{s^{10}}\\[1ex]
&\hskip 4ex
-\frac{20f^2f'f''^2f'''}{s^{10}}\Big) - \Lambda_0 \Big( \frac{1}{s} - \frac{\nu f'' f}{s^3}\Big)^2.
\end{aligned}
\right.
\end{equation}
The rank-3 operator in the formal series  $\rH\kk$  is given by
\begin{equation}
\label{eq:L3}
   \rH_{3}(z,\partial_{z})=\Big(-\frac{1}{s^2}+\frac{2\nu ff''}{s^4}-\frac{\nu^2f^2f''^2}{s^6}\Big)\Lambda_1,
\end{equation}
and the rank-4 operator can be written as
\begin{equation}
\label{eq:L4}
   \rH_{4}(z,\partial_{z})=
   \sum_{j=0}^4\rH_{4}^{(j)}(z)\partial_z^j,
   \quad\mbox{with}\quad
   \rH_{4}^{(4)}(z) = E\Big(\frac{4f^3f''}{s^8}+\frac{3f^4f''^2}{s^{10}} +\frac{f^2}{s^6}\Big)
\end{equation}
where the other terms $\rH_4^{(j)}(z)$ are smooth functions of $z$.

So,  $\rH_0(z,\partial_z)=\rH_0(z)$ is the multiplication by a function (which can be seen as a potential) and we check that $\rH_2$ is a selfadjoint operator of order $2$ on $H^1_0(\cI)$ with respect to the natural measure $\rd\cI = f(z)s(z)\,\rd z\,$:
\begin{equation}
\label{6E3}
   \big\langle \rH_2\eta,\etat \big\rangle_{\cI} =
   \int_{\cI} \Big(-\rH_2^{(2)}\!(z) \,\partial_z\eta\,\partial_z\etat
   + \rH_2^{(0)}\!(z)\,\eta\,\etat \Big)\,\rd\cI\,.
\end{equation}

We recall from \eqref{2E2} that the principal curvatures are $b^z_z=\frac{f''}{s^{3}}$ and $b^\varphi_\varphi=- \frac{1}{fs}$. Note that both are negative in the elliptic case.
\begin{remark}
(i) The function $\rH_0/E$ coincides with the square of the meridian curvature
\[
   \rH_0 = E\, (b^z_z)^2.
\]
(ii) There holds the following relation between $\rH^{(2)}_2$ and the principal curvatures
\begin{equation}
\label{6R2}
   -\rH^{(2)}_2 = 2E\,\frac{f^2}{s^2}\,b^z_z(b^\varphi_\varphi-b^z_z).
\end{equation}
(iii) Similarly
\begin{equation}
\label{6R4}
   \rH^{(4)}_4 = E\,\frac{f^4}{s^4}\,(b^\varphi_\varphi-3b^z_z)(b^\varphi_\varphi-b^z_z).
\end{equation}
\end{remark}

\subsection{High frequency analysis of the membrane operator in the elliptic case}
As mentioned above, we have to select one or several terms starting the series $\rH\kk$ that will play the role of an engine to work out a recurrence and allow to solve the formal series problem \eqref{scalreducfs}. In the parabolic case, this engine is $h^4\rH_4$. In the elliptic case, $\rH_0$ is the multiplication by the positive function $E(b^z_z)^2$. Its spectrum is essential and its bottom determines $\Lambda_0$
\begin{equation}
\label{eq:Lam0}
   \Lambda_0 = E\min_{z\in\overline\cI} (b^z_z)^2.
\end{equation}
We have to complete $\rH_0$ by further terms so that to obtain an operator with discrete spectrum close to the minimum energy $\Lambda_0$. This will be the case for the operator
\begin{equation}
\label{eq:H0+H2}
   \hM^k =  \rH_0+k^{-2}\rH_2
\end{equation}
if, cf.\ condition \eqref{eq:azi},
\begin{equation}
\label{eq:ellipH2}
   -\rH^{(2)}_2\ge0 \ \ \mbox{on}\ \ \overline\cI,\quad\mbox{i.e.}\quad
   |b^\varphi_\varphi| \ge |b^z_z| \ \ \mbox{on}\ \ \overline\cI,
\end{equation}
(use \eqref{6R2}), with strict inequalities for the values of $z$ where $\rH_0$ attains its minimum $\Lambda_0$.
It is interesting to note that the latter condition implies that, cf.\ \eqref{eq:Lam0} and \eqref{eq:essMk},
\[
   \min_{z\in\overline\cI} (b^z_z)^2 < \min_{z\in\overline\cI} (b^\varphi_\varphi)^2,
   \quad\mbox{i.e.}\quad
   \Lambda_0 < \min \sigma_{\sf ess}(\gM^k),
\]
which means that the expected limit at high frequency will be attained by eigenvalues 
below the essential spectrum.

\begin{remark}
We note that in the hyperbolic case, $f''>0$, so $b^z_z>0$. Hence the coefficient $-\rH^{(2)}_2$ is {\em always negative} and our analysis {\em never applies} in the hyperbolic case. Besides, in this case, $\Lambda_0$ is not the membrane high frequency limit, that is indeed $0$ (recall that the exponent in \eqref{order} is $\alpha = \frac23$ in hyperbolic case).
\end{remark}

From now on, we assume that \eqref{eq:ellipH2} holds and we discuss the lowest eigenpairs of the operators $\hM^k$ defined in \eqref{eq:H0+H2} and
\begin{equation}
\label{abegg}
\hK^k(\varepsilon) =  \rH_0+k^{-2}\rH_2 + \varepsilon^2k^4\rB_0,
   \quad\mbox{where}\quad \rB_0 = \frac{1}{3} \,\frac{E}{1-\nu^2} \, \frac{1}{f^4}
\end{equation}
in relation with properties of the ``potential'' $\rH_0$. For simplicity we denote
\begin{equation}
\label{eq:sigma}
   \rh(z) := -\rH^{(2)}_2(z),
\end{equation}
and consider successively the cases when $\rH_0$ has a non-degenerate minimum inside or on the boundary of the interval $\cI$, or when it is constant.

\subsection{Internal minimum of the potential (Gaussian case)}
\label{ss:gauss}
Besides \eqref{eq:ellipH2}, we assume that $\rH_0$ has a (unique) nondegenerate minimum in $z_0\in\cI$. Thus
\[
   \Lambda_0=\rH_0(z_0) \quad\mbox{and}\quad \partial_z^2\rH_0(z_0)>0.
\]
We assume moreover $\rh(z_0)>0$.

\subsubsection{High frequency analysis for the membrane operator}
Then the lowest eigenpairs of the membrane reduction  $\hK^k = \rH_0+k^{-2}\rH_2$  as $k \to \infty$ are driven by the harmonic oscillator
\begin{equation}
\label{6E10}
     -\rh(z_0) \,\partial_Z^2 + \frac{Z^2}{2}\, \partial_z^2\rH_0(z_0)\,.
\end{equation}
Here, the new homogenized variable $Z$ spans $\R$ and is linked to the physical variable $z$ by the relation
\begin{equation}
\label{6E11}
   Z =  \sqrt{k}\, (z-z_0)\,.
\end{equation}
This change of variable can be applied to the formal series reduction \eqref{eq:reduc} as follows: Let  $L[k] = \sum_{k\geq 0} k^{-n} L_n(z,\partial_z)$ be a formal series such that $L_n$ is an operator of order $n$. By Taylor expansion around $z_0$, we can expand for all $n$ the operator
$L_n(z,\partial_z) = \sum_{j \geq -n} k^{-j/2}\mathcal{L}_{n,j}(Z,\partial_Z)$. By reordering the powers of $k^{-j/2}$, we thus see that we can write
$$
   L[k] \equiv \mathcal{L}[k^{1/2}] = \sum_{n\geq 0} k^{-n/2} \mathcal{L}_n(Z,\partial_Z),
$$
where the operators $\mathcal{L}_n$ have polynomial coefficients in $Z$.
Applying this change of variable to the formal series reduction given in
Theorem \ref{reducform}, we obtain the new identity
\begin{equation}
\label{eq:redb}
   {\boldsymbol{\mathcal{M}}}[k^{1/2}]{\boldsymbol{\mathcal{V}}}[k^{1/2}] -
   \Lambda[k^{1/2}]\bA{\boldsymbol{\mathcal{V}}}[k^{1/2}]
   = \bV_{0} \circ (\mathcal{H}[k^{1/2}] - \Lambda[k^{1/2}] )\,,
\end{equation}
where ${\boldsymbol{\mathcal{M}}}[k^{1/2}]$, ${\boldsymbol{\mathcal{V}}}[k^{1/2}]$ and $\mathcal{H}[k^{1/2}]$ are the formal series induced by the formal series $\bM\kk$, $\bV\kk$ and $\rH\kk$ respectively. $\bV_{0}$ is still the embedding $\eta\mapsto(0,0,\eta)^\top$. We also agree that $\Lambda[k^{1/2}]$ is related with the old series $\Lambda_{\sf old}[k]=\sum_{n\ge0}k^{-n}\Lambda_{n,\sf old}$ by the identities $\Lambda_n=0$ if $n$ is odd, and $\Lambda_n=\Lambda_{n/2,\sf old}$ if $n$ is even.
Moreover, we calculate that
$$
   \cH_0= \rH_0(z_0),\quad  \cH_{1} = 0,\quad \mbox{and} \quad
   \cH_2 = -\rh(z_0) \,\partial_Z^2 + \frac{Z^2}{2}\, \partial_z^2\rH_0(z_0).
$$
Like for \eqref{scalreducfs}, the previous reduction leads to consider the formal series problem
\begin{equation}
\label{franck}
\cH[k^{1/2}] \eta[k^{1/2}] = \Lambda[k^{1/2}] \eta[k^{1/2}]\,.
\end{equation}
The first equation induced by this identity is $\cH_0\eta_0 = \Lambda_0\eta_0$, hence we have found again $\Lambda_0=\cH_0=\rH_0(z_0)$. Since for any $\eta$ we have now $\cH_0\eta=\Lambda_0\eta$, the next equations yield
\[
   \cH_1\eta_0 = \Lambda_1\eta_0\quad\mbox{and}\quad
   \cH_2\eta_0 = \Lambda_2\eta_0.
\]
Therefore $\Lambda_1=0$ (which is coherent with what was agreed in identity \eqref{eq:redb}) and $\eta_0$ is an eigenvector of the harmonic oscillator \eqref{6E10}.
The eigenvalues of this latter operator are
\begin{equation}
\label{eq:harmosc}
   (2\ell-1)\,\rc \,,\quad \ell=1,2,\ldots \quad\mbox{with}\quad
   \rc = \frac{1}{\sqrt2} \sqrt{ \rh(z_0) \,\partial_z^2\rH_0(z_0)}
\end{equation}
and the corresponding eigenvectors are Gaussian functions.
Taking $\eta_0(Z)$ as the first eigenmode ($\ell = 1$) we can construct the first terms $\eta_{1}$, $\eta_2, \ldots$ of the formal series problem \eqref{franck}. As the coefficients of the operators $\cH_{j}$ depend polynomially on $Z$, these terms are exponentially decreasing with respect to $Z$.
We can then define the pair $(\breve\Lambda{}^k,\breve\zetaz{}^k)$ by the formula
\begin{equation}
\label{grieg}
\begin{aligned}
   \breve\zetaz{}^k &= \chi( z ) \Big( \bV_0 \eta_0 +
   \sum_{1\le n+m\le6} k^{-(n+m)/2} \,{\boldsymbol{\mathcal{V}}}_n\eta_m \Big)  \big( \sqrt{k}(z - z_0) \big)  \\[-1ex]
   \breve\Lambda{}^k &=  \rH_0(z_0) + k^{-1} \rc\,,
\end{aligned}
\end{equation}
where $\chi\in C^\infty_0(\cI)$ is identically equal to $1$ in a neighborhood of $z_0$.
This pair is a quasimode for the full membrane operator $\gM^k$ as $k\to\infty$,
and we obtain that
\begin{equation}
\label{eq:memb6.3}
   \dist\big(\rH_0(z_0) + k^{-1} \rc\,,\,\gS(\gM^k)\big) \lesssim k^{-3/2},\quad k\to\infty.
\end{equation}
Note that in this case, the boundary conditions are automatically fulfilled as the quasimode constructed is localized near $z_0$.

\subsubsection{High frequency analysis for Koiter and Lam\'e operators}
Now we consider the operator $\hK^k(\varepsilon)$ defined in \eqref{abegg}.
We define its smallest eigenvalue $\lamone{\hK^k(\varepsilon)}$ and for each $\varepsilon>0$ small enough, look for $k(\varepsilon)$ such that $\lamone{\hK^k(\varepsilon)}$ is minimum.
Setting $\delta:=\varepsilon^2k^4$, we see that the operator $\hK^k(\varepsilon)$ has the form
\[
   W +  k^{-2}  \rH_2\quad\mbox{with}\quad W = \rH_0+\delta\rB_0.
\]
If $\delta$ is small enough, the function $W$ has the same property as $\rH_0$, i.e., it has a (unique) nondegenerate minimum. Let $z_0(\delta)$ be the point where this minimum is attained. By implicit function theorem, the correspondence $\delta\to z_0(\delta)$ is smooth for $\delta$ small enough and there holds, cf \eqref{eq:harmosc}
\[
\begin{aligned}
   \lamone{\hK^k(\varepsilon)} &=
   \rH_0(z_0(\delta))+\delta\rB_0(z_0(\delta)) +
   \frac{k^{-1}}{\sqrt2}
   \sqrt{ \rh (z_0(\delta)) \,\partial_z^2(\rH_0+\delta\rB_0)(z_0(\delta))}
   + \cO(k^{-3/2}) \\
\end{aligned}
\]
But
\[
\begin{aligned}
   &\rH_0(z_0(\delta)) = \rH_0(z_0) + \cO(\delta^2),\quad
   && \partial_z^2\rH_0(z_0(\delta)) = \partial_z^2\rH_0(z_0) + \cO(\delta),\\
   &\rB_0(z_0(\delta)) = \rB_0(z_0) + \cO(\delta),\quad
   &&  \rh(z_0(\delta)) = \rh(z_0) + \cO(\delta).
\end{aligned}
\]
Hence
\[
   \lamone{\hK^k(\varepsilon)} =
   \rH_0(z_0)+\delta\rB_0(z_0) +
   \frac{k^{-1}}{\sqrt2} \sqrt{ \rh(z_0) \,\partial_z^2\rH_0(z_0)}
   +\cO(\delta^2) + \cO(k^{-1}\delta) + \cO(k^{-3/2}).
\]
Let us set
\begin{equation}
\label{6E14bc}
   \rb=\rB_0(z_0)\quad\mbox{and}\quad \rc = \frac{1}{\sqrt2} \sqrt{ \rh(z_0) \,\partial_z^2\rH_0(z_0)}\,.
\end{equation}
So, replacing $\delta$ by its value $\varepsilon^2k^4$, we look for $k=k(\varepsilon)$ such that $ \varepsilon^2k^4 \rb + k^{-1} \rc$ is minimum and such that $\delta = \varepsilon^2 k^4$ is small\footnote{We check that $\delta = \mathcal{O}(\varepsilon^{2/5})$}. We homogenize the powers of $k$ by letting
$\gamma(\varepsilon) = k(\varepsilon)\, \varepsilon^{2/5}$,
and setting $\mu_1^{\hK}(\varepsilon) = \lamone{\hK^{k(\varepsilon)}(\varepsilon)}$ we find
\begin{equation}
\label{eq:6E16}
   k(\varepsilon) = \gamma\, \varepsilon^{-2/5}\quad\mbox{and}\quad
   \mu_1^{\hK}(\varepsilon) = \rH_0(z_0) + \am_1\,\varepsilon^{2/5} + \cO(\varepsilon^{3/5}),
\end{equation}
with the explicit constants $\gamma$ and $\am_1$:
\begin{equation}
\label{6E15}
   \gamma = \Big(\frac{\rc}{4\rb}\Big)^{1/5}\quad\mbox{and}\quad
   \am_1 = (4\rb \rc^{4})^{1/5} (1 + \frac14) \,.
\end{equation}
We find that the ratio $\rR$ of energies \eqref{eq:Rratio} is
\begin{equation}
\label{eq:Rratiogauss}
   \rR \simeq \frac{\varepsilon^2k^4\rb}{\rH_0(z_0) + \am_1\,\varepsilon^{2/5}}
   \simeq \frac{\rb}{\rH_0(z_0)} \Big(\frac{\rc}{4\rb}\Big)^{4/5}\,\varepsilon^{2/5}.
\end{equation}
Along the same lines as in the parabolic case, we convert the power law for $k$ \eqref{eq:6E16} into the power law $\varepsilon=(\gamma/k)^{5/2}$. We can then consider a formal series reduction as in \eqref{scriabine} and combine it with the change of variable $Z = \sqrt{k}( z - z_0)$.  The same analysis as before yields quasimodes 
$(\breve\Lambda{}^{k(\varepsilon)}, \breve\zetaz{}^{k(\varepsilon)} )$. Here $\breve\Lambda{}^{k(\varepsilon)}=\mu_1^{\hK}(\varepsilon)$ and $\breve\zetaz{}^{k(\varepsilon)}$ has a form similar to \eqref{grieg}, with $k=k(\varepsilon)$. 
Note that these quasimodes remain localized around $z_0$ and hence bending boundary layers do not show up as they did in the parabolic case.
We thus obtain
\begin{equation}
\label{eq:6F16}
   \dist\big(\mm_1(\varepsilon)\,,\, \gS(\gK^{k(\varepsilon)}(\varepsilon))\big) \lesssim \varepsilon^{3/5}
   \quad\mbox{with}\quad \mm_1(\varepsilon) = \rH_0(z_0) + \am_1\,\varepsilon^{2/5}.
\end{equation}

\begin{remark}
\label{rem:frank}
If $\rH_0$ attains its minimum $\am_0$ in a finite number of points $z_0^{(i)}$, we can construct quasimodes attached to each of these points of the same form as above, and with disjoint supports. The associated quantities obey to the same formulas as in \eqref{eq:6E16}-\eqref{eq:6F16}
\[
   k^{(i)}(\varepsilon) = \gamma^{(i)}\, \varepsilon^{-2/5}\quad\mbox{and}\quad
   \mm^{(i)}_1(\varepsilon) = \am_0 + \am_1^{(i)}\,\varepsilon^{2/5} ,
\]
with $\gamma^{(i)}$ and $\am_1^{(i)}$ defined by \eqref{6E15} with the values of quantities $\rb$ and $\rc$ at point $z_0^{(i)}$. Then $\mm_1(\varepsilon)=\min_i \mm^{(i)}_1(\varepsilon)$ and $k(\varepsilon)=k^{(i_0)}(\varepsilon)$ for $i_0$ such that the previous minimum is attained.
\end{remark}

\subsubsection{3D reconstruction and Rayleigh quotients}
As in the parabolic case,
we construct a three-component vector field on the surface $\cS$ by setting in normal coordinates
\[
   \zetaz^\varepsilon(z,\varphi) = e^{ik\varphi} \zetaz{}^k(z)
   \quad\mbox{with $\;\zetaz{}^k\;$ given in \eqref{grieg} \ and $\;k=\near{\varepsilon^{-2/5}\gamma_{\min}}$},
\]
and by setting
\[
   \bu^\varepsilon = \ocirc\rU\zetaz^\varepsilon.
\]
Since all traces of any order of $\zetaz^\varepsilon$ vanish on $\partial\cI$, there is no boundary corrector  and we are in case {\em (ii)} of Theorem \ref{th:rayl}. So we only have to estimate the behavior of the wave length  $ L = L^\varepsilon$ \eqref{eq:L} of $\zetaz^\varepsilon$ as $\varepsilon\to0$. We note the influence of:
\begin{itemize}
\item The profiles $\mathbf{G}_n(\sqrt{k} ( z - z_0 ))$ with $k \simeq \varepsilon^{-2/5}$,  that contribute a term in $1/\sqrt{k} \simeq \varepsilon^{1/5}$,
\item The azimuthal oscillation $e^{ik\varphi}$ that contributes a term in $k^{-1}\simeq\varepsilon^{2/5}$.
\end{itemize}
As a result we find in the nondegenerate parabolic case
\[
   L^\varepsilon  \gtrsim \varepsilon^{2/5}.
\]
So the assumptions of Theorem \ref{th:rayl} are uniformly satisfied for the family $(\zetaz^\varepsilon)_\varepsilon$ and the estimate \eqref{eq:rayl} reads now
\[
   \big\vert  Q^\varepsilon_{\aK}(\zetaz^\varepsilon) -
   Q^\varepsilon_{\aL}(\bu^\varepsilon)\big\vert
   \lesssim
   \varepsilon \,Q^\varepsilon_{\aK}(\zetaz^\varepsilon) \lesssim \varepsilon\,.
\]
Thus, we have exhibited a family of 3D displacements $\bu^\varepsilon$ with azimuthal frequency $k(\varepsilon)\equiv \varepsilon^{-2/5}\gamma$ such that
\[
  \big\vert  Q^\varepsilon_{\aL}(\bu^\varepsilon) - \mm_1(\varepsilon) \big\vert
  \lesssim \varepsilon^{3/5}\quad\mbox{with}\quad
  \mm_1(\varepsilon) = \rH_0(z_0) + \am_1\,\varepsilon^{2/5}.
\]
So we have proved the results summarized in the third line of Table \ref{tab:1}. 
In normal coordinate system, there holds:
\begin{equation}
\label{eq:qmgauss}
   \bu^\varepsilon\big|_{\cS}(z,\varphi) = e^{i\near{k(\varepsilon)}\varphi} \begin{pmatrix}
   0, \ 0, \  \eta_0 \big( \sqrt{k(\varepsilon)}(z - z_0) \big) 
   \end{pmatrix}^\top \ \mbox{mod.\ higher order terms as $\varepsilon\to0$,}
\end{equation}
with $\eta_0$ the first eigenvector of the harmonic oscillator.
So, the principal term of $\bu^\varepsilon$ displays a meridian concentration at scale $\sqrt{k(\varepsilon)}\sim \varepsilon^{-1/5}$. 
The numerical experiments (Model H, sect.\,\ref{modelH}) suggest that $(\mm_1(\varepsilon),\bu^\varepsilon)$ is indeed an approximation of the first 3D eigenpair.


\subsection{Minimum of the potential on the boundary (Airy case)}
\label{ss:airy}
We assume that $\rH_0$ attains its minimum at a point $z_0\in\partial\cI$ with $\partial_z\rH_0(z_0)\neq0$. Let us agree that $z_0$ is the left end of $\cI$, i.e.,  $z_0=z_-$, so that we have 
\[
   \Lambda_0 = \rH_0(z_0) \quad\mbox{and}\quad \partial_z\rH_0(z_0)>0.
\]
We still assume $\rh(z_0)>0$. The analysis is somewhat similar to the previous case, though a little more tricky.

\subsubsection{High frequency analysis for the membrane operator. Membrane boundary layers}
We meet the Airy-like operator
\begin{equation}
\label{6EAiry1}
   -\rh(z_0) \partial_Z^2+Z\,\partial_z\rH_0(z_0)
\end{equation}
on $H^1_0(\R^+)$ instead the harmonic oscillator \eqref{6E10}. The homogenized variable $Z$ is given by
\begin{equation}
\label{6EAiry2}
  Z =  (z-z_0)k^{2/3}.
\end{equation}
We can perform an analysis very similar to the previous case by doing a change of variable in the formal series reduction. This yields formal series problem in powers of $k^{-1/3}$ whose first terms are given by $\cH_0+k^{-2/3}\cH_2$ where $\cH_0=\rH_0(z_0)$ and $\cH_2$ is the operator \eqref{6EAiry1}.
The eigenvalues of the model operator \eqref{6EAiry1} are given by
\begin{equation}
\label{6EAiry3}
   \mathsf{z}_{\mathsf Airy}^{(\ell)} \, \big(\rh(z_0) \big)^{1/3}\, \big(\partial_z\rH_0(z_0)\big)^{2/3},
   \quad \ell=1,2,\ldots
\end{equation}
where $\mathsf{z}_{\mathsf Airy}^{(\ell)}$ is the $\ell$-th zero of the reverse Airy function $\mathsf{Ai}$.
We find that the first eigenvalue of the membrane reduction $\hM^k$ satisfies
\begin{equation}
\label{6EAiry4}
   \mu_1^{\hM}(k) = \rH_0(z_0) + k^{-2/3} \rc + \cO(k^{-1})\quad \mbox{with}\quad
   \rc=\mathsf{z}_{\mathsf Airy}^{(1)} \, \big( \rh(z_0)\big)^{1/3}\,
   \big(\partial_z\rH_0(z_0)\big)^{2/3}.
\end{equation}
Using the reconstruction operators ${\boldsymbol{\mathcal{V}}}[k^{1/3}]$ in the scaled variable allows to construct displacement $\breve\zetaz{}^k$ from an eigenfunction profile $\eta_0(Z)$ of the Airy operator.
However, the first terms of the reconstruction take the form:
$$
\begin{pmatrix}
   k^{-4/3}\zeta^k_z \\  k^{-1}\zeta^k_\varphi \\  \zeta^k_3
\end{pmatrix}
\quad\mbox{where}\quad
\begin{pmatrix}
   \zeta^k_z \\  \zeta^k_\varphi \\  \zeta^k_3
\end{pmatrix} =
\begin{pmatrix}
   f^2( b_{\varphi}^\varphi -(\nu+2) b_z^z)\partial_{Z}\eta_0  \\
   - if^2 ( b_\varphi^\varphi + \nu b_{z}^z) \eta_0 \\
   \eta_0
\end{pmatrix}
+\cO(k^{-1/3})\,.
$$
While we can impose $\eta_0(0) = 0$ to ensure that $\zeta_\varphi^k = 0$ at first order, we see that we have in general $\zeta_z^k \neq 0$. To construct a quasimode, we have to add new boundary layer terms to $\breve\zetaz{}^k$.

To determine such boundary layers near  $z_0=z_{-}$, we introduce the scaled variable $ Z = k d$ with
 $d = z - z_{-}$. Like already seen for the Koiter operator (sect.\,\ref{sss:bbl}), the formal series operator $\bM\kk$ is changed to a new formal series $\boldsymbol{\mathscr{M}}\kk$ whose first term is given by
$$
\boldsymbol{\mathscr{M}}_0 = \begin{pmatrix}
- \frac{1}{s^4} \partial_Z^2 + \frac{1-\nu}{2} (b_\varphi^\varphi)^2 & - \frac{1 + \nu}{2}i (b_\varphi^\varphi)^2 \partial_Z & \frac{1}{s^2}(b_z^z + \nu b_\varphi^\varphi)\partial_{Z}  \\
- \frac{1 + \nu}{2}i (b_\varphi^\varphi)^2 \partial_Z &- \frac{1 - \nu}{2}(b_\varphi^\varphi)^2 \partial_Z^2 + \frac{1}{f^4} & i\frac{1}{f^2}(b_\varphi^\varphi + \nu b_z^z) \\
-\frac{1}{s^2}(b_z^z + \nu b_\varphi^\varphi)\partial_{Z} & - i\frac{1}{f^2}(b_\varphi^\varphi + \nu b_z^z)   & (b_\varphi^\varphi + \nu b_z^z)^2
\end{pmatrix}
$$
where the quantities are evaluated in  $z_0=z_{-}$.
We can prove that this operator yields boundary layer profiles $\mathbf{G}(Z)$ exponentially decreasing with respect to $Z = k d$, and satisfying $G_z(0) = a_z$ for any given number $a_z$, which allows to compensate for the trace of the first term of $\zeta^k_z$.

We obtain a compound quasimode combining terms at scale $k^{2/3}d$ and terms at scale $kd$, and deduce in the end
\begin{equation}
\label{eq:memb6.4}
   \dist\big(\rH_0(z_0) + k^{-2/3} \rc\,,\,\gS(\gM^k)\big) \lesssim k^{-1},\quad k\to\infty.
\end{equation}

\subsubsection{High frequency analysis for Koiter operator}
The Koiter scalar reduction operator $\hK^k(\varepsilon)$, see \eqref{abegg} is still an Airy-like operator because the minimum of $\rH_0+\varepsilon^2k^4$ is still $z_0$ for $\varepsilon^2k^4$ small enough. We look for $k=k(\varepsilon)$ such that $ \varepsilon^2k^4 \rb + k^{-2/3} \rc$ is minimum. We homogenize $\varepsilon^2k^4$ with $k^{-2/3}$. We find
\begin{equation}
\label{6E16}
   k(\varepsilon) = \gamma\, \varepsilon^{-3/7}\quad\mbox{and}\quad
   \mu_1^{\hK}(\varepsilon) = \rH_0(z_0) + \am_1\,\varepsilon^{2/7} + \cO(\varepsilon^{3/7}),
\end{equation}
with the explicit constants $\gamma$ and $\am_1$, with $\rb=\rB_0(z_0)$ and $\rc$ defined in \eqref{6EAiry4}:
\begin{equation}
\label{6E15b}
   \gamma = \Big(\frac{\rc}{6\rb}\Big)^{3/14}\quad\mbox{and}\quad
   \am_1 = (6\rb \rc^{6})^{1/7} (1 + \frac16) \,.
\end{equation}
The ratio of energies $\rR$ \eqref{eq:Rratio} is equivalent to $\delta\,\varepsilon^{2/7}$ with an explicit constant $\delta$, compare with \eqref{eq:Rratiogauss}. 
Note that
in this case, two types of boundary layer terms are present: the one constructed above (membrane boundary layer) and the bending boundary layers terms associated with the Koiter operator, see sect.\,\ref{sss:bbl}. We obtain 
\[
   \dist\big(\mm_1(\varepsilon)\,,\, \gS(\gK^{k(\varepsilon)}(\varepsilon))\big) \lesssim \varepsilon^{3/7}
   \quad\mbox{with}\quad
   \mm_1(\varepsilon) = \rH_0(z_0) + \am_1\,\varepsilon^{2/7}.
\]
At this point, the reconstruction operator $\ocirc\rU$ is not precise enough to allow us to conclude as in the parabolic and Gaussian cases. Using more elaborate reconstruction as in \cite{Faou2} we would find a 3D vector field $\bu^\varepsilon$ with elastic energy $\simeq \mm_1(\varepsilon)$ and expression in normal coordinate system
\begin{equation}
\label{eq:qmairy}
   \bu^\varepsilon\big|_{\cS}(z,\varphi) = e^{i\near{k(\varepsilon)}\varphi} \begin{pmatrix}
   0, \ 0, \  \eta_0 \big( k(\varepsilon)^{2/3}(z - z_0) \big) 
   \end{pmatrix}^\top \ \mbox{mod.\ higher order terms as $\varepsilon\to0$,}
\end{equation}
with $\eta_0$ the first eigenvector of the Airy operator. The dominant meridian concentration scale is $ k(\varepsilon)^{2/3}\sim\varepsilon^{-2/7}$.
Numerical experiments (Model L, sect.\,\ref{modelL}) tend to confirm that the first eigenmode of the Lam\'e operator $\gL(\varepsilon)$ behaves like $\big(\mm_1(\varepsilon),\bu^\varepsilon)$.

\subsection{Constant potential (toroidal case)}
\label{ss:const}
Let us assume that $\rH_0$ is constant. We recall that $\rH_0=E(b^z_z)^2$. But $b^z_z$ coincides with the curvature of the arc $\cC$ of equation $r=f(z)$ in the meridian plane. So, $b^z_z$ is constant if and only if $\cC$ is a \emph{circular} arc.
Let $R$ be its radius and $(r_\circ,z_\circ)\in\R^2$ be its center.
Notice that the center of the circular arc may be at negative $r_\circ$. Then, in the elliptic case $f''<0$,
\begin{equation}
\label{eq:circ}
   f(z) = r_\circ + \sqrt{R^2-(z-z_\circ)^2},
\end{equation}
and the principal curvatures are given by
\begin{equation}
\label{eq:circb}
   b^z_z = -\frac{1}{R} \quad\mbox{and}\quad
   b^\varphi_\varphi(z) =
   -\frac{1}{R} \Big(1-\frac{r_\circ}{f(z)}\Big).
\end{equation}
So  in this case, we have 
\begin{equation}
\label{eq:circa}
   \rH_0 = \frac{E}{R^2} = \Lambda_0\quad\mbox{and}\quad
   \rh =-\rH^{(2)}_2 = - 2E\,\frac{f}{s^2}\,\frac{r_\circ}{R^2} \,.
\end{equation}
Now the Koiter scalar reduction operator $\hK^k(\varepsilon)$ is $\rH_0+k^{-2}\rH_2 + \varepsilon^2k^4\rB_0$ where $\rH_0$ is a constant function acting as a simple shift on the spectrum. In this case, no concentration occurs, and   we have simply to come back to the approach used for the parabolic case {\em mutatis mutandis}, with $k^{-2}$ instead of $k^{-4}$.
\subsubsection{Membrane scalar reduction}
We assume the sharp version of condition \eqref{eq:ellipH2} (strict inequalities) that ensures that $\rH_2$ has a compact resolvent and is semibounded from below. Thanks to \eqref{eq:circa}, we find that such condition is equivalent to
\begin{equation}
\label{eq:circc}
   r_\circ<0\,.
\end{equation}
Let $\Lambda_2$ be the first eigenvalue of $\rH_2$. Then the first eigenvalue of the membrane scalar reduction operator $\hM^k = \rH_0+k^{-2}\rH_2$ is $\Lambda_0+k^{-2}\Lambda_2$ and we can deduce that
\begin{equation}
\label{eq:memb6.5}
   \dist\big(\Lambda_0+k^{-2}\Lambda_2\,,\,\gS(\gM^k)\big) \lesssim k^{-5/2},\quad k\to\infty.
\end{equation}

\subsubsection{Koiter scalar reduction}
\begin{figure}[ht]
\includegraphics[scale=0.51]{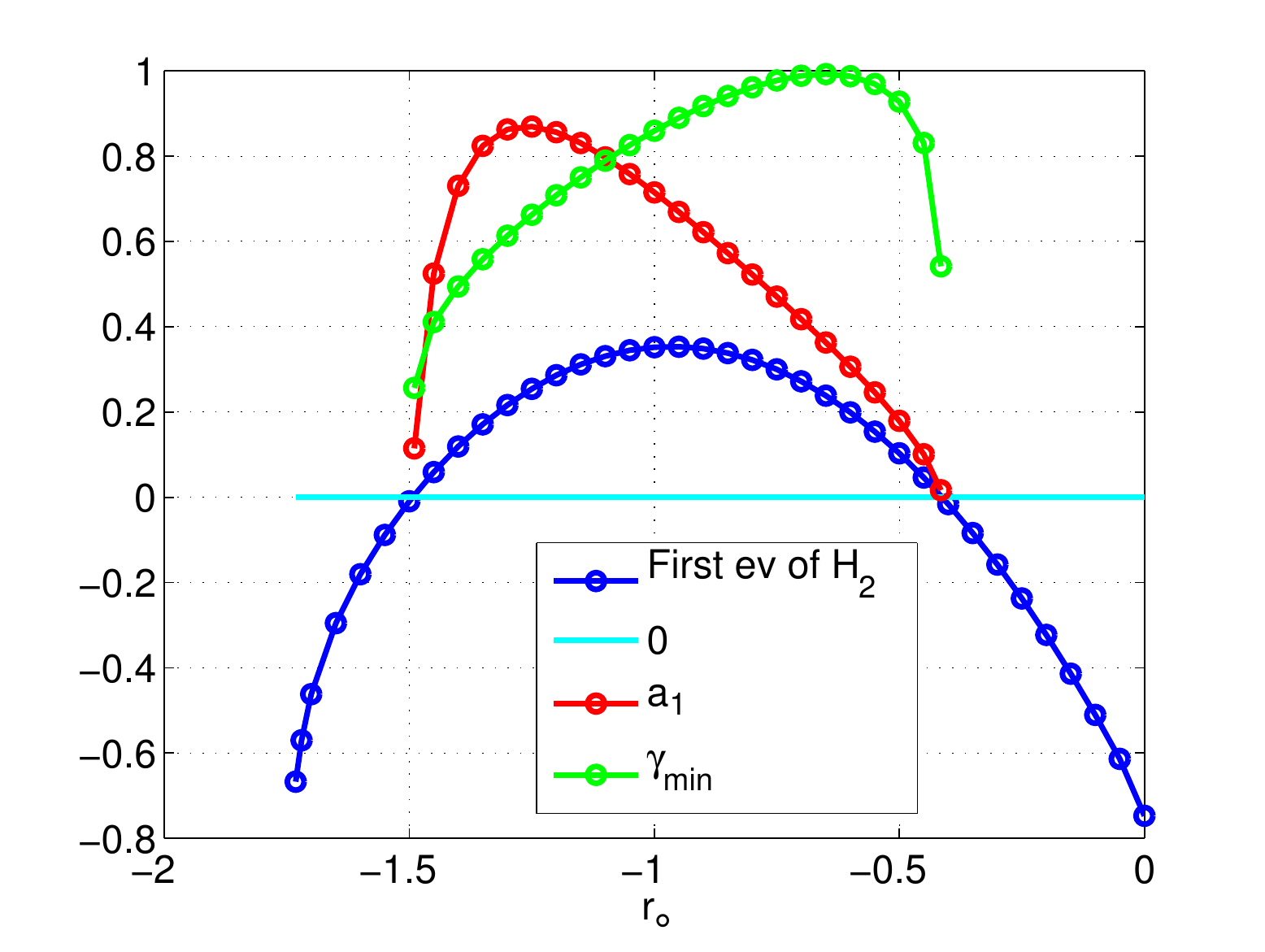}
   \caption{With $\cI=(-1,1)$, $R=2$ and $z_{\circ}=0$: 
   Quantities $\Lambda_2$, $\gamma_{\min}$ and $\am_1$ vs $r_{\circ}$.}
\label{VPH2opt}
\end{figure}

The operator $\hK^k(\varepsilon)=\rH_0+k^{-2}\rH_2 + \varepsilon^2k^4\rB_0$ is self-adjoint on $H^1_0(\cI)$. 
Its first eigenvalue is denoted by $\lamone{\hK^k(\varepsilon)}$. For any chosen $\varepsilon$ we look for $k(\varepsilon)\in\R^+$ realizing the minimum of $\lamone{\hK^k(\varepsilon)}$. We set
\begin{equation}
\label{5E5bis}
    \gamma(\varepsilon)  = k(\varepsilon) \;\varepsilon^{1/3}
\end{equation}
so that our operator becomes
\[
   \rH_0+\varepsilon^{2/3}\Big(\frac{1}{\gamma(\varepsilon)^2}\,\rH_2 +
   \gamma(\varepsilon)^4\,\rB_0\Big).
\]
Therefore $\gamma$ does not depend on $\varepsilon$. Let $\mu_1(\gamma)$ be the first eigenvalue of the operator
\begin{equation}
\label{5E6bis}
   \frac{1}{\gamma^2} \,\rH_2 + \gamma^4\, \rB_0\,.
\end{equation}
The function $\gamma\mapsto\mu_1(\gamma)$ is continuous. At this point we need the following extra assumption:
\begin{equation}
\label{eq:posit}
   \Lambda_2>0,\quad\mbox{i.e.}\quad \rH_2>0.
\end{equation}
Then the same argument as in the parabolic case allows to define $\gamma_{\min}$  as the (smallest) positive constant such that $\mu_1(\gamma)$ is minimum
\begin{equation}
\label{5E6bbis}
   \mu_1(\gamma_{\min}) = \min_{\gamma\in\R_+}\mu_1(\gamma) =: \am_1 \,.
\end{equation}
As an illustration of the non-trivial behavior of the quantities $\Lambda_2$, $\gamma_{\min}$ and $\am_1$, we plot them versus $r_{\circ}$ in Figure \ref{VPH2opt} (we choose $R=2$ and $z_{\circ}=0$).

Thus $k(\varepsilon)$ satisfies a power law that yields a formula for the minimal first eigenvalue $\mu_1^{\hK}(\varepsilon)$:
\begin{equation}
\label{5E7bis}
   k(\varepsilon) = \varepsilon^{-1/3} \gamma_{\min}\quad\mbox{and}\quad
   \mu_1^{\hK}(\varepsilon) =\mm_1(\varepsilon) = \rH_0 + \varepsilon^{2/3}\am_1\,,
\end{equation}
and after adding membrane and bending boundary layer terms as in the Airy case we arrive to
\[
   \dist\big(\mm_1(\varepsilon)\,,\, \gS(\gK^{k(\varepsilon)}(\varepsilon))\big) \lesssim \varepsilon
   \quad\mbox{with}\quad
   \mm_1(\varepsilon) = \rH_0 + \varepsilon^{2/3}\am_1.
\]
The ratio of energies $\rR$ \eqref{eq:Rratio} is equivalent to $\delta\,\varepsilon^{2/3}$. 
We note that, in contrast with the two previous cases (Gauss and Airy) when $\rH_0$ is not constant, the lower order term $\rH^{(0)}_2$ of the operator $\rH_2$ is involved in the asymptotics. Finally, like in the Airy case, we would need a more complete reconstruction operator 
to conclude the construction of 3D minimizers or 3D quasimodes. 
The resulting ``quasimode'' $\bu^\varepsilon$ has the same structure \eqref{eq:qmcyl} as in the cylindrical case. Moreover, 
numerical experiments (Model D, sect.\,\ref{modelD}) prove that, at least for some values of $r_\circ$, the lowest 3D eigenpairs fit the asymptotics $\lamone{\gL(\varepsilon)} \simeq \rH_0 + \varepsilon^{2/3}\am_1$ and $k(\varepsilon)\simeq\varepsilon^{-1/3} \gamma_{\min}$.

\section{Models}
\label{s:8}
We present in this section five models: cylinders (Model A: $f(z)=2$, $z\in(-1,1)$), cones (Model B: $f(z) = \frac32\,z-\frac12$, $z\in(-1,1)$), ``toroidal barrels'' (Model D: $f(z) = -1 + \sqrt{4 - z^2}$, $z\in(-1,1)$),  ``Gaussian barrels'' (Model H: $f(z) = 1 - \tfrac{z^2}{8} - \tfrac{z^4}{16}$, $z\in(-1,1)$), and ``Airy barrels'' (Model L: $f(z) = 1 - \tfrac{z^2}{8} - \tfrac{z^4}{16}$, $z\in(\frac12,\frac32)$), representing each of the five types that we could investigate from a theoretical point of view. We choose for all models
\[
   E=1 \quad\mbox{and}\quad \nu=0.3
\]
and perform 1D, 2D and 3D computations for each model. The 2D and 3D computations are performed with finite element codes ({\sc M\'elina}\footnote{M\'elina is an open source finite element library, see \href{https://anum-maths.univ-rennes1.fr/melina/}{https://anum-maths.univ-rennes1.fr/melina/}} for 2D and {\sc StressCheck}\footnote{Stress Check 9.0 is a trade mark of Engineering Software Research and Development, Inc., 
St. Louis, MO 63141, U.S.A.} for 3D) and for a finite set of values of $\varepsilon$ ranging from $0.2$ to $10^{-4}$ (in general this set contains the values $0.2$, $0.1$, and $5\cdot 10^{-j}$, $2\cdot 10^{-j}$, $10^{-j}$ for $j=2,3,4$). Let us mention that in our other paper \cite{ChDaFaYo16b} we present a more synthetic view of our theory together with a numerical study of two cases, a cylinder and an Airy barrel which coincide exactly with two models investigated in \cite{ArtioliBeiraoHakulaLovadina2009}. The agreement between our theory, our 2D-3D computations, and the computations presented in \cite{ArtioliBeiraoHakulaLovadina2009} is remarkable. Here we solve five different models to illustrate more completely the different cases pertaining to our approach.

The 1D calculations consist in computing the coefficients $\am_0$, $\am_1$ of \eqref{eq:mu1hK} and $\gamma$ of \eqref{eq:keps}. For Gaussian and Airy barrels we use our explicit formulas \eqref{6E15} and \eqref{6E15b}. For cones and toroidal barrels, we compute with {\sc Matlab} the spectrum of the one dimensional reduced operators \eqref{5E6} and \eqref{5E6bis} and optimize with respect to the parameter $\gamma$, whereas for cylinders, we calculate the eigenvalue $\mu^{\sf bilap}$ of the $\Delta^2$ operator by a semi-analytic method \cite[Ch.4]{beaudouin:tel}.

The 2D calculations solve the Lam\'e system $\gL^k$ at azimuthal frequency $k$ on meridian domains $\omega^\varepsilon$, see the corresponding variational formulations in Appendix \ref{app:B}. For each thickness parameter $\varepsilon$, any integer value of $k$ from $0$ to a certain cut-off frequency $k_{\max} = k_{\max}(\varepsilon)$ is used. The cut-off frequency $k_{\max}$ is determined so that we can observe a minimum for the first eigenvalue depending on $k$. This provides the numerical value $\tilde k(\varepsilon)$ for $k(\varepsilon)$. The domain $\omega^\varepsilon$ is meshed by curvilinear quadrilaterals of geometric degree 3. The meshes contain 2 elements in the thickness direction, and 8, 12 or 16 in the meridian direction. The polynomial interpolation degree of the FEM is 6 in each direction.

The 3D calculations solve the Lam\'e system $\gL$ on the 3D shells $\Omega^\varepsilon$. The azimuthal frequency $k(\varepsilon)$ is observed by counting the oscillations of the radial component \eqref{eq:rad} of the first eigenmode.

\begin{figure}[ht]
\includegraphics[scale=0.51]{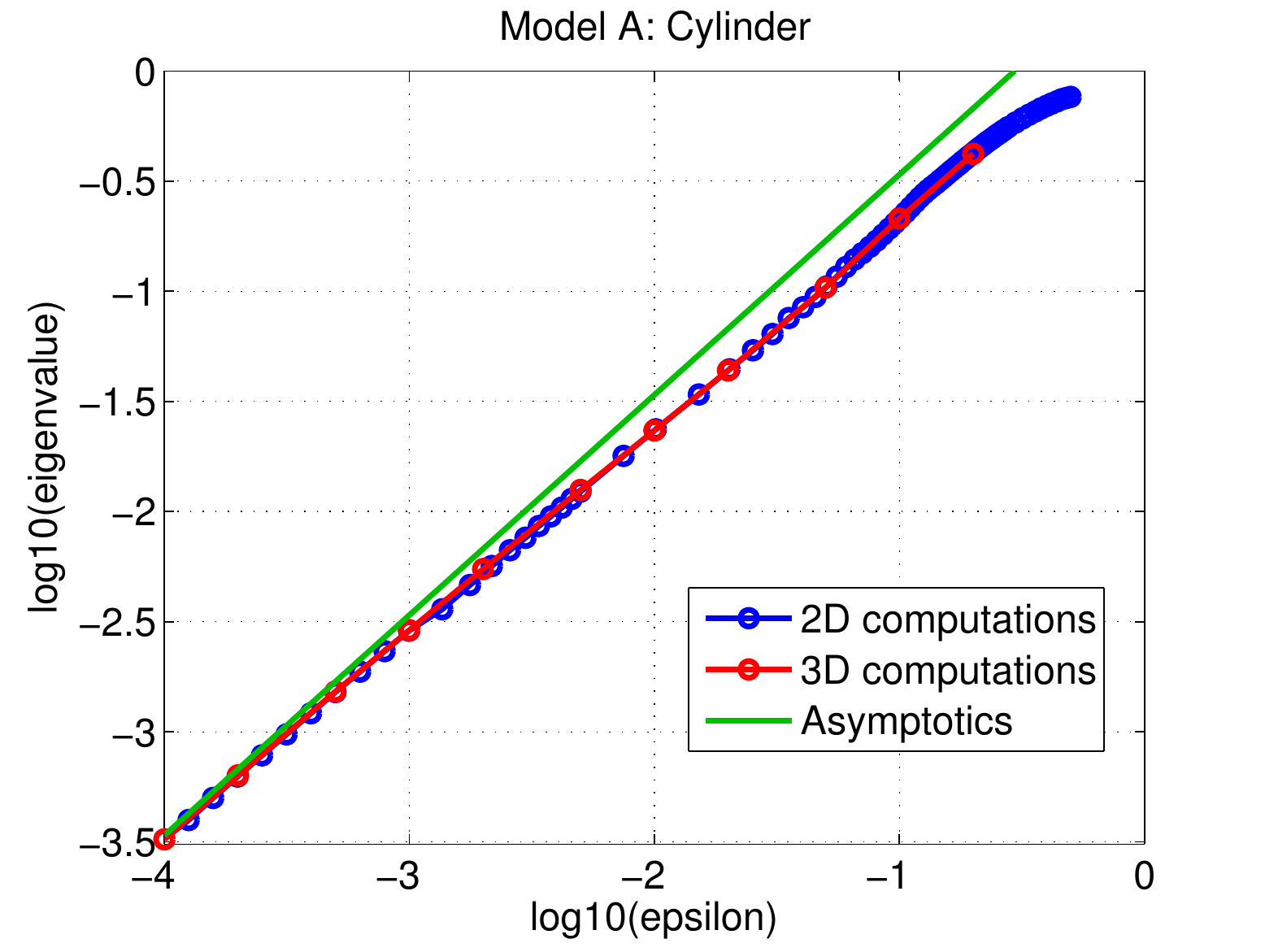}\hskip-1em
\includegraphics[scale=0.51]{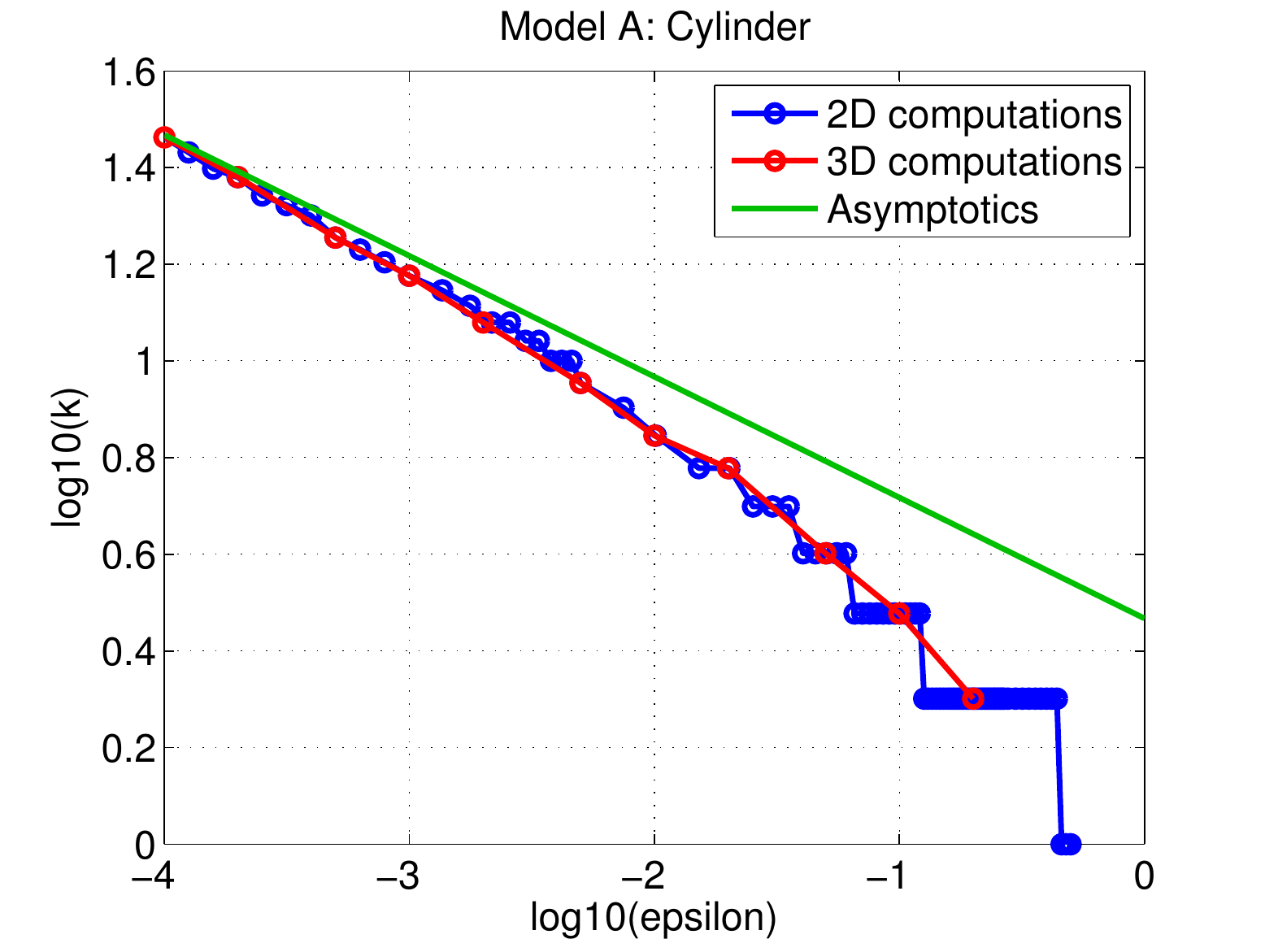}
   \caption{Model A: First eigenvalue $\tilde\lambda^{\varepsilon}_1$ and associated azimuthal frequency $k(\varepsilon)$.}
\label{f:A2}\end{figure}

\begin{figure}[ht]
\includegraphics[scale=0.51]{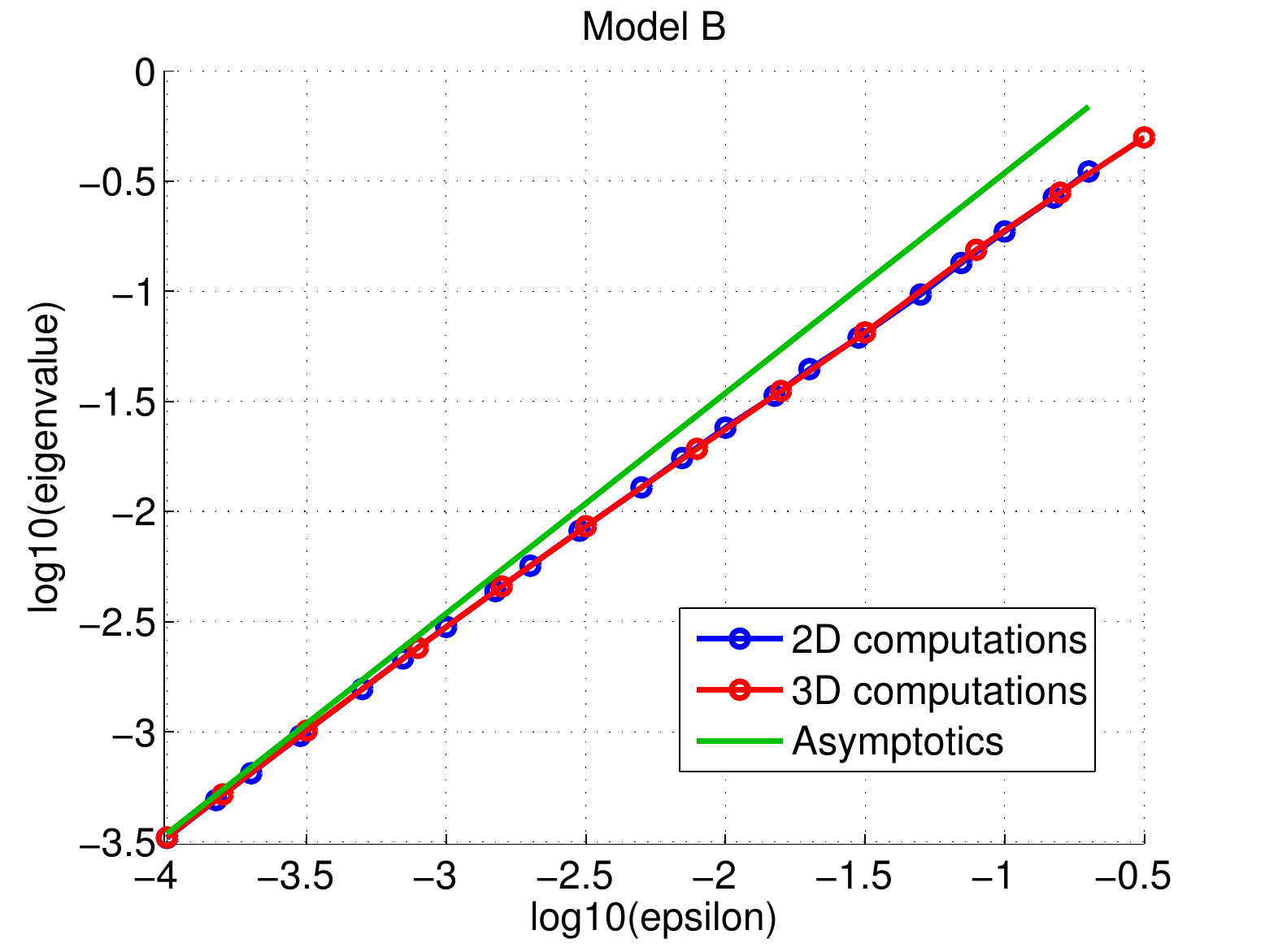}\hskip-1em
\includegraphics[scale=0.51]{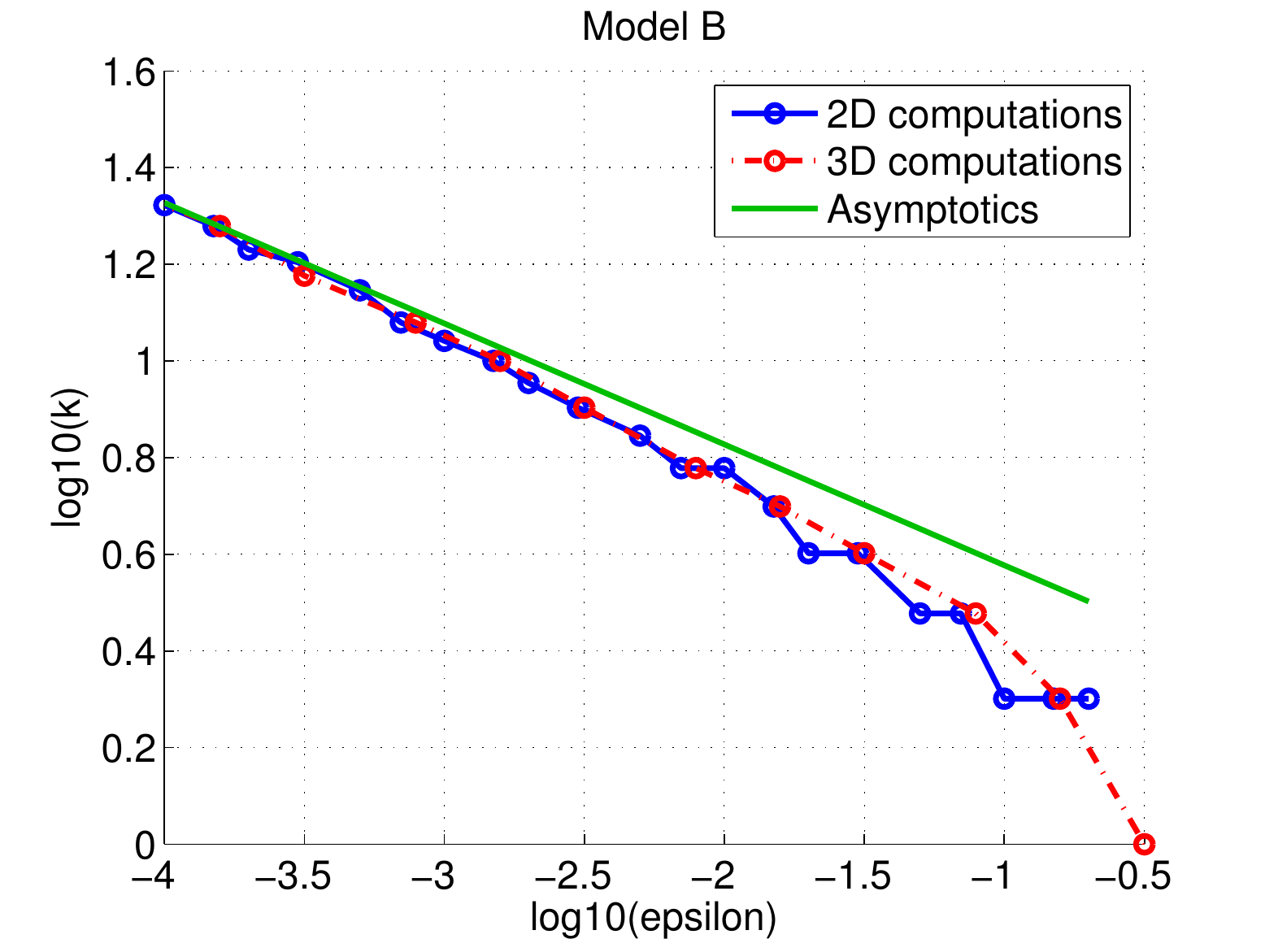}
   \caption{Model B: First eigenvalue $\tilde\lambda^{\varepsilon}_1$ and associated azimuthal frequency $k(\varepsilon)$.}
\label{f:B2}\end{figure}

\begin{figure}[ht]
\includegraphics[scale=0.51]{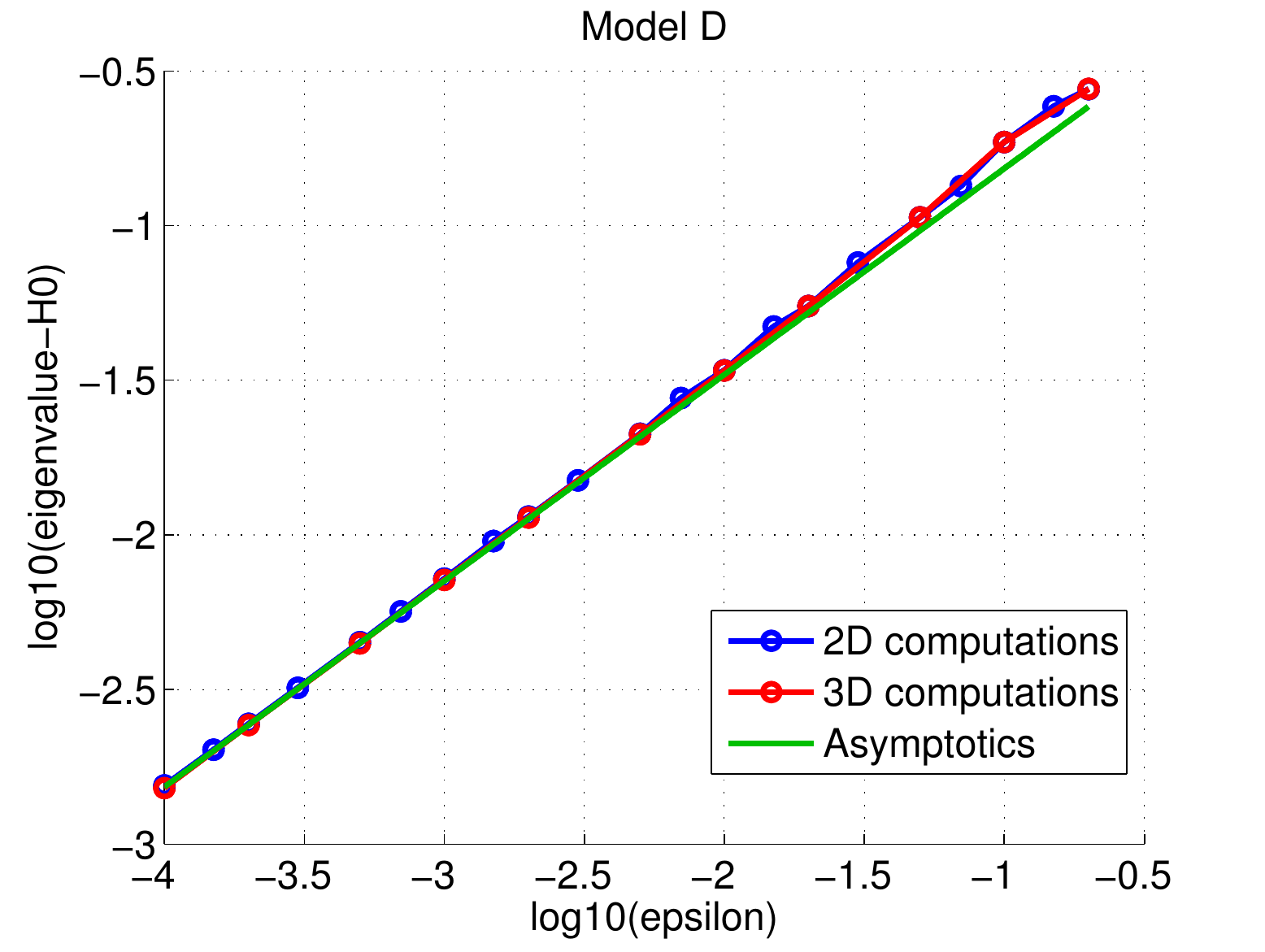}\hskip-1em
\includegraphics[scale=0.51]{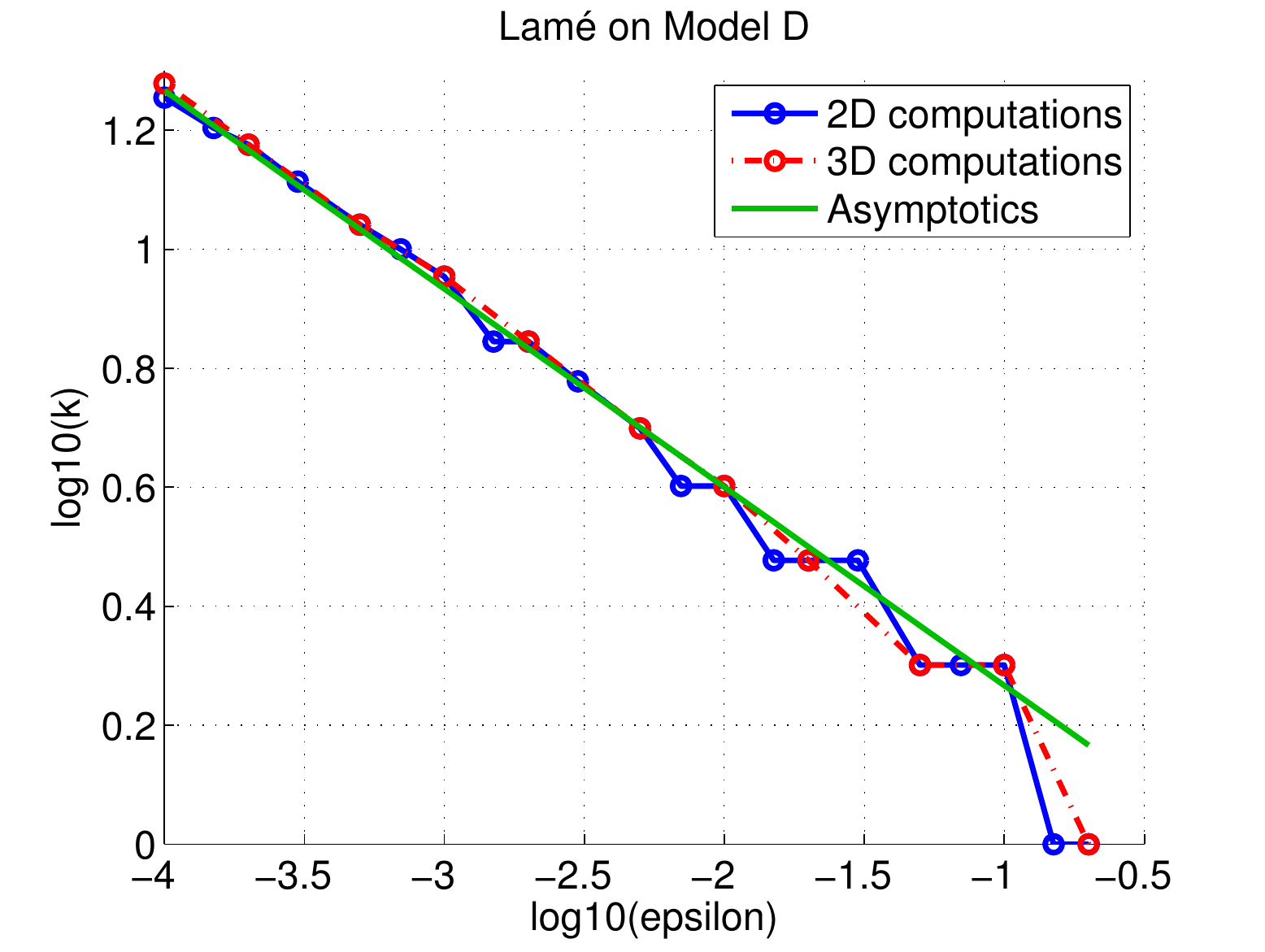}
   \caption{Model D: Difference $\tilde\lambda^{\varepsilon}_1-\Lambda_0$ and associated azimuthal frequency $k(\varepsilon)$.}
\label{f:D2}\end{figure}

\begin{figure}[ht]
\includegraphics[scale=0.51]{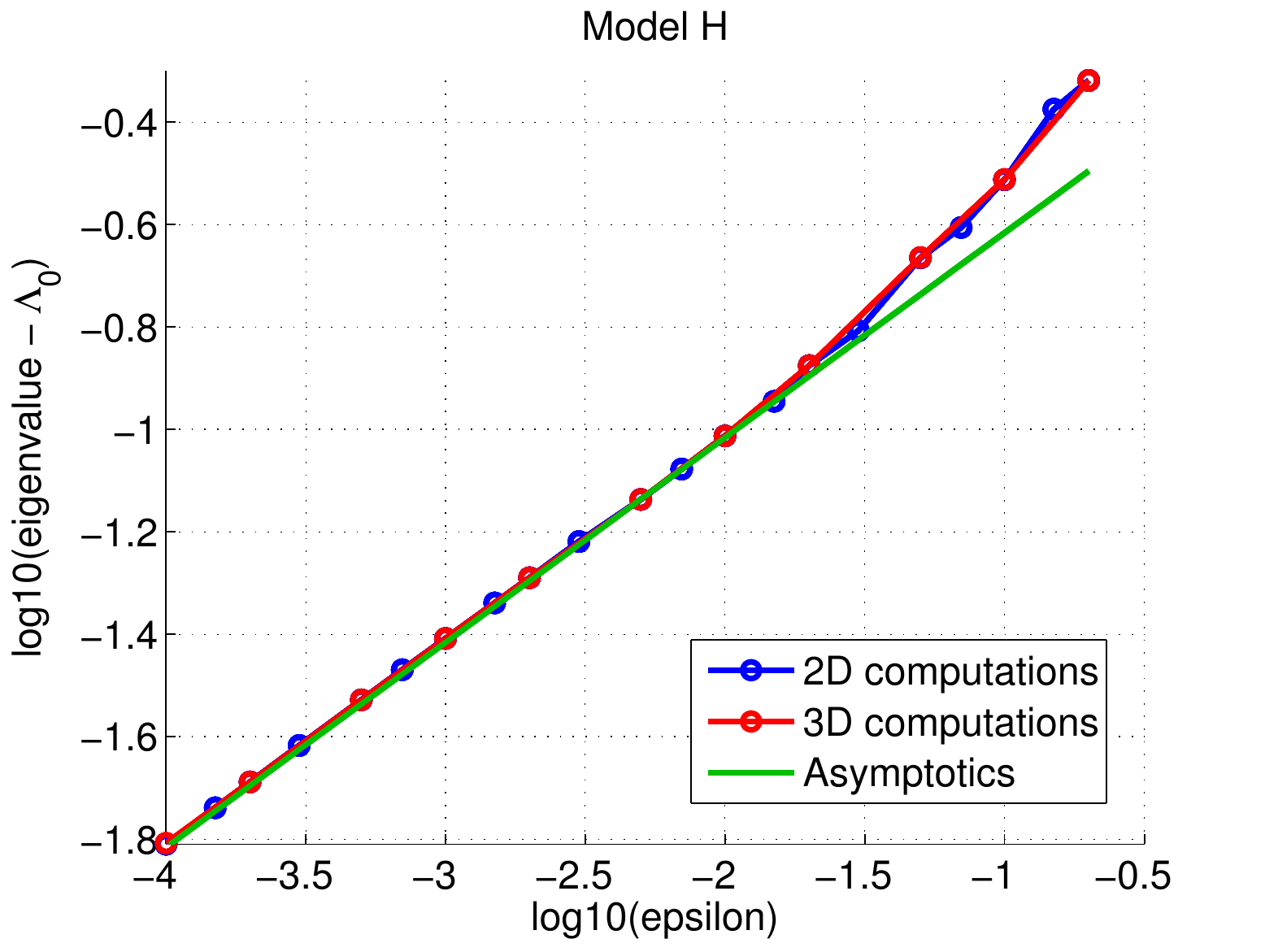}\hskip-1em
\includegraphics[scale=0.51]{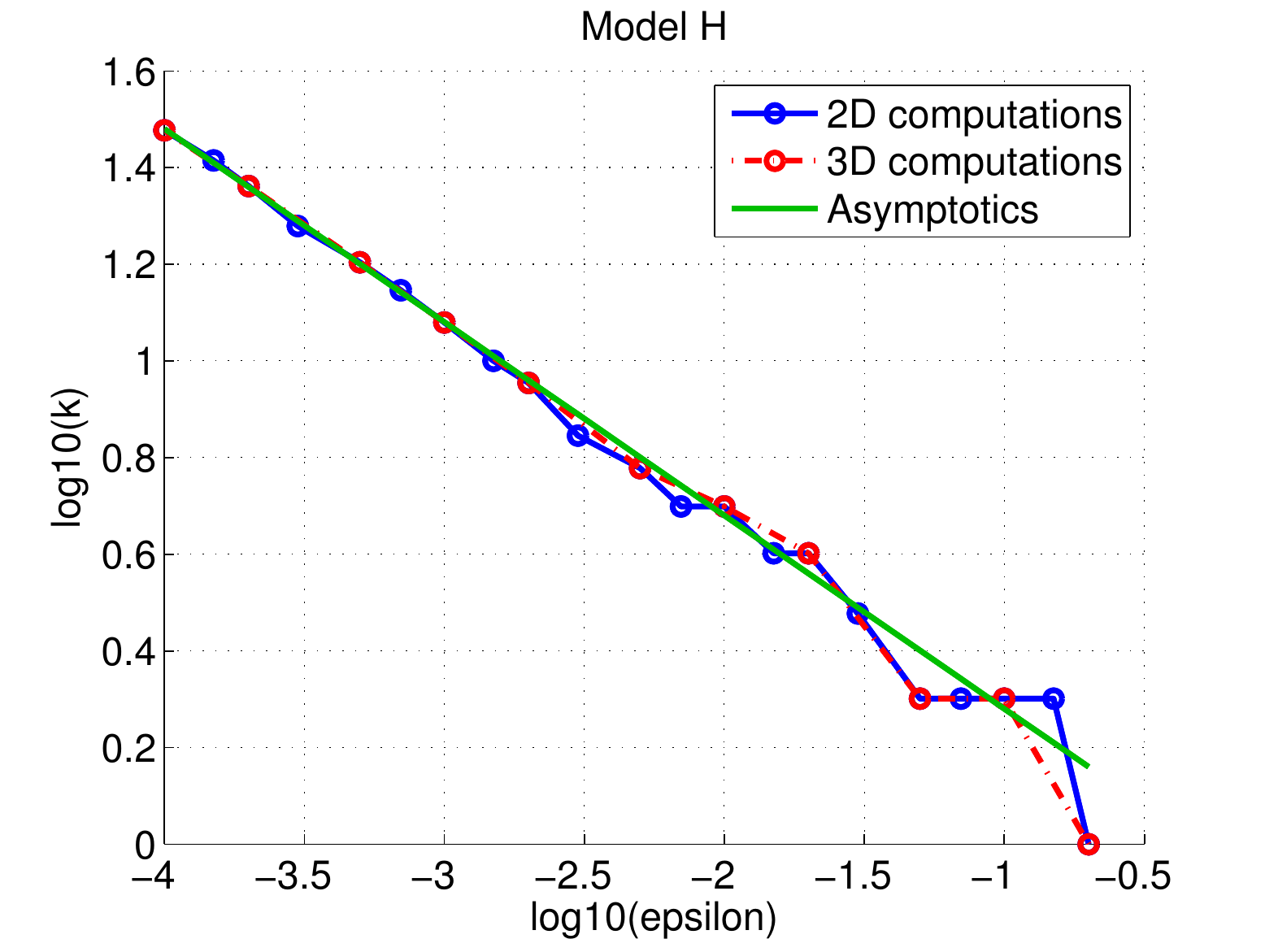}
   \caption{Model H: Difference $\tilde\lambda^{\varepsilon}_1-\Lambda_0$ and associated azimuthal frequency $k(\varepsilon)$.}
\label{f:H2}\end{figure}

\begin{figure}[ht]
\includegraphics[scale=0.35]{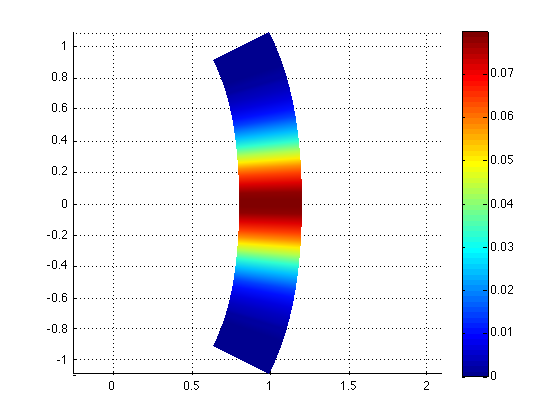}%
\includegraphics[scale=0.35]{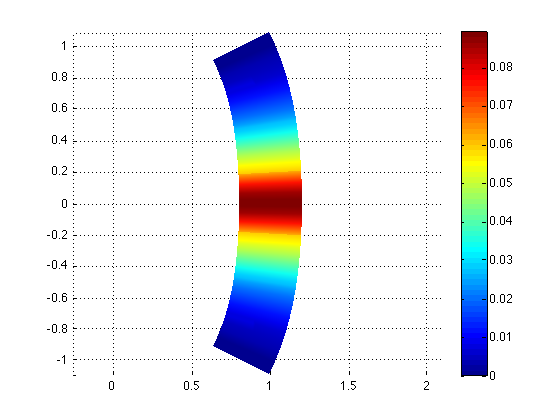}%
\includegraphics[scale=0.35]{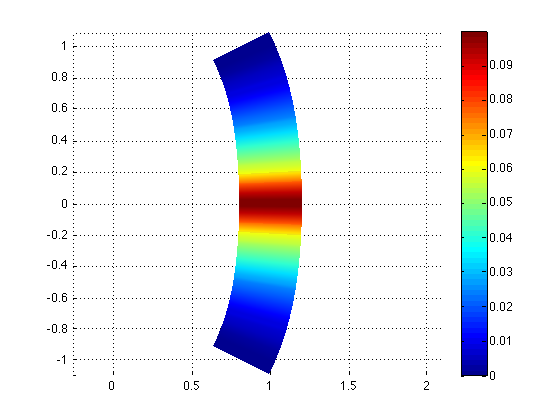}
   \caption{Model H: 2D first eigenmode (radial component) for $\varepsilon=10^{-3}$ and $k=12$, $\varepsilon=3\cdot10^{-4}$ and $k=19$, $\varepsilon=10^{-4}$ and $k=30$. Represented on $\omega^{0.2}$.}
\label{f:H4}\end{figure}

\begin{figure}[ht]
\includegraphics[scale=0.51]{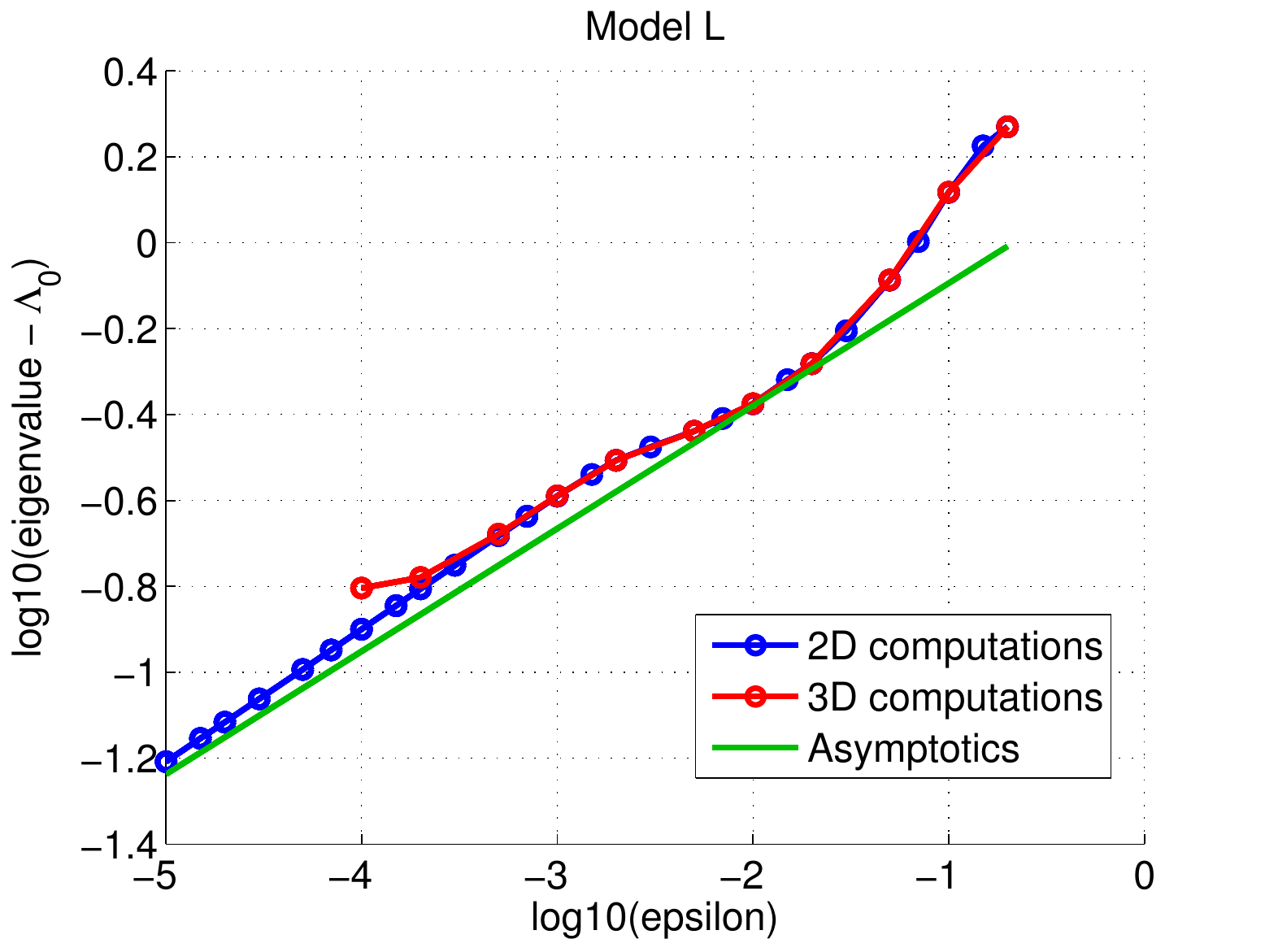}\hskip-1em
\includegraphics[scale=0.51]{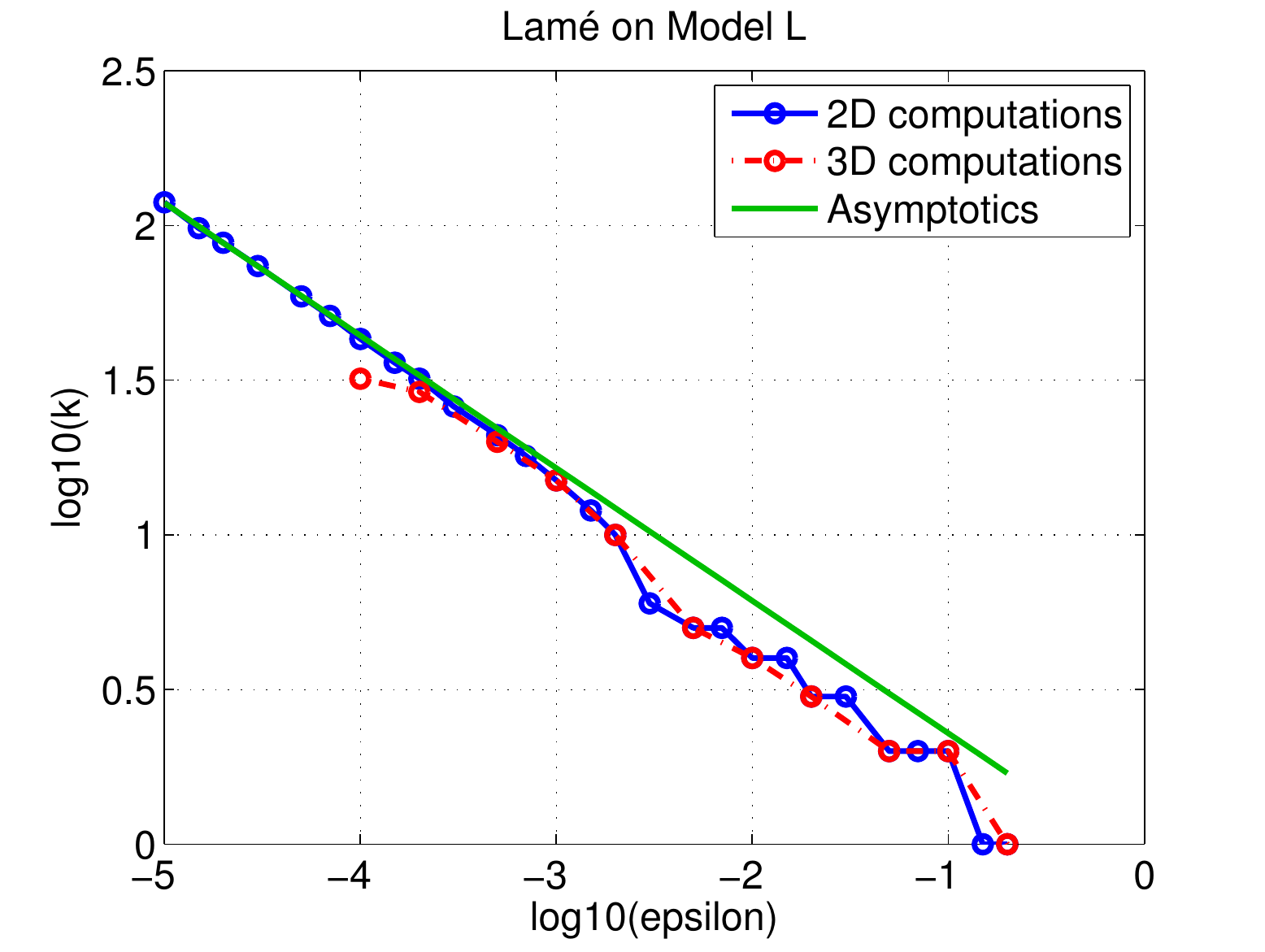}
   \caption{Model L: Difference $\tilde\lambda^{\varepsilon}_1-\Lambda_0$ and associated azimuthal frequency $k(\varepsilon)$.}
\label{f:L2}\end{figure}

\begin{figure}[ht]
\includegraphics[scale=0.35]{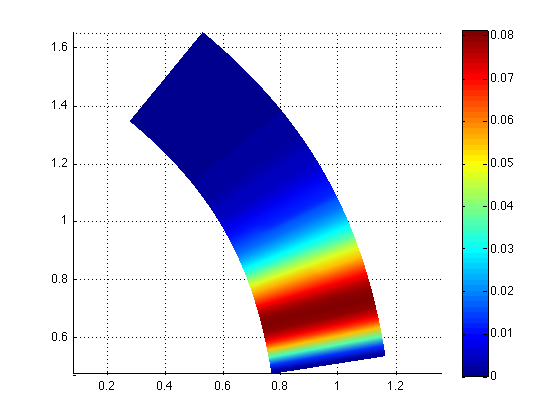}%
\includegraphics[scale=0.35]{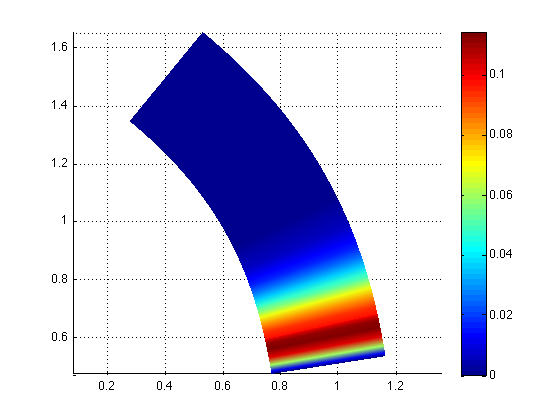}%
\includegraphics[scale=0.35]{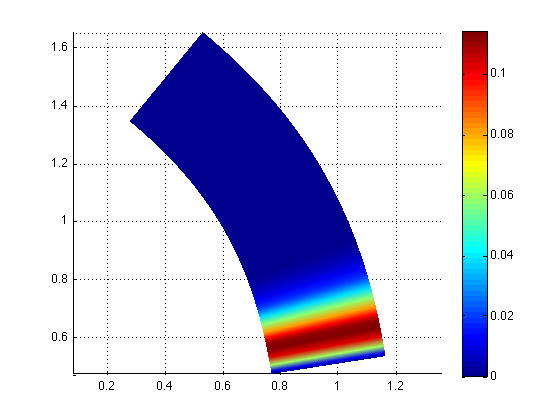}
   \caption{Model L: 2D first eigenmode (radial component) for $\varepsilon=10^{-3}$ and $k=15$, $\varepsilon=3\cdot10^{-4}$ and $k=26$, $\varepsilon=10^{-4}$ and $k=43$. Represented on $\omega^{0.2}$.}
\label{f:L4}\end{figure}

We represent in figures \ref{f:A1}, \ref{f:B1} and \ref{f:D1}, the  meridian domains $\omega^\varepsilon$ in the $(r,z)$ plane for models A, B and D, respectively. The curve $\cC$ is dotted. Figures \ref{f:H4} and \ref{f:L4} provide $\omega^\varepsilon$ for $\varepsilon=0.2$ for models H and L. Figures \ref{f:A2}, \ref{f:B2}, \ref{f:D2}, \ref{f:H2}, and \ref{f:L2} show the lowest computed eigenvalue $\tilde\lambda^{\varepsilon}_1$ and the associated azimuthal frequency $k(\varepsilon)$ versus $\varepsilon$ in loglog scale (in base 10). For elliptic models D, H and L, the difference $\tilde\lambda^{\varepsilon}_1-\Lambda_0$ is plotted. The 1D asymptotics is the line $\varepsilon\mapsto\mm_1(\varepsilon)$ \eqref{eq:quasi} (or $\varepsilon\mapsto\mm_1(\varepsilon)-\Lambda_0$ in elliptic models). 

\begin{figure}[ht]
\includegraphics[scale=0.225]{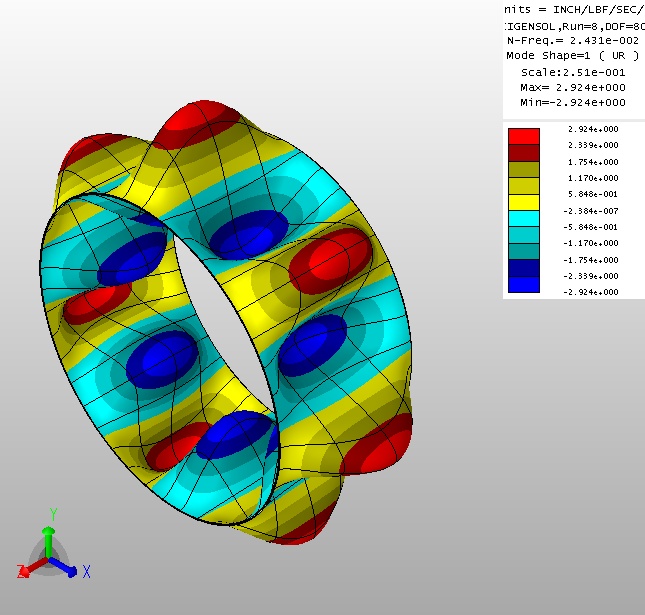}%
\includegraphics[scale=0.225]{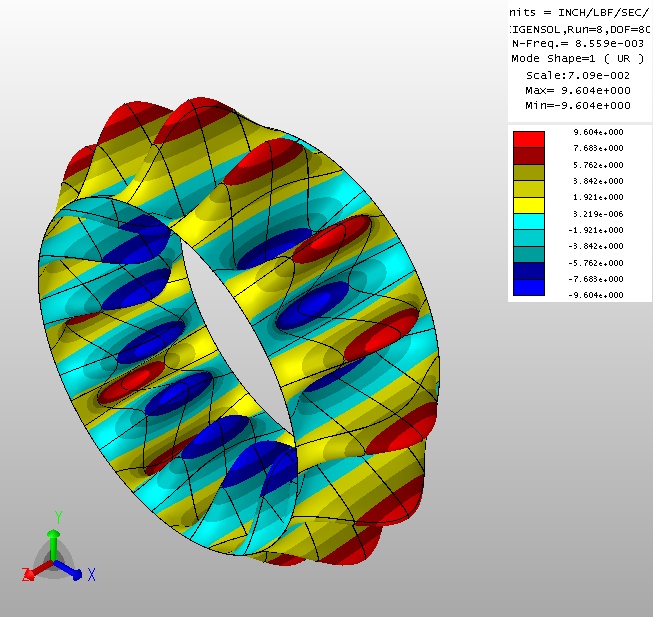}%
\includegraphics[scale=0.225]{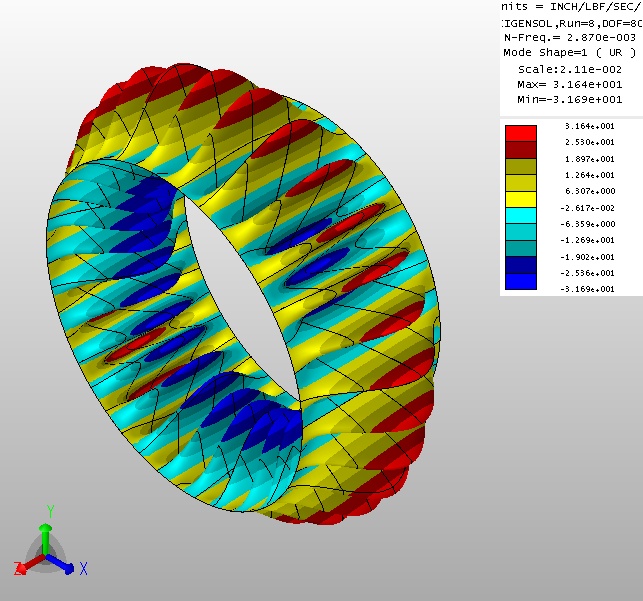}
   \caption{Model A: First eigenmode (radial component) for $\varepsilon=10^{-2}$, $10^{-3}$, $10^{-4}$.}
\label{f:A3}\end{figure}

\begin{figure}[ht]
\includegraphics[scale=0.225]{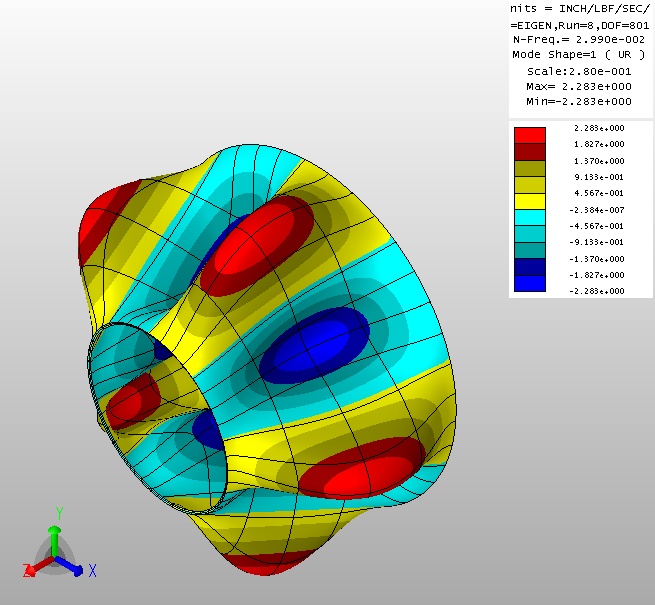}%
\includegraphics[scale=0.225]{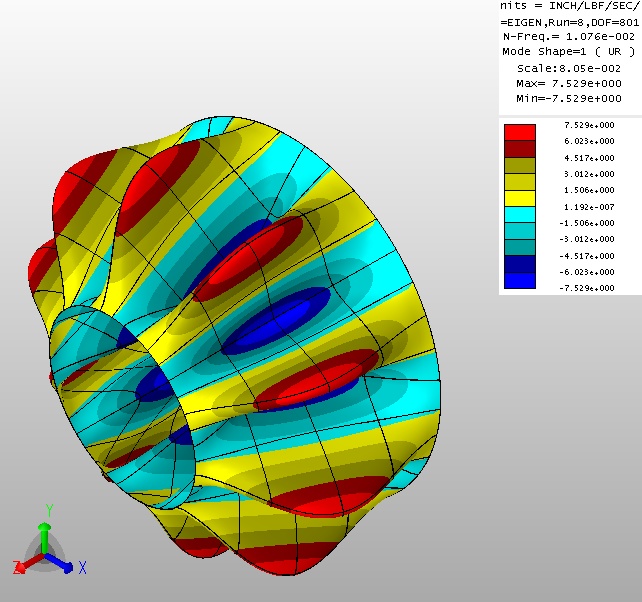}%
\includegraphics[scale=0.225]{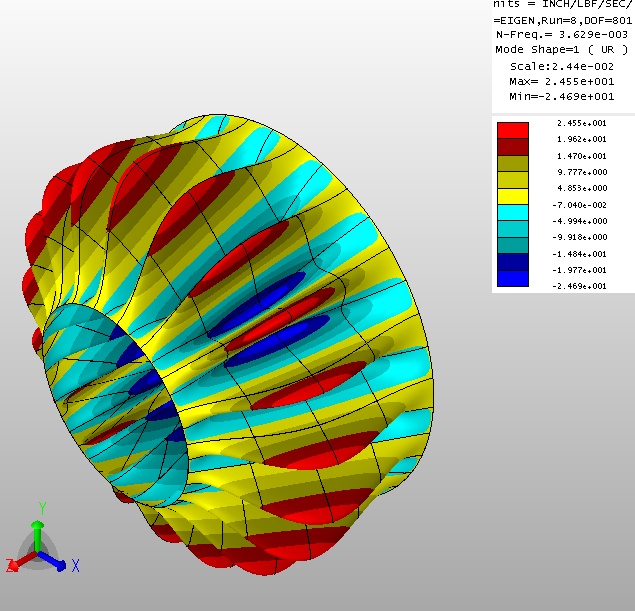}
   \caption{Model B: First eigenmode (radial component) for $\varepsilon=10^{-2}$, $10^{-3}$, $10^{-4}$.}
\label{f:B3}\end{figure}

\begin{figure}[ht]
\includegraphics[scale=0.225]{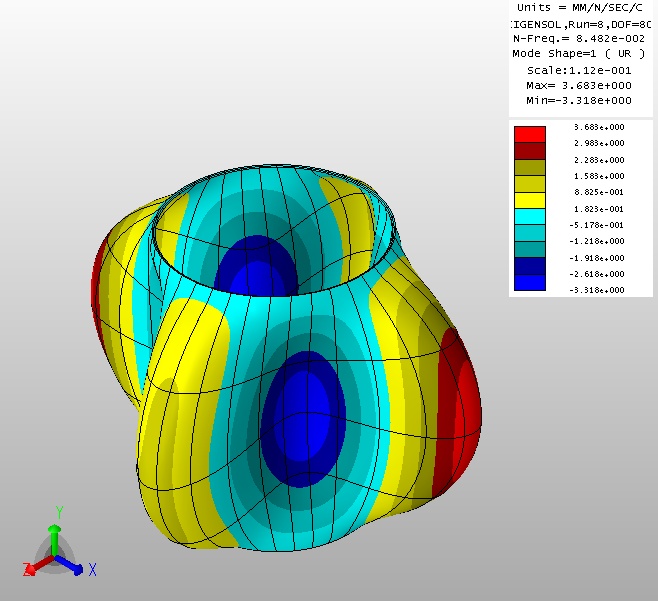}%
\includegraphics[scale=0.225]{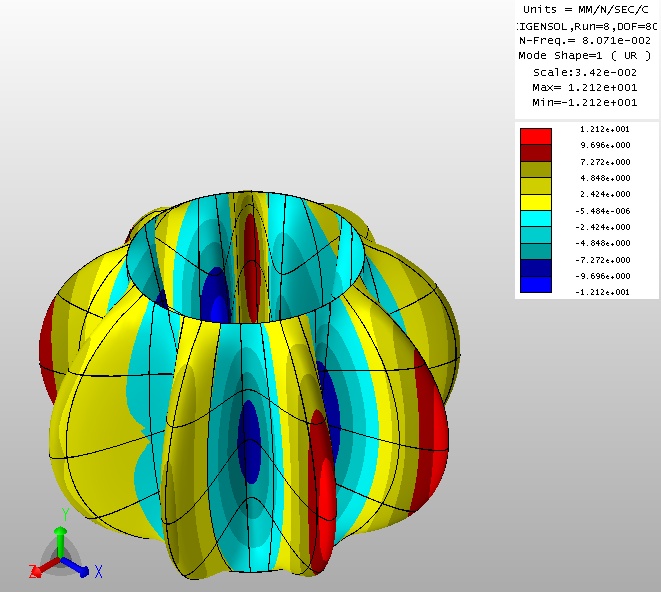}%
\includegraphics[scale=0.225]{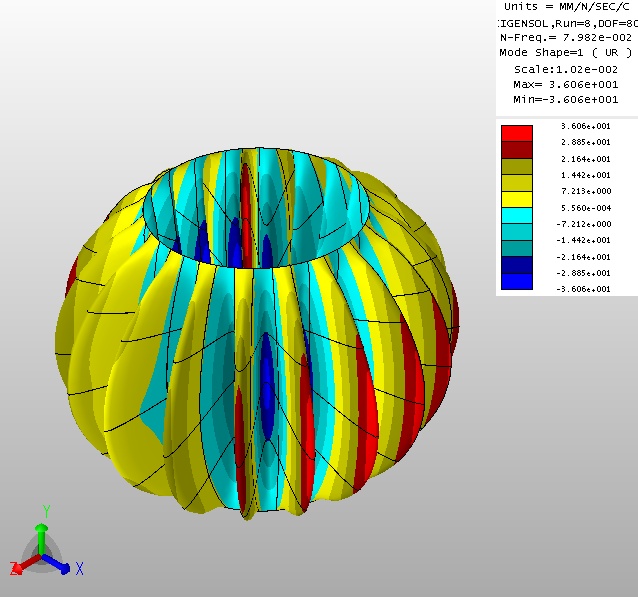}
   \caption{Model D: First eigenmode (radial component) for $\varepsilon=10^{-2}$, $10^{-3}$, $10^{-4}$.}
\label{f:D3}\end{figure}

\begin{figure}[ht]
\includegraphics[scale=0.225]{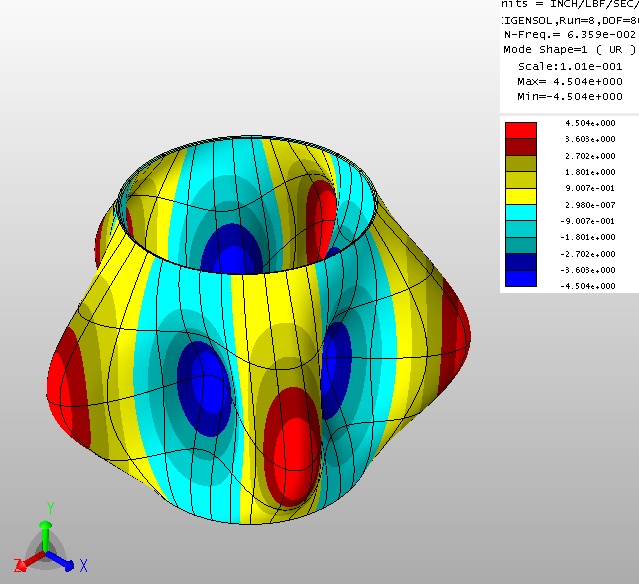}%
\includegraphics[scale=0.225]{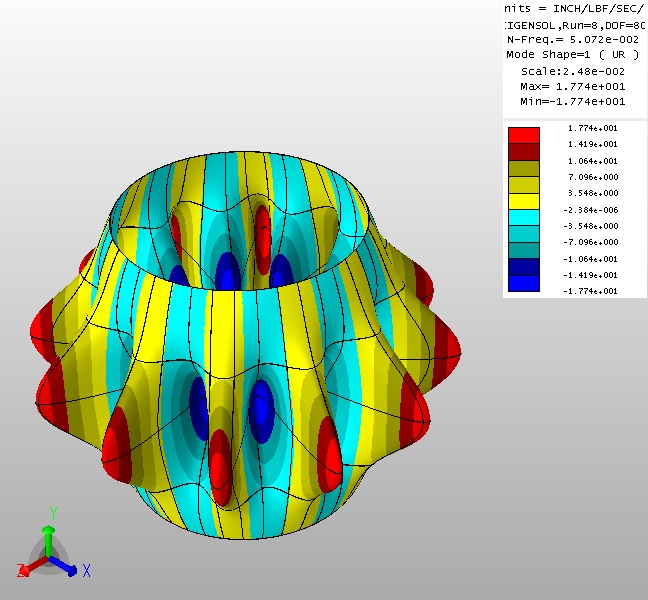}%
\includegraphics[scale=0.225]{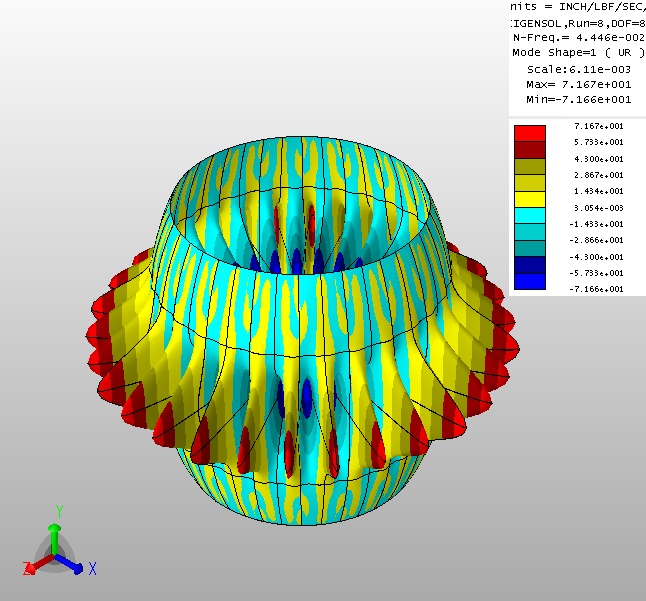}
   \caption{Model H: First eigenmode (radial component) for $\varepsilon=10^{-2}$, $10^{-3}$, $10^{-4}$.}
\label{f:H3}\end{figure}

\begin{figure}[ht]
\includegraphics[scale=0.225]{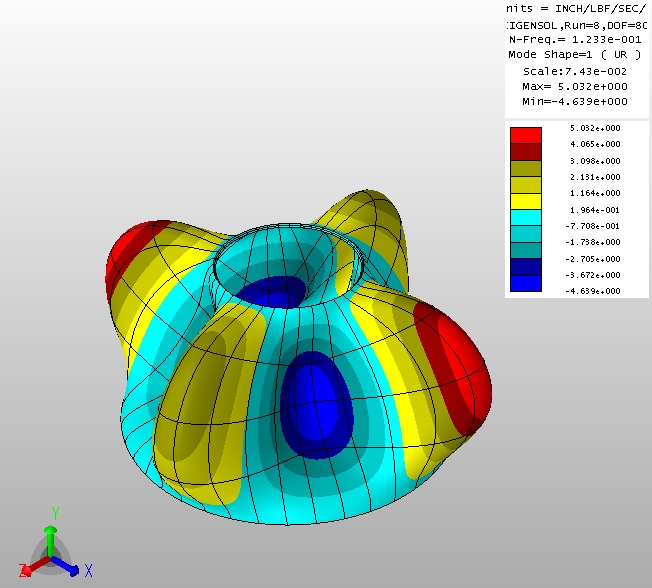}%
\includegraphics[scale=0.225]{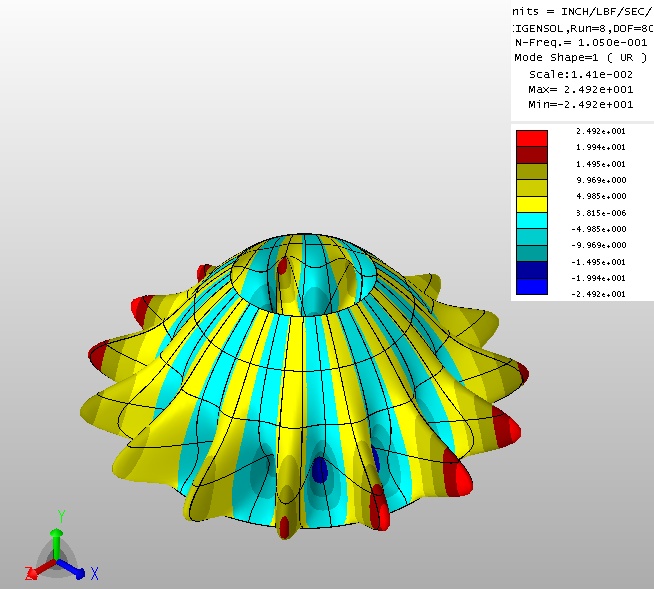}%
\includegraphics[scale=0.225]{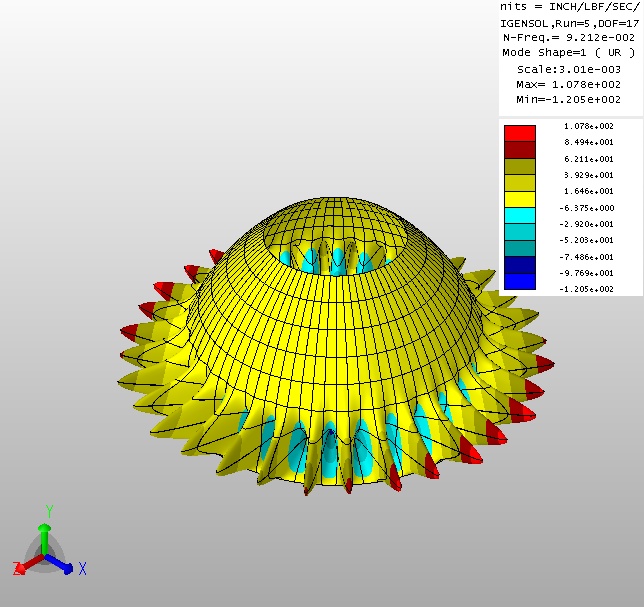}
   \caption{Model L: 3D first eigenmode (radial component) for $\varepsilon=10^{-2}$, $10^{-3}$, $10^{-4}$.}
\label{f:L3}\end{figure}

Figures \ref{f:A3}, \ref{f:B3}, \ref{f:D3}, \ref{f:H3}, and \ref{f:L3} show the radial component of the first 3D computed eigenvector $\tilde\bu^\varepsilon$ for three values of $\varepsilon$ and the five models, respectively. We note that this radial component has the same behavior as predicted for our quasimode $\bu^\varepsilon$, cf \eqref{eq:qmcyl} for models A, B, D, and \eqref{eq:qmgauss}, \eqref{eq:qmairy} for models H, L. Figures \ref{f:H4} and \ref{f:L4} are surface plots on the meridian domain $\omega^\varepsilon$ of the first 2D eigenvectors of the operator $\gL^{\near{k(\varepsilon)}}$ in Gaussian and Airy barrels, respectively, They clearly exhibit the meridian concentration of the modes as $\varepsilon$ decreases, cf the behavior in $z$-variable in \eqref{eq:qmgauss}, \eqref{eq:qmairy}. For visibility, they are scaled with respect to the width in order to be represented on the meridian domain with thickness parameter $\varepsilon=0.2$. The observable concentration scale is compatible with the theoretical scale induced from \eqref{eq:qmgauss}, \eqref{eq:qmairy}.
We summarize in Table \ref{tab:asy} the numerical values of the asymptotic quantities $\mm_1(\varepsilon)$ and $k(\varepsilon)$, as well as the observed asymptotics for the remainder $\tilde\lambda^\varepsilon_1-\mm_1(\varepsilon)$ for each of the five models. In Table \ref{tab:k} we list computed and theoretical values of $k(\varepsilon)$ for the four models B, D, H, and L (Model A is ommited because of its great similarity with Model B).

\begin{table}
{\small
\begin{tabular}{l@{\quad}ll@{\quad}l@{\quad}llllll}
\hline
Model       & $k(\varepsilon)$ & $\gamma$
            & $\mm_1(\varepsilon)$ & $\am_0$ & $\am_1$  & Remainder \\
\hline
\underline{\sc Parabolic} \\
A / Cylinder& $\gamma\varepsilon^{-1/4}$        & 2.9323
            & $\am_1\varepsilon$                & 0            & 3.3852
            & 15\,$\varepsilon^{3/2}$  (asymptotics)  \\
B / Cone    & $\gamma\varepsilon^{-1/4}$        & 2.1247
            & $\am_1\varepsilon$                & 0            & 3.4464
            & 15\,$\varepsilon^{3/2}$  (asymptotics)  \\
\underline{\sc Elliptic} \\
D / Toroidal& $\gamma\varepsilon^{-1/3}$        & 0.85935
            & $\am_0 + \am_1\varepsilon^{2/3}$  & 0.25000    & 0.71500
            &  0.2\,$\varepsilon$  (upper bound)\\
H / Gauss   & $\gamma\varepsilon^{-2/5}$        & 0.75901
            & $\am_0 + \am_1\varepsilon^{2/5}$  & 0.06250    & 0.60785
            & 0.05\,$\varepsilon^{3/5}$  (asymptotics) \\
L / Airy    & $\gamma\varepsilon^{-3/7}$        & 0.85141
            & $\am_0 + \am_1\varepsilon^{2/7}$  & 0.17804    & 1.55472
            & 2.9\,$\varepsilon^{4/7}$  (asymptotics) \\
\hline\\
\end{tabular}
}
\caption{Numerical values for asymptotic quantities $k(\varepsilon)$ and $\mm_1(\varepsilon)$. Observed asymptotics for the remainder $\tilde\lambda^\varepsilon_1-\mm_1(\varepsilon)$  (with $\tilde\lambda^\varepsilon_1$ obtained by 2D computations on the meridian domain $\omega^\varepsilon$).}
\label{tab:asy}
\end{table}

The formulas providing parameters $\gamma$, $\am_0$ and $\am_1$ are given in the following equations: \eqref{5E6b}--\eqref{5E10} for models A and B, \eqref{eq:circ}, \eqref{5E6bbis} for model D. Concerning models H and L, the function $\rH_0/E = \frac{f''^2}{s^6}$ is equal to $(\frac14+\frac34z^2)^2/(1+(\frac14z+\frac14z^3)^2)^{3}$ and reaches its minimum in the interior point $z_0=0$ for model H, and in the boundary point $z_0=0.5$ for model L. Formulas for $\gamma$ and $\am_1$ are given in \eqref{6E15} with \eqref{6E14bc}, and \eqref{6E15b} with \eqref{6EAiry4} for models H and L, respectively.

\begin{table}
\begin{tabular}{r@{\qquad}rr@{\qquad}rr@{\qquad}rr@{\qquad}rr}
\hline
               &  \multicolumn{2}{l}{Model B}&  \multicolumn{2}{l}{Model D}
               &  \multicolumn{2}{l}{Model H}&  \multicolumn{2}{l}{Model L} \\  
$\varepsilon$  & $\tilde k$ & $k$ & $\tilde k$ & $k$ & $\tilde k$ & $k$ & $\tilde k$ & $k$ \\               
\hline
0.20000  &   2  &   3.2  &   1  &   1.5  &   1  &   1.4  &   1  &   1.7\\ 
0.10000  &   2  &   3.8  &   2  &   1.8  &   2  &   1.9  &   2  &   2.3\\ 
0.05000  &   3  &   4.5  &   2  &   2.3  &   2  &   2.5  &   2  &   3.1\\ 
0.02000  &   4  &   5.6  &   3  &   3.2  &   4  &   3.6  &   3  &   4.6\\ 
0.01000  &   6  &   6.7  &   4  &   4.0  &   5  &   4.8  &   4  &   6.1\\ 
0.00500  &   7  &   8.0  &   5  &   5.0  &   6  &   6.3  &   5  &   8.2\\ 
0.00200  &   9  &  10.0  &   7  &   6.8  &   9  &   9.1  &  10  &  12.2\\ 
0.00100  &  11  &  11.9  &   9  &   8.6  &  12  &  12.0  &  15  &  16.4\\ 
0.00050  &  14  &  14.2  &  11  &  10.8  &  16  &  15.9  &  21  &  22.1\\ 
0.00020  &  17  &  17.9  &  15  &  14.7  &  23  &  22.9  &  32  &  32.7\\ 
0.00010  &  21  &  21.2  &  18  &  18.5  &  30  &  30.2  &  43  &  44.1\\ 
0.00005  &  25  &  25.3  &  24  &  23.3  &  40  &  39.9  &  59  &  59.3\\ 
\hline\\
\end{tabular}
\caption{For a collection of values of $\varepsilon$, observed azimuthal frequency $\tilde k=\tilde k(\varepsilon)$ versus theoretical value $k=k(\varepsilon)$ provided by our asymptotic formulas  for models B (cone), D (toroidal barrel), H (Gaussian barrel), and L (Airy barrel)}
\label{tab:k}
\end{table}

\newcommand\Cone[3]{
\framebox{\begin{minipage}{0.2\textwidth}
\def\thickness{#1}
\def\RR{#2}
\def\alph{#3}
\FPeval{\Rbot}{\RR-\alph}
\FPeval{\Rtop}{\RR+\alph}
     \figpt 1:(0, -2)
     \figpt 2:(0, 2)
     \figpt 3:(\Rbot, -1)
     \figpt 4:(\Rtop, 1)
     \figvectN 9 [3,4]
     \figvectU 10 [9]
     \figptstra 13 = 3,4 / \thickness, 10/
     \figptstra 23 = 3,4 / -\thickness, 10/

     \figdrawbegin{}
     \figdrawline [13,14,24,23,13]
     \figset (color={1 0 0}, width=1.5)
     \figdrawline [23,13]
     \figdrawline [24,14]
     \figset (color={0 0 1}, width=default, dash=2)
     \figdrawarrow [1,2]
     \figset (color={0 0.7 0}, dash=4)
     \figdrawline [3,4]

     \figdrawend
     \figvisu{\figbox}{}{%
     \figwrites 1 :{\em axis}(0.1)
     \figwrites 3 :{$\varepsilon=\thickness$}(0.5)
     }
   {\box\figbox}
\end{minipage}
}
}

\clearpage
\subsection{Model A: Cylindrical shells}
\label{modelA}
The midsurface parametrization, cf \eqref{5E1}, is given by
\[
   f(z) = R, \quad z\in(-1,1), \quad R=2.
\]

\begin{figure}[ht]
   \figinit{1cm}
   \Cone{0.2}{2}{0}
   \hskip 5mm
   \Cone{0.1}{2}{0}
   \hskip 5mm
   \Cone{0.05}{2}{0}
   \caption{Model A: Meridian domains for several values of $\varepsilon$.}
\label{f:A1}\end{figure}

\clearpage

\subsection{Model B: Conical shells}
\label{modelB}
The midsurface parametrization, cf \eqref{5E1}, is given by
\[
   f(z) = T z + R, \quad z\in(-1,1),\quad T=-0.5,\quad R=1.5.
\]

\begin{figure}[ht]
   \figinit{1cm}
   \Cone{0.2}{1.5}{-0.5}
   \hskip 5mm
   \Cone{0.1}{1.5}{-0.5}
   \hskip 5mm
   \Cone{0.05}{1.5}{-0.5}
   \caption{Model B: Meridian domains for several values of $\varepsilon$.}
\label{f:B1}\end{figure}

\newcommand\Torusplus[3]{
\framebox{\begin{minipage}{0.2\textwidth}
\def\thickness{#1}
\def\RR{#2}
\def\rz{#3}
\FPeval{\Rint}{\RR-\thickness}
\FPeval{\Rext}{\RR+\thickness}
\FPeval{\Centre}{\rz-\RR}
\FPeval{\unsurR}{1/\RR}
\FParcsin{\angleplus}{\unsurR}
\FPeval{\angleplus}{\angleplus/\FPpi}
\FPeval{\angleplus}{\angleplus*180}

     \figpt 1:(0, -2)
     \figpt 2:(0, 2)
     \figpt 5:(\Centre, 0)
     \figptcirc 3 :: 5;\RR (-\angleplus)
     \figptcirc 4 :: 5;\RR (\angleplus)
     \figptcirc 13 :: 5;\Rint (-\angleplus)
     \figptcirc 14 :: 5;\Rint (\angleplus)
     \figptcirc 23 :: 5;\Rext (-\angleplus)
     \figptcirc 24 :: 5;\Rext (\angleplus)
     \figdrawbegin{}
     \figdrawarccirc 5 ; \Rint (-\angleplus,\angleplus)
     \figdrawarccirc 5 ; \Rext (-\angleplus,\angleplus)
     \figset (color={1 0 0}, width=1.5)
     \figdrawline [23,13]
     \figdrawline [24,14]
     \figset (color={0 0 1}, width=default, dash=2)
     \figdrawarrow [1,2]
     \figset (color={0 0.7 0}, dash=4)
     \figdrawarccirc 5 ; \RR (-\angleplus,\angleplus)

     \figdrawend
     \figvisu{\figbox}{}{%
     \figwrites 1 :{\em axis}(0.1)
     \figwrites 3 :{$\quad\varepsilon=\thickness$}(0.5)
     \figset write(mark={$\Bl\bullet$})
     \figwritep[5]
     }
   {\box\figbox}
\end{minipage}
}
}

\clearpage

\subsection{Model D: Toroidal barrels}
\label{modelD}
The midsurface parametrization is, cf \eqref{eq:circ}
\[
   f(z) = r_\circ + \sqrt{R^2 - z^2}, \quad z\in(-1,1), \quad R=2\quad \mbox{and}\quad r_\circ=-1
\]

\begin{figure}[ht]
   \figinit{1cm}
   \Torusplus{0.2}{2}{1}
   \hskip 5mm
   \Torusplus{0.1}{2}{1}
   \hskip 5mm
   \Torusplus{0.05}{2}{1}
   \caption{Model D: Meridian domains for several values of $\varepsilon$.}
\label{f:D1}\end{figure}


\clearpage
\subsection{Model H: Gaussian barrel}
\label{modelH}
The midsurface parametrization is
\[
   f(z) = 1 - \tfrac{z^2}{8} - \tfrac{z^4}{16}, \quad z\in(-1,1).
\]

\clearpage
\subsection{Model L: Airy barrel}
\label{modelL}
The midsurface parametrization is
\[
   f(z) = 1 - \tfrac{z^2}{8} - \tfrac{z^4}{16}, \quad z\in(0.5,1.5).
\]

\section{Conclusion}
\label{s:9}
For five categories of clamped axisymmetric shells, we have exhibited a scalar 1D operator that determines the
asymptotic expansion of the azimuthal frequency $k(\varepsilon)$ of the first vibration mode, and a two-term
asymptotic expansion for $\mm_1(\varepsilon)=\am_0+\am_1\varepsilon^{\alpha_1}$ for the first eigenvalue of the full 3D Lam\'e system in the shell. These five categories are the cylinders and the trimmed cones (parabolic shells), as well as what we denote toroidal, Gaussian and Airy barrels (elliptic shells).
The most striking outcome of our analysis is the extremely good agreement of the three computation methods (3D, 2D and 1D) in all the five cases described above, strengthening the relevance of our constructions.
The presented methods demonstrate that the smallest eigenpairs for the Lam\'e system can be estimated for
specific shells by the reduced 1D model. Furthermore, the spreading or the concentration of the first eigenmode can be predicted accurately: The cases for which concentration occurs are the Gauss and Airy barrels and for those shells, the first eigenmode concentrates around a ring whose location $(f(z_0),z_0)$ is analytically known.

Another interesting observation is the comparison with the computations in \cite{ArtioliBeiraoHakulaLovadina2009}. The elliptic case that is considered there is $f(z)=1-\frac12 z^2$ on the interval $\cI=(-a,a)$ with $a=0.892668$. For this example we find
\[
   \sh(z) = \tfrac32 z^2 (1-\tfrac12 z^2) \quad\mbox{and}\quad
   \rH_0(z) = (1+z^2)^{-3}.
\]
So we see that we are in our admissible `Airy' case, predicting a behavior in $\varepsilon^{-3/7}$ for the azimuthal frequency $k(\varepsilon)$ of the first eigenvector. Noting that $\frac37\simeq 0.43$ and $\frac25=0.4$, we believe that this explains what have observed the authors \cite[p.55]{ArtioliBeiraoHakulaLovadina2009}:
{\em ``We also notice that in the elliptic case, $K_t$ {\em[$k(\varepsilon)$]} is probably growing slightly faster than exactly $t^{-2/5}$ {\em[$\varepsilon^{-2/5}$]}''}. The treatment in \cite{ChDaFaYo16b} of the same model as \cite{ArtioliBeiraoHakulaLovadina2009} confirms the asymptotics of $k(\varepsilon)$ and the very good argreement between eigenvalues of the Lam\'e operator and the eigenvalues of Koiter and Naghdi models.

The final question of interest is the overall validity of our approach. As noticed in the conclusion of \cite{ChDaFaYo16b}, the behavior at high angular frequency of the first membrane eigenvalue $\lamone{\gM^k}$ is of fundamental importance. For cylinders and trimmed cones, $\lamone{\gM^k}$ tend to $0$ as $k\to\infty$ and our asymptotic approach is always valid (at least for trimmed cones --- when the cone has a vertex we observe a similar behavior, but with a deteriorated accuracy of the 1D model). For barrels, the behavior of $\lamone{\gM^k}$ may happen to be more varied. Our approach is validated if the global infimum of $\lamone{\gM^k}$ is attained at infinity, and nowhere else. But, we have no {\em a priori} proof that this is the case. As visible for the Airy barrel of \cite{ArtioliBeiraoHakulaLovadina2009,ChDaFaYo16b}, $\lamone{\gM^k}$ has a local minimum at $k=0$, which causes axisymmetric modes to be dominant for moderately thin shells. Moreover, we have observed that for some narrow barrels, the global minimum of $\lamone{\gM^k}$ can sit at $k=0$.

\appendix

\section{High frequency reduction of the membrane operator}
\label{app:A}
Recall that the membrane operator is written as
$\bM \kk = k^2\bM_0 + k\bM_1 + \bM_2$.
Let us give elements of the proof of Theorem \ref{reducform}. We write the reduction formula \eqref{eq:reduc} in the form
\[
   \bM\kk \bV\kk = \bV_{0} \circ (\rH\kk - \Lambda\kk) + \Lambda\kk \bA\bV\kk\,.
\]
This formula holds in the sense of the formal series algebra: This means that it is equivalent to the collection of equations: For all $n \geq 0$,
\begin{equation}
\label{eq:recM}
   \bM_0 \bV_{n} + \bM_1 \bV_{n-1}  + \bM_2 \bV_{n-2}  =
   \bV_0 \circ (\rH_{n-2} - \Lambda_{n-2})+ \sum_{p + q = n - 2} \Lambda_{p} \bA  \bV_q.
\end{equation}
with the convention that $\bV_n$, $\rH_n$, and $\Lambda_n$ are $0$ for $n < 0$,
where $\bV_n$ and $\rH_n$ are the unknown coefficients of the formal series $\bV\kk = \sum_{n \geq 0} h^n \bV_n$ and $\rH\kk = \sum_{n \geq 0} h^n \rH_n$. Here the series $\Lambda\kk = \sum_{n \geq 0} h^n\Lambda_n$ is given.
The mass matrix $\bA$ is given by \eqref{eq:A}, and the operators $\bM_\ell$ by
\begin{equation}
\label{eq:M0A}
   \bM_0=\frac{E}{1-\nu^2}\begin{pmatrix}
   \frac{1-\nu}{2f^2s^2} & 0 & 0\\
   0 &\frac{1}{f^4}& 0\\
   0 & 0 & 0
\end{pmatrix} =:
\begin{pmatrix}
\rM_0^{zz} & 0 & 0 \\
0 & \rM_0^{\varphi\varphi} & 0 \\
0 & 0 & 0
\end{pmatrix}.
\end{equation}
\begin{equation}
\label{eq:M1}
\bM_1 =
i\,\frac{E}{1-\nu^2} \begin{pmatrix}
0 & -\frac{1+\nu}{2f^2s^2}\partial_{z} + \frac{2 f'}{f^3s^2} & 0\\[1ex]
\begin{matrix}
-\tfrac{1+\nu}{2f^2s^2}\partial_{z}\\
 +\Big(\frac{(\nu-3) f'}{2f^3s^2}+ \frac{(1+\nu)f'f''}{2f^2s^4}\Big)
 \end{matrix}
 & 0 &-\frac{1}{f^3s} +\frac{\nu f''}{f^2s^3}
 \\
0 &\frac{1}{f^3s}-\frac{\nu f''}{f^2s^3}  & 0
\end{pmatrix}
\end{equation}
and
\begin{equation}
\label{eq:M2}
\bM_2=\frac{E}{1-\nu^2}
\begin{pmatrix}
\begin{matrix}
-\frac{1}{s^4}\partial_{z}^2 +\Big(\tfrac{3f'f''}{s^6}-\tfrac{f'}{fs^4}\Big)\partial_{z} \\
+\tfrac{f''^2+f'f'''}{s^6}-\tfrac{4f'^2f''^2}{s^8} \\
-\tfrac{\nu f''}{fs^4}+\tfrac{(1+\nu)f'^2f''}{fs^6}+\tfrac{f'^2}{f^2s^4}
\end{matrix}
& 0 &
\begin{matrix}
\Big(\tfrac{f''}{s^5}-\tfrac{\nu}{fs^3}\Big)\partial_{z} \\
+\tfrac{f'}{f^2s^3}
+\frac{f'''}{s^5}\\
-\tfrac{3f'f''^2}{s^7}+\frac{f'f''}{fs^5}
\end{matrix}
 \\[1ex]
  0 &
  \begin{matrix}
-\tfrac{1-\nu}{2f^2s^2}\partial_{z}^2 \\
\!\!+\Big(\tfrac{(1-\nu)f'f''}{2f^2s^4}
+\tfrac{(1-\nu) f'}{2f^3s^2}\Big)\partial_{z} \!\!\\
 +\tfrac{ (1-\nu)f''}{f^3s^2}
 -\tfrac{(1-\nu)f'^2f''}{f^3s^4}
\end{matrix} &  0 \\[1ex]
\begin{matrix}
\Big(-\tfrac{f''}{s^5}+\tfrac{\nu}{fs^3}\Big)\partial_{z} \\
+\tfrac{f'f''^2}{s^7}+\tfrac{f'}{f^2s^3}
-\tfrac{2\nu f'f''}{fs^5}
\end{matrix}
& 0 &
\begin{matrix}
 \tfrac{f''^2}{s^6}
 +\tfrac{1}{f^2s^2}
-\tfrac{2\nu f''}{fs^4}
\end{matrix}
\end{pmatrix}.\hskip-1em
\end{equation}
Let us emphasize that the operators $\bM_1$ and $\bM_2$ have the following structure:
\begin{equation}
\label{eq:struct}
\bM_1 =
\begin{pmatrix}
0 & \rM_1^{z\varphi} & 0 \\
\rM_1^{\varphi z} & 0 & \rM_1^{\varphi 3} \\
0 & \rM_1^{3 \varphi} & 0
\end{pmatrix}
\quad\mbox{and}\quad
\bM_2 =
\begin{pmatrix}
\rM_2^{zz} & 0 & \rM_2^{z 3} \\
0 & \rM_2^{\varphi\varphi} & 0 \\
\rM_2^{3 z} & 0 & \rM_2^{33}
\end{pmatrix}.
\end{equation}

Let us now examine the collection of equation \eqref{eq:recM}. For $n = 0$, this equation reduces to
 $$
 \bM_{0}\bV_{0}=0,
 $$
 which is satisfied with the choice $\bV_{0}=(0,0,\mathrm{Id})^\top$. For $n = 1$, using the structure \eqref{eq:struct} of the operator $\bM_1$, the equation is $\bM_0 \bV_1 = - \bM_{1} \bV_0$ that can be written as the system
$$
\begin{pmatrix}
\rM_0^{zz}\rV_{1,z}  \\\rM_0^{\varphi\varphi} \rV_{1,\varphi} \\0
\end{pmatrix}
 =
\begin{pmatrix}
0 \\- \rM_1^{\varphi 3} \\ 0
\end{pmatrix}.
$$
Hence, solving this equation we find $\rV_{1,z} = 0$ and $\rV_{1,\varphi} = - (\rM_0^{\varphi\varphi})^{-1} \rM_1^{\varphi 3}$, i.e.,
\begin{equation}
\label{V1}
   \rV_{1,z}=0  \quad
   \rV_{1,\varphi}=\frac{if}{s}-\frac{i\nu f''f^2}{s^3}\,.
\end{equation}
The equation for $n = 2$ is written as
$$
\bM_0 \bV_{2}  = - \bM_1 \bV_{1}  - \bM_2 \bV_{0} + \bV_0 \circ (\rH_0 -  \Lambda_0) + \Lambda_0 \bA \bV_0\,.
$$
Using the structure \eqref{eq:struct} and the expressions of $\bV_0$ and $\bV_1$, it is equivalent to the system
$$
\left\{
\begin{array}{rcl}
\rM_0^{zz} \rV_{2,z} &=& - \rM_1^{z\varphi} \rV_{1,\varphi} - \rM_2^{z 3}  \\
\rM_0^{\varphi\varphi} \rV_{2,\varphi} &=& 0 \\
0 &=& - \rM_1^{3 \varphi}\rV_{1,\varphi} - \rM_2^{33} +  \rH_0 \,.
\end{array}
\right.
$$
The last equation of the previous system joint with \eqref{V1} gives the expression of the operator $\rH_0$
$$
\rH_0 = \rM_2^{33} - \rM_1^{3 \varphi} (\rM_0^{\varphi\varphi})^{-1} \rM_1^{\varphi 3},
$$
and using the first two equations, we can solve for $\bV_2$ by setting
$$
 \rV_{2,z} = (\rM_0^{zz})^{-1}( \rM_1^{z\varphi}(\rM_0^{\varphi\varphi})^{-1} \rM_1^{\varphi 3}  - \rM_2^{z 3} )
\quad\mbox{and}\quad
\rV_{2,\varphi} = (\rM_0^{\varphi\varphi})^{-1}a^{\varphi\varphi} \Lambda_0 .
$$
Thus we find that $\rH_0 = E \frac{f''^2}{s^6}$
and that the components of $\bV_2$ are
\begin{eqnarray}
\label{V2}
   \rV_{2,\varphi}&= & 0, \nonumber\\
   \rV_{2,z}&=&\Big(-\frac{f}{s}-\frac{(\nu+2)f^2f''}{s^3}\Big)\partial_{z}
   \nonumber\\ &&
   +\frac{f'}{s}+\frac{3(\nu+2)f^2f'f''^2}{s^5} -\frac{(\nu+2)f^2f'''}{s^3}-\frac{(2\nu+1)ff'f''}{s^3}\;.
\end{eqnarray}

Now let us assume that the operators $\rL_n$ and $\bV_{n+1}$ are constructed for $n \geq 1$. Then writing the equation \eqref{eq:recM} for $n +2$, we obtain the relation (using the fact that $-\bV_0\Lambda_n+\Lambda_n\bA\bV_0=0$)
$$
   \bM_0 \bV_{n+2} - \bV_0 \circ \rH_{n}  =  -\bM_1 \bV_{n+1}  - \bM_2 \bV_{n}
   +    \sum_{p=0}^{n-1} \Lambda_{p} \bA  \bV_{n-p}\,.
$$
This equation is equivalent to the system (using the fact that $\rV_{n,3} = 0$ for $n \geq 1$)
$$
\left\{
\begin{array}{rcl}
\rM_0^{\rho\rho} \rV_{n+2,z} &=&  - \rM_{1}^{z \varphi} \rV_{n+1,\varphi} - \rM_{2}^{z z} \rV_{n,z} + \sum_{p=0}^{n-1}  \Lambda_{p} a^{z z}  \rV_{n-p,z} \\[0.2ex]
\rM_0^{\varphi\varphi} \rV_{n+2,\varphi} &=&  - \rM_1^{\varphi z} \rV_{n+1,z} - \rM_{2}^{\varphi \varphi} \rV_{n,\varphi} + \sum_{p=0}^{n-1}  \Lambda_{p} a^{\varphi\varphi}  \rV_{n-p,\varphi} \\[0.2ex]
\rH_n &=& \rM_1^{3\varphi} \rV_{n+1,\varphi} + \rM_2^{3 z} \rV_{n,z}\,, \end{array}
\right.
$$
which gives the existence of the operators $\rV_{n+2,z}$, $\rV_{n+2,\varphi}$ and $\rH_n$. This shows the existence of the operators $\bV_{n}=(\rV_{n,z},\rV_{n,\varphi},0)^\top$. Moreover, we can check that $\bV_{n}$ is an operator of order $n-1$ and is polynomial in $\Lambda_j$, for $j\leq n-3$. The scalar operators $\rH_n$ are of order $n$, polynomial in $\Lambda_j$, for $j\leq n-2$.

Expressions of the operators $\rH_n$ for $n=0,\ldots,4$ are given by the formulas \eqref{6E1}--\eqref{eq:L4} and
the components of the operator $\bV_3$ are given by
\begin{align}
\label{V3}
\rV_{3,z} & =  0, \nonumber\\
\rV_{3,\varphi} & = 
   i\,\Big( -\frac{\nu f^3}{s^3}-\frac{(1+2\nu)f^4f''}{s^5}\Big) \partial_{z}^2 \nonumber\\
   & + \ i\,\Big(-\frac{(4\nu+6) f^3f'f''}{s^5} -\frac{(4\nu+2)f^4f'''}{s^5}+\frac{7(2\nu+1)f^4f'f''^2}{s^7}
     -\frac{f^2f'}{s^3}\Big) \partial_{z}\nonumber\\
   & + \ i\,\Big(\frac{(\nu^2+19\nu+19)f^3f'^2f''^2}{s^7}-\frac{(6\nu+6)f^3f'f'''}{s^5}
     -\frac{(5\nu+3)f^2f'^2f''}{s^5} \nonumber\\
   & \hskip 4ex - \frac{(36\nu+18)f^4f'^2f''^3}{s^9}+\frac{(20\nu+10)f^4f'f''f'''}{s^7}
      + \frac{ \nu f''f^2}{s^3}+  \frac{f'^2f}{s^3}\nonumber\\
	 & \hskip 4ex +\frac{ (6\nu+3)f^4f''^3}{s^7}-\frac{(2\nu+1)f^4f^{(4)}}{s^5}
	    -\frac{(\nu^2+\nu+1)f^3f''^2}{s^5}\Big)\nonumber\\
	 & + \ i\,\frac{1-\nu^2}{E}\Lambda_0\Big(\frac{f^3}{s}-\frac{\nu f^4f''}{s^3}\Big).
\end{align}

\section{Variational formulations on the meridian domain}
\label{app:B}
The 3D Lam\'e operator is independent of the azimuthal coordinate $\varphi$ if expressed in cylindrical components of displacements. The contravariant cylindrical components used in this paper are
\[
\begin{cases}
   \ru^r = \ru^{t_1}\cos\varphi + \ru^{t_2}\sin\varphi, &\mbox{(radial)}\\
   \ru^\varphi = -\ru^{t_1}\tfrac{1}{r}\sin\varphi + \ru^{t_2}\tfrac{1}{r}\cos\varphi,&\mbox{(azimuthal)}\\
   \ru^\tau = \ru^{t_3}&\mbox{(axial)}.
\end{cases}
\]
For a displacement $\bu$ defined on $\Omega^\varepsilon$, denote by $\widehat\bu_k=(\ru^r_k,\ru^\varphi_k,\ru^\tau_k)$ the Fourier coefficients of these components:
\[
   \ru^a_k(r,\tau) = \frac{1}{2\pi} \int_{0}^{2\pi} \ru^a(r,\varphi,\tau) \,e^{- i k\varphi} \,\rd \varphi,
   \quad a\in\{r,\varphi,\tau\},\quad k\in\Z,\quad (r,\tau)\in\omega^\varepsilon.
\]
The energy bilinear form $a^\varepsilon_{\aL}$ is decomposed in Fourier coefficients as follows
\[
   a^\varepsilon_{\aL}(\bu,\bv) = \sum_{k\in\Z} a^\varepsilon_{k}(\widehat\bu_k,\widehat\bv_k)\,.
\]
As soon as $k\neq0$, the energy bilinear form $a^\varepsilon_{k}$ at azimuthal frequency $k$ as non real coefficients. Nevertheless, by a simple change of components, the coefficients are back to real:
\[
	\begin{split}
	a^\varepsilon_{k}\big((\ru^r,i\ru^\varphi,&\ru^\tau),(\rv^r,-i\rv^\varphi,\rv^\tau)\big) =\\
	\frac{E}{1-\nu^2}\int_{\omega^\varepsilon}
	&\bigg\{\frac{(1-\nu)^2}{1-2\nu}\,
	\Big(r\partial_r \ru^r\partial_r \rv^r+r\partial_\tau \ru^\tau\partial_\tau \rv^\tau+\frac{k^2}{r^3}\ru^\varphi \rv^\varphi+\frac{1}{r}\ru^r\rv^r\Big)\\
	&+\frac{\nu(1-\nu)}{1-2\nu}\,  \Big[
	r\Big(\partial_r \ru^r\partial_\tau \rv^\tau+ \partial_\tau \ru^\tau\partial_r \rv^r\Big) + 
	\Big(\partial_r \ru^r\rv^r+\ru^r\partial_r \rv^r+\ru^r\partial_\tau \rv^\tau+\partial_\tau \ru^\tau\rv^r\Big)\Big]\\
	&+\frac{1-\nu}{2} \frac{k^2}{r} \Big(\ru^r\rv^r+\ru^\tau \rv^\tau\Big)\\
	&+\frac{1-\nu}{2} \frac{1}{r}\,\Big(\partial_r \ru^\varphi\partial_r\rv^\varphi-\frac{2}{r}\partial_r\ru^\varphi \rv^\varphi-\frac{2}{r}\ru^\varphi\partial_r\rv^\varphi+\frac{4}{r^2}\ru^\varphi\rv^\varphi+\partial_\tau \ru^\varphi\partial_\tau \rv^\varphi\Big)\\
	&+\frac{1-\nu}{2} \,r\Big(\partial_r \ru^\tau\partial_\tau \rv^r+\partial_\tau \ru^r\partial_r \rv^\tau+\partial_r \ru^\tau\partial_r \rv^\tau+\partial_\tau \ru^r\partial_\tau \rv^r\Big)\\
	&+k\bigg[
	\frac{(1-\nu)^2}{1-2\nu}\frac{1}{r^2}\Big(\ru^\varphi\rv^r+\ru^r\rv^\varphi\Big)
	+\frac{1-\nu}{r^2}\Big(\ru^\varphi \rv^r+\ru^r\rv^\varphi\Big)\\
	&+\frac{\nu(1-\nu)}{1-2\nu}  \frac{1}{r}\Big(\ru^\varphi\partial_r \rv^r+\partial_r \ru^r\rv^\varphi+\ru^\varphi\partial_\tau \rv^\tau+\partial_\tau \ru^\tau\rv^\varphi\Big)\\
	&-\frac{1-\nu}{2} \frac{1}{r}\Big(\ru^\tau\partial_\tau \rv^\varphi+\partial_\tau \ru^\varphi \rv^\tau+\ru^r\partial_r\rv^\varphi+\partial_r \ru^\varphi \rv^r\Big)
	\bigg]
	\bigg\} \ dr d\tau
	\end{split}
\]
The associate eigenproblem is: Find $\lambda$ and a nonzero $(\ru^r,\ru^\varphi,\ru^\tau)\in V(\omega^\varepsilon)$ such that for all $(\rv^r,\rv^\varphi,\rv^\tau)\in V(\omega^\varepsilon)$
\[
   a^\varepsilon_{k}\big((\ru^r,i\ru^\varphi,\ru^\tau),(\rv^r,-i\rv^\varphi,\rv^\tau)\big) =
   \lambda \int_{\omega^\varepsilon} [\ru^{r}\rv^{r}+\frac{1}{r^2}\ru^{\varphi}\rv^{\varphi}+\ru^{\tau}\rv^{\tau}]\ r\rd r\rd\tau.
\]
Here the variational space $V(\omega^\varepsilon)$ corresponds to $V(\Omega^\varepsilon)$
\[
V(\omega^\varepsilon) := \{\widehat\bu = (\ru^r,\ru^\varphi,\ru^\tau) \in H^1(\omega^\varepsilon)^3 \, , \quad
 \widehat\bu = 0 \quad \mbox{on}\quad \partial_0\omega^\varepsilon\}.
\]
Note that the eigenvalues of $a^\varepsilon_{k}$ are the same as those of $a^\varepsilon_{-k}$ because of the identity
\[
   a^\varepsilon_{-k}\big((\ru^r,i\ru^\varphi,\ru^\tau),(\rv^r,-i\rv^\varphi,\rv^\tau)\big) =
   a^\varepsilon_{k}\big((\ru^r,-i\ru^\varphi,\ru^\tau),(\rv^r,i\rv^\varphi,\rv^\tau)\big).
\]

\bibliographystyle{siam}
\bibliography{shell}

\end{document}